\documentclass[leqno,openany,draft]{book}

\usepackage{makeidx}
\makeindex

\def\diam{\mathop{\rm diam}}
\def\dist{\mathop{\rm dist}}
\def\supp{\mathop{\rm supp}}
\def\Lip{\mathop{\rm Lip}}

\newtheorem{theorem}{Theorem}
\newtheorem{lemma}[theorem]{Lemma}
\newtheorem{proposition}[theorem]{Proposition}
\newtheorem{definition}[theorem]{Definition}
\newtheorem{corollary}[theorem]{Corollary}

\newcommand{\begintheorem}{\addtocounter{equation}{1}\begin{theorem}}
\newcommand{\beginlemma}{\addtocounter{equation}{1}\begin{lemma}}
\newcommand{\beginproposition}{\addtocounter{equation}{1}\begin{proposition}}
\newcommand{\begindefinition}{\addtocounter{equation}{1}\begin{definition}}
\newcommand{\begincorollary}{\addtocounter{equation}{1}\begin{corollary}}

\begin{document}

\frontmatter

\title{Some aspects of analysis related to \\ $p$-adic numbers, 2}

\author{Stephen Semmes \\
        Rice University}

\date{}

\maketitle

\chapter*{Preface}

        Some aspects of analysis involving fields with absolute
value functions are discussed, which includes the real or complex
numbers with their standard absolute values, as well as ultrametric
situations like the $p$-adic numbers.

\tableofcontents

\mainmatter

\chapter{Basic notions}
\label{basic notions}

\section{Metrics and ultrametrics}
\label{metrics, ultrametrics}

        Remember that a \emph{metric}\index{metrics} on a set $M$
is a nonnegative real-valued function $d(x, y)$ defined for $x, y \in M$
that satisfies the following three conditions.  First,
\begin{equation}
\label{d(x, y) = 0 if and only if x = y}
        d(x, y) = 0 \hbox{ if and only if } x = y.
\end{equation}
Second, $d(x, y)$ should be symmetric in $x$ and $y$, so that
\begin{equation}
\label{d(x, y) = d(y, x)}
        d(x, y) = d(y, x)
\end{equation}
for every $x, y \in M$.  Third,
\begin{equation}
\label{d(x, z) le d(x, y) + d(y, z)}
        d(x, z) \le d(x, y) + d(y, z)
\end{equation}
for every $x, y, z \in M$, which is to say that $d(\cdot, \cdot)$
satisfies the triangle inequality\index{triangle inequality} on $M$.
If
\begin{equation}
\label{d(x, z) le max(d(x, y), d(y, z))}
        d(x, z) \le \max(d(x, y), d(y, z))
\end{equation}
for every $x, y, z \in M$, then $d(\cdot, \cdot)$ is said to be an
\emph{ultrametric}\index{ultrametrics} on $M$.  Of course, the
ultrametric version of the triangle inequality (\ref{d(x, z) le
  max(d(x, y), d(y, z))}) automatically implies the ordinary triangle
inequality (\ref{d(x, z) le d(x, y) + d(y, z)}).  It is easy to see
that the discrete metric\index{discrete metric} on $M$ is an
ultrametric, which is defined by putting $d(x, y)$ equal to $1$ when
$x \ne y$ and to $0$ otherwise.

        Let $(M, d(x, y))$ be a metric space,\index{metric spaces}
so that $d(x, y)$ is a metric on a set $M$.  The open and closed 
balls\index{open balls}\index{closed balls} centered at a point $x \in
M$ with radius $r > 0$ are defined as usual by
\begin{equation}
\label{B(x, r) = {z in M : d(x, z) < r}}
        B(x, r) = \{z \in M : d(x, z) < r\}
\end{equation}
and
\begin{equation}
\label{overline{B}(x, r) = {z in M : d(x, z) le r}}
        \overline{B}(x, r) = \{z \in M : d(x, z) \le r\},
\end{equation}
respectively.  It is sometimes convenient to allow $r = 0$ in
(\ref{overline{B}(x, r) = {z in M : d(x, z) le r}}), so that
(\ref{overline{B}(x, r) = {z in M : d(x, z) le r}}) reduces to
$\{x\}$.  If $y$ is any element of $B(x, r)$, then $t = r - d(x, y) >
0$, and one can check that
\begin{equation}
\label{B(y, t) subseteq B(x, r)}
        B(y, t) \subseteq B(x, r),
\end{equation}
using the triangle inequality.  Similarly, if $y \in \overline{B}(x,
r)$, then $t = r - d(x, y) \ge 0$, and
\begin{equation}
\label{overline{B}(y, t) subseteq overline{B}(x, r)}
        \overline{B}(y, t) \subseteq \overline{B}(x, r),
\end{equation}
by the triangle inequality.  One can also define collections of
open and closed subsets of $M$ in the standard way, to get a topology
on $M$ determined by the metric.  Open balls are open sets with
respect to this topology, by (\ref{B(y, t) subseteq B(x, r)}),
and it is well known that closed balls are closed sets too.

        If $d(\cdot, \cdot)$ is an ultrametric on $M$, then it is
easy to see that (\ref{B(y, t) subseteq B(x, r)}) holds with $t = r$.
More precisely,
\begin{equation}
\label{B(x, r) = B(y, r)}
        B(x, r) = B(y, r)
\end{equation}
for every $x, y \in M$ with $d(x, y) < r$, because each of the two
balls is contained in the other, by the same argument.  Similarly,
(\ref{overline{B}(y, t) subseteq overline{B}(x, r)}) holds with $t =
r$ when $d(\cdot, \cdot)$ is an ultrametric on $M$, and in fact
\begin{equation}
\label{overline{B}(x, r) = overline{B}(y, r)}
        \overline{B}(x, r) = \overline{B}(y, r)
\end{equation}
for every $x, y \in M$ with $d(x, y) \le r$.  It follows that closed
balls of positive radius are also open subsets of $M$ in this case.
One can check that open balls in $M$ are closed sets too when
$d(\cdot, \cdot)$ is an ultrametric, which is equivalent to saying
that the complement of an open ball is an open set in this situation.

        Let us continue to suppose that $d(\cdot, \cdot)$ is an
ultrametric on $M$, and let $x$, $y$, $z$ be elements of $M$.  If
$d(x, y) \le d(y, z)$, then (\ref{d(x, z) le max(d(x, y), d(y, z))})
implies that
\begin{equation}
\label{d(x, z) le d(y, z)}
        d(x, z) \le d(y, z).
\end{equation}
If $d(x, y) < d(y, z)$, then
\begin{equation}
\label{d(y, z) le max(d(y, x), d(x, z))}
        d(y, z) \le \max(d(y, x), d(x, z))
\end{equation}
implies that
\begin{equation}
\label{d(y, z) le d(x, z)}
        d(y, z) \le d(x, z).
\end{equation}
It follows that
\begin{equation}
\label{d(x, z) = d(y, z)}
        d(x, z) = d(y, z)
\end{equation}
when $d(x, y) < d(y, z)$, by combining (\ref{d(x, z) le d(y, z)})
and (\ref{d(y, z) le d(x, z)}).

\section{Quasimetrics}
\label{quasimetrics}

        Let $M$ be a set, and let $d(x, y)$ be a nonnegative real-valued
function defined for $x, y \in M$ that satisfies the first two
requirements (\ref{d(x, y) = 0 if and only if x = y}) and (\ref{d(x,
  y) = d(y, x)}) of a metric in Section \ref{metrics, ultrametrics}.
If there is a real number $C \ge 1$ such that
\begin{equation}
\label{d(x, z) le C (d(x, y) + d(y, z))}
        d(x, z) \le C \, (d(x, y) + d(y, z))
\end{equation}
for every $x, y, z \in M$, then $d(\cdot, \cdot)$ is said to be a
\emph{quasimetric}\index{quasimetrics} on $M$.  Equivalently,
$d(\cdot, \cdot)$ is a quasimetric on $M$ if there is a real number
$C' \ge 1$ such that
\begin{equation}
\label{d(x, z) le C' max(d(x, y), d(y, z))}
        d(x, z) \le C' \, \max(d(x, y), d(y, z))
\end{equation}
for every $x, y, z \in M$.  More precisely, (\ref{d(x, z) le C (d(x,
  y) + d(y, z))}) implies (\ref{d(x, z) le C' max(d(x, y), d(y, z))})
with $C'$ taken to be $2 \, C$, and (\ref{d(x, z) le C' max(d(x, y),
  d(y, z))}) implies (\ref{d(x, z) le C (d(x, y) + d(y, z))}) with $C$
taken to be $C'$.  Of course, (\ref{d(x, z) le C (d(x, y) + d(y, z))})
is the same as (\ref{d(x, z) le d(x, y) + d(y, z)}) when $C = 1$, and
(\ref{d(x, z) le C' max(d(x, y), d(y, z))}) is the same as (\ref{d(x,
  z) le max(d(x, y), d(y, z))}) when $C' = 1$.  If $d(\cdot, \cdot)$
is a quasimetric on $M$, then one can define open and closed balls in
$M$ with respect to $d(\cdot, \cdot)$ in the same way as before, as in
(\ref{B(x, r) = {z in M : d(x, z) < r}}) and (\ref{overline{B}(x, r) =
  {z in M : d(x, z) le r}}).  One can also use open balls in $M$ to
define a collection of open subsets of $M$ in the standard way, which
leads to a topology on $M$.  However, open balls in $M$ are not
necessarily open sets in this situation, and one should be a bit
careful about some other differences as well.  We shall not be too
concerned with quasimetrics here, but the terminology will sometimes
be convenient.

        Suppose that $d(x, y)$ is a quasimetric on a set $M$ that
satisfies (\ref{d(x, z) le C' max(d(x, y), d(y, z))}) for some $C' \ge
1$, and that $a$ is a positive real number.  It is easy to see that
$d(x, y)^a$ is also a quasimetric on $M$ under these conditions,
because
\begin{equation}
\label{d(x, z)^a le (C')^a max(d(x, y)^a, d(y, z)^a)}
        d(x, z)^a \le (C')^a \, \max(d(x, y)^a, d(y, z)^a)
\end{equation}
for every $x, y, z \in M$.  In particular, if $d(x, y)$ is an
ultrametric on $M$, then $d(x, y)^a$ is an ultrametric on $M$ for
every $a > 0$, since one can take $C' = 1$ in (\ref{d(x, z)^a le
  (C')^a max(d(x, y)^a, d(y, z)^a)}).  If $d(x, y)$ is any quasimetric
on $M$ and $a > 0$, then the open ball centered at a point in $M$ with
radius $r > 0$ with respect to $d(x, y)$ is the same as the open ball
centered at the same point in $M$ with radius $r^a$ with respect to
$d(x, y)^a$.  This implies that the topology on $M$ associated to
$d(x, y)^a$ is the same as the topology associated to $d(x, y)$ for
every $a > 0$.

        If $0 < a \le 1$, then it is well known that
\begin{equation}
\label{(r + t)^a le r^a + t^a}
        (r + t)^a \le r^a + t^a
\end{equation}
for every $r, t \ge 0$, which is the same as saying that
\begin{equation}
\label{r + t le (r^a + t^a)^{1/a}}
        r + t \le (r^a + t^a)^{1/a}.
\end{equation}
Indeed,
\begin{equation}
\label{max(r, t) le (r^a + t^a)^{1/a}}
        \max(r, t) \le (r^a + t^a)^{1/a}
\end{equation}
for every $r, t \ge 0$, which implies that
\begin{equation}
\label{r + t le max(r, t)^{1 - a} (r^a + t^a) le ... = (r^a + t^a)^{1/a}}
 \quad  r + t \le \max(r, t)^{1 - a} \, (r^a + t^a)
               \le (r^a + t^a)^{((1 - a) / a) + 1} = (r^a + t^a)^{1/a},
\end{equation}
as desired.  If $d(x, y)$ is a quasimetric on a set $M$ that satisfies
(\ref{d(x, z) le C (d(x, y) + d(y, z))}) for some $C \ge 1$, then it
follows that
\begin{equation}
\label{d(x, z)^a le C^a (d(x, y)^a + d(y, z)^a)}
        d(x, z)^a \le C^a \, (d(x, y)^a + d(y, z)^a)
\end{equation}
for every $x, y, z \in M$ when $0 < a \le 1$.  In particular, if $d(x,
y)$ is a metric on $M$, then $d(x, y)^a$ is also a metric on $M$ when
$0 < a \le 1$, since we can take $C = 1$ in (\ref{d(x, z)^a le C^a
  (d(x, y)^a + d(y, z)^a)}).

        If $a \ge 1$, then $f(r) = r^a$ is a convex function on $[0, \infty)$,
which implies that
\begin{equation}
\label{(r + t)^a = 2^a (r/2 + t/2)^a le ... = 2^{a - 1} (r^a + t^a)}
        (r + t)^a = 2^a \, (r/2 + t/2)^a \le 2^a \, (r^a / 2 + t^a / 2)
                                           = 2^{a - 1} \, (r^a + t^a)
\end{equation}
for every $r, t \ge 0$.  If $d(x, y)$ is again a quasimetric on a set $M$
that satisfies (\ref{d(x, z) le C (d(x, y) + d(y, z))}) for some $C \ge 1$,
then we get that
\begin{equation}
\label{d(x, z)^a le 2^{a - 1} C^a (d(x, y)^a + d(y, z)^a)}
        d(x, z)^a \le 2^{a - 1} \, C^a \, (d(x, y)^a + d(y, z)^a)
\end{equation}
for every $x, y, z \in M$ when $a \ge 1$.  This gives another way to
see that $d(x, y)^a$ is a quasimetric on $M$ for every $a > 0$ when
$d(x, y)$ is a quasimetric on $M$, using (\ref{d(x, z)^a le C^a (d(x,
  y)^a + d(y, z)^a)}) and (\ref{d(x, z)^a le 2^{a - 1} C^a (d(x, y)^a
  + d(y, z)^a)}) instead of (\ref{d(x, z)^a le (C')^a max(d(x, y)^a,
  d(y, z)^a)}).  If $M$ is the real line ${\bf R}$ and $d(x, y)$ is
the standard metric on ${\bf R}$, then it is easy to see that $d(x, y)^a$
is not a metric on ${\bf R}$ for any $a > 1$.

\section{Absolute value functions}
\label{absolute value functions}

        Let $k$ be a field.  A nonnegative real-valued function $|x|$
defined for $x \in k$ is said to be an \emph{absolute value
  function}\index{absolute value functions} on $k$ if it satisfies the
following three conditions.  First,
\begin{equation}
\label{|x| = 0 if and only if x = 0}
        |x| = 0 \hbox{ if and only if } x = 0.
\end{equation}
Second, $|\cdot|$ should be multiplicative on $k$, which is to say
that
\begin{equation}
\label{|x y| = |x| |y|}
        |x \, y| = |x| \, |y|
\end{equation}
for every $x, y \in k$.  Third, $|\cdot|$ should satisfy the triangle
inequality on $k$, in the sense that
\begin{equation}
\label{|x + y| le |x| + |y|}
        |x + y| \le |x| + |y|
\end{equation}
for every $x, y \in k$.

        Of course, the standard absolute value of a real number $x$
is defined by putting $|x| = x$ when $ x \ge 0$ and $|x| = - x$ when
$x \le 0$.  This satisfies the three conditions mentioned in the
preceding paragraph, and hence defines an absolute value function on
the field ${\bf R}$ of real numbers.\index{R@${\bf R}$} Similarly, the
usual absolute value or modulus of a complex number defines an
absolute value function on the field ${\bf C}$ of complex
numbers.\index{C@${\bf C}$}  If $k$ is any field, then the \emph{trivial
absolute value function}\index{trivial absolute value function}
is defined on $k$ by putting $|0| = 0$ and
\begin{equation}
\label{|x| = 1 for every x in k with x ne 0}
        |x| = 1 \hbox{ for every } x \in k \hbox{ with } x \ne 0.
\end{equation}
It is easy to see that this satisfies the three conditions in the
previous paragraph as well.

        Let $|\cdot|$ be an absolute value function on a field $k$.
As usual, we let $0$ and $1$ denote the additive and multiplicative
identity elements of $k$, respectively, as well as their counterparts
in ${\bf R}$, and it should always be clear from the context which is
intended.  Note that $|1| > 0$, since $1 \ne 0$ in $k$, by definition
of a field.  Because $1^2 = 1$ in $k$, we get that $|1|^2 = |1^2| =
|1|$, by (\ref{|x y| = |x| |y|}), and hence
\begin{equation}
\label{|1| = 1}
        |1| = 1.
\end{equation}
If $x \in k$ satisfies $x^n = 1$ for some positive integer $n$, then
it follows that
\begin{equation}
\label{|x|^n = |x^n| = |1| = 1}
        |x|^n = |x^n| = |1| = 1,
\end{equation}
so that $|x| = 1$.

        If $x \in k$, then the additive inverse of $x$ in $k$ is denoted
$-x$, as usual.  In particular, $-1$ denotes the additive inverse of $1$
in $k$, and it is easy to see that
\begin{equation}
\label{(-1) x = -x}
        (-1) \, x = -x
\end{equation}
for every $x \in k$.  Applying this to $x = -1$, we get that $(-1)^2 = 1$,
and hence
\begin{equation}
\label{|-1| = 1}
        |-1| = 1,
\end{equation}
by (\ref{|x|^n = |x^n| = |1| = 1}).  It follows that
\begin{equation}
\label{|-x| = |x|}
        |-x| = |x|
\end{equation}
for every $x \in k$, by (\ref{|x y| = |x| |y|}), (\ref{(-1) x = -x}),
and (\ref{|-1| = 1}).  If $x \in k$ and $x \ne 0$, then $x$ has a
multiplicative inverse $x^{-1}$ in $k$, and
\begin{equation}
\label{|x^{-1}| = |x|^{-1}}
        |x^{-1}| = |x|^{-1}
\end{equation}
by (\ref{|x y| = |x| |y|}) and (\ref{|1| = 1}).

        If we put
\begin{equation}
\label{d(x, y) = |x - y|}
        d(x, y) = |x - y|
\end{equation}
for each $x, y \in k$, then $d(x, y)$ is symmetric in $x$ and $y$, by
(\ref{|-x| = |x|}).  Thus (\ref{d(x, y) = |x - y|}) defines a metric
on $k$, since (\ref{d(x, y) = 0 if and only if x = y}) and (\ref{d(x,
  z) le d(x, y) + d(y, z)}) in Section \ref{metrics, ultrametrics}
follow from (\ref{|x| = 0 if and only if x = 0}) and (\ref{|x + y| le
  |x| + |y|}).  Let us say that $|\cdot|$ is an \emph{ultrametric
  absolute value function}\index{ultrametric absolute value functions}
on $k$ if it satisfies
\begin{equation}
\label{|x + y| le max(|x|, |y|)}
        |x + y| \le \max(|x|, |y|)
\end{equation}
for every $x, y \in k$, which implies (\ref{|x + y| le |x| + |y|}).
In this case, (\ref{d(x, y) = |x - y|}) is an ultrametric on $k$,
because (\ref{d(x, z) le max(d(x, y), d(y, z))}) in Section
\ref{metrics, ultrametrics} follows from (\ref{|x + y| le max(|x|,
  |y|)}).  The trivial absolute value function on $k$ is an
ultrametric absolute value function, for which the corresponding
ultrametric (\ref{d(x, y) = |x - y|}) is the discrete metric.

        As another class of examples, let $p$ be a prime number,
and let us recall the definition of the \emph{$p$-adic absolute
  value}\index{p-adic absolute value@$p$-adic absolute value} $|x|_p$
of a rational number $x$.  Of course, $|0|_p = 0$.  Otherwise, if $x
\ne 0$, then $x$ can be expressed as $p^j \, (a / b)$, where $a$, $b$,
and $j$ are integers, $a, b \ne 0$, and neither $a$ nor $b$ is
divisible by $p$.  In this case, we put
\begin{equation}
\label{|x|_p = p^{-j}}
        |x|_p = p^{-j},
\end{equation}
and one can check that this defines an ultrametric absolute value
function on the field ${\bf Q}$\index{Q@${\bf Q}$} of rational
numbers.  The corresponding ultrametric
\begin{equation}
\label{d_p(x, y) = |x - y|_p}
        d_p(x, y) = |x - y|_p
\end{equation}
is known as the \emph{$p$-adic metric}\index{p-adic metric@$p$-adic metric}
on ${\bf Q}$.

        If $|\cdot|$ is any absolute value function on a field $k$,
then the corresponding metric (\ref{d(x, y) = |x - y|}) determines a
topology on $k$ in the usual way.  It is easy to see that addition and
multiplication define continuous mappings from $k \times k$ into $k$
under these conditions, using the product topology on $k \times k$
associated to this topology.  Similarly, one can check that
\begin{equation}
\label{x mapsto x^{-1}}
        x \mapsto x^{-1}
\end{equation}
defines a continuous mapping from $k \setminus \{0\}$ into itself in
this situation.  The proofs of these statements are analogous to
standard arguments for real and complex numbers.

        Let $|\cdot|$ be an ultrametric absolute value function on a
field $k$, and let $x, y \in k$ be given.  If $|x| \le |y|$, then
\begin{equation}
\label{|x + y| le max(|x|, |y|) = |y|}
        |x + y| \le \max(|x|, |y|) = |y|.
\end{equation}
If $|x| < |y|$, then
\begin{equation}
\label{|y| = |-x + (x + y)| le max(|-x|, |x + y|) = max(|x|, |x + y|)}
        |y| = |-x + (x + y)| \le \max(|-x|, |x + y|) = \max(|x|, |x + y|)
\end{equation}
implies that
\begin{equation}
\label{|y| le |x + y|}
        |y| \le |x + y|.
\end{equation}
It follows that
\begin{equation}
\label{|x + y| = |y|}
        |x + y| = |y|
\end{equation}
when $|x| < |y|$, by combining (\ref{|x + y| le max(|x|, |y|) = |y|})
and (\ref{|y| le |x + y|}).  Of course, this can also be considered as
a special case of the remarks at the end of Section \ref{metrics,
  ultrametrics}.

\section{Completeness}
\label{completeness}

        Remember that a metric space $(M, d(x, y))$ is said to be
\emph{complete}\index{complete metric spaces} if every Cauchy sequence
of elements of $M$ converges to an element of $M$.  If $M$ is not
complete, then it is well known that $M$ can be completed, in the
sense that there is an isometric embedding of $M$ onto a dense subset
of a complete metric space.  Such a completion of $M$ is also unique
up to isometric equivalence.  If $d(x, y)$ is an ultrametric on $M$,
then it is not difficult to see that the completion of $M$ will also
be an ultrametric space too.

        Let $k$ be a field, let $|\cdot|$ be an absolute value function
on $k$, and consider the corresponding metric (\ref{d(x, y) = |x - y|}).
If $k$ is not already complete with respect to this metric, then $k$ can
be completed as a metric space, as in the preceding paragraph.  One
can check that the absolute value function and field operations on $k$
can also be extended to the completion in a natural way.  More precisely,
the absolute value function on $k$ is the same as the distance to $0$
with respect to the corresponding metric, and so its extension to the
completion of $k$ is already included in the metric on the completion
of $k$.  Addition and multiplication on $k$ can be extended as mappings
from $k \times k$ into $k$ to analogous mappings for the completion of $k$,
with the appropriate continuity properties.  Nonzero elements of the
completion of $k$ have multiplicative inverses in the completion of $k$,
so that the completion of $k$ becomes a field.  The extension of the
absolute value function on $k$ to the completion of $k$ is an absolute
value function on the completion of $k$, which corresponds to the metric
on the completion of $k$ in the same way as before.  If $|\cdot|$ is an
ultrametric absolute value function on $k$, then its extension to the
completion of $k$ is also an ultrametric absolute value function.

        Of course, the real and complex numbers are already complete with
respect to their standard Euclidean metrics.  The set ${\bf Q}$ of
rational numbers is not complete with respect to the standard
Euclidean metric, and its completion as a metric space corresponds to
the real line with the standard Euclidean metric.  More precisely, the
completion of ${\bf Q}$ as a field with the standard absolute value
function corresponds to the field ${\bf R}$ of real numbers with the
standard absolute value function, as in the preceding paragraph.  One
can also show that ${\bf Q}$ is not complete with respect to the
$p$-adic metric for any prime number $p$.  The completion of ${\bf Q}$
with respect to the $p$-adic metric leads to the field ${\bf
  Q}_p$\index{Q_p@${\bf Q}_p$} of \emph{$p$-adic
  numbers}.\index{p-adic numbers@$p$-adic numbers}  The extension of
the $p$-adic absolute value and metric to ${\bf Q}_p$ are denoted
$|\cdot|_p$ and $d_p(\cdot, \cdot)$, as before.  Because the
$p$-adic absolute value is an ultrametric absolute value function
on ${\bf Q}$, its extension to ${\bf Q}_p$ is an ultrametric absolute
value function as well.

        If $p$ is a prime number, $y \in {\bf Q}$, and $y \ne 0$, then
$|y|_p$ is an integer power of $p$, by construction.  Similarly, if
$y \in {\bf Q}_p$ and $y \ne 0$, then $|y|_p$ is still an integer
power of $p$.  This follows from the construction of ${\bf Q}_p$ as
the completion of ${\bf Q}$ with respect to the $p$-adic metric,
and it can also be derived from the analogous statement for ${\bf Q}$
and the fact that ${\bf Q}$ is dense in ${\bf Q}_p$, using (\ref{|x + y| = |y|})
in Section \ref{absolute value functions}.

\section{Quasimetric absolute value functions}
\label{quasimetric absolute value functions}

        Let $k$ be a field again, and let $|\cdot|$ be a nonnegative
real-valued function on $k$ that satisfies (\ref{|x| = 0 if and only
  if x = 0}) and (\ref{|x y| = |x| |y|}) in Section \ref{absolute
  value functions}.  As before, this implies that $|\cdot|$ satisfies
(\ref{|1| = 1}), (\ref{|x|^n = |x^n| = |1| = 1}), (\ref{|-1| = 1}),
(\ref{|-x| = |x|}), and (\ref{|x^{-1}| = |x|^{-1}}) in Section
\ref{absolute value functions}.  Let us say that $|\cdot|$ is a
\emph{quasimetric absolute value function}\index{quasimetric absolute
  value functions} on $k$ if there is a real number $C \ge 1$ such that
\begin{equation}
\label{|x + y| le C (|x| + |y|)}
        |x + y| \le C \, (|x| + |y|)
\end{equation}
for every $x, y \in k$.  Equivalently, $|\cdot|$ is a quasimetric
absolute value function on $k$ if there is a $C' \ge 1$ such that
\begin{equation}
\label{|x + y| le C' max(|x|, |y|)}
        |x + y| \le C' \, \max(|x|, |y|)
\end{equation}
for every $x, y \in k$.  As in Section \ref{quasimetrics}, (\ref{|x +
  y| le C (|x| + |y|)}) implies (\ref{|x + y| le C' max(|x|, |y|)})
with $C'$ taken to be $2 \, C$, and (\ref{|x + y| le C' max(|x|,
  |y|)}) implies (\ref{|x + y| le C (|x| + |y|)}) with $C$ taken to be
$C'$.  Note that (\ref{|x + y| le |x| + |y|}) in Section \ref{absolute
  value functions} is the same as (\ref{|x + y| le C (|x| + |y|)})
with $C = 1$, and that (\ref{|x + y| le max(|x|, |y|)}) in Section
\ref{absolute value functions} is the same as (\ref{|x + y| le C'
  max(|x|, |y|)}) with $C' = 1$.  If $|\cdot|$ is a quasimetric
absolute value function on $k$, then (\ref{d(x, y) = |x - y|}) is a
quasimetric on $k$, where (\ref{d(x, z) le C (d(x, y) + d(y, z))}) and
(\ref{d(x, z) le C' max(d(x, y), d(y, z))}) in Section \ref{metrics,
  ultrametrics} correspond exactly to (\ref{|x + y| le C (|x| + |y|)})
and (\ref{|x + y| le C' max(|x|, |y|)}), respectively.

        If $|x|$ is a quasimetric absolute value function on $k$, then
$|x|^a$ is also a quasimetric absolute value function on $k$ for every
positive real number $a$.  More precisely, if $|x|$ satisfies (\ref{|x
  + y| le C' max(|x|, |y|)}) for some $C' \ge 1$, then
\begin{equation}
\label{|x + y|^a le (C')^a max(|x|^a, |y|^a)}
        |x + y|^a \le (C')^a \, \max(|x|^a, |y|^a)
\end{equation}
for every $x, y \in k$.  In particular, if $|x|$ is an ultrametric
absolute value function on $k$, then $|x|^a$ is an ultrametric
absolute value function on $k$ for every $a > 0$, since one can take
$C' = 1$ in (\ref{|x + y|^a le (C')^a max(|x|^a, |y|^a)}).  Similarly,
if $|x|$ satisfies (\ref{|x + y| le C (|x| + |y|)}) for some $C \ge
1$, then
\begin{equation}
\label{|x + y|^a le C^a (|x|^a + |y|^a)}
        |x + y|^a \le C^a \, (|x|^a + |y|^a)
\end{equation}
for every $x, y \in k$ when $0 < a \le 1$, by (\ref{(r + t)^a le r^a +
  t^a}) in Section \ref{quasimetrics}.  If $a \ge 1$, then we get that
\begin{equation}
\label{|x + y|^a le 2^{a - 1} C^a (|x|^a + |y|^a)}
        |x + y|^a \le 2^{a - 1} \, C^a \, (|x|^a + |y|^a)
\end{equation}
for every $x, y \in k$, by (\ref{(r + t)^a = 2^a (r/2 + t/2)^a le
  ... = 2^{a - 1} (r^a + t^a)}) in Section \ref{quasimetrics}.  If
$|x|$ is an absolute value function on $k$, then it follows from
(\ref{|x + y|^a le C^a (|x|^a + |y|^a)}) that $|x|^a$ is an absolute
value function on $k$ when $0 < a \le 1$, since one can take $C = 1$.
If $|x|$ is the standard absolute value function on ${\bf Q}$ or ${\bf
  R}$, then $|x|^a$ is not an absolute value function for any $a > 1$.

        Let $|\cdot|$ be a nonnegative real-valued function on a field
$k$ that satisfies (\ref{|x| = 0 if and only if x = 0}) and
(\ref{|x y| = |x| |y|}) in Section \ref{absolute value functions}
again.  If $|\cdot|$ satisfies (\ref{|x + y| le C' max(|x|, |y|)})
for some $C' \ge 1$, then it follows that
\begin{equation}
\label{|1 + z| le C'}
        |1 + z| \le C'
\end{equation}
for every $z \in k$ with $|z| \le 1$, by (\ref{|1| = 1}) in Section
\ref{absolute value functions}.  Conversely, one can check that
(\ref{|1 + z| le C'}) implies (\ref{|x + y| le C' max(|x|, |y|)}),
using (\ref{|x y| = |x| |y|}) in Section \ref{absolute value
  functions}.  More precisely, (\ref{|x + y| le C' max(|x|, |y|)}) is
trivial when $x = 0$ or $y = 0$, and so one may as well suppose that
$x, y \ne 0$.  If $|x| \le |y|$, then one can apply (\ref{|1 + z| le
  C'}) to $z = x \, y^{-1}$ to get (\ref{|x + y| le C' max(|x|,
  |y|)}).  Similarly, if $|y| \le |x|$, then (\ref{|x + y| le C'
  max(|x|, |y|)}) follows from (\ref{|1 + z| le C'}) applied to $z =
x^{-1} \, y$.  Thus $|\cdot|$ is a quasimetric absolute value function
on $k$ if and only if it satisfies (\ref{|1 + z| le C'}) for some $C'
\ge 1$.  This corresponds to Definition 1.1 on p12 of \cite{c}, with
different terminology.

        If $|\cdot|$ is an absolute value function on a field $k$,
then $|\cdot|$ satisfies (\ref{|x + y| le C' max(|x|, |y|)}) with $C' = 2$.
Conversely, if $|\cdot|$ is a quasimetric absolute value function
on $k$ that satisfies (\ref{|x + y| le C' max(|x|, |y|)}) with $C' = 2$,
then it can be shown that $|\cdot|$ is an absolute value function on $k$.
See Lemma 1.2 on p13 of \cite{c}.  If $|x|$ is a quasimetric absolute
value function on $k$ that satisfies (\ref{|x + y| le C' max(|x|, |y|)})
for some $C' \ge 1$, then it follows that $|x|^a$ is an absolute value
function on $k$ for every $a > 0$ such that $(C')^a \le 2$, by
(\ref{|x + y|^a le (C')^a max(|x|^a, |y|^a)}).  In particular, this
condition holds for all sufficiently small $a > 0$.

\section{The archimedian property}
\label{archimedian property}

        Let $k$ be a field, and let ${\bf Z}_+$\index{Z_+@${\bf Z}_+$}
denote the set of positive integers, as usual.  If $x \in k$ and $n$
is a positive integer, then we let $n \cdot x$ denote the sum of $n$
$x$'s in $k$.  It is easy to see that
\begin{equation}
\label{n_1 cdot (n_2 cdot x) = (n_1 n_2) cdot x}
        n_1 \cdot (n_2 \cdot x) = (n_1 \, n_2) \cdot x
\end{equation}
for every $n_1, n_2 \in {\bf Z}_+$ and $x \in k$, and that
\begin{equation}
\label{(n cdot x) y = n cdot (x y)}
        (n \cdot x) \, y = n \cdot (x \, y)
\end{equation}
for every $n \in {\bf Z}_+$ and $x, y \in k$.  In particular,
\begin{equation}
\label{n cdot x = (n cdot 1) x}
        n \cdot x = (n \cdot 1) \, x
\end{equation}
for every $n \in {\bf Z}_+$ and $x \in k$, where $1$ refers to the
multiplicative identity element in $k$.  Similarly,
\begin{equation}
\label{(n_1 cdot x) (n_2 cdot y) = (n_1 n_2) cdot (x y)}
        (n_1 \cdot x) \, (n_2 \cdot y) = (n_1 \, n_2) \cdot (x \, y)
\end{equation}
for every $n_1, n_2 \in {\bf Z}_+$ and $x, y \in k$, which can be
verified directly, or using (\ref{n_1 cdot (n_2 cdot x) = (n_1 n_2)
  cdot x}) and (\ref{(n cdot x) y = n cdot (x y)}).

        Let $|\cdot|$ be an absolute value function on $k$.  Observe that
\begin{equation}
\label{|n cdot 1| le n}
        |n \cdot 1| \le n
\end{equation}
for every $n \in {\bf Z}_+$, by (\ref{|1| = 1}) in Section
\ref{absolute value functions}.  If $|\cdot|$ is an ultrametric
absolute value function on $k$, then
\begin{equation}
\label{|n cdot 1| le 1}
        |n \cdot 1| \le 1
\end{equation}
for every $n \in {\bf Z}_+$.  Of course,
\begin{equation}
\label{|(n_1 n_2) cdot 1| = ... = |n_1 cdot 1| |n_2 cdot 1|}
        |(n_1 \, n_2) \cdot 1| = |(n_1 \cdot 1) (n_2 \cdot 1)|
                               = |n_1 \cdot 1| \, |n_2 \cdot 1|
\end{equation}
for every $n_1, n_2 \in {\bf Z}_+$, which implies that
\begin{equation}
\label{|n^j cdot 1| = |n cdot 1|^j}
        |n^j \cdot 1| = |n \cdot 1|^j
\end{equation}
for every $j, n \in {\bf Z}_+$.  If $|n \cdot 1| > 1$ for some
$n \in {\bf Z}_+$, then it follows that (\ref{|n^j cdot 1| = |n cdot 1|^j})
tends to infinity as $j \to \infty$.

        An absolute value function $|\cdot|$ on a field $k$ is said
to be \emph{archimedian}\index{archimedian property} if the nonnegative
real numbers of the form $|n \cdot 1|$ with $n \in {\bf Z}_+$ do not
have a finite upper bound.  Otherwise, $|\cdot|$ is said to be
\emph{non-archimedian},\index{non-archimedian property} which means
that there is a positive real number $A$ such that
\begin{equation}
\label{|n cdot 1| le A}
        |n \cdot 1| \le A
\end{equation}
for every $n \in {\bf Z}_+$.  Equivalently, $|\cdot|$ is archimedian
if $|n \cdot 1| > 1$ for some $n \in {\bf Z}_+$, and $|\cdot|$ is
non-archimedian if (\ref{|n cdot 1| le A}) holds with $A = 1$, by the
remarks in the preceding paragraph.  Ultrametric absolute value
functions are obviously non-archimedian, as in (\ref{|n cdot 1| le
  1}).  Conversely, if an absolute value function $|\cdot|$ on $k$
satisfies (\ref{|n cdot 1| le 1}) for every $n \in {\bf Z}_+$, then it
can be shown that $|\cdot|$ is an ultrametric absolute value function
on $k$.  See Lemma 1.5 on p16 of \cite{c}, or Theorem 2.2.2 on p28 of
\cite{fg}.  This also works for quasimetric absolute value functions,
using an analogous argument, or by reducing to the case of ordinary
absolute value functions, as mentioned at the end of the preceding
section.

        If $k$ has positive characteristic, then there are only finitely
many elements of $k$ of the form $n \cdot 1$ for some $n \in {\bf Z}_+$.
This implies that any absolute value function on $k$ is non-archimedian,
and hence an ultrametric absolute value function.

\section{Topological equivalence}
\label{topological equivalence}

        Let $k$ be a field, and let $|\cdot|_1$ and $|\cdot|_2$ be
absolute value functions on $k$.  Also let
\begin{equation}
\label{d_1(x, y) = |x - y|_1 and d_2(x, y) = |x - y|_2}
        d_1(x, y) = |x - y|_1 \quad\hbox{and}\quad d_2(x, y) = |x - y|_2
\end{equation}
be the corresponding metrics on $k$, as in (\ref{d(x, y) = |x - y|})
in Section \ref{absolute value functions}.  Let us say that
$|\cdot|_1$ and $|\cdot|_2$ are \emph{equivalent}\index{equivalent
  absolute value functions} as absolute value functions on $k$ if
there is a positive real number $a$ such that
\begin{equation}
\label{|x|_2 = |x|_1^a}
        |x|_2 = |x|_1^a
\end{equation}
for every $x \in k$.  Of course, this implies that
\begin{equation}
\label{d_2(x, y) = d_1(x, y)^a}
        d_2(x, y) = d_1(x, y)^a
\end{equation}
for every $x, y \in k$, and hence that $d_1(x, y)$ and $d_2(x, y)$
determine the same topology on $k$.  Conversely, if $|\cdot|_1$ and
$|\cdot|_2$ are topologically equivalent, in the sense that $d_1(x,
y)$ and $d_2(x, y)$ determine the same topology on $k$, then it can be
shown that $|\cdot|_1$ and $|\cdot|_2$ are equivalent in this sense.
See Lemma 3.2 on p20 of \cite{c}, or Lemma 3.1.2 on p42 of \cite{fg}.
Part of the proof is to observe that the open unit ball in $k$ with
respect to an absolute value function can be described topologically
as the set of $x \in k$ such that $x^j \to 0$ as $j \to \infty$.  Thus
topological equivalence of the absolute value functions implies that
they have the same open unit balls in $k$, one can show that this
implies that the absolute value functions are equivalent in the sense
of (\ref{|x|_2 = |x|_1^a}).

        Of course, the trivial absolute value function on $k$ corresponds
to the discrete metric on $k$, and hence the discrete topology.
Conversely, if the topology on $k$ determined by the metric associated
to an absolute value function $|\cdot|$ is the discrete topology, then
$|\cdot|$ is the trivial absolute value function on $k$.  This follows
from the characterization of topological equivalence mentioned in the
preceding paragraph, but one can also check it more directly.  More
precisely, if the topology determined by the metric associated to
$|\cdot|$ is discrete, then the open unit ball in $k$ with respect to
$|\cdot|$ contains only $0$, because of the topological description of
the open unit ball in the previous paragraph.  It is easy to see that
this implies that $|\cdot|$ is the trivial absolute value function on
$k$, using (\ref{|x^{-1}| = |x|^{-1}}) in Section \ref{absolute value
  functions}.

        Let $|\cdot|_1$ and $|\cdot|_2$ be absolute value functions
on $k$ again, and suppose for the moment that the topology on $k$
determined by $d_1(x, y)$ is at least as strong as the topology
determined by $d_2(x, y)$.  This means that every open set in $k$ with
respect to $d_2(x, y)$ is also an open set with respect to $d_1(x,
y)$, and hence that any sequence of elements of $k$ that converges to
$0$ with respect to $d_1(x, y)$ also converges to $0$ with respect to
$d_2(x, y)$.  It follows that the open unit ball in $k$ with respect
to $|\cdot|_1$ is contained in the open unit ball with respect to
$|\cdot|_2$, by the topological description of the open unit ball
mentioned earlier.  Of course, this holds automatically when
$|\cdot|_1$ is the trivial absolute value function on $k$.

        If $|\cdot|_1$ is not the trivial absolute value function on $k$,
and if the open unit ball in $k$ with respect to $|\cdot|_1$ is
contained in the open unit ball with respect to $|\cdot|_2$, then
$|\cdot|_1$ and $|\cdot|_2$ are equivalent on $k$.  This is Lemma 3.1
on p18 of \cite{c}.  More precisely, one can show that the open unit
balls in $k$ with respect to $|\cdot|_1$ and $|\cdot|_2$ are also the
same under these conditions, and then the rest of the proof is the
same as before.  It follows that $|\cdot|_1$ and $|\cdot|_2$ are
equivalent on $k$ when $|\cdot|_1$ is not the trivial absolute value
function and the topology on $k$ determined by $d_1(x, y)$ is at least
as strong as the topology determined by $d_2(x, y)$, by the remarks
in the previous paragraph.

\section{Ostrowski's theorems}
\label{ostrowski's theorems}

        Let $|\cdot|$ be a nontrivial absolute value function on the
field ${\bf Q}$ of rational numbers.  A famous theorem of
Ostrowski\index{Ostrowski's theorems} states that $|\cdot|$ is either
equivalent to the standard Euclidean absolute value on ${\bf Q}$, or
to the $p$-adic absolute value on ${\bf Q}$ for some prime number $p$.
See Theorem 2.1 on p16 of \cite{c}, or Theorem 3.1.3 on p44 of
\cite{fg}.  More precisely, $|\cdot|$ is equivalent to the standard
Euclidean absolute value on ${\bf Q}$ exactly when $|\cdot|$ has the
archimedian property on ${\bf Q}$.  If $|\cdot|$ is a nontrivial
absolute value function on ${\bf Q}$ that is non-archimedian, then
$|n| \le 1$ for every $n \in {\bf Z}_+$, and $|n| < 1$ for some $n \in
{\bf Z}_+$.  If $p$ is the smallest positive integer such that $|p| <
1$, then $p > 1$, and it is easy to see that $p$ has to be a prime
number, because of the multiplicative property of absolute value
functions.  Under these conditions, one can show that $|\cdot|$ is
equivalent to the $p$-adic absolute value on ${\bf Q}$ for this prime
number $p$.

        Suppose now that $k$ is a field of characteristic $0$, and that
$|\cdot|$ is an absolute value function on $k$.  It is well known that there
is a natural embedding of ${\bf Q}$ into $k$ under these conditions,
so that $|\cdot|$ induces an absolute value function on ${\bf Q}$.
Note that $|\cdot|$ has the archimedian property on $k$ if and only if
the induced absolute value function has the archimedian property on
${\bf Q}$.  In this case, Ostrowki's theorem implies that the induced
absolute value function on ${\bf Q}$ is equivalent to the standard
Euclidean absolute value on ${\bf Q}$.  If $k$ is also complete with
respect to the metric associated to $|\cdot|$, then the natural
embedding of ${\bf Q}$ into $k$ can be extended continuously to an
embedding of ${\bf R}$ into $k$.

        Similarly, $|\cdot|$ is non-archimedian on $k$ if and only if
the induced absolute value function on ${\bf Q}$ is non-archimedian.
In this case, if the induced absolute value function on ${\bf Q}$
is also nontrivial, then Ostrowski's theorem implies that the induced
absolute value function on ${\bf Q}$ is equivalent to the $p$-adic
absolute value for some prime number $p$.  If $k$ is complete with
respect to the metric associated to $|\cdot|$, then the natural embedding
of ${\bf Q}$ into $k$ can be extended continuously to an embedding of
${\bf Q}_p$ into $k$.

        Let $|\cdot|$ be an archimedian absolute value function on a
field $k$, so that $k$ has characteristic $0$, as in Section
\ref{archimedian property}.  If $k$ is complete with respect to the
corresponding metric, then another famous theorem of Ostrowski implies
that $k$ is isomorphic to either the real or complex numbers, where
$|\cdot|$ corresponds to an absolute value function on ${\bf R}$ or
${\bf C}$ that is equivalent to the standard one.  See Theorem 1.1 on
p33 of \cite{c}.  As before, the natural embedding of ${\bf Q}$ into
$k$ extends continuously to an embedding of ${\bf R}$ into $k$ under
these conditions, and the first possibility is that this embedding is
surjective.  Otherwise, Ostrowki's theorem implies that $k$ is
isomorphic to the complex numbers, where this embedding of ${\bf R}$
into $k$ corresponds exactly to the standard embedding of ${\bf R}$
into ${\bf C}$.

\section{Discrete absolute value functions}
\label{discrete absolute value functions}

        Let $k$ be a field, and let $|\cdot|$ be an absolute value
function on $k$.  If $|\cdot|$ is nontrivial on $k$, then there is an
$x \in k$ such that $x \ne 0$ and $|x| \ne 1$.  More precisely, either
there is a $y \in k$ such that $y \ne 0$ and $|y| < 1$, or there is a
$z \in k$ such that $|z| > 1$.  In fact there are both such a $y$ and
$z$ in $k$, since each type of element of $k$ can be obtained from the
other by taking the multiplicative inverse.  By taking powers of such
elements of $k$, one can get nonzero elements of $k$ whose absolute
value is arbitrarily large or small.

        If $|\cdot|$ is any absolute value function on $k$, then
\begin{equation}
\label{{|x| : x in k, x ne 0}}
        \{|x| : x \in k, \, x \ne 0\}
\end{equation}
is a subgroup of the multiplicative group ${\bf R}_+$\index{R_+@${\bf
    R}_+$} of positive real numbers.  Let us say that $|\cdot|$ is
\emph{discrete}\index{discrete absolute value functions} on $k$ if
there is a positive real number $\rho < 1$ such that
\begin{equation}
\label{|x| le rho}
        |x| \le \rho
\end{equation}
for every $x \in k$ with $|x| < 1$.  Equivalently, this means that
\begin{equation}
\label{|x| ge 1/rho}
        |x| \ge 1/\rho
\end{equation}
for every $x \in k$ with $|x| > 1$, by applying (\ref{|x| le rho}) to
$1/x$.  This is also the same as saying that $1$ is not a limit point
of (\ref{{|x| : x in k, x ne 0}}) with respect to the standard metric
on ${\bf R}$.  One can check that this implies that (\ref{{|x| : x in
    k, x ne 0}}) has no limit points in ${\bf R}_+$, although $0$ is a
limit point of (\ref{{|x| : x in k, x ne 0}}) in ${\bf R}$ when
$|\cdot|$ is nontrivial, as in the previous paragraph.  Of course, the
trivial absolute value function on any field is discrete.  If $p$ is a
prime number, then the $p$-adic absolute value function is discrete on
$k = {\bf Q}$ or ${\bf Q}_p$, with $\rho = 1/p$.

        Suppose for the moment that $|\cdot|$ is an archimedian
absolute value function on $k$.  This implies that $k$ has
characteristic $0$, as in Section \ref{archimedian property}, so that
there is a natural embedding of ${\bf Q}$ into $k$.  We have also seen
that the induced absolute value function on ${\bf Q}$ is archimedian
under these conditions, and hence that the induced absolute value
function on ${\bf Q}$ is a positive power of the standard absolute
value function on ${\bf Q}$, by Ostrowski's theorem.  The standard
absolute value function on ${\bf Q}$ is obviously not discrete,
which means that $|\cdot|$ is not discrete on $k$.  This shows that
discrete absolute value functions are always non-archimedian.

        Let $|\cdot|$ be an absolute value function on a field $k$
again, and put
\begin{equation}
\label{rho_1 = sup {|x| : x in k, |x| < 1}}
        \rho_1 = \sup \{|x| : x \in k, \, |x| < 1\},
\end{equation}
so that $0 \le \rho_1 \le 1$.  Thus $\rho_1 < 1$ if and only if
$|\cdot|$ is discrete on $k$, and $\rho_1 = 0$ if and only if
$|\cdot|$ is the trivial absolute value function on $k$.  Suppose now
that $|\cdot|$ is discrete and nontrivial on $k$, so that $0 < \rho_1
< 1$.  It is not too difficult to check that there is an $x_1 \in k$
such that
\begin{equation}
\label{|x_1| = rho_1}
        |x_1| = \rho_1
\end{equation}
under these conditions.  More precisely, $\rho_1$ is an element of the
closure of (\ref{{|x| : x in k, x ne 0}}) in ${\bf R}_+$ when
$|\cdot|$ is nontrivial on $k$, and (\ref{{|x| : x in k, x ne 0}}) has
no limit points in ${\bf R}_+$ when $|\cdot|$ is discrete on $k$.
This implies that $\rho_1$ is an element of (\ref{{|x| : x in k, x ne
    0}}), so that $\rho_1$ can be expressed as in (\ref{|x_1| =
  rho_1}).  It follows that
\begin{equation}
\label{|x_1^j| = |x_1|^j = rho_1^j}
        |x_1^j| = |x_1|^j = \rho_1^j
\end{equation}
for each $j \in {\bf Z}$.  If $w$ is any nonzero element of $k$, then
one can also verify that there is a $j \in {\bf Z}$ such that
\begin{equation}
\label{|w| = rho_1^j}
        |w| = \rho_1^j,
\end{equation}
so that (\ref{{|x| : x in k, x ne 0}}) consists exactly of integer
powers of $\rho_1$.  Otherwise, if $|w|$ lies between two successive
powers of $\rho_1$, then one can use a suitable power of $x_1$ to
reduce to the case where $|w|$ lies strictly between $\rho_1$ and $1$,
contradicting the definition of $\rho_1$.

\section{Nonnegative sums}
\label{nonnegative sums}

        Let $X$ be a nonempty set, and let $f(x)$ be a nonnegative
real-valued function on $X$.  Of course, if $X$ has only finitely
many elements, then the sum
\begin{equation}
\label{sum_{x in X} f(x)}
        \sum_{x \in X} f(x)
\end{equation}
can be defined in the usual way.  Otherwise, (\ref{sum_{x in X} f(x)})
is defined as a nonnegative extended real number to be the supremum of
the sums
\begin{equation}
\label{sum_{x in A} f(x)}
        \sum_{x \in A} f(x)
\end{equation}
over all nonempty finite subsets $A$ of $X$.  If the finite sums
(\ref{sum_{x in A} f(x)}) have an upper bound in ${\bf R}$, so that
(\ref{sum_{x in X} f(x)}) is finite, then $f$ is said to be
\emph{summable}\index{summable functions} on $X$.  It is sometimes
convenient to permit sums over the empty set, in which case the sum is
interpreted as being equal to $0$.  It will also sometimes be
convenient to consider sums of nonnegative extended real numbers,
which can be handled with straightforward adjustments.  The main point
is that a sum is automatically infinite when any of the terms being
summed is infinite.

        If $X = {\bf Z}_+$, then it is customary to look at the infinite series
\begin{equation}
\label{sum_{j = 1}^infty f(j)}
        \sum_{j = 1}^\infty f(j)
\end{equation}
as the limit of the corresponding sequence
\begin{equation}
\label{sum_{j = 1}^n f(j)}
        \sum_{j = 1}^n f(j)
\end{equation}
of partial sums.  Of course, the partial sums (\ref{sum_{j = 1}^n
  f(j)}) are monotonically increasing in $n$ when $f$ is a nonnegative
real-valued function on ${\bf Z}_+$.  If the partial sums (\ref{sum_{j
    = 1}^n f(j)}) are bounded, then they converge to their supremum as
a sequence of real numbers, and otherwise they tend to $+\infty$ as $n
\to \infty$ in the usual sense.  It is easy to see that the supremum
of the partial sums (\ref{sum_{j = 1}^n f(j)}) over $n \in {\bf Z}_+$
is the same as the supremum of the sums (\ref{sum_{x in A} f(x)}) over
finite subsets $A$ of $X = {\bf Z}_+$, since any such set $A$ is
contained in the set $\{1, 2, 3, \ldots, n\}$ for some $n \in {\bf Z}_+$.
This implies that this interpretation of the sum (\ref{sum_{j = 1}^infty f(j)})
is equivalent to (\ref{sum_{x in X} f(x)}) when $X = {\bf Z}_+$.

        If $X$ is any countably infinite set, then one can enumerate the
elements of $X$ by a sequence to reduce to the case of ordinary
infinite series again.  Note that the resulting value of the sum does
not depend on the choice of the enumeration of the elements of $X$.
This follows from the fact that any rearrangement of an infinite
series of nonnegative real numbers has the same sum as the initial
series.  The definition of (\ref{sum_{x in X} f(x)}) as the supremum
of all finite subsums (\ref{sum_{x in A} f(x)}) is already invariant
under any permutation of the elements of $X$, by construction.

        Let $f$ be a nonnegative real-valued summable function on any
set $X$, and let $\epsilon > 0$ be given.  Observe that
\begin{equation}
\label{f(x) ge epsilon}
        f(x) \ge \epsilon
\end{equation}
for at most finitely many $x \in X$, since otherwise the sums
(\ref{sum_{x in A} f(x)}) would be unbounded.  More precisely, the
number of $x \in X$ such that (\ref{f(x) ge epsilon}) holds is less
than or equal to $1/\epsilon$ times (\ref{sum_{x in X} f(x)}).  It
follows that there are only finitely or countably many $x \in X$ such
that $f(x) > 0$, by taking $\epsilon = 1/n$ for each $n \in {\bf
  Z}_+$.  Thus the sum (\ref{sum_{x in X} f(x)}) can always be reduced
to a finite sum or an infinite series when it is finite.

        If $f$, $g$ are nonnegative real-valued functions on any set
$X$, then one can check that
\begin{equation}
\label{sum_{x in X} (f(x) + g(x)) = sum_{x in X} f(x) + sum_{x in X} g(x)}
        \sum_{x \in X} (f(x) + g(x)) = \sum_{x \in X} f(x) + \sum_{x \in X} g(x),
\end{equation}
where the right side of (\ref{sum_{x in X} (f(x) + g(x)) = sum_{x in
    X} f(x) + sum_{x in X} g(x)}) is interpreted as being $+\infty$
when either of the two sums is infinite.  Similarly,
\begin{equation}
\label{sum_{x in X} a f(x) = a sum_{x in X} f(x)}
        \sum_{x \in X} a \, f(x) = a \, \sum_{x \in X} f(x)
\end{equation}
for every nonnegative real-valued function $f$ on $X$ and positive
real number $a$, where the right side of (\ref{sum_{x in X} a f(x) = a
  sum_{x in X} f(x)}) is interpreted as being infinite when
(\ref{sum_{x in X} f(x)}) is infinite.  If $a = 0$, then the right
side of (\ref{sum_{x in X} a f(x) = a sum_{x in X} f(x)}) should be
interpreted as being equal to $0$, even when $f$ is not summable on
$X$.

        If $f$ is a nonnegative real-valued function on a set $X$ and
$B$, $C$ are disjoint subsets of $X$, then
\begin{equation}
\label{sum_{x in B cup C} f(x) = sum_{x in B} f(x) + sum_{x in C} f(x)}
 \sum_{x \in B \cup C} f(x) = \sum_{x \in B} f(x) + \sum_{x \in C} f(x).
\end{equation}
This can be verified directly, or derived from (\ref{sum_{x in X}
  (f(x) + g(x)) = sum_{x in X} f(x) + sum_{x in X} g(x)}).

        Let $f$ be a nonnegative real-valued summable function on a
set $X$ again, and let $\epsilon > 0$ be given.  By definition of
(\ref{sum_{x in X} f(x)}), there is a finite set $A(\epsilon) \subseteq X$
such that
\begin{equation}
\label{sum_{x in X} f(x) < sum_{x in A(epsilon)} f(x) + epsilon}
        \sum_{x \in X} f(x) < \sum_{x \in A(\epsilon)} f(x) + \epsilon.
\end{equation}
Equivalently, this means that
\begin{equation}
\label{sum_{x in X setminus A(epsilon)} f(x) < epsilon}
        \sum_{x \in X \setminus A(\epsilon)} f(x) < \epsilon,
\end{equation}
by (\ref{sum_{x in B cup C} f(x) = sum_{x in B} f(x) + sum_{x in C} f(x)}).

\section{Nonnegative sums, continued}
\label{nonnegative sums, continued}

        Let $I$ be a nonempty set, and let $E_j$ be a set for each
$j \in I$.  Suppose that
\begin{equation}
\label{E_j cap E_l = emptyset}
        E_j \cap E_l = \emptyset
\end{equation}
for every $j, l \in I$ with $j \ne l$, and put
\begin{equation}
\label{E = bigcup_{j in I} E_j}
        E = \bigcup_{j \in I} E_j.
\end{equation}
If $f$ is a nonnegative real-valued function on $E$, then
\begin{equation}
\label{sum_{x in E} f(x) = sum_{j in I} (sum_{x in E_j} f(x))}
 \sum_{x \in E} f(x) = \sum_{j \in I} \Big(\sum_{x \in E_j} f(x)\Big).
\end{equation}
More precisely, if
\begin{equation}
\label{sum_{x in E_j} f(x) = +infty}
        \sum_{x \in E_j} f(x) = +\infty
\end{equation}
for any $j \in I$, then the sum over $I$ on the right side of
(\ref{sum_{x in E} f(x) = sum_{j in I} (sum_{x in E_j} f(x))}) is
automatically interpreted as being infinite, as mentioned in the
previous section.  Otherwise, the right side of (\ref{sum_{x in E}
  f(x) = sum_{j in I} (sum_{x in E_j} f(x))}) is a sum of nonnegative
real numbers, which may still be infinite.

        In order to prove (\ref{sum_{x in E} f(x) = sum_{j in I}
(sum_{x in E_j} f(x))}), let us first observe that
\begin{equation}
\label{sum_{j in I_1} (sum_{x in E_j} f(x)) le sum_{x in E} f(x)}
 \sum_{j \in I_1} \Big(\sum_{x \in E_j} f(x)\Big) \le \sum_{x \in E} f(x)
\end{equation}
for every finite set $I_1 \subseteq I$, by (\ref{sum_{x in B cup C} f(x)
  = sum_{x in B} f(x) + sum_{x in C} f(x)}).  This implies that
\begin{equation}
\label{sum_{j in I} (sum_{x in E_j} f(x)) le sum_{x in E} f(x)}
 \sum_{j \in I} \Big(\sum_{x \in E_j} f(x)\Big) \le \sum_{x \in E} f(x),
\end{equation}
by taking the supremum over all finite subsets $I_1$ of $I$.  To get
the opposite inequality, it suffices to verify that
\begin{equation}
\label{sum_{x in A} f(x) le sum_{j in I} (sum_{x in E_j} f(x))}
 \sum_{x \in A} f(x) \le \sum_{j \in I} \Big(\sum_{x \in E_j} f(x)\Big)
\end{equation}
for every finite set $A \subseteq E$.  Of course, if $A$ is any
finite subset of $E$, then there is a finite set $I_1 \subseteq I$
such that
\begin{equation}
\label{A subseteq bigcup_{j in I_1} E_j}
        A \subseteq \bigcup_{j \in I_1} E_j.
\end{equation}
In this case, it is enough to take the sum over $j \in I_1$ on the
right side of (\ref{sum_{x in A} f(x) le sum_{j in I} (sum_{x in E_j}
  f(x))}).

        Now let $Y$ and $Z$ be nonempty sets, and let $f(y, z)$ be
a nonnegative real-valued function on their Cartesian product $Y
\times Z$.  Thus
\begin{equation}
\label{sum_{y in Y} f(y, z)}
        \sum_{y \in Y} f(y, z)
\end{equation}
can be defined as a nonnegative extended real number for each $z \in Z$
as before, and similarly
\begin{equation}
\label{sum_{z in Z} f(y, z)}
        \sum_{z \in Z} f(y, z)
\end{equation}
can be defined as an nonnegative extended real number for each $y \in Y$.
The iterated sums
\begin{equation}
\label{sum_{z in Z} (sum_{y in Y} f(y, z))}
        \sum_{z \in Z} \Big(\sum_{y \in Y} f(y, z)\Big)
\end{equation}
and
\begin{equation}
\label{sum_{y in Y} (sum_{z in Z} f(y, z))}
        \sum_{y \in Y} \Big(\sum_{z \in Z} f(y, z)\Big)
\end{equation}
are also defined as nonnegative extended real numbers, as well as the sum
\begin{equation}
\label{sum_{(y, z) in Y times Z} f(y, z)}
        \sum_{(y, z) \in Y \times Z} f(y, z)
\end{equation}
taken over $Y \times Z$ directly.  Under these conditions, the sums
(\ref{sum_{z in Z} (sum_{y in Y} f(y, z))}), (\ref{sum_{y in Y}
  (sum_{z in Z} f(y, z))}), and (\ref{sum_{(y, z) in Y times Z} f(y,
  z)}) are all equal.  More precisely, the equality of either
(\ref{sum_{z in Z} (sum_{y in Y} f(y, z))}) or (\ref{sum_{y in Y}
  (sum_{z in Z} f(y, z))}) with (\ref{sum_{(y, z) in Y times Z} f(y,
  z)}) follows from (\ref{sum_{x in E} f(x) = sum_{j in I} (sum_{x in
    E_j} f(x))}).

\chapter{Norms on vector spaces}
\label{norms on vector spaces}

        Throughout this chapter, we let $k$ be a field, and $|\cdot|$
be an absolute value function on $k$.

\section{Norms and ultranorms}
\label{norms, ultranorms}

        Let $V$ be a vector space over $k$.  A nonnegative real-valued
function $N$ on $V$ over $k$ is said to be a \emph{norm}\index{norms}
on $V$ if it satisfies the following three conditions.  First,
\begin{equation}
\label{N(v) = 0 if and only if v = 0}
        N(v) = 0 \hbox{ if and only if } v = 0.
\end{equation}
Second,
\begin{equation}
\label{N(t v) = |t| N(v)}
        N(t \, v) = |t| \, N(v)
\end{equation}
for every $t \in k$ and $v \in V$.  Third,
\begin{equation}
\label{N(v + w) le N(v) + N(w)}
        N(v + w) \le N(v) + N(w)
\end{equation}
for every $v, w \in V$.  If
\begin{equation}
\label{N(v + w) le max(N(v), N(w))}
        N(v + w) \le \max(N(v), N(w))
\end{equation}
for every $v, w \in V$, then $N$ is said to be an
\emph{ultranorm}\index{ultranorms} on $V$.  Of course, (\ref{N(v + w)
  le max(N(v), N(w))}) automatically implies (\ref{N(v + w) le N(v) +
  N(w)}), and it is easy to see that $|\cdot|$ must be an ultrametric
absolute value function on $k$ if there is an ultranorm on $V$ and $V
\ne \{0\}$.

        If $N$ is a norm on $V$, then
\begin{equation}
\label{d(v, w) = N(v - w)}
        d(v, w) = N(v - w)
\end{equation}
defines a metric on $V$, which is an ultrametric when $N$ is an
ultranorm on $V$.  Remember that $|-1| = 1$, as in (\ref{|-1| = 1}) in
Section \ref{absolute value functions}, which together with (\ref{N(t
  v) = |t| N(v)}) implies that (\ref{d(v, w) = N(v - w)}) is symmetric
in $v$ and $w$.  If $|\cdot|$ is the trivial absolute value function
on $k$, then the \emph{trivial ultranorm}\index{trivial ultranorm} on
$V$ is defined by putting $N(0) = 0$ and
\begin{equation}
\label{N(v) = 1 for every v in V with v ne 0}
        N(v) = 1 \hbox{ for every } v \in V \hbox{ with } v \ne 0.
\end{equation}
It is easy to see that this is an ultranorm on $V$, for which the
corresponding metric (\ref{d(v, w) = N(v - w)}) is the discrete
metric.

        Let $n$ be a positive integer, and let $k^n$ be the space of
$n$-tuples of elements of $k$, considered as a vector space over $k$
with respect to coordinatewise addition and scalar multiplication.
Note that
\begin{equation}
\label{N_0(v) = max(|v_1|, ldots, |v_n|)}
        N_0(v) = \max(|v_1|, \ldots, |v_n|)
\end{equation}
defines a norm on $k^n$, which is an ultranorm when $|\cdot|$ is an
ultrametric absolute value function on $k$.  Let $d_0(v, w)$ be the
metric associated to $N_0$ as in (\ref{d(v, w) = N(v - w)}).  By
construction, an open or closed ball in $k^n$ with respect to $d_0(v,
w)$ is the same as the Cartesian product of $n$ open or closed balls
in $k$ of the same radius with respect to the metric associated to
$|\cdot|$, respectively.  In particular, the topology on $k^n$
determined by $d_0(v, w)$ is the same as the product topology
corresponding to the topology on $k$ determined by the metric
associated to $|\cdot|$.  If $k$ is complete with respect to the
metric associated to $|\cdot|$, then it is easy to see that $k^n$ is
complete with respect to $d_0(v, w)$.  If $|\cdot|$ is the trivial
absolute value function on $k$, then $N_0$ is the trivial ultranorm on
$k^n$.

        Let $V$ be any vector space over $k$ again, and let $N$ be a
norm on $V$.  If $V$ is not already complete with respect to the
associated metric (\ref{d(v, w) = N(v - w)}), then one can pass to its
completion as a metric space in the usual way, as in Section
\ref{completeness}.  By standard arguments, the vector space
operations and the norm $N$ can be extended to the completion of $V$,
in such a way that the completion of $V$ also becomes a vector space
over $k$, and so that the extension of the norm $N$ to the completion
of $V$ is a norm on the completion of $V$ as a vector space over $k$
too.  Of course, if $k$ is not complete with respect to the metric
associated to the absolute value function $|\cdot|$, then one can pass
to its completion as well, as in Section \ref{completeness}.  If $V$
is already complete, and $k$ is not complete, then one can extend
scalar multiplication on $V$ to the completion of $k$, so that $V$
becomes a vector space over the completion of $k$, and so that $N$ is
a norm on $V$ as a vector space over the completion of $k$.

\section{The supremum norm}
\label{supremum norm}

        Let $X$ be a nonempty set, and let $(M, d(\cdot, \cdot))$
be a metric space.  Remember that a subset of $M$ is said to be
bounded\index{bounded sets} if it is contained in a ball of finite
radius, and that a function $f$ on $X$ with values in $M$ is said to
be bounded\index{bounded functions} if $f(X)$ is a bounded set in $M$.
Let $B(X, M)$\index{B(X, M)@$B(X, M)$} be the space of bounded
functions on $X$ with values in $M$.  If $f, g \in B(X, M)$, then
\begin{equation}
\label{d(f(x), g(x))}
        d(f(x), g(x))
\end{equation}
is a bounded nonnegative real-valued function on $X$, so that
\begin{equation}
\label{sup_{x in X} d(f(x), g(x))}
        \sup_{x \in X} d(f(x), g(x))
\end{equation}
is defined as a nonnegative real number.  It is well known that
(\ref{sup_{x in X} d(f(x), g(x))}) defines a metric on the space of
bounded functions on $X$ with values in $M$, called the \emph{supremum
  metric}.\index{supremum metric} If $d(\cdot, \cdot)$ is an
ultrametric on $M$, then it is easy to see that the supremum metric is
also an ultrametric on $B(X, M)$.  If $M$ is complete with respect to
$d(\cdot, \cdot)$, then $B(X, M)$ is complete with respect to the
supremum metric, by standard arguments.

        Now let $V$ be a vector space over $k$, and let $N$ be a norm
on $V$, so that the remarks in the preceding paragraph can be applied
to $M = V$, with the metric associated to $N$.  Of course, a
$V$-valued function $f$ on $X$ is bounded if and only if $N(f(x))$ is
bounded as a nonnegative real-valued function on $X$.  Let us use the
notation $\ell^\infty(X, V)$\index{l^infinity(X, V)@$\ell^\infty(X, V)$}
for the space of bounded $V$-valued functions $f$ on $X$, which is a
vector space over $k$ with respect to pointwise addition and scalar
multiplication.  It is easy to see that
\begin{equation}
\label{||f||_infty = ||f||_{ell^infty(X, V)} = sup_{x in X} N(f(x))}
        \|f\|_\infty = \|f\|_{\ell^\infty(X, V)} = \sup_{x \in X} N(f(x))
\end{equation}
defines a norm on $\ell^\infty(X, V)$, known as the \emph{supremum
  norm},\index{supremum norm} for which the corresponding metric is
the supremum metric.  If $N$ is an ultranorm on $V$, then the supremum
norm is an ultranorm on $\ell^\infty(X, V)$ too.

        A $V$-valued function $f$ on $X$ is said to \emph{vanish at 
infinity}\index{vanishing at infinity} on $X$ if for each $\epsilon > 0$,
\begin{equation}
\label{N(f(x)) < epsilon}
        N(f(x)) < \epsilon
\end{equation}
for all but finitely many $x \in X$.  It is easy to see that this
implies that $f$ is bounded on $V$, by taking $\epsilon = 1$.  Thus
the collection $c_0(X, V)$\index{c_0(X, V)@$c_0(X, V)$} of $V$-valued
functions on $X$ that vanish at infinity is contained in
$\ell^\infty(X, V)$, and is in fact a linear subspace of
$\ell^\infty(X, V)$.  One can also check that $c_0(X, V)$ is a closed
set in $\ell^\infty(X, V)$ with respect to the supremum metric.
Note that $f$ vanishes at infinity on $X$ if and only if $N(f(x))$
vanishes at infinity on $X$ as a real-valued function on $X$.

        The \emph{support}\index{supports of functions} of a $V$-valued
function $f$ on $X$ is defined to be the subset $\supp f$ of $X$
given by
\begin{equation}
\label{supp f = {x in X : f(x) ne 0}}
        \supp f = \{x \in X : f(x) \ne 0\}.
\end{equation}
The collection $c_{00}(X, V)$ \index{c_00(X, V)@$c_{00}(X, V)$} of
$V$-valued functions $f$ on $X$ such that $\supp f$ has only finitely
many elements is a linear subspace of $c_0(X, V)$.  If $f$ is a
$V$-valued function on $X$ that vanishes at infinity, then $f$ can be
approximated by elements of $c_{00}(X, V)$ with respect to the
supremum norm, so that $c_0(X, V)$ is the same as the closure of
$c_{00}(X, V)$ in $\ell^\infty(X, V)$.  Observe that
\begin{equation}
\label{supp f = bigcup_{j = 1}^infty {x in X : N(f(x)) ge 1/j}}
        \supp f = \bigcup_{j = 1}^\infty \{x \in X : N(f(x)) \ge 1/j\}
\end{equation}
has only finitely or countably many elements when $f \in c_0(X, V)$,
because
\begin{equation}
\label{{x in X : N(f(x)) ge 1/j}}
        \{x \in X : N(f(x)) \ge 1/j\}
\end{equation}
has only finitely many elements for each positive integer $j$.  If $N$
is the trivial ultranorm on $V$, then a $V$-valued function $f$ on $X$
vanishes at infinity only when $\supp f$ has only finitely many
elements.

\section{$\ell^r$ Norms}
\label{ell^r norms}

        Let $V$ be a vector space over $k$ again, and let $N$ be a
norm on $V$.  Also let $X$ be a nonempty set, and let $r$ be a
positive real number.  A $V$-valued function $f$ on $X$ is said to
\emph{$r$-summable}\index{rsummable functions@$r$-summable functions}
if
\begin{equation}
\label{N(f(x))^r}
        N(f(x))^r
\end{equation}
is summable as a nonnegative real-valued function on $X$, as in
Section \ref{nonnegative sums}.  If $f$ is $r$-summable with $r = 1$,
then we may simply say that $f$ is summable\index{summable functions}
on $X$.  The space of $V$-valued $r$-summable functions on $X$ is
denoted $\ell^r(X, V)$.\index{l^r(X, V)@$\ell^r(X, V)$}

        If $f$, $g$ are $V$-valued functions on $X$, then
\begin{equation}
\label{N(f(x) + g(x))^r le (N(f(x)) + N(g(x)))^r}
         N(f(x) + g(x))^r \le (N(f(x)) + N(g(x)))^r
\end{equation}
for every $r > 0$ and $x \in X$, because of the triangle inequality
(\ref{N(v + w) le N(v) + N(w)}) for $N$.  This implies that
\begin{equation}
\label{N(f(x) + g(x))^r le N(f(x))^r + N(g(x))^r}
        N(f(x) + g(x))^r \le N(f(x))^r + N(g(x))^r
\end{equation}
for every $x \in X$ when $r \le 1$, by (\ref{(r + t)^a le r^a + t^a})
in Section \ref{quasimetrics}.  Similarly,
\begin{equation}
\label{N(f(x) + g(x))^r le 2^{r - 1} (N(f(x))^r + N(g(x))^r)}
        N(f(x) + g(x))^r \le 2^{r - 1} \, (N(f(x))^r + N(g(x))^r)
\end{equation}
for every $x \in X$ when $r \ge 1$, by (\ref{(r + t)^a = 2^a (r/2 +
  t/2)^a le ... = 2^{a - 1} (r^a + t^a)}) in Section
\ref{quasimetrics}.  In both cases, we can take the sum over $x \in
X$, to get that $f + g$ is $r$-summable when $f$ and $g$ are both
$r$-summable.  Of course, if $f(x)$ is $r$-summable on $X$ and $t \in
k$, then $t \, f(x)$ is $r$-summable on $X$ too, so that $\ell^r(X,
V)$ is a vector space with respect to pointwise addition and scalar
multiplication.

        Put
\begin{equation}
\label{||f||_r = ||f||_{ell^r(X, V)} = (sum_{x in X} N(f(x))^r)^{1/r}}
        \|f\|_r = \|f\|_{\ell^r(X, V)} = \Big(\sum_{x \in X} N(f(x))^r\Big)^{1/r}
\end{equation}
for each $f \in \ell^r(X, V)$.  Thus $\|f\|_r$ is a nonnegative real
number that is equal to $0$ exactly when $f = 0$, and
\begin{equation}
\label{||t f||_r = |t| ||f||_r}
        \|t \, f\|_r = |t| \, \|f\|_r
\end{equation}
for every $t \in k$ and $f \in \ell^r(X, V)$.  If $r \ge 1$, then
\begin{equation}
\label{||f + g||_r le ||f||_r + ||g||_r}
        \|f + g\|_r \le \|f\|_r + \|g\|_r
\end{equation}
for every $f, g \in \ell^r(X, V)$.  This is well known when $V = {\bf
  R}$ with the standard absolute value function, and otherwise one can
reduce to that case using (\ref{N(f(x) + g(x))^r le (N(f(x)) +
  N(g(x)))^r}).  If $r \le 1$, then
\begin{equation}
\label{||f + g||_r^r le ||f||_r^r + ||g||_r^r}
        \|f + g\|_r^r \le \|f\|_r^r + \|g\|_r^r
\end{equation}
for every $f, g \in \ell^r(X, V)$, as one can see by summing
(\ref{N(f(x) + g(x))^r le N(f(x))^r + N(g(x))^r}) over $x \in X$.

        Suppose for the moment that $N$ is an ultranorm on $V$,
and let $f$, $g$ be $V$-valued functions on $X$ again.  In this case,
we have that
\begin{eqnarray}
\label{N(f(x) + g(x))^r le ... le N(f(x))^r + N(g(x))^r}
        N(f(x) + g(x))^r & \le & \max(N(f(x)), N(g(x)))^r \\
                          & = &\max(N(f(x))^r, N(g(x))^r) \nonumber \\
                         & \le & N(f(x))^r + N(g(x))^r \nonumber
\end{eqnarray}
for every $r > 0$ and $x \in X$.  This implies that (\ref{||f +
  g||_r^r le ||f||_r^r + ||g||_r^r}) holds for every $r > 0$ and $f, g
\in \ell^r(X, V)$, by summing (\ref{N(f(x) + g(x))^r le ... le
  N(f(x))^r + N(g(x))^r}) over $x \in X$.  Note that (\ref{||f +
  g||_r^r le ||f||_r^r + ||g||_r^r}) automatically implies (\ref{||f +
  g||_r le ||f||_r + ||g||_r}) when $r \ge 1$, by (\ref{(r + t)^a le
  r^a + t^a}) in Section \ref{quasimetrics}.

        It follows from (\ref{||f + g||_r le ||f||_r + ||g||_r}) that
$\|\cdot\|_r$ defines a norm on $\ell^r(X, V)$ when $r \ge 1$,
which determines a metric on $\ell^r(X, V)$ in the usual way.  If $r <
1$, then (\ref{||f + g||_r^r le ||f||_r^r + ||g||_r^r}) implies that
\begin{equation}
\label{||f - g||_r^r}
        \|f - g\|_r^r
\end{equation}
defines a metric on $\ell^r(X, V)$.  Note that every $r$-summable
$V$-valued function $f$ on $X$ is bounded, with
\begin{equation}
\label{||f||_infty le ||f||_r}
        \|f\|_\infty \le \|f\|_r.
\end{equation}
More precisely, such a function $f$ vanishes at infinity on $X$, so
that
\begin{equation}
\label{ell^r(X, V) subseteq c_0(X, V)}
        \ell^r(X, V) \subseteq c_0(X, V)
\end{equation}
for every $r > 0$.  If $V$ is complete with respect to the metric
associated to $N$, then one can check that $\ell^r(X, V)$ is complete
with respect to the corresponding metric for every $r > 0$, by
standard arguments.

        Suppose that $f \in \ell^q(X, V)$ for some $q > 0$, with
$q \le r$.  Thus $f$ is bounded on $X$, as in the preceding paragraph,
and
\begin{equation}
\label{N(f(x))^r le ||f||_infty^{r - q} N(f(x))^q le ||f||_q^{r - q} N(f(x))^q}
 N(f(x))^r \le \|f\|_\infty^{r - q} \, N(f(x))^q \le \|f\|_q^{r - q} \, N(f(x))^q
\end{equation}
for every $x \in X$.  This implies that $f$ is also $r$-summable on
$X$, with
\begin{equation}
\label{||f||_r^r le ||f||_q^{r - q} ||f||_q^q = ||f||_q^r}
        \|f\|_r^r \le \|f\|_q^{r - q} \, \|f\|_q^q = \|f\|_q^r,
\end{equation}
by summing (\ref{N(f(x))^r le ||f||_infty^{r - q} N(f(x))^q le
  ||f||_q^{r - q} N(f(x))^q}) over $x \in X$.  Equivalently,
\begin{equation}
\label{||f||_r le ||f||_q}
        \|f\|_r \le \|f\|_q
\end{equation}
for every $f \in \ell^q(X, V)$ when $q \le r$.

        Of course, a $V$-valued function on $X$ with finite support is
$r$-summable for every $r > 0$, so that
\begin{equation}
\label{c_{00}(X, V) subseteq ell^r(X, V)}
        c_{00}(X, V) \subseteq \ell^r(X, V)
\end{equation}
for every $r > 0$.  It is not too difficult to check that $c_{00}(X,
V)$ is dense in $\ell^r(X, V)$ for every $r > 0$, with respect to the
appropriate metric, as before.  The main point is that if $f \in
\ell^r(X, V)$, then for each $\epsilon > 0$ there is a finite set
$A(\epsilon) \subseteq X$ such that
\begin{equation}
\label{sum_{x in X setminus A(epsilon)} N(f(x))^r < epsilon}
        \sum_{x \in X \setminus A(\epsilon)} N(f(x))^r < \epsilon,
\end{equation}
as in (\ref{sum_{x in X setminus A(epsilon)} f(x) < epsilon}) in
Section \ref{nonnegative sums}.  If $N$ is the trivial ultranorm on
$V$, then every $r$-summable function on $X$ has finite support.

\section{Bounded linear mappings}
\label{bounded linear mappings}

        Let $V$, $W$ be vector spaces over $k$ equipped with norms
$N_V$, $N_W$, respectively.  A linear mapping $T$ from $V$ into $W$
is said to be \emph{bounded}\index{bounded linear mappings} if
there is a nonnegative real number $C$ such that
\begin{equation}
\label{N_W(T(v)) le C N_V(v)}
        N_W(T(v)) \le C \, N_V(v)
\end{equation}
for every $v \in V$.  This implies that
\begin{equation}
\label{N_W(T(u) - T(v)) = N_W(T(u - v)) le C N_V(u - v)}
        N_W(T(u) - T(v)) = N_W(T(u - v)) \le C \, N_V(u - v)
\end{equation}
for every $u, v \in V$, and hence that $T$ is uniformly continuous as
a mapping from $V$ into $W$, with respect to the metrics associated to
their norms.  It is easy to see that the collection $\mathcal{BL}(V,
W)$\index{BL(V, W)@$\mathcal{BL}(V, W)$} of bounded linear mappings
from $V$ into $W$ is a vector space with respect to pointwise addition
and scalar multiplication, as usual.

        As a simple class of examples, let $X$ be a nonempty set, and
consider the vector space $c_{00}(X, k)$ of $k$-vaued functions with
finite support on $X$.  Also let $\|f\|_\infty$ be the corresponding
supremum norm on $c_{00}(X, k)$, and let $\|f\|_1$ be the $\ell^1$
norm on $c_{00}(X, k)$, as in previous two sections.  Thus
\begin{equation}
\label{||f||_infty le ||f||_1}
        \|f\|_\infty \le \|f\|_1
\end{equation}
for every $f \in c_{00}(X, k)$, as in (\ref{||f||_infty le ||f||_r}),
which implies that the identity operator $I$ on $c_{00}(X, k)$ is
bounded as a linear mapping from $c_{00}(X, k)$ equipped with the
$\ell^1$ norm into $c_{00}(X, k)$ equipped with the supremum norm.  Of
course, there is an analogous statement for the standard inclusion of
$\ell^1(X, k)$ in $c_0(X, k)$.  However, if $X$ has infinitely many
elements, then there is no $C < \infty$ such that
\begin{equation}
\label{||f||_1 le C ||f||_infty}
        \|f\|_1 \le C \, \|f\|_\infty
\end{equation}
for every $f \in c_{00}(X, k)$, which means that the identity operator
$I$ is not bounded as a linear mapping from $c_{00}(X, k)$ equipped
with the supremum norm into $c_{00}(X, k)$ equipped with the $\ell^1$
norm in this case.  If $|\cdot|$ is the trivial absolute value
function on $k$, then the corresponding supremum norm $\|f\|_\infty$
on $c_{00}(X, k)$ is the same as the trivial ultranorm on $c_{00}(X,
k)$, and $c_{00}(X, k)$ is the same as $c_0(X, k)$ and $\ell^1(X, k)$.
The topology on $c_{00}(X, k)$ determined by the metric associated to
the $\ell^1$ norm is the discrete topology, which is the same as the
topology determined by the metric associated to the supremum norm.

        Let $V$ and $W$ be arbitrary vector spaces over $k$ again,
equipped with norms $N_V$ and $N_W$, and let $T$ be a linear mapping
from $V$ into $W$.  Suppose that there is a positive real number
$r$ and a nonnegative real number $A$ such that
\begin{equation}
\label{N_W(T(v)) le A}
        N_W(T(v)) \le A
\end{equation}
for every $v \in V$ with $N_V(v) < r$.  In particular, this condition
holds when $T$ is continuous at $0$ as a mapping from $V$ into $W$,
with respect to the topologies determined by the metrics associated to
the corresponding norms.  If $|\cdot|$ is not the trivial absolute
value function on $k$, then it is easy to see that $T$ has to be
bounded as a linear mapping from $V$ into $W$.  This does not always
work when $|\cdot|$ is the trivial absolute value function on $k$, as
in the preceding paragraph.

        If $T$ is a bounded linear mapping from $V$ into $W$,
then the \emph{operator norm}\index{operator norm} of $T$ is defined by
\begin{equation}
\label{||T||_{op} = ||T||_{op, VW} = inf {C ge 0 : ... holds}}
        \|T\|_{op} = \|T\|_{op, VW} = \inf \{C \ge 0 :
                      \hbox{(\ref{N_W(T(v)) le C N_V(v)}) holds}\}.
\end{equation}
It is easy to see that (\ref{N_W(T(v)) le C N_V(v)}) holds with $C =
\|T\|_{op}$, so that
\begin{equation}
\label{N_W(T(v)) le ||T||_{op} N_V(v)}
        N_W(T(v)) \le \|T\|_{op} \, N_V(v)
\end{equation}
for every $v \in V$.  One can also check that $\|T\|_{op}$ does define
a norm on $\mathcal{BL}(V, W)$, which is an ultranorm when $N_W$ is an
ultranorm on $W$.  If $W$ is complete with respect to the metric
associated to the norm $N_W$, then one can verify that
$\mathcal{BL}(V, W)$ is complete with respect to the metric associated
to the operator norm, by standard arguments.

        In some situations, the operator norm may be defined by
\begin{equation}
\label{sup {N_W(T(v)) : v in V, N_V(v) le 1}}
        \sup \{N_W(T(v)) : v \in V, \, N_V(v) \le 1\},
\end{equation}
which is clearly less than or equal to (\ref{||T||_{op} = ||T||_{op,
    VW} = inf {C ge 0 : ... holds}}).  Similarly, one might consider
\begin{equation}
\label{sup {r^{-1} N_W(T(v)) : v in V, N_V(v) le r}}
        \sup \{r^{-1} \, N_W(T(v)) : v \in V, \, N_V(v) \le r\}
\end{equation}
for any positive real number $r$, which is also automatically less
than or equal to (\ref{||T||_{op} = ||T||_{op, VW} = inf {C ge 0 :
    ... holds}}).  Thus (\ref{||T||_{op} = ||T||_{op, VW} = inf {C ge
    0 : ... holds}}) is equal to (\ref{sup {r^{-1} N_W(T(v)) : v in V,
    N_V(v) le r}}) for some $r > 0$ when (\ref{N_W(T(v)) le C N_V(v)})
holds with $C$ equal to (\ref{sup {r^{-1} N_W(T(v)) : v in V, N_V(v)
    le r}}).  In particular, this happens for every $r > 0$ when $k =
{\bf R}$ or ${\bf C}$ with the standard absolute value function, as
one can see using scalar multiplication.

        Suppose that $|\cdot|$ is an absolute value function on a
field $k$, and that $|\cdot|$ is not discrete on $k$, in the sense
described in Section \ref{discrete absolute value functions}.  This
means that there are $t \in k$ such that $|t| \ne 1$ and $|t|$ is as
close as one wants to $1$.  In this case, one can check that the
values of $|\cdot|$ on $k$ are dense in the set of nonnegative real
numbers with respect to the standard metric on the real line, using
integer powers of these elements $t$ of $k$.  This permits one to show
that (\ref{||T||_{op} = ||T||_{op, VW} = inf {C ge 0 : ... holds}}) is
equal to (\ref{sup {r^{-1} N_W(T(v)) : v in V, N_V(v) le r}}) for
every $r > 0$, using scalar multiplication again.

        Let $|\cdot|$ be any absolute value function on a field $k$
again.   If for each $v \in V$ there is a $t \in k$ such that
\begin{equation}
\label{N_V(v) = |t|}
        N_V(v) = |t|,
\end{equation}
then the same type of argument can be used to show that
(\ref{||T||_{op} = ||T||_{op, VW} = inf {C ge 0 : ... holds}}) is
equal to (\ref{sup {N_W(T(v)) : v in V, N_V(v) le 1}}).  Note that
(\ref{sup {r^{-1} N_W(T(v)) : v in V, N_V(v) le r}}) is always the
same as its counterpart with $r$ replaced by $r \, |t|$ for any $t \in
k$ with $t \ne 0$.  If $|\cdot|$ is any nontrivial absolute value
function on $k$, then (\ref{||T||_{op} = ||T||_{op, VW} = inf {C ge 0
    : ... holds}}) is less than or equal to a constant multiple of
(\ref{sup {r^{-1} N_W(T(v)) : v in V, N_V(v) le r}}) for every $r >
0$, where the constant depends only on $|\cdot|$.  Of course, these
types of arguments using scalar multiplication do not work so well
when $|\cdot|$ is the trivial absolute value function on $k$.

        Suppose that $E$ is a dense linear subspace of $V$, with
respect to the metric on $V$ associated to the norm $N_V$.  Also let
$T$ be a bounded linear mapping from $E$ into $W$, using the
restriction of $N_V$ to $E$.  Thus $T$ is uniformly continuous as a
mapping from $E$ into $W$ with respect to the corresponding metrics,
as before.  If $W$ is complete, then a well-known fact about metric
spaces implies that $T$ has a unique extension to a uniformly
continuous mapping from $V$ into $W$.  In this situation, the
extension of $T$ is a bounded linear mapping from $V$ into $W$, with
the same operator norm as $T$ has on $E$.

\section{Infinite series}
\label{infinite series}

        Let $V$ be a vector space over $k$, and let $N$ be a norm on $V$.
As usual, an infinite series $\sum_{j = 1}^\infty a_j$ with terms in
$V$ is said to \emph{converge}\index{convergence of infinite series}
if the corresponding sequence
\begin{equation}
\label{s_n = sum_{j = 1}^n a_j}
        s_n = \sum_{j = 1}^n a_j
\end{equation}
of partial sums converges to an element of $V$.  More precisely, this
means that $\{s_n\}_{n = 1}^\infty$ converges to an element of $V$
with respect to the metric $d(\cdot, \cdot)$ associated to $N$, as in
(\ref{d(v, w) = N(v - w)}) in Section \ref{norms, ultranorms}.  In
this case, the value of the sum $\sum_{j = 1}^\infty a_j$ is defined
to be the limit of the sequence $\{s_n\}_{n = 1}^\infty$.

        Note that the sequence $\{s_n\}_{n = 1}^\infty$ of partial sums
is a Cauchy sequence in $V$ with respect to $d(\cdot, \cdot)$ if and
only if for each $\epsilon > 0$ there is an $L \ge 1$ such that
\begin{equation}
\label{N(sum_{j = l + 1}^n a_j) < epsilon}
        N\Big(\sum_{j = l + 1}^n a_j\Big) < \epsilon
\end{equation}
for every $n > l \ge L$.  In particular, this implies that $\{a_j\}_{j
  = 1}^\infty$ converges to $0$ in $V$, by taking $n = l + 1$.  Of
course, $\{s_n\}_{n = 1}^\infty$ is a Cauchy sequence in $V$ when
$\{s_n\}_{n = 1}^\infty$ converges in $V$, and the converse holds when
$V$ is complete with respect to $d(\cdot, \cdot)$.

        If
\begin{equation}
\label{sum_{j = 1}^infty N(a_j)}
        \sum_{j = 1}^\infty N(a_j)
\end{equation}
converges as an infinite series of nonnegative real numbers, then we
say that $\sum_{j = 1}^\infty a_j$ \emph{converges
  absolutely}.\index{absolute convergence} It is easy to see that this
implies that $\{s_n\}_{n = 1}^\infty$ is a Cauchy sequence in $V$,
because
\begin{equation}
\label{N(sum_{j = l + 1}^n a_j) le sum_{j = l + 1}^n N(a_j)}
        N\Big(\sum_{j = l + 1}^n a_j\Big) \le \sum_{j = l + 1}^n N(a_j)
\end{equation}
for each $n > l \ge 1$.  If $V$ is complete, then it follows that
$\sum_{j = 1}^\infty a_j$ converges in $V$, in which case we also have that
\begin{equation}
\label{N(sum_{j = 1}^infty a_j) le sum_{j = 1}^infty N(a_j)}
        N\Big(\sum_{j = 1}^\infty a_j\Big) \le \sum_{j = 1}^\infty N(a_j).
\end{equation}

        Similarly, if $N$ is an ultranorm on $V$, then
\begin{equation}
\label{N(sum_{j = l + 1}^n a_j) le max_{l + 1 le j le n} N(a_j)}
        N\Big(\sum_{j = l + 1}^n a_j\Big) \le \max_{l + 1 \le j \le n} N(a_j)
\end{equation}
for every $n > l \ge 1$.  This implies that $\{s_n\}_{n = 1}^\infty$
is a Cauchy sequence in $V$ when $\{a_j\}_{j = 1}^\infty$ converges to
$0$.  If $V$ is complete, then it follows that $\sum_{j = 1}^\infty
a_j$ converges in $V$ under these conditions, and that
\begin{equation}
\label{N(sum_{j = 1}^infty a_j) le max_{j ge 1} N(a_j)}
        N\Big(\sum_{j = 1}^\infty a_j\Big) \le \max_{j \ge 1} N(a_j).
\end{equation}
Note that the maximum on the right side of (\ref{N(sum_{j = 1}^infty
  a_j) le max_{j ge 1} N(a_j)}) is attained when $N(a_j) \to 0$ as $j
\to \infty$, as a sequence of nonnegative real numbers.

\section{Generalized convergence}
\label{generalized convergence}

        Let $V$ be a vector space over $k$ again, equipped with a norm
$N$.  Also let $X$ be a nonempty set, and let $f$ be a function on $X$
with values in $V$.  If $A$ is a finite subset of $X$, then the sum
\begin{equation}
\label{sum_{x in A} f(x), 2}
        \sum_{x \in A} f(x)
\end{equation}
can be defined in the usual way, where (\ref{sum_{x in A} f(x), 2}) is
interpreted as being equal to $0$ when $A = \emptyset$.  Of course,
the collection of finite subsets of $X$ is partially ordered by
inclusion.  In fact, the collection of finite subsets of $X$ is a
directed system, because any two finite subsets $A_1$, $A_2$ of $X$
are contained in the finite subset $A_1 \cup A_2$ of $X$.  Thus the
family of finite sums (\ref{sum_{x in A} f(x), 2}) may be considered
as a net of elements of $V$ indexed by the collection of finite
subsets of $X$, and so the convergence\index{generalized convergence}
of the sum
\begin{equation}
\label{sum_{x in X} f(x), 2}
        \sum_{x \in X} f(x)
\end{equation}
may be defined in terms of the convergence of this net in $V$.  More
precisely, this net converges to an element $v$ of $V$ if for every
$\epsilon > 0$ there is a finite set $A(\epsilon) \subseteq X$ such
that
\begin{equation}
\label{N(sum_{x in A} f(x) - v) < epsilon}
        N\Big(\sum_{x \in A} f(x) - v\Big) < \epsilon
\end{equation}
for every finite set $A \subseteq X$ such that $A(\epsilon) \subseteq
A$.  It is easy to see that the limit $v$ of this net is unique when
it exists, in which case the value of the sum (\ref{sum_{x in X} f(x),
  2}) is defined to be $v$.  If $X$ has only finitely many elements,
then this reduces to the usual definition of the sum (\ref{sum_{x in
    X} f(x), 2}).

        Similarly, the net of finite sums (\ref{sum_{x in A} f(x), 2})
is a Cauchy net in $V$ if for each $\epsilon > 0$ there is a
finite set $A_1(\epsilon) \subseteq X$ such that
\begin{equation}
\label{N(sum_{x in A} f(x) - sum_{x in A'} f(x)) < epsilon}
        N\Big(\sum_{x \in A} f(x) - \sum_{x \in A'} f(x)\Big) < \epsilon
\end{equation}
for any two finite sets $A, A' \subseteq X$ such that $A_1(\epsilon)
\subseteq A, A'$.  If the net of finite sums (\ref{sum_{x in A} f(x),
  2}) converges in $V$, then it is easy to see that is a Cauchy net.
This follows from the triangle inequality, with $A_1(\epsilon)$ taken
to be the set $A(\epsilon/2)$ in the definition of convergence of the
net.  If $N$ is an ultranorm on $V$, then one can take $A_1(\epsilon)
= A(\epsilon)$.

        As a variant of this, let us say that the sum 
(\ref{sum_{x in X} f(x), 2}) satisfies the \emph{generalized Cauchy 
criterion}\index{generalized Cauchy criterion} if for each $\epsilon > 0$
there is a finite subset $A_0(\epsilon)$ of $X$ such that
\begin{equation}
\label{N(sum_{x in B} f(x)) < epsilon}
        N\Big(\sum_{x \in B} f(x)\Big) < \epsilon
\end{equation}
for every finite set $B \subseteq X$ that satisfies $A_0(\epsilon)
\cap B = \emptyset$.  If the net of finite sums (\ref{sum_{x in A}
  f(x), 2}) is a Cauchy net, then the sum (\ref{sum_{x in X} f(x), 2})
satisfies the generalized Cauchy criterion, with $A_0(\epsilon) =
A_1(\epsilon)$ for each $\epsilon > 0$.  More precisely, if $B
\subseteq X$ is a finite set such that $A_1(\epsilon) \cap B =
\emptyset$, then we can take $A = A_1(\epsilon) \cup B$ and $A' =
A_1(\epsilon)$ in (\ref{N(sum_{x in A} f(x) - sum_{x in A'} f(x)) <
  epsilon}) to get (\ref{N(sum_{x in B} f(x)) < epsilon}).
Conversely, if the sum (\ref{sum_{x in X} f(x), 2}) satisfies the
generalized Cauchy criterion, then the net of finite sums (\ref{sum_{x
    in A} f(x), 2}) is a Cauchy net, with $A_1(\epsilon) =
A_0(\epsilon/2)$ for each $\epsilon > 0$.  To see this, let $\epsilon
> 0$ be given, and let $A, A' \subseteq X$ be finite sets such that
$A_0(\epsilon/2) \subseteq A, A'$.  Put $B = A \setminus (A \cap A')$
and $B' = A' \setminus (A \cap A')$, so that
\begin{equation}
\label{sum_{x in A} f(x) - sum_{x in A'} f(x) = ...}
        \sum_{x \in A} f(x) - \sum_{x \in A'} f(x)
          = \sum_{x \in B} f(x) - \sum_{x \in B'} f(x),
\end{equation}
and hence
\begin{equation}
\label{N(sum_{x in A} f(x) - sum_{x in A'} f(x)) le ...}
        N\Big(\sum_{x \in A} f(x) - \sum_{x \in A'} f(x)\Big)
            \le N\Big(\sum_{x \in B} f(x)\Big) + N\Big(\sum_{x \in B'} f(x)\Big).
\end{equation}
By hypothesis, both terms on the right side of (\ref{N(sum_{x in A}
  f(x) - sum_{x in A'} f(x)) le ...}) are less than $\epsilon / 2$,
since $B$ and $B'$ are disjoint from $A_0(\epsilon/2)$.  This implies
that (\ref{N(sum_{x in A} f(x) - sum_{x in A'} f(x)) < epsilon})
holds, as desired.  If $N$ is an ultranorm on $V$, then one can take
$A_1(\epsilon) = A_0(\epsilon)$ for every $\epsilon > 0$, by an
analogous argument.

        Suppose that $f$ is summable as a $V$-valued function on $X$,
as in Section \ref{ell^r norms}, and let us check the sum (\ref{sum_{x
    in X} f(x), 2}) satisfies the generalized Cauchy criterion.
Remember that the summability of $f$ on $X$ means that $N(f(x))$ is
summable as a nonnegative real-valued function on $X$, as in Section
\ref{nonnegative sums}.  This implies that for each $\epsilon > 0$
there is a finite subset $A_0(\epsilon)$ of $X$ such that
\begin{equation}
\label{sum_{x in X setminus A_0(epsilon)} N(f(x)) < epsilon}
        \sum_{x \in X \setminus A_0(\epsilon)} N(f(x)) < \epsilon,
\end{equation}
as in (\ref{sum_{x in X setminus A(epsilon)} f(x) < epsilon}) in
Section \ref{nonnegative sums}.  If $B$ is a finite subset of $X
\setminus A_0(\epsilon)$, then it follows that
\begin{equation}
\label{N(sum_{x in B} f(x)) le sum_{x in B} N(f(x)) le ... < epsilon}
        N\Big(\sum_{x \in B} f(x)\Big) \le \sum_{x \in B} N(f(x))
          \le \sum_{x \in X \setminus A_0(\epsilon)} N(f(x)) < \epsilon,
\end{equation}
as desired.  Similarly, if $N$ is an ultranorm on $V$, and if $f$ is a
$V$-valued function on $X$ that vanishes at infinity, then the sum
(\ref{sum_{x in X} f(x), 2}) satisfies the generalized Cauchy
criterion.  In this case, we can simply take
\begin{equation}
\label{A_0(epsilon) = {x in X : N(f(x)) ge epsilon}}
        A_0(\epsilon) = \{x \in X : N(f(x)) \ge \epsilon\}
\end{equation}
for each $\epsilon > 0$, which has only finitely many elements by
hypothesis.  If $B$ is a finite subset of $X \setminus A_0(\epsilon)$,
then
\begin{equation}
\label{N(sum_{x in B} f(x)) le max_{x in B} N(f(x)) < epsilon}
        N\Big(\sum_{x \in B} f(x)\Big) \le \max_{x \in B} N(f(x)) < \epsilon,
\end{equation}
as desired.

\section{Generalized convergence, continued}
\label{generalized convergence, continued}

        Let $V$ be a vector space over $k$ equipped with a norm $N$
again, and let $f$ be a $V$-valued function on a nonempty set $X$ such
that the sum (\ref{sum_{x in X} f(x), 2}) satisfies the generalized
Cauchy criterion.  This implies that $f$ vanishes at infinity on $X$,
by applying (\ref{N(sum_{x in B} f(x)) < epsilon}) to sets $B$ with
exactly one element.  In particular, it follows that the support of
$f$ has only finitely or countably many elements, as in Section
\ref{supremum norm}.  Of course, if the support of $f$ has only
finitely many elements, then the sum (\ref{sum_{x in X} f(x), 2}) can
be defined in the usual way.

        Otherwise, let $\{x_j\}_{j = 1}^\infty$ be a sequence of distinct
elements of $X$ such that every element of the support of $f$ is of
the form $x_j$ for some $j$.  It is easy to see that the infinite series
\begin{equation}
\label{sum_{j = 1}^infty f(x_j)}
        \sum_{j = 1}^\infty f(x_j)
\end{equation}
satisfies the usual Cauchy criterion under these conditions, as in
Section \ref{infinite series}, so that the partial sums of this series
form a Cauchy sequence in $V$.  If $V$ is complete with respect to the
metric associated to $N$, then it follows that (\ref{sum_{j = 1}^infty
  f(x_j)}) converges as an infinite series in $V$.  Using the
generalized Cauchy criterion for the sum (\ref{sum_{x in X} f(x), 2})
again, one can check that the net of finite sums (\ref{sum_{x in A}
  f(x), 2}) converges in $V$ to the same value of the sum (\ref{sum_{j
    = 1}^infty f(x_j)}), as in Section \ref{generalized convergence}.

        If $f$ is summable on $X$, then (\ref{sum_{x in X} f(x), 2})
satisfies the generalized Cauchy criterion, as in the previous
section.  In this case, it is easy to see that
\begin{equation}
\label{N(sum_{x in X} f(x)) le sum_{x in X} N(f(x))}
        N\Big(\sum_{x \in X} f(x)\Big) \le \sum_{x \in X} N(f(x))
\end{equation}
when $V$ is complete with respect to the metric associated to $N$, so
that the net of finite sums (\ref{sum_{x in A} f(x), 2}) converges in
$V$.  Similarly, if $N$ is an ultranorm on $V$ and $f$ vanishes at
infinity on $X$, then we have seen that (\ref{sum_{x in X} f(x), 2})
satisfies the generalized Cauchy criterion again.  If $V$ is also
complete with respect to the ultranorm associated to $N$, so that
(\ref{sum_{x in X} f(x), 2}) is defined as an element of $V$, then we
have that
\begin{equation}
\label{N(sum_{x in X} f(x)) le max_{x in X} N(f(x))}
        N\Big(\sum_{x \in X} f(x)\Big) \le \max_{x \in X} N(f(x)).
\end{equation}
Note that the maximum on the right side of (\ref{N(sum_{x in X} f(x))
  le max_{x in X} N(f(x))}) is attained when $f$ vanishes at infinity
on $X$.

        If $f$ is a nonnegative real-valued function on $X$ which
is summable, then it is easy to see that the net of finite sums
(\ref{sum_{x in A} f(x), 2}) converges to their supremum.  Thus the
definition of the sum (\ref{sum_{x in X} f(x), 2}) in Section
\ref{nonnegative sums} is equivalent to the definition of the sum in
Section \ref{generalized convergence} in this situation.  Similarly,
if $f$ is a nonnegative real-valued function on $X$ that is not
summable, then the finite sums (\ref{sum_{x in A} f(x), 2}) tend to
$+\infty$ in a suitable sense.  If $f$ is a real or complex-valued
summable function on $X$, then $f$ can be expressed as a linear
combination of summable nonnegative real-valued functions on $X$.
This gives another way to look at the convergence of the net of finite
sums (\ref{sum_{x in A} f(x), 2}) in this case.

        An infinite series with terms in $V$ may be considered as a
sum over $X = {\bf Z}_+$, to which the earlier discussion applies.  If
the corresponding net of all finite subsums of such a sum over ${\bf
  Z}_+$ converges in $V$, then the partial sums of any rearrangement
of the series converges to the same value.  Similarly, if a sum over
${\bf Z}_+$ satisfies the generalized Cauchy criterion, as in Section
\ref{generalized convergence}, then the sequence of partial sums of
any rearrangement of the series is a Cauchy sequence.  In this case,
the convergence of any of these Cauchy sequences implies the
convergence of the whole net of finite sums to the same value, as
before.

        Let $X$ be any nonempty set again, and let $f$ be a $V$-valued
function on $X$ such that the sum (\ref{sum_{x in X} f(x), 2})
satisfies the generalized Cauchy criterion, as in the previous
section.  If $E$ is any nonempty subset of $X$, then it is easy to see
that the restriction of $f$ to $E$ has the same property, so that
\begin{equation}
\label{sum_{x in E} f(x)}
        \sum_{x \in E} f(x)
\end{equation}
satisfies the generalized Cauchy criterion as a sum over $E$.  If $V$
is complete with respect to the metric associated to $N$, then it
follows that the net of all finite subsums of (\ref{sum_{x in E}
  f(x)}) converges to an element of $V$ for each $E \subseteq X$, as
before.  If $E_1$ and $E_2$ are pairwise-disjoint subsets of $X$, then
one can also check that
\begin{equation}
\label{sum_{x in E_1 cup E_2} f(x) = sum_{x in E_1} f(x) + sum_{x in E_2} f(x)}
 \sum_{x \in E_1 \cup E_2} f(x) = \sum_{x \in E_1} f(x) + \sum_{x \in E_2} f(x)
\end{equation}
under these conditions.

        Suppose now that $f$ is a $V$-valued function on a nonempty set $X$
such that the sum (\ref{sum_{x in X} f(x), 2}) does not satisfies the
generalized Cauchy criterion.  This means that there is an $\epsilon >
0$ such that for each finite set $A \subseteq X$ there is a finite set
$B \subseteq X \setminus A$ that satisfies
\begin{equation}
\label{N(sum_{x in B} f(x)) ge epsilon}
        N\Big(\sum_{x \in B} f(x)\Big) \ge \epsilon.
\end{equation}
If $X = {\bf Z}_+$, then one can use this to find a rearrangement of
the infinite series $\sum_{j = 1}^\infty f(j)$ for which the
corresponding sequence of partial sums does not form a Cauchy
sequence.  Alternatively, if $X = {\bf Z}_+$, then one can use this to
find a strictly increasing sequence $\{j_l\}_{l = 1}^\infty$ of
positive integers such that the sequence of partial sums
\begin{equation}
\label{sum_{l = 1}^n f(j_l)}
        \sum_{l = 1}^n f(j_l)
\end{equation}
does not form a Cauchy sequence in $V$.

\section{Bounded finite sums}
\label{bounded finite sums}

        Let $V$ be a vector space over $k$ with a norm $N$ again,
and let $X$ be a nonemnpty set.  Let us say that a $V$-valued function
$f$ on $X$ has \emph{bounded finite sums}\index{bounded finite sums}
if the sums
\begin{equation}
\label{sum_{x in A} f(x), 3}
        \sum_{x \in A} f(x)
\end{equation}
over all finite subsets $A$ of $X$ have bounded norm in $V$.  More
precisely, this means that there is a nonnegative real number $C$
depending on $f$ such that
\begin{equation}
\label{N(sum_{x in A} f(x)) le C}
        N\Big(\sum_{x \in A} f(x)\Big) \le C
\end{equation}
for every finite set $A \subseteq X$.  It is easy to see that the
space $BFS(X, V)$\index{BFS(X, V)@$BFS(X, V)$} of $V$-valued functions
on $X$ with bounded finite sums is a vector space with respect to
pointwise addition and scalar multiplication, and that
\begin{eqnarray}
\label{||f||_{BFS} = sup {N(sum_{x in A} f(x)) : A is a finite subset of X}}
        \|f\|_{BFS} & = & \|f\|_{BFS(X, V)} \\
 & = & \sup \bigg\{N\Big(\sum_{x \in A} f(x)\Big)
                         : A \hbox{ is a finite subset of } X\bigg\} \nonumber
\end{eqnarray}
defines a norm on $BFS(X, V)$.  In fact, $BFS(X, V)$ is a linear
subspace of $\ell^\infty(X, V)$, and
\begin{equation}
\label{||f||_infty le ||f||_{BFS}}
        \|f\|_\infty \le \|f\|_{BFS}
\end{equation}
for every $f \in BFS(X, V)$, as one can see by restricting one's attention
to subsets $A$ of $X$ with exactly one element.  If $N$ is an ultranorm
on $V$, then every bounded $V$-valued function on $X$ has bounded finite
sums, so that $BFS(X, V)$ is the same as $\ell^\infty(X, V)$, and
\begin{equation}
\label{||f||_infty = ||f||_{BFS}}
        \|f\|_\infty = \|f\|_{BFS}
\end{equation}
for every $f \in \ell^\infty(X, V)$.  If $V$ is complete with respect
to the metric associated to $N$, then one can check that $BFS(X, V)$
is complete with respect to the metric associated to the BFS norm,
by standard arguments.

        Let $f$ be a $V$-valued function on $X$ for which the sum
\begin{equation}
\label{sum_{x in X} f(x), 3}
        \sum_{x \in X} f(x)
\end{equation}
satisfies the generalized Cauchy criterion, and let us check that $f$
has bounded finite sums.  To do this, let $A_0(1)$ be a finite subset
of $X$ such that (\ref{N(sum_{x in B} f(x)) < epsilon}) in Section
\ref{generalized convergence} holds with $\epsilon = 1$ for every
finite set $B \subseteq X$ that is disjoint from $A_0(1)$.  If $A$ is
any finite subset of $X$, then
\begin{equation}
\label{sum_{x in A} f(x) = ...}
        \sum_{x \in A} f(x) = \sum_{x \in A \cap A_0(1)} f(x)
                               + \sum_{x \in A \setminus A_0(1)} f(x),
\end{equation}
so that
\begin{equation}
\label{N(sum_{x in A} f(x)) le ...}
 N\Big(\sum_{x \in A} f(x)\Big) \le N\Big(\sum_{x \in A \cap A_0(1)} f(x)\Big)
                               + N\Big(\sum_{x \in A \setminus A_0(1)} f(x)\Big).
\end{equation}
Because $A \setminus A_0(1)$ is disjoint from $A_0(1)$, we get that
\begin{equation}
\label{N(sum_{x in A} f(x)) < sum_{x in A_0(1)} N(f(x)) + 1}
        N\Big(\sum_{x \in A} f(x)\Big) < \sum_{x \in A_0(1)} N(f(x)) + 1,
\end{equation}
using the triangle inequality to estimate the first term on the right
side of (\ref{N(sum_{x in A} f(x)) le ...}).  This shows that the sums
(\ref{sum_{x in A} f(x), 2}) have bounded norm, and of course one could
get a better estimate when $N$ is an ultranorm on $V$.

        Let $GCC(X, V)$\index{GCC(X, V)@$GCC(X, V)$} be the space of
$V$-valued functions $f$ on $X$ such that (\ref{sum_{x in X} f(x), 3})
satisfies the generalized Cauchy criterion, as in Section
\ref{generalized convergence}.  The argument in the preceding
paragraph implies that
\begin{equation}
\label{GCC(X, V) subseteq BFS(X, V)}
        GCC(X, V) \subseteq BFS(X, V),
\end{equation}
and it is easy to see that $GCC(X, V)$ is a linear subspace of $BFS(X,
V)$.  One can also check that $GCC(X, V)$ is a closed set in $BFS(X,
V)$, with respect to the metric on $BFS(X, V)$ associated to the BFS
norm.  Of course,
\begin{equation}
\label{c_{00}(X, V) subseteq GCC(X, V)}
        c_{00}(X, V) \subseteq GCC(X, V),
\end{equation}
and in fact $GCC(X, V)$ is the same as the closure of $c_{00}(X, V)$
in $BFS(X, V)$ with respect to the BFS norm.  Note that
\begin{equation}
\label{GCC(X, V) subseteq c_0(X, V)}
        GCC(X, V) \subseteq c_0(X, V),
\end{equation}
as mentioned at the beginning of Section \ref{generalized convergence,
  continued}, and that
\begin{equation}
\label{GCC(X, V) = c_0(X, V)}
        GCC(X, V) = c_0(X, V)
\end{equation}
when $N$ is an ultranorm on $V$, as indicated near the end of Section
\ref{generalized convergence}.

        If $N$ is any norm on $V$ and $f$ is a $V$-valued summable
function on $X$, then
\begin{equation}
\label{N(sum_{x in A} f(x)) le sum_{x in A} N(f(x)) le ||f||_1}
        N\Big(\sum_{x \in A} f(x)\Big) \le \sum_{x \in A} N(f(x)) \le \|f\|_1
\end{equation}
for every finite set $A \subseteq X$.  Thus $f$ has bounded finite
sums on $X$, and
\begin{equation}
\label{||f||_{BFS} le ||f||_1}
        \|f\|_{BFS} \le \|f\|_1.
\end{equation}
More precisely,
\begin{equation}
\label{ell^1(X, V) subseteq GCC(X, V)}
        \ell^1(X, V) \subseteq GCC(X, V),
\end{equation}
as in Section \ref{generalized convergence}.  Alternatively, we have
seen that $c_{00}(X, V)$ is dense in $\ell^1(X, V)$ with respect to
the $\ell^1$ norm, as in Section \ref{ell^r norms}.  This implies that
every $V$-valued summable function $f$ on $X$ can be approximated by
functions with finite support with respect to the BFS norm, by
(\ref{||f||_{BFS} le ||f||_1}), so that $f \in GCC(X, V)$, as in the
preceding paragraph.

        Let $f$ be a real-valued function on $X$ with bounded finite sums
with respect to the standard absolute value function on ${\bf R}$, so that
\begin{equation}
\label{|sum_{x in A} f(x)| le C}
        \biggl|\sum_{x \in A} f(x)\biggr| \le C
\end{equation}
for some nonnegative real number $C$ and every finite set $A \subseteq
X$.  This implies that
\begin{equation}
\label{sum_{x in A} |f(x)| le 2 C}
        \sum_{x \in A} |f(x)| \le 2 \, C
\end{equation}
for every finite set $A \subseteq X$, by applying (\ref{|sum_{x in A}
  f(x)| le C}) to the subsets of $A$ consisting of $x \in A$ such that
$f(x) \ge 0$ or $f(x) \le 0$, respectively.  It follows that $f$ is a
summable function on $X$ under these conditions, with
\begin{equation}
\label{sum_{x in X} |f(x)| le 2 C}
        \sum_{x \in X} |f(x)| \le 2 \, C.
\end{equation}
Similarly, if $f$ is a complex-valued function on $X$ with bounded
finite sums with respect to the standard absolute value function on
${\bf C}$, then one can apply the previous remarks to the real and
imaginary parts of $f$, to get that $f$ is summable on $X$.

        Let $V$ be a vector space over $k$ with a norm $N$ again,
and suppose that $V$ is complete with respect to the corresponding
metric.  If $f \in GCC(X, V)$, then the net of finite sums
(\ref{sum_{x in A} f(x), 3}) converges in $V$, as in Section
\ref{generalized convergence, continued}.  The value of the sum
(\ref{sum_{x in X} f(x), 3}) is defined to be the limit of this net,
which satisfies
\begin{equation}
\label{N(sum_{x in X} f(x)) le ||f||_{BFS}}
        N\Big(\sum_{x \in X} f(x)\Big) \le \|f\|_{BFS}.
\end{equation}
Thus
\begin{equation}
\label{f mapsto sum_{x in X} f(x)}
        f \mapsto \sum_{x \in X} f(x)
\end{equation}
defines a bounded linear mapping from $GCC(X, V)$ into $V$, using the
restriction of the BFS norm to $GCC(X, V)$.  More precisely, it is
easy to see that the operator norm of (\ref{f mapsto sum_{x in X}
  f(x)}) is equal to $1$, by considering $V$-valued functions $f(x)$
on $X$ that are equal to $0$ at all but one point in $X$.

        Of course, the sum (\ref{sum_{x in X} f(x), 3}) can be defined
in the usual way when $f$ has finite support on $X$, so that (\ref{f
  mapsto sum_{x in X} f(x)}) may be considered initially as a linear
mapping from $c_{00}(X, V)$ into $V$.  In this case, (\ref{N(sum_{x in
    X} f(x)) le ||f||_{BFS}}) follows directly from the definition of
the BFS norm, so that (\ref{f mapsto sum_{x in X} f(x)}) is a bounded
linear mapping with respect to the BFS norm on $c_{00}(X, V)$.  If $V$
is complete, then this mapping has a unique extension to a bounded
linear mapping from the closure of $c_{00}(X, V)$ in $BFS(X, V)$ into
$V$, by the remarks at the end of Section \ref{bounded linear
  mappings}.  We have also seen that the closure of $c_{00}(X, V)$ in
$BFS(X, V)$ is the same as $GCC(X, V)$.  This gives another way to
look at (\ref{f mapsto sum_{x in X} f(x)}) as a bounded linear mapping
from $GCC(X, V)$ into $V$, with respect to the BFS norm on $GCC(X,
V)$.

\section{Sums of sums}
\label{sums of sums}

        Let $V$ be a vector space over $k$ equipped with a norm $N$,
and let us suppose throughout this section that $V$ is complete with
respect to the associated metric.  Also let $X$ be a nonempty set,
and let $f$ be a $V$-valued function on $X$ such that
\begin{equation}
\label{sum_{x in X} f(x), 4}
        \sum_{x \in X} f(x)
\end{equation}
satisfies the generalized Cauchy criterion, as in Section
\ref{generalized convergence}.  If $E$ is any subset of $X$, then it
follows that
\begin{equation}
\label{sum_{x in E} f(x), 2}
        \sum_{x \in E} f(x)
\end{equation}
satisfies the generalized Cauchy criterion too, as mentioned in
Section \ref{generalized convergence, continued}.  This implies that
the net of all finite subsums of (\ref{sum_{x in E} f(x), 2})
converges in $V$, because $V$ is complete, as in Section
\ref{generalized convergence, continued} again.  We also have that
\begin{equation}
\label{N(sum_{x in E} f(x)) le ||f||_{BFS(E, V)} le ||f||_{BFS(X, V)}}
 N\Big(\sum_{x \in E} f(x)\Big) \le \|f\|_{BFS(E, V)} \le \|f\|_{BFS(X, V)}
\end{equation}
for every $E \subseteq X$, as in (\ref{N(sum_{x in X} f(x)) le
  ||f||_{BFS}}) in the previous section.  Here $\|f\|_{BFS(X, V)}$ is
the usual BFS norm of $f$ on $X$, as in (\ref{||f||_{BFS} = sup
  {N(sum_{x in A} f(x)) : A is a finite subset of X}}), and
$\|f\|_{BFS(E, V)}$ refers to the BFS norm of the restriction of $f$
to $E$.  Of course, if $f$ is summable on $X$, then the restriction of
$f$ to any set $E \subseteq X$ is summable, and
\begin{equation}
\label{N(sum_{x in E} f(x)) le sum_{x in E} N(f(x)) le sum_{x in X} N(f(x))}
        N\Big(\sum_{x \in E} f(x)\Big) \le \sum_{x \in E} N(f(x))
                                       \le \sum_{x \in X} N(f(x)),
\end{equation}
as in (\ref{N(sum_{x in X} f(x)) le sum_{x in X} N(f(x))}) in Section
\ref{generalized convergence, continued}.  If $N$ is an ultranorm
on $V$, then it suffices to ask that $f$ vanish at infinity on $X$,
as in Section \ref{generalized convergence}, which implies that
the restriction of $f$ to any set $E \subseteq X$ vanishes at
infinity on $E$.  In this case, we have that
\begin{equation}
\label{N(sum_{x in E} f(x)) le max_{x in E} N(f(x)) le max_{x in X} N(f(x))}
        N\Big(\sum_{x \in E} f(x)\Big) \le \max_{x \in E} N(f(x))
                                       \le \max_{x \in X} N(f(x))
\end{equation}
for each $E \subseteq X$, as in (\ref{N(sum_{x in X} f(x)) le max_{x
    in X} N(f(x))}) in Section \ref{generalized convergence,
  continued}.

        Let $I$ be a nonempty set, and let $\{E_j\}_{j \in I}$ be a
family of pairwise-disjoint subsets of $X$.  Thus
\begin{equation}
\label{a(j) = sum_{x in E_j} f(x)}
        a(j) = \sum_{x \in E_j} f(x)
\end{equation}
is defined as an element of $V$ for each $j \in I$, as in the
preceding paragraph.  If $j_1, \ldots, j_n$ are finitely many distinct
elements of $I$, then
\begin{equation}
\label{sum_{l = 1}^n a(j_l) = sum_{x in bigcup_{l = 1}^n E_{j_l}} f(x)}
        \sum_{l = 1}^n a(j_l) = \sum_{x \in \bigcup_{l = 1}^n E_{j_l}} f(x),
\end{equation}
as in (\ref{sum_{x in E_1 cup E_2} f(x) = sum_{x in E_1} f(x) + sum_{x
    in E_2} f(x)}) in Section \ref{generalized convergence,
  continued}.  It follows that
\begin{equation}
\label{N(sum_{l = 1}^n a(j_l)) le ... le ||f||_{BFS(X, V)}}
        N\Big(\sum_{l = 1}^n a(j_l)\Big)
 \le \|f\|_{BFS\big(\bigcup_{l = 1}^n E_{j_l}, V\big)} \le \|f\|_{BFS(X, V)},
\end{equation}
as in (\ref{N(sum_{x in E} f(x)) le ||f||_{BFS(E, V)} le ||f||_{BFS(X,
    V)}}).  This implies that $a$ has bounded finite sums on $I$, with
\begin{equation}
\label{||a||_{BFS(I, V)} le ... le ||f||_{BFS(X, V)}}
        \|a\|_{BFS(I, V)} \le \|f\|_{BFS\big(\bigcup_{j \in I} E_j, V\big)}
                         \le \|f\|_{BFS(X, V)}.
\end{equation}
One can also check that
\begin{equation}
\label{sum_{j in I} a(j)}
        \sum_{j \in I} a(j)
\end{equation}
satisfies the generalized Cauchy criterion under these conditions,
using (\ref{N(sum_{l = 1}^n a(j_l)) le ... le ||f||_{BFS(X, V)}}) and
the analogous property of $f$.  More precisely, one can verify that
\begin{equation}
\label{sum_{j in I} a(j) = sum_{x in bigcup_{j in I} E_j} f(x)}
        \sum_{j \in I} a(j) = \sum_{x \in \bigcup_{j \in I} E_j} f(x),
\end{equation}
by considering approximations of the various sums by finite subsums,
and where the right side of (\ref{sum_{j in I} a(j) = sum_{x in
    bigcup_{j in I} E_j} f(x)}) is defined as in the previous
paragraph.  If $f$ is summable on $X$, then it is easy to see that $a$
is summable on $I$, with
\begin{equation}
\label{||a||_{ell^1(I, V)} le ... le ||f||_{ell^1(X, V)}}
        \|a\|_{\ell^1(I, V)} \le \sum_{j \in I} \sum_{x \in E_j} N(f(x))
                     \le \|f\|_{\ell^1(X, V)}.
\end{equation}
Similarly, if $N$ is an ultranorm on $V$, and $f$ vanishes at infinity
on $X$, then one can check directly that $a$ vanishes at infinity
on $I$, and that
\begin{equation}
\label{||a||_{ell^infty(I, V)} le ... le ||f||_{ell^infty(X, V)}}
        \|a\|_{\ell^\infty(I, V)} \le \max_{j \in I} \max_{x \in E_j} N(f(x))
                                \le \|f\|_{\ell^\infty(X, V)}.
\end{equation}

        Note that (\ref{a(j) = sum_{x in E_j} f(x)}) defines a linear
mapping from $f \in GCC(X, V)$ into the vector space of $V$-valued
functions $a$ on $I$.  More precisely, this is a bounded linear
mapping from $GCC(X, V)$ into $BFS(I, V)$ with respect to the $BFS$
norms on $X$ and $I$, by (\ref{||a||_{BFS(I, V)} le ... le
  ||f||_{BFS(X, V)}}).  If $f$ has finite support in $X$, then $a$ has
finite support in $I$, because the $E_j$'s are pairwise disjoint.
Using this, it is easy to see that this mapping sends $GCC(X, V)$ into
$GCC(I, V)$, since $GCC(X, V)$ is the closure of $c_{00}(X, V)$ in
$BFS(X, V)$, and similarly for $GCC(I, V)$.  Clearly (\ref{sum_{j in
    I} a(j) = sum_{x in bigcup_{j in I} E_j} f(x)}) holds when $f$ has
finite support in $X$, which implies that (\ref{sum_{j in I} a(j) =
  sum_{x in bigcup_{j in I} E_j} f(x)}) holds for every $f \in GCC(X,
V)$, since $c_{00}(X, V)$ is dense in $GCC(X, V)$ with respect to the
BFS norm.  One can also look at (\ref{a(j) = sum_{x in E_j} f(x)}) as
defining a bounded linear mapping from $\ell^1(X, V)$ into $\ell^1(I,
V)$, by (\ref{||a||_{ell^1(I, V)} le ... le ||f||_{ell^1(X, V)}}).  If
$N$ is an ultranorm on $V$, then (\ref{a(j) = sum_{x in E_j} f(x)})
can be used to initially define a bounded linear mapping from $c_0(X,
V)$ into $\ell^\infty(I, V)$ with respect to the $\ell^\infty$ norms
on $X$ and $I$, by (\ref{||a||_{ell^infty(I, V)} le ... le
  ||f||_{ell^infty(X, V)}}).  As before, this mapping sends $c_{00}(X,
V)$ into $c_{00}(I, V)$, which implies that it sends $c_0(X, V)$ into
$c_0(I, V)$ in this case.

        Suppose now that $X = Y \times Z$ is the Cartesian product
of two nonempty sets $Y$ and $Z$.  If $f$ is a $V$-valued function on
$X$ such that (\ref{sum_{x in X} f(x), 4}) satisfies the generalized
Cauchy criterion, as before, then
\begin{equation}
\label{sum_{y in Y} f(y, z), 2}
        \sum_{y \in Y} f(y, z)
\end{equation}
satisfies the generalized Cauchy criterion for each $z \in Z$, and
\begin{equation}
\label{sum_{z in Z} f(y, z), 2}
        \sum_{z \in Z} f(y, z)
\end{equation}
satisfies the generalized Cauchy criterion for each $y \in Y$.  Here
$f(y, z)$ refers to the value of $f$ at $x = (y, z) \in Y \times Z$
for each $y \in Y$ and $z \in Z$, so that (\ref{sum_{y in Y} f(y, z), 2})
and (\ref{sum_{z in Z} f(y, z), 2}) are simply sums of $f$ over subsets
of $Y \times Z$.  Thus the sum (\ref{sum_{y in Y} f(y, z), 2}) is defined
as an element of $V$ for each $z \in Z$, because $V$ is complete, and
similarly (\ref{sum_{z in Z} f(y, z), 2}) is defined as an element of $V$
for each $y \in Y$.  We also have that
\begin{equation}
\label{sum_{z in Z} (sum_{y in Y} f(y, z)), 2}
        \sum_{z \in Z} \Big(\sum_{y \in Y} f(y, z)\Big)
\end{equation}
satisfies the generalized Cauchy criterion as a sum over $z \in Z$, and that
\begin{equation}
\label{sum_{y in Y} (sum_{z in Z} f(y, z)), 2}
        \sum_{y \in Y} \Big(\sum_{z \in Z} f(y, z)\Big)
\end{equation}
satisfies the generalized Cauchy criterion as a sum over $y \in Y$.
Both of these statements may be considered as instances of the
analogous statement for (\ref{sum_{j in I} a(j)}) discussed earlier.
Using (\ref{sum_{j in I} a(j) = sum_{x in bigcup_{j in I} E_j} f(x)}),
we get that (\ref{sum_{z in Z} (sum_{y in Y} f(y, z)), 2}) and
(\ref{sum_{y in Y} (sum_{z in Z} f(y, z)), 2}) are both equal to
\begin{equation}
\label{sum_{(y, z) in Y times Z} f(y, z), 2}
        \sum_{(y, z) \in Y \times Z} f(y, z).
\end{equation}
If $f$ is summable on $X$, then all of these sums are sums of summable
functions on the corresponding sets, as before.  If $N$ is an
ultranorm on $X$ and $f$ vanishes at infinity on $X$, then one can
check directly that these sums are sums of functions that vanish at
infinity on the corresponding sets.

\section{Finite-dimensional vector spaces}
\label{finite-dimensional vector spaces}

        Let $n$ be a positive integer, so that $k^n$ may be considered
as a vector space over $k$, as in Section \ref{norms, ultranorms}.
Also let $W$ be a vector space over $k$, and let $T$ be a linear
mapping from $k^n$ into $W$.  The standard basis vectors
$e(1), \ldots, e(n)$ in $k^n$ can be defined in the usual way, so that
the $l$th coordinate of $e(j)$ is equal to $1$ when $j = l$, and to
$0$ otherwise.  Thus each $v = (v_1, \ldots, v_n) \in k^n$ may be
expressed as
\begin{equation}
\label{v = sum_{j = 1}^n v_j e(j)}
        v = \sum_{j = 1}^n v_j \, e(j),
\end{equation}
which implies that
\begin{equation}
\label{T(v) = sum_{j = 1}^n v_j T(e(j))}
        T(v) = \sum_{j = 1}^n v_j \, T(e(j)).
\end{equation}
Let $N_0$ be the norm on $k^n$ defined in (\ref{N_0(v) = max(|v_1|,
  ldots, |v_n|)}), and let $N_W$ be a norm on $W$.  Observe that
\begin{eqnarray}
\label{N_W(T(v)) le ... le (sum_{j = 1}^n N_W(e(j))) N_0(v)}
 N_W(T(v)) \le \sum_{j = 1}^n N_W(v_j \, e(j))
           & = & \sum_{j = 1}^n |v_j| \, N_W(e(j)) \\
            & \le & \Big(\sum_{j = 1}^n N_W(e(j))\Big) \, N_0(v) \nonumber
\end{eqnarray}
for every $v \in k^n$, by (\ref{T(v) = sum_{j = 1}^n v_j T(e(j))}).
This implies that $T$ is a bounded linear mapping from $k^n$ into $W$,
with
\begin{equation}
\label{||T||_{op} le sum_{j = 1}^n N_W(e(j))}
        \|T\|_{op} \le \sum_{j = 1}^n N_W(e(j)),
\end{equation}
as in Section \ref{bounded linear mappings}.  If $N_W$ is an ultranorm
on $W$, then we get that
\begin{eqnarray}
\label{N_W(T(v)) le ... le (max_{1 le j le n} N_W(e(j))) N_0(v)}
 N_W(T(v)) \le \max_{1 \le j \le n} N_W(v_j \, e(j))
           & = & \max_{1 \le j \le n} (|v_j| \, N_W(e(j))  \\
           & \le & \Big(\max_{1 \le j \le n} N_W(e(j))\Big) \, N_0(v) \nonumber
\end{eqnarray}
for every $v \in k^n$, and hence
\begin{equation}
\label{||T||_{op} le max_{1 le j le n} N_W(e(j))}
        \|T\|_{op} \le \max_{1 \le j \le n} N_W(e(j)).
\end{equation}
More precisely,
\begin{equation}
\label{||T||_{op} = max_{1 le j le n} N_W(e(j))}
        \|T\|_{op} = \max_{1 \le j \le n} N_W(e(j))
\end{equation}
in this case, because equality holds in (\ref{N_W(T(v)) le ... le
  (max_{1 le j le n} N_W(e(j))) N_0(v)}) with $v = e(l)$ for some $l$.
Note that $N_0$ satisfies the condition indicated in (\ref{N_V(v) =
  |t|}) in Section \ref{bounded linear mappings}, which implies that
the operator norm of $T$ can be given as in (\ref{sup {N_W(T(v)) : v
    in V, N_V(v) le 1}}) in that section.

        If $N$ is any norm on $k^n$, then there is a positive real
number $C_1$ such that
\begin{equation}
\label{N(v) le C_1 N_0(v)}
        N(v) \le C_1 \, N_0(v)
\end{equation}
for every $v \in k^n$.  This follows from (\ref{N_W(T(v)) le ... le
  (sum_{j = 1}^n N_W(e(j))) N_0(v)}) applied to $W = k^n$, $N_W = N$,
and $T$ equal to the identity mapping on $k^n$.  If $k$ is complete
with respect to the metric associated to $|\cdot|$, then one can show
that there is also a positive real number $C_2$ such that
\begin{equation}
\label{N_0(v) le C_2 N(v)}
        N_0(v) \le C_2 \, N(v)
\end{equation}
for every $v \in k^n$.  See Lemma 2.1 on p116 of \cite{c}, or Theorem
5.2.1 on p137 of \cite{fg}.  This implies that the topology on $k^n$
determined by the metric associated to $N$ is the same as the topology
determined by the metric associated to $N_0$.

        As in Section \ref{norms, ultranorms}, the topology on $k^n$
determined by the metric associated to $N_0$ is the same as the
product topology corresponding to the topology on $k$ determined by
the metric associated to the absolute value function.  If $k$ is
locally compact, then $k^n$ is locally compact with respect to this
topology as well.  Of course, if $|\cdot|$ is the trivial absolute
value function on $k$, then $N_0$ is the trivial ultranorm on $k^n$,
and the corresponding topologies are discrete.  Suppose for the moment
that $|\cdot|$ is not the trivial absolute value function on $k$, and
that $k$ is locally compact with respect to the topology determined by
the metric associated to $|\cdot|$.  In this case, it is easy to see
that all closed balls in $k$ are compact, so that closed and bounded
subsets of $k$ are compact.  It follows that closed balls in $k^n$
with respect to the metric associated to $N_0$ are compact too, by
Tychonoff's theorem.  This implies that closed and bounded subsets of
$k^n$ are compact too.  Note that $k$ is complete when $k$ is locally
compact, as one can check using the fact that compact metric spaces
are complete.

        It is a bit simpler to show (\ref{N_0(v) le C_2 N(v)}) when
$|\cdot|$ is nontrivial on $k$ and $k$ is locally compact, so that
\begin{equation}
\label{{v in k^n : N_0(v) = 1}}
        \{v \in k^n : N_0(v) = 1\}
\end{equation}
is a compact subset of $k^n$.  This also uses the fact that $N$ is
continuous as a real-valued function on $k^n$ with respect to the
metric associated to $N_0$, which can be derived from (\ref{N(v) le
  C_1 N_0(v)}).  If $k$ is only asked to be complete with respect to
the metric associated to $|\cdot|$, then one can use induction on $n$
to prove (\ref{N_0(v) le C_2 N(v)}).  The base case $n = 1$ is easy,
and when $n \ge 2$ the induction hypothesis implies that a condition
like (\ref{N_0(v) le C_2 N(v)}) holds on $k^{n - 1} \times \{0\}$.  In
particular, it follows that $k^{n - 1} \times \{0\}$ is complete with
respect to the metric associated to $N$, because $k$ is complete, by
hypothesis.  This implies that $k^{n - 1} \times \{0\}$ is a closed
subset of $k^n$ with respect to the metric associated to $N$, by
standard arguments.  If $e(n)$ is the $n$th standard basis vector in
$k^n$, as before, then it follows that there is a positive lower bound
for the distances between $e(n)$ and elements of $k^{n - 1} \times
\{0\}$ with respect to the metric associated to $N$.  Equivalently,
this means that $|v_n|$ is bounded by a constant multiple of $N(v)$
for each $v \in k^n$.  This permits a condition like (\ref{N_0(v) le
  C_2 N(v)}) to be obtained on $k^n$ from an analogous condition on
$k^{n - 1} \times \{0\}$.

\section{$q$-Norms}
\label{q-norms}

        Let $V$ be a vector space over $k$, and let $q$ be a positive
real number.  Also let $N$ be a nonnegative real-valued function on
$V$ that satisfies the same positivity and homogeneity conditions as
for a norm, as in (\ref{N(v) = 0 if and only if v = 0}) and (\ref{N(t
  v) = |t| N(v)}) in Section \ref{norms, ultranorms}.  Let us say that
$N$ is a \emph{$q$-norm}\index{qnorms@$q$-norms} on $V$ if
\begin{equation}
\label{N(v + w)^q le N(v)^q + N(w)^q}
        N(v + w)^q \le N(v)^q + N(w)^q
\end{equation}
for every $v, w \in V$.  Of course, (\ref{N(v + w)^q le N(v)^q +
  N(w)^q}) is the same as the usual triangle inequality (\ref{N(v + w)
  le N(v) + N(w)}) when $q = 1$, so that a $1$-norm is the same as a
norm.  If $N$ is an ultranorm on $V$, then
\begin{equation}
\label{N(v + w)^q le max(N(v), N(w))^q le N(v)^q + N(w)^q}
        N(v + w)^q \le \max(N(v), N(w))^q \le N(v)^q + N(w)^q
\end{equation}
for every $v, w \in V$ and $q > 0$, so that $N$ is a $q$-norm on $V$
for every $q > 0$.

        Note that (\ref{N(v + w)^q le N(v)^q + N(w)^q}) is equivalent
to asking that
\begin{equation}
\label{N(v + w) le (N(v)^q + N(w)^q)^{1/q}}
        N(v + w) \le (N(v)^q + N(w)^q)^{1/q}
\end{equation}
for every $v, w \in V$.  Clearly
\begin{equation}
\label{max(N(v), N(w)) le (N(v)^q + N(w)^q)^{1/q}}
        \max(N(v), N(w)) \le (N(v)^q + N(w)^q)^{1/q}
\end{equation}
for every $v, w \in V$ and $q > 0$, which is the same as the second
step in (\ref{N(v + w)^q le max(N(v), N(w))^q le N(v)^q + N(w)^q}).
We also have that
\begin{equation}
\label{N(v)^q + N(w)^q le 2 max(N(v)^q, N(w)^q)}
        N(v)^q + N(w)^q \le 2 \, \max(N(v)^q, N(w)^q)
\end{equation}
for every $v, w \in V$ and $q > 0$, so that
\begin{eqnarray}
\label{(N(v)^q + N(w)^q)^{1/q} le 2^{1/q} max(N(v), N(w))}
        (N(v)^q + N(w)^q)^{1/q} \le 2^{1/q} \, \max(N(v), N(w)).
\end{eqnarray}
It follows from (\ref{max(N(v), N(w)) le (N(v)^q + N(w)^q)^{1/q}}) and
(\ref{(N(v)^q + N(w)^q)^{1/q} le 2^{1/q} max(N(v), N(w))}) that
\begin{equation}
\label{lim_{q to infty} (N(v)^q + N(w)^q)^{1/q} = max(N(v), N(w))}
        \lim_{q \to \infty} (N(v)^q + N(w)^q)^{1/q} = \max(N(v), N(w))
\end{equation}
for every $v, w \in V$, since $2^{1/q} \to 1$ as $q \to \infty$.  Thus
one might interpret (\ref{N(v + w) le (N(v)^q + N(w)^q)^{1/q}}) as
being the ultrametric version (\ref{N(v + w) le max(N(v), N(w))}) of
the triangle inequality when $q = \infty$.

        If $0 < q \le r < \infty$, then
\begin{equation}
\label{N(v)^r + N(w)^r le max(N(v), N(w))^{r - q} (N(v)^q + N(w)^q)}
        N(v)^r + N(w)^r \le \max(N(v), N(w))^{r - q} \, (N(v)^q + N(w)^q)
\end{equation}
for every $v, w \in V$.  This implies that
\begin{equation}
\label{N(v)^r + N(w)^r le ... = (N(v)^q + N(w)^q)^{r/q}}
\qquad  N(v)^r + N(w)^r \le (N(v)^q + N(w)^q)^{(r - q)/q + 1}
                             = (N(v)^q + N(w)^q)^{r/q},
\end{equation}
for every $v, w \in V$, using (\ref{max(N(v), N(w)) le (N(v)^q +
  N(w)^q)^{1/q}}) in the first step.  Thus
\begin{equation}
\label{(N(v)^r + N(w)^r)^{1/r} le (N(v)^q + N(w)^q)^{1/q}}
        (N(v)^r + N(w)^r)^{1/r} \le (N(v)^q + N(w)^q)^{1/q}
\end{equation}
for every $v, w \in V$ when $q \le r$, by taking the $r$th root of
both sides of (\ref{N(v)^r + N(w)^r le ... = (N(v)^q +
  N(w)^q)^{r/q}}).  This could also be derived from (\ref{(r + t)^a le
  r^a + t^a}) or (\ref{r + t le (r^a + t^a)^{1/a}}) in Section
\ref{quasimetrics}, or from (\ref{||f||_r le ||f||_q}) in Section
\ref{ell^r norms}.  If $N$ is an $r$-norm on $V$, then it follows that
$N$ is a $q$-norm on $V$ too when $q \le r$.

        Suppose for the moment that
\begin{equation}
\label{|t + t'|^q le |t|^q + |t'|^q}
        |t + t'|^q \le |t|^q + |t'|^q
\end{equation}
for some $q > 0$ and every $t, t' \in k$, so that $|t|^q$ is also an
absolute value function on $k$.  If $N$ is a $q$-norm on $V$, then
\begin{equation}
\label{N(v)^q}
        N(v)^q
\end{equation}
may be considered as a norm on $V$ with respect to $|t|^q$ as an
absolute value function on $k$.  More precisely, the homogeneity
condition (\ref{N(t v) = |t| N(v)}) for $N$ with respect to $|t|$ on
$k$ implies that $N(v)^q$ satisfies the analogous condition with
respect to $|t|^q$ on $k$.  Similarly, (\ref{N(v + w)^q le N(v)^q +
  N(w)^q}) is the same as the standard triangle inequality for
$N(v)^q$ as a norm on $V$ with respect to $|t|^q$.

        Suppose now that $V \ne \{0\}$, and that $N$ is a $q$-norm
on $V$ for some $q > 0$.  Let $u$ be a nonzero element of $V$, and let
$t$, $t'$ be arbitrary elements of $k$.  If we apply (\ref{N(v + w)^q
  le N(v)^q + N(w)^q}) to $v = t \, u$ and $w = t' \, u$, then we get
that
\begin{eqnarray}
\label{|t + t'|^q N(u)^q = ... = |t|^q N(u)^q + |t'|^q N(u)^q}
        |t + t'|^q \, N(u)^q = N(t \, u + t' \, u)^q
                          & \le &  N(t \, u)^q + N(t' \, u)^q \\
                      & = & |t|^q \, N(u)^q + |t'|^q \, N(u)^q, \nonumber
\end{eqnarray}
using also the homogeneity property (\ref{N(t v) = |t| N(v)}) of $N$
with respect to $|\cdot|$.  This shows that (\ref{|t + t'|^q le |t|^q
  + |t'|^q}) holds under these conditions, since $N(u) > 0$.

        If $N$ is a norm on $V$, then
\begin{equation}
\label{N(v - w)}
        N(v - w)
\end{equation}
defines a metric on $V$, as in (\ref{d(v, w) = N(v - w)}) in Section
\ref{norms, ultranorms}.  Similarly, if $N$ is a $q$-norm on $V$, then
\begin{equation}
\label{N(v - w)^q}
        N(v - w)^q
\end{equation}
defines a metric on $V$, which is the same as (\ref{N(v - w)}) when $q
= 1$.  If $N$ is a $q$-norm on $V$ and $|t|^q$ is an absolute value
function on $k$, then $N(v)^q$ may be considered as a norm on $V$ with
respect to $|t|^q$ on $k$, as before, and (\ref{N(v - w)^q}) is the
same as the metric associated to this norm.  Of course, if $N$ is a
$q$-norm on $V$ and $q \ge 1$, then $N$ is a norm on $V$ too, so that
(\ref{N(v - w)}) is a metric on $V$ as well, which determines the same
topology on $V$ as (\ref{N(v - w)^q}).  If $q < 1$, then (\ref{N(v -
  w)}) is at least a quasimetric on $V$, as in Section
\ref{quasimetrics}.

        Suppose that $N$ is a $q$-norm on $V \ne \{0\}$ for some $q > 0$,
and let $X$ be a nonempty set.  Also let $f$ be a $V$-valued function
on $X$ such that
\begin{equation}
\label{N(f(x))^q}
        N(f(x))^q
\end{equation}
is summable as a nonnegative real-valued function on $X$, as in
Section \ref{nonnegative sums}.  As before, $|t|^q$ is an absolute
value function on $k$ under these conditions, and $N(v)^q$ may be
considered as a norm on $V$ with respect to $|t|^q$ on $k$.  Thus the
summability of (\ref{N(f(x))^q}) on $X$ is the same as the summability
of $f$ as a $V$-valued function on $X$ with respect to $N(v)^q$ as a
norm on $V$ with respect to $|t|^q$ on $k$.  This implies that
$\sum_{x \in X} f(x)$ satisfies the generalized Cauchy condition with
respect to $N(v)^q$ as a norm on $V$ with respect to $|t|^q$ on $k$,
as in Section \ref{generalized convergence}.

\section{$\ell^r$ Norms, continued}
\label{ell^r norms, continued}

        Let $V \ne \{0\}$ be a vector space over $k$ again, and let
$N$ be a $q$-norm on $V$ for some $q > 0$.  It follows that $|t|^q$
is an absolute value function on $k$ too, as in the previous section,
and that $N(v)^q$ may be considered as a norm on $V$ with respect to
$|t|^q$ on $k$.  Also let $X$ be a nonempty set, and let $r$ be a
positive real number.  As in Section \ref{ell^r norms}, a $V$-valued
function $f$ on $X$ is said to be
\emph{$r$-summable}\index{rsummable functions@$r$-summable functions}
with respect to $N$ on $V$ if
\begin{equation}
\label{N(f(x))^r, 2}
        N(f(x))^r
\end{equation}
is summable as a nonnegative real-valued function on $X$.  Let us
denote the space of $r$-summable $V$-valued functions on $X$ by
$\ell^r(X, V)$,\index{l^r(X, V)@$\ell^r(X, V)$} as before, or by
$\ell^r_N(X, V)$,\index{l^r_N(X, V)@$\ell^r_N(X, V)$} to indicate the
role of $N$.

        Of course,
\begin{equation}
\label{N(f(x))^r = (N(f(x))^q)^{r/q}}
        N(f(x))^r = (N(f(x))^q)^{r/q}
\end{equation}
for every $x \in X$, so that $f$ is $r$-summable with respect to $N$
on $V$ if and only if $f$ is $(r/q)$-summable with respect to $N(v)^q$
as a norm on $V$ with respect to $|t|^q$ on $k$.  Thus
\begin{equation}
\label{ell^r_N(X, V) = ell^{r/q}_{N^q}(X, V)}
        \ell^r_N(X, V) = \ell^{r/q}_{N^q}(X, V),
\end{equation}
where $N^q$ is considered as a norm on $V$ with respect to $|t|^q$ on
$k$ on the right side of (\ref{ell^r_N(X, V) = ell^{r/q}_{N^q}(X,
  V)}).  The discussion in Section \ref{ell^r norms} implies that the
right side of (\ref{ell^r_N(X, V) = ell^{r/q}_{N^q}(X, V)}) is a
vector space with respect to pointwise addition and scalar
multiplication of $V$-valued functions on $X$, so that the same
conclusion holds for the left side of (\ref{ell^r_N(X, V) =
  ell^{r/q}_{N^q}(X, V)}).

        Put
\begin{equation}
\label{||f||_r = ... = (sum_{x in X} N(f(x))^r)^{1/r}}
        \|f\|_r = \|f\|_{\ell^r(X, V)} = \|f\|_{\ell^r_N(X, V)}
                     = \Big(\sum_{x \in X} N(f(x))^r\Big)^{1/r}
\end{equation}
for each $f \in \ell^r_N(X, V)$, as in (\ref{||f||_r = ||f||_{ell^r(X,
    V)} = (sum_{x in X} N(f(x))^r)^{1/r}}).  It is easy to see that
this satisfies the usual positivity and homogeneity conditions for a
norm, because of the corresponding properties of $N$.  Note that
\begin{equation}
\label{||f||_{ell^{r/q}_{N^q}(X, V)} = ... = (sum_{x in X} N(f(x))^r)^{q/r}}
 \|f\|_{\ell^{r/q}_{N^q}(X, V)} = \Big(\sum_{x \in X} (N(f(x))^q)^{r/q}\Big)^{q/r}
                           = \Big(\sum_{x \in X} N(f(x))^r\Big)^{q/r}
\end{equation}
for every $f$ in (\ref{ell^r_N(X, V) = ell^{r/q}_{N^q}(X, V)}), so that
\begin{equation}
\label{||f||_{ell^{r/q}_{N^q}(X, V)} = (||f||_{ell^r_N(X, V)})^q}
        \|f\|_{\ell^{r/q}_{N^q}(X, V)} = \big(\|f\|_{\ell^r_N(X, V)}\big)^q.
\end{equation}

        Suppose for the moment that $q = 1$, so that $N$ is a norm
on $V$.  If $r \ge 1$, then (\ref{||f||_r = ... = (sum_{x in X}
  N(f(x))^r)^{1/r}}) defines a norm on $\ell^r_N(X, V)$, by (\ref{||f
  + g||_r le ||f||_r + ||g||_r}).  Similarly, if $0 < r \le 1$, then
(\ref{||f||_r = ... = (sum_{x in X} N(f(x))^r)^{1/r}}) defines an
$r$-norm on $V$, by (\ref{||f + g||_r^r le ||f||_r^r + ||g||_r^r}).
If $N$ is an ultranorm on $V$, then (\ref{||f + g||_r^r le ||f||_r^r +
  ||g||_r^r}) holds for every $r > 0$, as mentioned in Section
\ref{ell^r norms}.  This implies that (\ref{||f||_r = ... = (sum_{x in
    X} N(f(x))^r)^{1/r}}) is an $r$-norm on $\ell^r_N(X, V)$ for
every $r > 0$ in this case.

        Now let $N$ be a $q$-norm on $V$ for some $q > 0$ again,
so that $N(v)^q$ is a norm on $V$ with respect to $|t|^q$ on $k$.  It
follows that (\ref{||f||_{ell^{r/q}_{N^q}(X, V)} = ... = (sum_{x in X}
  N(f(x))^r)^{q/r}}) is a norm on $\ell^{r/q}_{N^q}(X, V)$ when $r/q
\ge 1$, and that (\ref{||f||_{ell^{r/q}_{N^q}(X, V)} = ... = (sum_{x
    in X} N(f(x))^r)^{q/r}}) is an $(r/q)$-norm on
$\ell^{r/q}_{N^q}(X, V)$ when $r/q \le 1$, as in the preceding
paragraph.  More precisely, we still use $|t|^q$ as the absolute value
function on $k$ for these two statements, but the main point is the
version of the triangle inequality that we get.  If $r/q \ge 1$,
then we have that
\begin{equation}
\label{||f + g||_{ell^{r/q}_{N^q}(X, V)} le ...}
 \|f + g\|_{\ell^{r/q}_{N^q}(X, V)} \le \|f\|_{\ell^{r/q}_{N^q}(X, V)}
                                         + \|g\|_{\ell^{r/q}_{N^q}(X, V)}
\end{equation}
for every $f, g \in \ell^{r/q}_{N^q}(X, V)$.  If $r/q \le 1$, then
\begin{equation}
\label{||f + g||_{ell^{r/q}_{N^q}(X, V)}^{r/q} le ...}
 \|f + g\|_{\ell^{r/q}_{N^q}(X, V)}^{r/q} \le \|f\|_{\ell^{r/q}_{N^q}(X, V)}^{r/q}
                                              + \|g\|_{\ell^{r/q}_{N^q}(X, V)}^{r/q}
\end{equation}
for every $f, g \in \ell^{r/q}_{N^q}(X, V)$.

        These two statements can be reformulated in terms of
(\ref{||f||_r = ... = (sum_{x in X} N(f(x))^r)^{1/r}}), using
(\ref{ell^r_N(X, V) = ell^{r/q}_{N^q}(X, V)}) and
(\ref{||f||_{ell^{r/q}_{N^q}(X, V)} = (||f||_{ell^r_N(X, V)})^q}).
If $r \ge q$, then (\ref{||f + g||_{ell^{r/q}_{N^q}(X, V)} le ...})
implies that
\begin{equation}
\label{||f + g||_{ell^r_N(X, V)}^q le ...}
 \|f + g\|_{\ell^r_N(X, V)}^q \le \|f\|_{\ell^r_N(X, V)}^q + \|g\|_{\ell^r_N(X, V)}^q
\end{equation}
for every $f, g \in \ell^r_N(X, V)$, so that (\ref{||f||_r = ... =
  (sum_{x in X} N(f(x))^r)^{1/r}}) defines a $q$-norm on $\ell^r_N(X,
V)$.  If $r \le q$, then (\ref{||f + g||_{ell^{r/q}_{N^q}(X, V)}^{r/q}
  le ...}) implies that
\begin{equation}
\label{||f + g||_{ell^r_N(X, V)}^r le ...}
 \|f + g\|_{\ell^r_N(X, V)}^r \le \|f\|_{\ell^r_N(X, V)}^r + \|g\|_{\ell^r_N(X, V)}^r
\end{equation}
for every $f, g \in \ell^r_N(X, V)$, so that (\ref{||f||_r = ... =
  (sum_{x in X} N(f(x))^r)^{1/r}}) defines an $r$-norm on $\ell^r_N(X,
V)$.

        Let us take $\ell^\infty(X, V) = 
\ell^\infty_N(X, V)$\index{l^infinity(X, V)@$\ell^\infty(X, V)$}
to\index{l^infinity_N(X, V)@$\ell^\infty_N(X, V)$} be the space of
$V$-valued functions $f$ on $X$ that are bounded, in the sense that
\begin{equation}
\label{N(f(x))}
        N(f(x))
\end{equation}
is bounded as a nonnegative real-valued function on $X$.  Of course,
this is the same as saying that $N(f(x))^q$ is bounded on $X$, so that
(\ref{ell^r_N(X, V) = ell^{r/q}_{N^q}(X, V)}) also holds when $r =
\infty$.  In particular, $\ell^\infty_N(X, V)$ is a vector space with
respect to pointwise addition and scalar multiplication, as in Section
\ref{supremum norm}.  If we put
\begin{equation}
\label{||f||_infty = ... = sup_{x in X} N(f(x))}
        \|f\|_\infty = \|f\|_{\ell^\infty(X, V)} = \|f\|_{\ell^\infty_N(X, V)}
                                             = \sup_{x \in X} N(f(x))
\end{equation}
for every $f \in \ell^\infty_N(X, V)$, as in (\ref{||f||_infty =
  ||f||_{ell^infty(X, V)} = sup_{x in X} N(f(x))}), then
(\ref{||f||_{ell^{r/q}_{N^q}(X, V)} = (||f||_{ell^r_N(X, V)})^q})
holds when $r = \infty$ too.  It is easy to see that (\ref{||f||_infty
  = ... = sup_{x in X} N(f(x))}) defines a $q$-norm on
$\ell^\infty_N(X, V)$ under these conditions, directly from the
definitions, or using (\ref{||f||_{ell^{r/q}_{N^q}(X, V)} =
  (||f||_{ell^r_N(X, V)})^q}) with $r = \infty$ to reduce to the case
of norms.

\chapter{Additional examples and results}
\label{additional examples, results}

\section{Cauchy products}
\label{cauchy products}

        Let $\sum_{j = 0}^\infty a_j$ and $\sum_{l = 0}^\infty b_l$
be infinite series with terms in a field $k$, and put
\begin{equation}
\label{c_n = sum_{j = 0}^n a_j b_{n - j}}
        c_n = \sum_{j = 0}^n a_j \, b_{n - j}
\end{equation}
for each nonnegative integer $n$.  The infinite series $\sum_{n =
  0}^\infty c_n$ is known as the \emph{Cauchy product}\index{Cauchy
  products} of the series $\sum_{j = 0}^\infty a_j$ and $\sum_{l =
  0}^\infty b_l$, and it is easy to see that
\begin{equation}
\label{sum_{n = 0}^infty c_n = (sum_{j = 0}^infty a_j) (sum_{l = 0}^infty b_l)}
        \sum_{n = 0}^\infty c_n = \Big(\sum_{j = 0}^\infty a_j\Big)
                                 \, \Big(\sum_{l = 0}^\infty b_l\Big)
\end{equation}
formally.  In particular, this holds when $a_j = 0$ for all but
finitely many $j \ge 0$ and $b_l = 0$ for all but finitely many $l \ge
0$, in which case $c_n = 0$ for all but finitely many $n$.  We can
look at this in terms of the discussion in Section \ref{sums of sums},
with
\begin{equation}
\label{X = ({bf Z}_+ cup {0}) times ({bf Z}_+ cup {0})}
        X = ({\bf Z}_+ \cup \{0\}) \times ({\bf Z}_+ \cup \{0\}).
\end{equation}
Put
\begin{equation}
\label{E_n = {(j, l) in X : j + l = n}}
        E_n = \{(j, l) \in X : j + l = n\}
\end{equation}
for each nonnegative integer $n$, so that $E_n$ is a finite set with
exactly $n + 1$ elements for each $n \ge 0$, the $E_n$'s are pairwise
disjoint, and
\begin{equation}
\label{bigcup_{n = 0}^infty E_n}
        X = \bigcup_{n = 0}^\infty E_n.
\end{equation}
If $f \in c_{00}(X, k)$, then it follows that
\begin{equation}
\label{sum_{(j, l) in X} f(j, l) = ...}
 \sum_{(j, l) \in X} f(j, l) = \sum_{n = 0}^\infty
                               \Big(\sum_{(j, l) \in E_n} f(j, l)\Big).
\end{equation}
Let $f$ be the $k$-valued function on $X$ defined by
\begin{equation}
\label{f(j, l) = a_j b_l}
        f(j, l) = a_j \, b_l
\end{equation}
for each $j, l \ge 0$, so that
\begin{equation}
\label{sum_{(j, l) in E_n} f(j, l) = c_n}
        \sum_{(j, l) \in E_n} f(j, l) = c_n
\end{equation}
for every $n \ge 0$.  If $a_j = 0$ for all but finitely many $j$, and
$b_l = 0$ for all but finitely many $l$, then $f \in c_{00}(X, k)$,
and
\begin{equation}
\label{sum_{(j, l) in X} f(j, l) = ..., 2}
        \sum_{(j, l) \in X} f(j, l) = \Big(\sum_{j = 0}^\infty a_j\Big)
                                     \, \Big(\sum_{l = 0}^\infty b_l\Big).
\end{equation}
This corresponds to summing $f(j, l)$ over $j$ and $l$ separately.
Thus (\ref{sum_{n = 0}^infty c_n = (sum_{j = 0}^infty a_j) (sum_{l =
    0}^infty b_l)}) follows from (\ref{sum_{(j, l) in X} f(j, l) =
  ...}), (\ref{sum_{(j, l) in E_n} f(j, l) = c_n}), and (\ref{sum_{(j,
    l) in X} f(j, l) = ..., 2}) under these conditions.

        Suppose for the moment that $k = {\bf R}$, and the $a_j$'s
and $b_l$'s are nonnegative real numbers for each $j, l \ge 0$.  This
implies that the $c_n$'s are nonnegative real numbers for each $n \ge
0$, and that each of the three sums in (\ref{sum_{n = 0}^infty c_n =
  (sum_{j = 0}^infty a_j) (sum_{l = 0}^infty b_l)}) is defined as a
nonnegative extended real number.  In this case, it is well known and
not difficult to check that (\ref{sum_{n = 0}^infty c_n = (sum_{j =
    0}^infty a_j) (sum_{l = 0}^infty b_l)}) always holds, with
suitable interpretations when the sums are infinite.  More precisely,
the right side of (\ref{sum_{n = 0}^infty c_n = (sum_{j = 0}^infty
  a_j) (sum_{l = 0}^infty b_l)}) should be interpreted as being equal
to $0$ whenever one of the factors is equal to $0$, even if the other
factor is infinite, and otherwise the right side of (\ref{sum_{n =
    0}^infty c_n = (sum_{j = 0}^infty a_j) (sum_{l = 0}^infty b_l)})
should be interpreted as being infinite when one of the factors is
infinite and the other is positive.  This may be considered as a
consequence of the discussion in Section \ref{nonnegative sums,
  continued} for nonnegative real-valued functions, using the
interpretations just mentioned for the right side of (\ref{sum_{(j, l)
    in X} f(j, l) = ..., 2}).

        Suppose now that $k = {\bf R}$ or ${\bf C}$, with the standard
absolute value function.  Note that
\begin{equation}
\label{|c_n| le sum_{j = 0}^n |a_j| |b_{n - j}|}
        |c_n| \le \sum_{j = 0}^n |a_j| \, |b_{n - j}|
\end{equation}
for each $n \ge 0$, and hence that
\begin{equation}
\label{sum_{n = 0}^infty |c_n| le ...}
        \sum_{n = 0}^\infty |c_n|
 \le \sum_{n = 0}^\infty \Big(\sum_{j = 0}^n |a_j| \, |b_{n - j}|\Big)
   = \Big(\sum_{j = 0}^\infty |a_j|\Big) \, \Big(\sum_{l = 0}^\infty |b_l|\Big),
\end{equation}
with suitable interpretations for the right side of (\ref{sum_{n =
    0}^infty |c_n| le ...}), as in the preceding paragraph.  This
implies that $\sum_{n = 0}^\infty c_n$ converges absolutely when
$\sum_{j = 0}^\infty a_j$ and $\sum_{l = 0}^\infty b_l$ converge
absolutely, in which case one can check that (\ref{sum_{n = 0}^infty
  c_n = (sum_{j = 0}^infty a_j) (sum_{l = 0}^infty b_l)}) holds, by
approximating the various sums by finite sums.  This can also be seen
as a consequence of the discussion in Section \ref{sums of sums} for
summable functions, using the fact that (\ref{f(j, l) = a_j b_l}) is
summable on (\ref{X = ({bf Z}_+ cup {0}) times ({bf Z}_+ cup {0})}).
Alternatively, one can reduce to the analogous statement for
nonnegative real numbers, by expressing $\sum_{j = 0}^\infty a_j$ and
$\sum_{l = 0}^\infty b_l$ as linear combinations of convergent series
of nonnegative real numbers.

        If $k$ is any field with an ultrametric absolute value function
$|\cdot|$, then
\begin{equation}
\label{|c_n| le max_{0 le j le n} (|a_j| |b_{n - j}|)}
        |c_n| \le \max_{0 \le j \le n} (|a_j| \, |b_{n - j}|)
\end{equation}
for each $n \ge 0$.  Using this, it is easy to see that $\{c_n\}_{n =
  0}^\infty$ converges to $0$ in $k$ when $\{a_j\}_{j = 0}^\infty$ and
$\{b_l\}_{l = 0}^\infty$ converge to $0$ in $k$.  If $k$ is complete
with respect to the metric associated to $|\cdot|$, then it follows
that the corresponding infinite series converge in $k$.  In this
situation, one can check that (\ref{sum_{n = 0}^infty c_n = (sum_{j =
    0}^infty a_j) (sum_{l = 0}^infty b_l)}) holds, by approximating
the various sums by finite sums again.  As before, this can also be
derived from the discussion in Section \ref{sums of sums}, using the
fact that (\ref{f(j, l) = a_j b_l}) vanishes at infinity on (\ref{X =
  ({bf Z}_+ cup {0}) times ({bf Z}_+ cup {0})}) when $\{a_j\}_{j =
  0}^\infty$ and $\{b_l\}_{l = 0}^\infty$ converge to $0$ in $k$.

\section{Formal power series}
\label{formal power series}

        Let $k_0$ be a field, and let $T$ be an indeterminate.
By a \emph{formal power series}\index{formal power series} in $T$
with coefficients in $k_0$ we mean an expression of the form
\begin{equation}
\label{f(T) = sum_{j = 0}^infty f_j T^j}
        f(T) = \sum_{j = 0}^\infty f_j \, T^j,
\end{equation}
where $f_j \in k_0$ for each nonnegative integer $j$.  The collection
of all formal power series in $T$ with coefficients in $k_0$ is
denoted $k_0[[T]]$,\index{k_0[[t]]@$k_0[[T]]$} as usual.  More
precisely, the elements of $k_0[[T]]$ correspond to sequences
$\{f_j\}_{j = 0}^\infty$ of elements of $k_0$, or equivalently to
functions from the set ${\bf Z}_+ \cup \{0\}$ of nonnegative integers
into $k_0$.  Thus $k_0[[T]]$ may be defined as the collection of such
sequences, or equivalently as the space of $k_0$-valued functions on
${\bf Z}_+ \cup \{0\}$.  However, it is often more convenient to
represent elements of $k_0[[T]]$ as in (\ref{f(T) = sum_{j = 0}^infty
  f_j T^j}).  Note that $k_0[[T]]$ is a vector space over $k_0$ with
respect to termwise addition and scalar multiplication.

        Let $f(T)$ and $g(T)$ be elements of $k_0[[T]]$, where
$f(T)$ is as in (\ref{f(T) = sum_{j = 0}^infty f_j T^j}), and
similarly
\begin{equation}
\label{g(T) = sum_{l = 0}^infty g_l T^l}
        g(T) = \sum_{l = 0}^\infty g_l \, T^l
\end{equation}
for some $g_l \in k_0$.  The product of $f(T)$ and $g(T)$ can be defined
formally in the usual way, with
\begin{equation}
\label{T^j T^l = T^{j + l}}
        T^j \, T^l = T^{j + l}
\end{equation}
for all $j, l \ge 0$.  This means that
\begin{equation}
\label{f(T) g(T) = sum_{j = 0}^infty sum_{l = 0}^infty f_j g_l T^{j + l}}
 f(T) \, g(T) = \sum_{j = 0}^\infty \sum_{l = 0}^\infty f_j \, g_l \, T^{j + l},
\end{equation}
where there are only finitely many terms involving $T^n$ for each
nonnegative integer $n$.  Collecting these terms, we get that
\begin{equation}
\label{f(T) g(T) = sum_{n = 0}^infty (f g)_n T^n}
        f(T) \, g(T) = \sum_{n = 0}^\infty (f \, g)_n \, T^n,
\end{equation}
where
\begin{equation}
\label{(f g)_n = sum_{j = 0}^n f_j g_{n - j}}
        (f \, g)_n = \sum_{j = 0}^n f_j \, g_{n - j}
\end{equation}
for each $n \ge 0$.  Of course, this is the same as the Cauchy
product, discussed in the previous section.  More precisely, one can
use (\ref{f(T) g(T) = sum_{n = 0}^infty (f g)_n T^n}) and (\ref{(f
  g)_n = sum_{j = 0}^n f_j g_{n - j}}) as the official definition of
multiplication on $k_0[[T]]$, which makes sense directly at the level
of the corresponding sequences of coefficients in $k_0$.  One can also
check that this makes $k_0[[T]]$ into a commutative ring, and in fact
an algebra over $k_0$.

        A formal polynomial $f(T)$ in $T$ with coefficients in $k_0$
may be considered as a formal power series (\ref{f(T) = sum_{j =
    0}^infty f_j T^j}) such that $f_j = 0$ for all but finitely many
$j$.  Thus the collection $k_0[T]$\index{k_0[T]@$k_0[T]$} of formal
polynomials in $T$ may be considered as a subalgebra of $k_0[[T]]$.
Similarly, $k_0$ may be identified with the subalgebra of $k_0[T]$
consisting of power series (\ref{f(T) = sum_{j = 0}^infty f_j T^j})
such that $f_j = 0$ for every $j \ge 1$.  With this identification,
the multiplicative identity element $1$ in $k_0$ corresponds to $T^0$
in $k_0[T]$, which is the multiplicative identity element in
$k_0[[T]]$.

        If $f(T)$ is a nonzero formal power series in $T$ with
coefficients in $k_0$, then let $n(f(T))$ be the smallest nonnegative
integer $n$ such that
\begin{equation}
\label{f_n ne 0}
        f_n \ne 0.
\end{equation}
If $g(T)$ is another nonzero formal power series in $T$ with
coefficients in $k_0$, then it is easy to see that $f(T) \, g(T) \ne
0$ too, and that
\begin{equation}
\label{n(f(T) g(T)) = n(f(T)) + n(g(T))}
        n(f(T) \, g(T)) = n(f(T)) + n(g(T)).
\end{equation}
Let us extend $n(f(T))$ to the case where $f(T) = 0$ by putting $n(0)
= +\infty$, so that (\ref{n(f(T) g(T)) = n(f(T)) + n(g(T))}) holds
with the usual interpretations for every $f(T), g(T) \in k_0[[T]]$.
Note that
\begin{equation}
\label{n(a f(T)) = n(f(T))}
        n(a \, f(T)) = n(f(T))
\end{equation}
for every $f(T) \in k_0[[T]]$ and $a \in k_0$ with $a \ne 0$, which
may be considered as a special case of (\ref{n(f(T) g(T)) = n(f(T)) +
  n(g(T))}).  We also have that
\begin{equation}
\label{n(f(T) + g(T)) ge min(n(f(T)), n(g(T)))}
        n(f(T) + g(T)) \ge \min(n(f(T)), n(g(T)))
\end{equation}
for every $f(T), g(T) \in k_0[[T]]$, with the usual interpretations
for infinite values of $n(\cdot)$.

        Let $r$ be a positive real number strictly less than $1$,
and put
\begin{equation}
\label{|f(T)| = |f(T)|_r = r^{n(f(T))}}
        |f(T)| = |f(T)|_r = r^{n(f(T))}
\end{equation}
for every $f(T) \in k_0[[T]]$ with $f(T) \ne 0$, and $|f(T)| = 0$ when
$f(T) = 0$.  Thus
\begin{equation}
\label{|f(T) + g(T)| le max(|f(T)|, |g(T)|)}
        |f(T) + g(T)| \le \max(|f(T)|, |g(T)|)
\end{equation}
for every $f(T), g(T) \in k_0[[T]]$, by (\ref{n(f(T) + g(T)) ge
  min(n(f(T)), n(g(T)))}).  Similarly,
\begin{equation}
\label{|f(T) g(T)| = |f(T)| |g(T)|}
        |f(T) \, g(T)| = |f(T)| \, |g(T)|
\end{equation}
for every $f(T), g(T) \in k_0[[T]]$, by (\ref{n(f(T) g(T)) = n(f(T)) +
  n(g(T))}).  In particular,
\begin{equation}
\label{|a f(T)| = |f(T)|}
        |a \, f(T)| = |f(T)|
\end{equation}
for every $f(T) \in k_0[[T]]$ and $a \in k_0$ with $a \ne 0$, as in
(\ref{n(a f(T)) = n(f(T))}).  It follows that $|f(T)|$ defines an
ultranorm on $k_0[[T]]$ as a vector space over $k_0$, and using the
trivial absolute value function on $k_0$.

        This implies that
\begin{equation}
\label{|f(T) - g(T)|}
        |f(T) - g(T)|
\end{equation}
defines an ultrametric on $k_0[[T]]$, which determines a topology on
$k_0[[T]]$ in the usual way.  As before, $k_0[[T]]$ can be identified
with the Cartesian product of a family of copies of $k_0$ indexed by
${\bf Z}_+ \cup \{0\}$.  It is easy to see that the topology on
$k_0[[T]]$ determined by (\ref{|f(T) - g(T)|}) corresponds to the
product topology on this Cartesian product, using the discrete
topology on $k_0$.  If $k_0$ has only finitely many elements, then
$k_0[[T]]$ is compact with respect to this topology, by Tychonoff's
theorem.  Note that $k_0[T]$ is dense in $k_0[[T]]$ with respect to
this topology, for any $k_0$.

        If $\alpha$ is any positive real number, then
\begin{equation}
\label{|f(T)|_r^alpha = |f(T)|_{r^alpha}}
        |f(T)|_r^\alpha = |f(T)|_{r^\alpha}
\end{equation}
for every $f(T) \in k_0[[T]]$, by the definition (\ref{|f(T)| =
  |f(T)|_r = r^{n(f(T))}}) of $|f(T)|_r$.  The corresponding
ultrametric
\begin{equation}
\label{|f(T) - g(T)|_r^alpha = |f(T) - g(T)|_{r^alpha}}
        |f(T) - g(T)|_r^\alpha = |f(T) - g(T)|_{r^\alpha}
\end{equation}
determines the same topology on $k_0[[T]]$ for every $\alpha > 0$, as
in Section \ref{quasimetrics}.  Of course, this also follows from the
description of this topology on $k_0[[T]]$ as the product topology
associated to the discrete topology on $k_0$, as in the preceding
paragraph.

        Let $\{f_l(T)\}_{l = 1}^\infty$ be a sequence of elements
of $k_0[[T]]$, with
\begin{equation}
\label{f_l(T) = sum_{j = 0}^infty f_{j,l} T^j}
        f_l(T) = \sum_{j = 0}^\infty f_{j,l} \, T^j
\end{equation}
for each $l \ge 1$, and let $f(T)$ be another element of $k_0[[T]]$,
as in (\ref{f(T) = sum_{j = 0}^infty f_j T^j}).  One can check that
$\{f_l(T)\}_{l = 1}^\infty$ converges to $f(T)$ with respect to
the ultrametric (\ref{|f(T) - g(T)|}) if and only if for each
$j \ge 0$ we have that
\begin{equation}
\label{f_{j, l} = f_j}
        f_{j, l} = f_j
\end{equation}
for all sufficiently large $l$, depending on $j$.  Similarly,
$\{f_l(T)\}_{l = 1}^\infty$ is a Cauchy sequence in $k_0[[T]]$ with
respect to (\ref{|f(T) - g(T)|}) if and only if for each $j \ge 0$,
$f_{j, l}$ is constant in $l$ for sufficiently large $l$, depending on
$j$.  It follows that every Cauchy sequence in $k_0[[T]]$ with respect
to (\ref{|f(T) - g(T)|}) converges to an element of $k_0[[T]]$, so
that $k_0[[T]]$ is complete as a metric space with respect to
(\ref{|f(T) - g(T)|}).

\section{Geometric series}
\label{geometric series}

        Let $k$ be a field, with an absolute value function $|\cdot|$.
It is well known and easy to check that
\begin{equation}
\label{(1 - x) sum_{j = 0}^n x^j = 1 - x^{n + 1}}
        (1 - x) \, \sum_{j = 0}^n x^j = 1 - x^{n + 1}
\end{equation}
for every $x \in k$ and nonnegative integer $n$, where $x^0$ is
interpreted as being equal to $1$, as usual.  This implies that
\begin{equation}
\label{sum_{j = 0}^n x^j = frac{1 - x^{n + 1}}{1 - x}}
        \sum_{j = 0}^n x^j = \frac{1 - x^{n + 1}}{1 - x}
\end{equation}
for every $n \ge 0$ when $x \ne 1$, and hence that
\begin{equation}
\label{lim_{n to infty} sum_{j = 0}^n x^j = frac{1}{1 - x}}
        \lim_{n \to \infty} \sum_{j = 0}^n x^j = \frac{1}{1 - x}
\end{equation}
when $|x| < 1$, because $|x^{n + 1}| = |x|^{n + 1} \to 0$ as $n \to
\infty$.  Thus the geometric series\index{geometric series} $\sum_{j =
  0}^\infty x^j$ converges in $k$ when $|x| < 1$, with sum equal to
$1/(1 - x)$, as usual.

        Consider the case where $k = {\bf Q}$, equipped with the
$p$-adic absolute value $|\cdot|_p$ for some prime number $p$.
If $w \in {\bf Z}$, then
\begin{equation}
\label{|p w|_p le 1/p < }
        |p \, w|_p \le 1/p < 1,
\end{equation}
and we get that
\begin{equation}
\label{lim_{n to infty} sum_{j = 0}^n (p w)^j = frac{1}{1 - p w}}
        \lim_{n \to \infty} \sum_{j = 0}^n (p \, w)^j = \frac{1}{1 - p \, w},
\end{equation}
as in (\ref{lim_{n to infty} sum_{j = 0}^n x^j = frac{1}{1 - x}}).  Of
course, the limit in (\ref{lim_{n to infty} sum_{j = 0}^n (p w)^j =
  frac{1}{1 - p w}}) is taken with respect to the $p$-adic metric on
${\bf Q}$.

        Now let $k_0$ be a field, let $T$ be an indeterminate, and
let $k_0[[T]]$ be the algebra of formal power series in $T$ with
coefficients in $k_0$, as in the previous section.  Also let $r \in
(0, 1)$ be given, and let $|f(T)|$ be as defined in (\ref{|f(T)| =
  |f(T)|_r = r^{n(f(T))}}).  Note that $a(T) \in k_0[[T]]$ satisfies
\begin{equation}
\label{|a(T)| < 1}
        |a(T)| < 1
\end{equation}
if and only if the constant term in $a(T)$ is equal to $0$, which is
the same as saying that
\begin{equation}
\label{a(T) = T b(T)}
        a(T) = T \, b(T)
\end{equation}
for some $b(T) \in k_0[[T]]$.  This implies that
\begin{equation}
\label{a(T)^l = T^l b(T)^l}
        a(T)^l = T^l \, b(T)^l
\end{equation}
for each positive integer $l$, and hence that the $T^j$ term in
$a(T)^l$ is equal to $0$ when $j < l$.  We also have that
\begin{equation}
\label{(1 - a(T)) sum_{l = 0}^n a(T)^l = 1 - a(T)^{n + 1}}
        (1 - a(T)) \, \sum_{l = 0}^n a(T)^l = 1 - a(T)^{n + 1}
\end{equation}
for each nonnegative integer $n$, as in (\ref{(1 - x) sum_{j = 0}^n
  x^j = 1 - x^{n + 1}}).

        Observe that the $T^j$ term in
\begin{equation}
\label{sum_{l = 0}^n a(T)^l = sum_{l = 0}^n T^l b(T)^l}
        \sum_{l = 0}^n a(T)^l = \sum_{l = 0}^n T^l \, b(T)^l
\end{equation}
does not depend on $n$ when $n \ge j$, so that one can define
\begin{equation}
\label{sum_{l = 0}^infty a(T)^l = sum_{l = 0}^infty T^l b(T)^l}
        \sum_{l = 0}^\infty a(T)^l = \sum_{l = 0}^\infty T^l \, b(T)^l
\end{equation}
as an element of $k_0[[T]]$ in the obvious way.  Equivalently, the
sequence of partial sums (\ref{sum_{l = 0}^n a(T)^l = sum_{l = 0}^n
  T^l b(T)^l}) converges in $k_0[[T]]$ with respect respect to the
topology described in the previous section.  One can also check that
\begin{equation}
\label{(1 - a(T)) sum_{l = 0}^infty a(T)^l = 1}
        (1 - a(T)) \, \sum_{l = 0}^\infty a(T)^l = 1
\end{equation}
under these conditions, using (\ref{(1 - a(T)) sum_{l = 0}^n a(T)^l =
  1 - a(T)^{n + 1}}).  Thus $1 - a(T)$ has a multiplicative inverse in
$k_0[[T]]$ when $a(T) \in k_0[[T]]$ satisfies (\ref{|a(T)| < 1}),
which is the same as (\ref{a(T) = T b(T)}).

\section{Formal Laurent series}
\label{formal laurent series}

        Let $k_0$ be a field again, and let $T$ be an indeterminate.
By a \emph{formal Laurent series}\index{formal Laurent series} in $T$
with coefficients in $k_0$ we mean an expression of the form
\begin{equation}
\label{f(T) = sum_{j = -infty}^infty f_j T^j}
        f(T) = \sum_{j = -\infty}^\infty f_j \, T^j,
\end{equation}
with $f_j \in k_0$ for each integer $j$.  As before, such a formal
Laurent series $f(T)$ is supposed to correspond exactly to a
doubly-infinite sequence $\{f_j\}_{j = -\infty}^\infty$ of elements of
$k_0$, or equivalently to a $k_0$-valued function on the set ${\bf
  Z}$\index{Z@${\bf Z}$} of integers.  In particular, the space of
these series is a vector space over $k_0$ in an obvious way, with
respect to termwise addition and scalar multiplication.

        Let $k_0((T))$\index{k_0((T))@$k_0((T))$} be the space of
formal Laurent series $f(T)$ in $T$ with coefficients in $k_0$
such that $f_j = 0$ for all but finitely many negative integers $j$.
This is a linear subspace of the vector space of all formal Laurent
series in $T$ with coefficients in $k_0$, as in the preceding
paragraph.  An element of $k((T))$ may be expressed as
\begin{equation}
\label{f(T) = sum_{j = n}^infty f_j T^j}
        f(T) = \sum_{j = n}^\infty f_j \, T^j
\end{equation}
for some integer $n$, where it is understood that $f_j = 0$ when $j <
n$.  It is sometimes convenient to use the notation
\begin{equation}
\label{f(T) = sum_{j >> -infty} f_j T^j}
        f(T) = \sum_{j >> -\infty} f_j \, T^j,
\end{equation}
as on p27 of \cite{c}, to indicate that $f_j = 0$ for all but finitely
many negative integers $j$, without specifying an integer $n$ as in
(\ref{f(T) = sum_{j = n}^infty f_j T^j}).

        Let $f(T), g(T) \in k_0((T))$ be given, where $f(T)$ is as
in (\ref{f(T) = sum_{j >> -infty} f_j T^j}), and similarly
\begin{equation}
\label{g(T) = sum_{l >> -infty} g_l T^l}
        g(T) = \sum_{l >> -\infty} g_l \, T^l
\end{equation}
for some $g_l \in k_0$.  As in Section \ref{formal power series}, the
product of $f(T)$ and $g(T)$ can be defined formally by
\begin{equation}
\label{f(T) g(T) = sum_{j >> -infty} sum_{l >> -infty} f_j g_l T^{j + l}}
 f(T) \, g(T) = \sum_{j >> -\infty} \, \sum_{l >> -\infty} f_j \, g_l \, T^{j + l},
\end{equation}
since there are only finitely many terms involving $T^n$ for any
integer $n$.  More precisely, if we collect the terms involving $T^n$,
then we get that
\begin{equation}
\label{f(T) g(T) = sum_{n >> -infty} (f g)_n T^n}
        f(T) \, g(T) = \sum_{n >> -\infty} (f \, g)_n \, T^n,
\end{equation}
where
\begin{equation}
\label{(f g)_n = sum_{j + l = n atop j, l >> -infty} f_j g_l}
        (f \, g)_n = \sum_{j + l = n \atop j, l >> -\infty} f_j \, g_l
\end{equation}
for each $n \in {\bf Z}$.  Note that (\ref{(f g)_n = sum_{j + l = n
    atop j, l >> -infty} f_j g_l}) is indeed a sum over finitely many
$j, l \in {\bf Z}$ for every $n \in {\bf Z}$, and that (\ref{(f g)_n =
  sum_{j + l = n atop j, l >> -infty} f_j g_l}) is equal to $0$ for
all but finitely many negative integers $n$, so that (\ref{f(T) g(T) =
  sum_{n >> -infty} (f g)_n T^n}) is an element of $k_0((T))$.  As
before, we use (\ref{f(T) g(T) = sum_{n >> -infty} (f g)_n T^n}) and
(\ref{(f g)_n = sum_{j + l = n atop j, l >> -infty} f_j g_l}) as the
official definition of multiplication on $k_0((T))$, which makes sense
directly at the level of the corresponding sequences of coefficients
in $k_0$.

        It is not difficult to check that $k_0((T))$ is a commutative
ring with respect to this definition of multiplication, and in fact a
commutative algebra over $k_0$.  Let us identify each $f(T) \in k_0[[T]]$
with an element of $k_0((T))$, by putting $f_j = 0$ when $j < 0$.
This makes $k_0[[T]]$ a subalgebra of $k_0((T))$, and thus $k_0$ and
$k_0[T]$ can be identified with subalgebras of $k_0((T))$ as well.
In particular, the multiplicative identity element $1$ in $k_0$ 
corresponds to the multiplicative identity element in $k_0((T))$.

        If $f(T) \in k_0((T))$ and $f(T) \ne 0$, then $f(T)$ can be
expressed as
\begin{equation}
\label{f(T) = c T^n (1 - T b(T))}
        f(T) = c \, T^n \, (1 - T \, b(T))
\end{equation}
for some $c \in k_0$ with $c \ne 0$, $n \in {\bf Z}$, and $b(T) \in
k_0[[T]]$.  Remember that $1 - T \, b(T)$ has a multiplicative inverse
in $k_0[[T]]$, as in the previous section.  This implies that $f(T)$
has a multiplicative inverse in $k_0((T))$, given by
\begin{equation}
\label{f(T)^{-1} = c^{-1} T^{-n} (1 - T b(T))^{-1}}
        f(T)^{-1} = c^{-1} \, T^{-n} \, (1 - T \, b(T))^{-1},
\end{equation}
so that $k_0((T))$ is a field.

        Let $n(f(T))$ be the unique integer $n$ as in
(\ref{f(T) = c T^n (1 - T b(T))}) when $f(T) \in k_0((T))$ and
$f(T) \ne 0$, which is the same as saying that $f_n \ne 0$ and
$f_j = 0$ for every $j < n$.  Let us also put $n(0) = +\infty$,
so that this definition of $n(f(T))$ extends the one for $f(T) \in
k_0[[T]]$ in Section \ref{formal power series}.  It is easy to see
that (\ref{n(f(T) g(T)) = n(f(T)) + n(g(T))}), (\ref{n(a f(T)) = n(f(T))}),
and (\ref{n(f(T) + g(T)) ge min(n(f(T)), n(g(T)))}) continue to hold
for $f(T), g(T) \in k_0((T))$, with the usual interpretations
for infinite values of $n(\cdot)$.  Of course, $k_0[[T]]$ corresponds
exactly to the set of $f(T) \in k_0((T))$ such that $n(f(T)) \ge 0$.

        As in Section \ref{formal power series} again, we let $r$ be
a positive real number strictly less than $1$, and put
\begin{equation}
\label{|f(T)| = |f(T)|_r = r^{n(f(T))}, 2}
        |f(T)| = |f(T)|_r = r^{n(f(T))}
\end{equation}
when $f(T) \in k_0((T))$ and $f(T) \ne 0$, and $|0| = 0$.  This
extension of $|f(T)|$ to $f(T) \in k_0((T))$ continues to satisfy
(\ref{|f(T) + g(T)| le max(|f(T)|, |g(T)|)}), (\ref{|f(T) g(T)| =
  |f(T)| |g(T)|}), and (\ref{|a f(T)| = |f(T)|}), by the analogues of
(\ref{n(f(T) g(T)) = n(f(T)) + n(g(T))}), (\ref{n(a f(T)) = n(f(T))}),
and (\ref{n(f(T) + g(T)) ge min(n(f(T)), n(g(T)))}) for $f(T), g(T)
\in k_0((T))$.  It follows that $|f(T)|$ defines an ultrametric
absolute value function on the field $k_0((T))$, whose restriction to
$k_0$ is the trivial absolute value function on $k_0$.  Note that
(\ref{|f(T)|_r^alpha = |f(T)|_{r^alpha}}) also continues to hold for
every $\alpha > 0$ and $f(T) \in k_0((T))$.

        Put
\begin{eqnarray}
\label{T^n k_0[[T]] = ... = {g(t) in k_0((T)) : |g(T)| le r^n}}
        T^n \, k_0[[T]] & = & \{T^n \, f(T) : f(T) \in k_0[[T]]\} \\
                        & = & \{g(t) \in k_0((T)) : |g(T)| \le r^n\} \nonumber
\end{eqnarray}
for each $n \in k_0$.  If $k_0$ has only finitely many elements, then
we have seen that $k_0[[T]]$ is compact with respect to the topology
determined by the ultrametric associated to $|f(T)|$, as in Section
\ref{formal power series}.  Similarly, (\ref{T^n k_0[[T]] = ... =
  {g(t) in k_0((T)) : |g(T)| le r^n}}) is compact for each $n \in {\bf
  Z}$ in this case.  This implies that closed and bounded subsets of
$k_0((T))$ are compact, because every bounded subset of $k_0((T))$ is
contained in (\ref{T^n k_0[[T]] = ... = {g(t) in k_0((T)) : |g(T)| le
    r^n}}) for some $n \in {\bf Z}$.

        Let $k_0$ be any field again, and let $\{f_l(T)\}_{l = 1}^\infty$
be a sequence of elements of $k_0((T))$.  Thus for each $l \ge 1$,
$f_l(T)$ may be expressed as
\begin{equation}
\label{f_l(T) = sum_{j >> -infty} f_{j, l} T^j}
        f_l(T) = \sum_{j >> -\infty} f_{j, l} \, T^j,
\end{equation}
where $f_{j, l} \in k_0$.  This sequence is bounded in $k_0((T))$ with
respect to the absolute value function $|\cdot|$ defined earlier if
and only if there is an $n \in {\bf Z}$ such that $f_l(T)$ is an
element of (\ref{T^n k_0[[T]] = ... = {g(t) in k_0((T)) : |g(T)| le
    r^n}}) for each $l \ge 1$.  Equivalently, this means that there is
an $n \in {\bf Z}$ such that $f_l(T)$ can be expressed as
\begin{equation}
\label{f_l(T) = sum_{j = n}^infty f_{j, l} T^j}
        f_l(T) = \sum_{j = n}^\infty f_{j, l} \, T^j
\end{equation}
for each $l \ge 1$.

        Let $\{f_l(T)\}_{l = 1}^\infty$ be a sequence of elements
of $k_0((T))$ again, as in (\ref{f_l(T) = sum_{j >> -infty} f_{j, l}
  T^j}), and let $f(T)$ be another element of $k_0((T))$, as in
(\ref{f(T) = sum_{j >> -infty} f_j T^j}).  One can check that
$\{f_l(T)\}_{l = 1}^\infty$ converges to $f(T)$ with respect to the
ultrametric associated to the absolute value function defined earlier
if and only if $\{f_l(T)\}_{l = 1}^\infty$ is a bounded sequence in
$k_0((T))$, and for each $j \in {\bf Z}$ we have that
\begin{equation}
\label{f_{j, l} = f_j, 2}
        f_{j, l} = f_j
\end{equation}
for all sufficiently large $l$, depending on $j$.  Similarly,
$\{f_l(T)\}_{l = 1}^\infty$ is a Cauchy sequence in $k_0((T))$ with
respect to the ultrametric associated to the absolute value function
defined earlier if and only if $\{f_l(T)\}_{l = 1}^\infty$ is a
bounded sequence in $k_0((T))$, and for each $j \in {\bf Z}$, $f_{j,
  l}$ is constant in $l$ for sufficiently large $l$, depending on $j$.
In particular, it follows from this that every Cauchy sequence in
$k_0((T))$ converges to an element of $k_0((T))$, so that $k_0((T))$
is complete with respect to the ultrametric associated to the absolute
value function defined earlier.

\section{$p$-Adic integers}
\label{p-adic integers}

        Let $p$ be a prime number, and let $|\cdot|_p$ be the $p$-adic
absolute value on ${\bf Q}$, as in Section \ref{absolute value functions}.
Thus every integer $x$ satisfies
\begin{equation}
\label{|x|_p le 1}
        |x|_p \le 1,
\end{equation}
which implies that (\ref{|x|_p le 1}) also holds for every $x \in {\bf
  Q}$ in the closure of the set ${\bf Z}$ of integers with respect to
the $p$-adic metric.  Conversely, suppose that $y \in {\bf Q}$
satisfies $|y|_p \le 1$, and let us check that $y$ is in the closure
of ${\bf Z}$ in ${\bf Q}$ with respect to the $p$-adic metric.  By
definition of the $p$-adic absolute value, $y = a / b$ for some $a, b
\in {\bf Z}$ such that $b \ne 0$ and $b$ is not an integer multiple of
$p$.  It follows that there are $c, w \in {\bf Z}$ such that
\begin{equation}
\label{b c = 1 - p w}
        b \, c = 1 - p \, w,
\end{equation}
because the integers modulo $p$ form a field, and hence
\begin{equation}
\label{y = frac{a}{b} = frac{a c}{b c} = frac{a c}{1 - p w}}
        y = \frac{a}{b} = \frac{a \, c}{b \, c} = \frac{a \, c}{1 - p \, w}.
\end{equation}
We have seen that $1/(1 - p \, w)$ can be expressed as the limit of a
sequence of integers with respect to the $p$-adic metric, as in
(\ref{lim_{n to infty} sum_{j = 0}^n (p w)^j = frac{1}{1 - p w}}) in
Section \ref{geometric series}.  This implies that $y$ has the same
property, as desired.

        The set ${\bf Z}_p$\index{Z_p@${\bf Z}_p$} of \emph{$p$-adic
integers}\index{p-adic integers@$p$-adic integers} is defined to
be the closed unit ball in ${\bf Q}_p$ with respect to the $p$-adic metric,
which is to say that
\begin{equation}
\label{{bf Z}_p = {x in {bf Q}_p : |x|_p le 1}}
        {\bf Z}_p = \{x \in {\bf Q}_p : |x|_p \le 1\}.
\end{equation}
Clearly ${\bf Z} \subseteq {\bf Z}_p$, which implies that ${\bf Z}_p$
contains the closure of ${\bf Z}$ in ${\bf Q}_p$ with respect to the
$p$-adic metric.  In fact, one can check that ${\bf Z}_p$ is equal to
the closure of ${\bf Z}$ in ${\bf Q}_p$, using the remarks in the
preceding paragraph.  More precisely, one can first verify that ${\bf
  Q} \cap {\bf Z}_p$ is dense in ${\bf Z}_p$, because ${\bf Q}$ is
dense in ${\bf Q}_p$, and using the ultrametric version of the
triangle inequality.  The discussion in the previous paragraph implies
that ${\bf Z}$ is dense in ${\bf Q} \cap {\bf Z}_p$, and hence that
${\bf Z}$ is dense in ${\bf Z}_p$, as desired.

        Put
\begin{equation}
\label{p^j {bf Z} = {p^j x : x in {bf Z}}}
        p^j \, {\bf Z} = \{p^j \, x : x \in {\bf Z}\}
\end{equation}
and
\begin{equation}
\label{p^j {bf Z}_p = ... = {y in {bf Q}_p : |y|_p le p^{-j}}}
        p^j \, {\bf Z}_p = \{p^j \, x : x \in {\bf Z}_p\}
                         = \{y \in {\bf Q}_p : |y|_p \le p^{-j}\}
\end{equation}
for each $j \in {\bf Z}$.  It is easy to see that $p^j \, {\bf Z}_p$
is the closure of $p^j \, {\bf Z}$ in ${\bf Q}_p$ for each $j \in {\bf
  Z}$, because of the statement for $j = 0$ discussed in the preceding
paragraph.  Note that $p^j \, {\bf Z}$ is a subgroup of ${\bf Q}$ with
respect to addition for each $j \in {\bf Z}$, and similarly $p^j \,
{\bf Z}_p$ is a subgroup of ${\bf Q}_p$ with respect to addition for
each $j \in {\bf Z}$.  Of course, ${\bf Z}$ is a subring of ${\bf Q}$,
and one can check that ${\bf Z}_p$ is a subring of ${\bf Q}_p$ too.
If $j$ is a nonnegative integer, then $p^j \, {\bf Z}$ is an
ideal in ${\bf Z}$, and $p^j \, {\bf Z}_p$ is an ideal in ${\bf Z}_p$.
This implies that the quotients
\begin{equation}
\label{{bf Z} / p^j {bf Z}}
        {\bf Z} / p^j \, {\bf Z}
\end{equation}
and
\begin{equation}
\label{{bf Z}_p / p^j {bf Z}_p}
        {\bf Z}_p / p^j \, {\bf Z}_p
\end{equation}
are defined as commutative rings when $j \ge 0$.  The obvious
inclusion of ${\bf Z}$ in ${\bf Z}_p$ leads to a ring homomorphism
from (\ref{{bf Z} / p^j {bf Z}}) into (\ref{{bf Z}_p / p^j {bf Z}_p})
for each nonnegative integer $j$, since $p^j \, {\bf Z}$ is contained
in $p^j \, {\bf Z}_p$.  This homomorphism is actually injective for
each $j \ge 0$, because
\begin{equation}
\label{{bf Z} cap (p^j {bf Z}_p) = p^j {bf Z}}
        {\bf Z} \cap (p^j \, {\bf Z}_p) = p^j \, {\bf Z}
\end{equation}
for every nonnegative integer $j$, as one can verify directly from the
definitions.  One can also check that this homomorphism from (\ref{{bf
    Z} / p^j {bf Z}}) into (\ref{{bf Z}_p / p^j {bf Z}_p}) is
surjective for each $j \ge 0$, using the fact that ${\bf Z}$ is dense
in ${\bf Z}_p$.  Thus this homomorphism from (\ref{{bf Z} / p^j {bf
    Z}}) into (\ref{{bf Z}_p / p^j {bf Z}_p}) is an isomorphism for
each $j \ge 0$.

        It follows that (\ref{{bf Z}_p / p^j {bf Z}_p}) has exactly
$p^j$ elements for every nonnegative integer $j$, so that ${\bf Z}_p$
can be expressed as the union of $p^j$ pairwise-disjoint closed balls
of radius $p^{-j}$ in ${\bf Q}_p$ for each $j \ge 0$.  This implies
that ${\bf Z}_p$ is compact with respect to the $p$-adic metric,
because ${\bf Z}_p$ is closed and totally bounded in ${\bf Q}_p$, and
${\bf Q}_p$ is complete.  Similarly, $p^l \, {\bf Z}_p$ is a compact
subset of ${\bf Q}_p$ for every $l \in {\bf Z}$.  Of course, every bounded
subset of ${\bf Q}_p$ is contained in $p^l \, {\bf Z}_p$ for some
$l \in {\bf Z}$, and hence every closed and bounded subset of ${\bf Q}_p$
is compact.

\section{Radius of convergence}
\label{radius of convergence}

        Let $k$ be a field with an absolute value function $|\cdot|$,
and suppose that $k$ is complete with respect to the metric associated
to $|\cdot|$.  Also let $a_0, a_1, a_2, a_3, \ldots$ be a sequence of
elements of $k$, and consider the corresponding formal power series
\begin{equation}
\label{f(X) = sum_{j = 0}^infty a_j X^j}
        f(X) = \sum_{j = 0}^\infty a_j \, X^j,
\end{equation}
where $X$ is an indeterminate.  As in \cite{c, fg}, we use upper-case
letters like $X$ for indeterminates, and lower-case letters like $x$
for elements of $k$ or other fields.  If $x \in k$, then we can
consider the convergence of the power series\index{power series}
\begin{equation}
\label{sum_{j = 0}^infty a_j x^j}
        \sum_{j = 0}^\infty a_j \, x^j,
\end{equation}
where $x^0$ is interpreted as being the multiplicative identity
element $1$ in $k$, as usual.  Of course, if (\ref{sum_{j = 0}^infty
  a_j x^j}) converges for some $x \in k$, then
\begin{equation}
\label{lim_{j to infty} a_j x^j = 0}
        \lim_{j \to \infty} a_j \, x^j = 0
\end{equation}
in $k$, and hence
\begin{equation}
\label{{a_j x^j}_{j = 0}^infty is a bounded sequence}
        \{a_j \, x^j\}_{j = 0}^\infty \hbox{ is a bounded sequence}
\end{equation}
in $k$.

        If
\begin{equation}
\label{{|a_j| t^j}_{j = 0}^infty is a bounded sequence}
        \{|a_j| \, t^j\}_{j = 0}^\infty \hbox{ is a bounded sequence}
\end{equation}
in ${\bf R}$ for some nonnegative real number $t$, then
\begin{equation}
\label{lim_{j to infty} |a_j| r^j = 0}
        \lim_{j \to \infty} |a_j| \, r^j = 0
\end{equation}
for every nonnegative real number $r < t$, and in fact
\begin{equation}
\label{sum_{j = 0}^infty |a_j| r^j}
        \sum_{j = 0}^\infty |a_j| \, r^j
\end{equation}
converges in ${\bf R}$ when $0 \le r < t$.  Let $\rho$ be the supremum
of the set of $r \ge 0$ such that (\ref{lim_{j to infty} |a_j| r^j =
  0}) holds, which automatically includes $r = 0$.  As usual, $\rho$
is taken to be $+\infty$ when (\ref{lim_{j to infty} |a_j| r^j = 0})
holds for arbitrarily large $r$.  Equivalently, $\rho$ may be defined
as the supremum of the set of $t \ge 0$ such that (\ref{{|a_j| t^j}_{j
    = 0}^infty is a bounded sequence}) holds, or as the supremum of
the $r \ge 0$ such that (\ref{sum_{j = 0}^infty |a_j| r^j}) converges.
We can also characterize $\rho$ as the unique nonnegative extended
real number such that (\ref{sum_{j = 0}^infty |a_j| r^j}) converges
when $0 \le r < \rho$, and (\ref{{|a_j| t^j}_{j = 0}^infty is a
  bounded sequence}) does not hold for any $t > \rho$.  It follows
that (\ref{sum_{j = 0}^infty a_j x^j}) converges absolutely for each
$x \in k$ with $|x| < \rho$, and that (\ref{sum_{j = 0}^infty a_j
  x^j}) does not converge for any $x \in k$ with $|x| > \rho$.  Of
course, $\rho$ is known as the \emph{radius of
  convergence}\index{radius of convergence} of the formal power series
(\ref{f(X) = sum_{j = 0}^infty a_j X^j}).  It is well known that
\begin{equation}
\label{rho = (limsup_{j to infty} |a_j|^{1/j})^{-1}}
        \rho = \Big(\limsup_{j \to \infty} |a_j|^{1/j}\Big)^{-1},
\end{equation}
with the usual conventions that $1/0 = +\infty$ and $1/+\infty = 0$.

        Let us suppose for the rest of the section that either
\begin{equation}
\label{k = {bf R} or {bf C}, with the standard absolute value function}
        k = {\bf R} \hbox{ or } {\bf C},
              \hbox{ with the standard absolute value function,}
\end{equation}
or that
\begin{eqnarray}
\label{|cdot| is an ultrametric absolute value function on a field k, ...}
 && |\cdot| \hbox{ is an ultrametric absolute value function on a field } k, \\
 && \hbox{and $k$ is complete with respect to the associated ultrametric.}
                                                                  \nonumber
\end{eqnarray}
If $k = {\bf R}$ or ${\bf C}$, and (\ref{sum_{j = 0}^infty |a_j| r^j})
converges for some $r \ge 0$, then (\ref{sum_{j = 0}^infty a_j x^j})
converges absolutely for every $x \in k$ with $|x| \le r$, and the
sequence of partial sums
\begin{equation}
\label{sum_{j = 0}^n a_j x^j}
        \sum_{j = 0}^n a_j \, x^j
\end{equation}
converges uniformly to the sum (\ref{sum_{j = 0}^infty a_j x^j}) on
the closed ball
\begin{equation}
\label{overline{B}(0, r) = {x in k : |x| le r}}
        \overline{B}(0, r) = \{x \in k : |x| \le r\}.
\end{equation}
Similarly, if $k$ is as in (\ref{|cdot| is an ultrametric absolute
  value function on a field k, ...}), and (\ref{lim_{j to infty} |a_j|
  r^j = 0}) holds for some $r \ge 0$, then (\ref{sum_{j = 0}^infty a_j
  x^j}) converges in $k$ for every $x \in k$ with $|x| \le r$, and the
partial sums (\ref{sum_{j = 0}^n a_j x^j}) converge to the whole sum
(\ref{sum_{j = 0}^infty a_j x^j}) uniformly on $\overline{B}(0, r)$.
In both cases, it follows that (\ref{sum_{j = 0}^infty a_j x^j})
defines a continuous $k$-valued function on $\overline{B}(0, r)$.
Using this, one can check that (\ref{sum_{j = 0}^infty a_j x^j})
defines a continuous $k$-valued function on the open ball
\begin{equation}
\label{B(0, rho) = {x in k : |x| < rho}}
        B(0, \rho) = \{x \in k : |x| < \rho\},
\end{equation}
where $\rho$ is the radius of convergence, as in the preceding
paragraph.

        Let
\begin{equation}
\label{sum_{l = 0}^infty b_l x^l}
        \sum_{l = 0}^\infty b_l \, x^l
\end{equation}
be another power series with coefficients in $k$, and let
\begin{equation}
\label{c_n = sum_{j = 0}^n a_j b_{n - j}, 2}
        c_n = \sum_{j = 0}^n a_j \, b_{n - j}
\end{equation}
be the Cauchy product of the coefficients of (\ref{sum_{j = 0}^infty
  a_j x^j}) and (\ref{sum_{l = 0}^infty b_l x^l}) for each $n \ge 0$,
as in (\ref{c_n = sum_{j = 0}^n a_j b_{n - j}}) in Section \ref{cauchy
  products}.  Thus
\begin{equation}
        c_n \, x^n = \sum_{j = 0}^n (a_j x^j) \, (b_{n - j} \, x^{n - j})
\end{equation}
for each $n \ge 0$, so that
\begin{equation}
\label{sum_{n = 0}^infty c_n x^n}
        \sum_{n = 0}^\infty c_n \, x^n
\end{equation}
is the same as the Cauchy product of (\ref{sum_{j = 0}^infty a_j x^j})
and (\ref{sum_{l = 0}^infty b_l x^l}).  If $k = {\bf R}$ or ${\bf C}$,
and (\ref{sum_{j = 0}^infty a_j x^j}) and (\ref{sum_{l = 0}^infty b_l
  x^l}) converge absolutely for some $x \in k$, then it follows that
(\ref{sum_{n = 0}^infty c_n x^n}) also converges absolutely, and satisfies
\begin{equation}
\label{sum_{n = 0}^infty c_n x^n = ...}
 \sum_{n = 0}^\infty c_n \, x^n = \Big(\sum_{j = 0}^\infty a_j \, x^j\Big)
                                 \, \Big(\sum_{l = 0}^\infty b_l \, x^l\Big),
\end{equation}
as in Section \ref{cauchy products}.  Similarly, if $k$ is as in
(\ref{|cdot| is an ultrametric absolute value function on a field k,
  ...}), and (\ref{sum_{j = 0}^infty a_j x^j}) and (\ref{sum_{l =
    0}^infty b_l x^l}) converge in $k$ for some $x \in k$, then
(\ref{sum_{n = 0}^infty c_n x^n}) converges in $k$ too, and satisfies
(\ref{sum_{n = 0}^infty c_n x^n = ...}).

\section{Compositions}
\label{compositions}

        Let $k$ be a field, and let
\begin{equation}
\label{f(x) = sum_{j = 0}^infty a_j x^j}
        f(x) = \sum_{j = 0}^\infty a_j \, x^j
\end{equation}
and
\begin{equation}
\label{g(y) = sum_{l = 0}^infty b_l y^l}
        g(y) = \sum_{l = 0}^\infty b_l \, y^l
\end{equation}
be power series with coefficients in $k$.  We would like to consider
the composition
\begin{equation}
\label{f(g(y)) = sum_{j = 0}^infty a_j g(y)^j}
        f(g(y)) = \sum_{j = 0}^\infty a_j \, g(y)^j
\end{equation}
of these two series, at least formally for the moment.  Put
\begin{equation}
\label{E_j = ({bf Z}_+ cup {0})^j}
        E_j = ({\bf Z}_+ \cup \{0\})^j
\end{equation}
for each $j \in {\bf Z}_+$, which is the $j$th Cartesian power of
${\bf Z}_+ \cup \{0\}$, consisting of $j$-tuples $\alpha = (\alpha_1,
\ldots, \alpha_j)$ nonnegative integers.  Also put
\begin{equation}
\label{d_j(alpha) = alpha_1 + alpha_2 + cdots + alpha_j}
        d_j(\alpha) = \alpha_1 + \alpha_2 + \cdots + \alpha_j
\end{equation}
and
\begin{equation}
\label{beta_j(alpha) = b_{alpha_1} b_{alpha_2} cdots b_{alpha_j}}
        \beta_j(\alpha) = b_{\alpha_1} \, b_{\alpha_2} \cdots b_{\alpha_j}
\end{equation}
for each $\alpha \in E_j$.  Thus
\begin{equation}
\label{g(y)^j = sum_{alpha in E_j} beta_j(alpha) y^{d_j(alpha)}}
        g(y)^j = \sum_{\alpha \in E_j} \beta_j(\alpha) \, y^{d_j(\alpha)}
\end{equation}
for each $j \in {\bf Z}_+$, at least formally, and in particular this
holds for every $y \in k$ when $b_l = 0$ for all but finitely many
$l$, so that the sum on the right side of (\ref{g(y)^j = sum_{alpha in
    E_j} beta_j(alpha) y^{d_j(alpha)}}) reduces to a finite sum.

        It will be convenient to take $E_0$ to be a set with exactly
one element not in $E_j$ for any $j \in {\bf Z}_+$, so that the
$E_j$'s are pairwise disjoint for all $j \ge 0$.  Put
\begin{equation}
\label{E = bigcup_{j = 0}^infty E_j}
        E = \bigcup_{j = 0}^\infty E_j,
\end{equation}
and let $\phi$ be the $k$-valued function on $E$ defined by
\begin{equation}
\label{phi(alpha) = a_j beta_j(alpha)}
        \phi(\alpha) = a_j \, \beta_j(\alpha)
\end{equation}
for each $\alpha \in E_j$ when $j \ge 1$, and $\phi = a_0$ on $E_0$.
Similarly, let $d$ be the function on $E$ with values in ${\bf Z}_+
\cup \{0\}$ defined by
\begin{equation}
\label{d = d_j on E_j}
        d = d_j \quad\hbox{on } E_j
\end{equation}
for each $j \ge 0$, with $d_0 = 0$ on $E_0$.  Combining (\ref{f(g(y))
  = sum_{j = 0}^infty a_j g(y)^j}) and (\ref{g(y)^j = sum_{alpha in
    E_j} beta_j(alpha) y^{d_j(alpha)}}), we get that
\begin{equation}
\label{f(g(y)) = ... = sum_{alpha in E} phi(alpha) y^{d(alpha)}}
 f(g(y)) = a_0 + \sum_{j = 1}^\infty a_j \, \Big(\sum_{\alpha \in E_j}
                                         \beta_j(\alpha) \, y^{d_j(\alpha)}\Big)
         = \sum_{\alpha \in E} \phi(\alpha) \, y^{d(\alpha)},
\end{equation}
at least formally.  In particular, if $a_j = 0$ for all but finitely
many $j$, and $b_l = 0$ for all but finitely many $l$, then $\phi \in
c_{00}(E, k)$, and (\ref{f(g(y)) = ... = sum_{alpha in E} phi(alpha)
  y^{d(alpha)}}) holds for every $y \in k$.

        Put
\begin{equation}
\label{A_n = {alpha in E : d(alpha) = n} = ...}
        A_n = \{\alpha \in E : d(\alpha) = n\}
            = \bigcup_{j = 0}^\infty \{\alpha \in E_j : d_j(\alpha) = n\}
\end{equation}
for each nonnegative integer $n$, so that the $A_n$'s are pairwise
disjoint and
\begin{equation}
\label{E = bigcup_{n = 0}^infty A_n}
        E = \bigcup_{n = 0}^\infty A_n.
\end{equation}
If we put
\begin{equation}
\label{c_n = sum_{alpha in A_n} phi(alpha)}
        c_n = \sum_{\alpha \in A_n} \phi(\alpha)
\end{equation}
for each $n \ge 0$, then we get that
\begin{equation}
\label{f(g(y)) = ... = sum_{n = 0}^infty c_n y^n}
 f(g(y)) = \sum_{n = 0}^\infty \Big(\sum_{\alpha \in A_n} \phi(\alpha)
                                                         \, y^{d(\alpha)}\Big) 
         = \sum_{n = 0}^\infty c_n \, y^n,
\end{equation}
at least formally, by (\ref{f(g(y)) = ... = sum_{alpha in E}
  phi(alpha) y^{d(alpha)}}).  As before, if $a_j = 0$ for all but
finitely many $j$, and $b_l = 0$ for all but finitely many $l$, then
$\phi \in c_{00}(E, k)$, and (\ref{f(g(y)) = ... = sum_{n = 0}^infty
  c_n y^n}) holds for all $y \in k$.

        Suppose for the moment that $k = {\bf R}$ or ${\bf C}$, with
the standard absolute value function, and that
\begin{equation}
\label{sum_{j = 0}^infty |a_j| r^j, 2}
        \sum_{j = 0}^\infty |a_j| \, r^j
\end{equation}
converges for some $r \ge 0$.  Suppose also that 
\begin{equation}
\label{sum_{l = 0}^infty |b_l| t^l le r}
        \sum_{l = 0}^\infty |b_l| \, t^l \le r
\end{equation}
for some $t > 0$, so that for each $y \in k$ with $|y| \le t$ the
series in (\ref{g(y) = sum_{l = 0}^infty b_l y^l}) converges
absolutely and satisfies
\begin{equation}
\label{|g(y)| le r}
        |g(y)| \le r.
\end{equation}
It follows that the series in (\ref{f(g(y)) = sum_{j = 0}^infty a_j
  g(y)^j}) converges absolutely when $|y| \le t$ too.  Observe that
\begin{equation}
\label{sum_{alpha in E_j} |beta_j(alpha)| t^{d_j(alpha)} = ... le r^j}
        \sum_{\alpha \in E_j} |\beta_j(\alpha)| \, t^{d_j(\alpha)}
             = \Big(\sum_{l = 0}^\infty |b_l| \, t^l\Big)^j \le r^j
\end{equation}
for each $j \ge 1$, and hence that
\begin{equation}
\label{sum_{alpha in E} |phi(alpha)| t^{d(alpha)} le ... }
 \quad \sum_{\alpha \in E} |\phi(\alpha)| \, t^{d(\alpha)}
         = |a_0| + \sum_{j = 1}^\infty \Big(\sum_{\alpha \in E_j} |a_j| \, 
                                    |\beta_j(\alpha)| \, t^{d_j(\alpha)}\Big) 
          \le \sum_{j = 0}^\infty |a_j| \, r^j.
\end{equation}
If $y \in k$ and $|y| \le t$, then we get that
\begin{equation}
\label{beta_j(alpha) y^{d_j(alpha)} is summable on E_j}
        \beta_j(\alpha) \, y^{d_j(\alpha)} \hbox{ is summable on } E_j
\end{equation}
for each $j \ge 1$, and that
\begin{equation}
\label{phi(alpha) y^{d(alpha)} is summable on E}
        \phi(\alpha) \, y^{d(\alpha)} \hbox{ is summable on } E.
\end{equation}
Using (\ref{beta_j(alpha) y^{d_j(alpha)} is summable on E_j}), it is
easy to see that (\ref{g(y)^j = sum_{alpha in E_j} beta_j(alpha)
  y^{d_j(alpha)}}) holds for each $j$, as before.  More precisely, the
sum on the right side of (\ref{g(y)^j = sum_{alpha in E_j}
  beta_j(alpha) y^{d_j(alpha)}}) may be treated as an iterated sum
over each of the $j$ factors of ${\bf Z}_+ \cup \{0\}$ in $E_j$, which
can be evaluated using (\ref{g(y) = sum_{l = 0}^infty b_l y^l}) in
each coordinate.  Similarly, (\ref{f(g(y)) = ... = sum_{alpha in E}
  phi(alpha) y^{d(alpha)}}) holds under these conditions, by the
remarks in Section \ref{sums of sums}.

        Using (\ref{E = bigcup_{n = 0}^infty A_n}), we can rearrange
the sum in (\ref{sum_{alpha in E} |phi(alpha)| t^{d(alpha)} le ... }),
to get that
\begin{equation}
\label{sum_{n = 0}^infty (sum_{alpha in A_n} |phi(alpha)| t^{d(alpha)}) = ...}
 \sum_{n = 0}^\infty \Big(\sum_{\alpha \in A_n} |\phi(\alpha)| \, t^{d(\alpha)}\Big)
               = \sum_{\alpha \in E} |\phi(\alpha)| \, t^{d(\alpha)}
                  \le \sum_{j = 0}^\infty |a_j| \, r^j.
\end{equation}
Equivalently, this means that
\begin{equation}
\label{sum_{n = 0}^infty (sum_{alpha in A_n} |phi(alpha)|) t^n le ...}
 \sum_{n = 0}^\infty \Big(\sum_{\alpha \in A_n} |\phi(\alpha)|\Big) \, t^n
                                     \le \sum_{j = 0}^\infty |a_j| \, r^j,
\end{equation}
by the definition (\ref{A_n = {alpha in E : d(alpha) = n} = ...}) of
$A_n$.  Because $t > 0$, it follows that
\begin{equation}
\label{phi(alpha) is summable on A_n}
        \phi(\alpha) \hbox{ is summable on } A_n
\end{equation}
for each $n \ge 0$, so that the sum in (\ref{c_n = sum_{alpha in A_n}
  phi(alpha)}) is defined for each $n \ge 0$.  If $|y| \le t$, then
(\ref{phi(alpha) y^{d(alpha)} is summable on E}) permits us to go from
(\ref{f(g(y)) = ... = sum_{alpha in E} phi(alpha) y^{d(alpha)}}) to
(\ref{f(g(y)) = ... = sum_{n = 0}^infty c_n y^n}), as in Section
\ref{sums of sums}.  More precisely, the sum on the right side of
(\ref{f(g(y)) = ... = sum_{n = 0}^infty c_n y^n}) converges absolutely
when $|y| \le t$, and the value of the sum is equal to $f(g(y))$.
Note that
\begin{equation}
\label{|c_n| le sum_{alpha in A_n} |phi(alpha)|}
        |c_n| \le \sum_{\alpha \in A_n} |\phi(\alpha)|
\end{equation}
for each $n \ge 0$, by the definition (\ref{c_n = sum_{alpha in A_n}
  phi(alpha)}) of $c_n$.  Thus (\ref{sum_{n = 0}^infty (sum_{alpha in
    A_n} |phi(alpha)|) t^n le ...}) implies that
\begin{equation}
\label{sum_{n = 0}^infty |c_n| t^n le sum_{j = 0}^infty |a_j| r^j}
 \sum_{n = 0}^\infty |c_n| \, t^n \le \sum_{j = 0}^\infty |a_j| \, r^j,
\end{equation}
and in particular that the left side of (\ref{sum_{n = 0}^infty |c_n|
  t^n le sum_{j = 0}^infty |a_j| r^j}) converges under these
conditions.

        Now let $k$ be any field with an ultrametric absolute value
function $|\cdot|$ such that $k$ is complete with respect to the
ultrametric associated to $|\cdot|$.  Suppose that
\begin{equation}
\label{lim_{j to infty} |a_j| r^j = 0, 2}
        \lim_{j \to \infty} |a_j| \, r^j = 0
\end{equation}
for some $r \ge 0$, and that $t > 0$ satisfies
\begin{equation}
\label{lim_{l to infty} |b_l| t^l = 0}
        \lim_{l \to \infty} |b_l| \, t^l = 0
\end{equation}
and
\begin{equation}
\label{max_{l ge 0} |b_l| t^l le r}
        \max_{l \ge 0} |b_l| \, t^l \le r.
\end{equation}
If $y \in k$ and $|y| \le t$, then (\ref{lim_{l to infty} |b_l| t^l =
  0}) implies that the series in (\ref{g(y) = sum_{l = 0}^infty b_l
  y^l}) defining $g(y)$ converges in $k$, and (\ref{max_{l ge 0} |b_l|
  t^l le r}) implies that
\begin{equation}
\label{|g(y)| le r, 2}
        |g(y)| \le r.
\end{equation}
It follows that the series in (\ref{f(g(y)) = sum_{j = 0}^infty a_j
  g(y)^j}) converges in $k$ as well under these conditions, by
(\ref{lim_{j to infty} |a_j| r^j = 0, 2}).  Using (\ref{lim_{l to
    infty} |b_l| t^l = 0}), one can check that
\begin{equation}
\label{|beta_j(alpha)| t^{d_j(alpha)} vanishes at infinity on E_j}
        |\beta_j(\alpha)| \, t^{d_j(\alpha)} \hbox{ vanishes at infinity on } E_j
\end{equation}
for each $j \ge 1$.  Moreover,
\begin{equation}
\label{max_{alpha in E_j} |beta_j(alpha)| t^{d_j(alpha)} = ... le r^j}
        \max_{\alpha \in E_j} |\beta_j(\alpha)| \, t^{d_j(\alpha)}
             = \Big(\max_{l \ge 0} (|b_l| \, t^l)\Big)^j \le r^j
\end{equation}
for each $j \ge 1$, by (\ref{max_{l ge 0} |b_l| t^l le r}).  Thus
\begin{equation}
\label{max_{alpha in E_j} |phi(alpha)| t^{d(alpha)} = ... le |a_j| r^j}
        \max_{\alpha \in E_j} |\phi(\alpha)| \, t^{d(\alpha)}
          = |a_j| \, \max_{\alpha \in E_j} |\beta_j(\alpha)| \, t^{d_j(\alpha)}
             \le |a_j| \, r^j
\end{equation}
for each $j \ge 1$, which tends to $0$ as $j \to \infty$, by
(\ref{lim_{j to infty} |a_j| r^j = 0, 2}).  This implies that
\begin{equation}
\label{|phi(alpha)| t^{d(alpha)} vanishes at infinity on E}
  |\phi(\alpha)| \, t^{d(\alpha)} \hbox{ vanishes at infinity on } E,
\end{equation}
using (\ref{|beta_j(alpha)| t^{d_j(alpha)} vanishes at infinity on
  E_j}) to get that the restriction of $|\phi(\alpha)| \,
t^{d(\alpha)}$ to $E_j$ vanishes at infinity on $E_j$ for each $j \ge
1$.  If $y \in k$ and $|y| \le t$, then we obtain that
\begin{equation}
\label{beta_j(alpha) y^{d_j(alpha)} vanishes at infinity on E_j}
        \beta_j(\alpha) \, y^{d_j(\alpha)} \hbox{ vanishes at infinity on } E_j
\end{equation}
for each $j \ge 1$, by (\ref{|beta_j(alpha)| t^{d_j(alpha)} vanishes
  at infinity on E_j}), and that
\begin{equation}
\label{phi(alpha) y^{d(alpha)} vanishes at infinity on E}
        \phi(\alpha) \, y^{d(\alpha)} \hbox{ vanishes at infinity on } E,
\end{equation}
by (\ref{|phi(alpha)| t^{d(alpha)} vanishes at infinity on E}).
Hence
\begin{equation}
\label{sum_{alpha in E_j} beta_j(alpha) y^{d_j(alpha)} satisfies the gcc}
        \sum_{\alpha \in E_j} \beta_j(\alpha) \, y^{d_j(\alpha)}
                     \hbox{ satisfies the generalized Cauchy criterion}
\end{equation}
for each $j \ge 1$, as in Section \ref{generalized convergence}, and
similarly
\begin{equation}
\label{sum_{alpha in E} phi(alpha) y^{d(alpha)} satisfies the gcc}
        \sum_{\alpha \in E} \phi(\alpha) \, y^{d(\alpha)}
                       \hbox{ satisfies the generalized Cauchy criterion}.
\end{equation}
This means that these sums can be defined as elements of $k$, as in
Section \ref{generalized convergence, continued}, because $k$ is
complete.  As before, (\ref{g(y)^j = sum_{alpha in E_j} beta_j(alpha)
  y^{d_j(alpha)}}) and (\ref{f(g(y)) = ... = sum_{alpha in E}
  phi(alpha) y^{d(alpha)}}) hold under these conditions, by the
remarks in Section \ref{sums of sums}.

        Of course, (\ref{|phi(alpha)| t^{d(alpha)} vanishes at infinity on E})
implies that the restriction of $\phi(\alpha) \, t^{d(\alpha)}$ to
$\alpha \in A_n$ vanishes at infinity on $A_n$ for each $n \ge 0$.
This implies that
\begin{equation}
        \phi(\alpha) \hbox{ vanishes at infinity on } A_n
\end{equation}
for each $n \ge 0$, since $d(\alpha) = n$ for every $\alpha \in A_n$,
by the definition (\ref{A_n = {alpha in E : d(alpha) = n} = ...}) of
$A_n$, and $t > 0$.  It follows that
\begin{equation}
\label{sum_{alpha in A_n} phi(alpha) satisfies the gcc}
        \sum_{\alpha \in A_n} \phi(\alpha)
                     \hbox{ satisfies the generalized Cauchy criterion}
\end{equation}
for each $n \ge 0$, as in Section \ref{generalized convergence,
  continued}, so that (\ref{c_n = sum_{alpha in A_n} phi(alpha)}) is
well defined for each $n \ge 0$.  If $y \in k$ and $|y| \le t$, then
(\ref{sum_{alpha in E} phi(alpha) y^{d(alpha)} satisfies the gcc})
permits us to go from (\ref{f(g(y)) = ... = sum_{alpha in E}
  phi(alpha) y^{d(alpha)}}) to (\ref{f(g(y)) = ... = sum_{n = 0}^infty
  c_n y^n}) again, as in Section \ref{sums of sums}.  More precisely,
this means that the sum on the right side of (\ref{f(g(y)) = ... =
  sum_{n = 0}^infty c_n y^n}) converges in $k$ when $|y| \le t$,
and that the value of the sum is equal to $f(g(y))$.  Note that
\begin{equation}
\label{|c_n| le max_{alpha in A_n} |phi(alpha)|}
        |c_n| \le \max_{\alpha \in A_n} |\phi(\alpha)|
\end{equation}
for each $n \ge 0$, by the definition (\ref{c_n = sum_{alpha in A_n}
  phi(alpha)}) of $c_n$ and the ultrametric version of the triangle
inequality.  Thus
\begin{equation}
\label{|c_n| t^n le ... = max_{alpha in A_n} |phi(alpha)| t^{d(alpha)}}
        |c_n| \, t^n \le \max_{\alpha \in A_n} |\phi(\alpha)| \, t^n
                    = \max_{\alpha \in A_n} |\phi(\alpha)| \, t^{d(\alpha)}
\end{equation}
for each $n \ge 0$, using the definition (\ref{A_n = {alpha in E :
    d(alpha) = n} = ...}) of $A_n$ in the second step.  In particular,
\begin{equation}
\label{lim_{n to infty} |c_n| t^n = 0}
        \lim_{n \to \infty} |c_n| \, t^n = 0,
\end{equation}
because of (\ref{|phi(alpha)| t^{d(alpha)} vanishes at infinity on
  E}), and because the $A_n$'s are pairwise-disjoint subsets of $E$.

\section{Compositions, continued}
\label{compositions, continued}

        Let $k$ be a field, and let
\begin{equation}
\label{f(X) = sum_{j = 0}^infty a_j X^j, 2}
        f(X) = \sum_{j = 0}^\infty a_j \, X^j
\end{equation}
and
\begin{equation}
\label{g(Y) = sum_{l = 0}^infty b_l Y^l}
        g(Y) = \sum_{l = 0}^\infty b_l \, Y^l
\end{equation}
be formal power series with coefficients in $k$.  As in the previous
section, we would like to consider the composition
\begin{equation}
\label{f(g(Y)) = sum_{j = 0}^infty a_j g(Y)^j}
        f(g(Y)) = \sum_{j = 0}^\infty a_j \, g(Y)^j
\end{equation}
of these two series, at least formally.  If $E_j$, $d_j(\alpha)$, and
$\beta_j(\alpha)$ are as in (\ref{E_j = ({bf Z}_+ cup {0})^j}),
(\ref{d_j(alpha) = alpha_1 + alpha_2 + cdots + alpha_j}), and
(\ref{beta_j(alpha) = b_{alpha_1} b_{alpha_2} cdots b_{alpha_j}}),
respectively, then we have that
\begin{equation}
\label{g(Y)^j = sum_{alpha in E_j} beta_j(alpha) Y^{d_j(alpha)}}
        g(Y)^j = \sum_{\alpha \in E_j} \beta_j(\alpha) \, Y^{d_j(\alpha)}
\end{equation}
for each $j \in {\bf Z}_+$, as in (\ref{g(y)^j = sum_{alpha in E_j}
  beta_j(alpha) y^{d_j(alpha)}}).  More precisely, put
\begin{equation}
\label{A_{j, n} = {alpha in E_j : d_j(alpha) = n}}
        A_{j, n} = \{\alpha \in E_j : d_j(\alpha) = n\}
\end{equation}
for each $j \in {\bf Z}_+$ and nonnegative integer $n$, so that the
$A_{j, n}$'s are pairwise-disjoint finite subsets of $E_j$ such that
\begin{equation}
\label{E_j = bigcup_{n = 0}^infty A_{j, n}}
        E_j = \bigcup_{n = 0}^\infty A_{j, n}.
\end{equation}
Thus
\begin{equation}
\label{c_{j, n} = sum_{alpha in A_{j, n}} beta_j(alpha)}
        c_{j, n} = \sum_{\alpha \in A_{j, n}} \beta_j(\alpha)
\end{equation}
is defined as a finite sum of elements of $k$ for each $j \in {\bf
  Z}_+$ and $n \ge 0$, and (\ref{g(Y)^j = sum_{alpha in E_j}
  beta_j(alpha) Y^{d_j(alpha)}}) may be interpreted as saying that
\begin{equation}
\label{g(Y)^j = sum_{n = 0}^infty c_{j, n} Y^n}
        g(Y)^j = \sum_{n = 0}^\infty c_{j, n} \, Y^n
\end{equation}
for each $j \in {\bf Z}_+$, as formal power series in $Y$.

        If $E$, $\phi$, and $d$ are as in (\ref{E = bigcup_{j = 0}^infty E_j}),
(\ref{phi(alpha) = a_j beta_j(alpha)}), and (\ref{d = d_j on E_j}),
respectively, then we get that
\begin{equation}
\label{f(g(Y)) = ... = sum_{alpha in E} phi(alpha) Y^{d(alpha)}}
 f(g(Y)) = a_0 + \sum_{j = 1}^\infty \Big(\sum_{\alpha \in E_j} a_j \,
                                     \beta_j(\alpha) \, Y^{d_j(\alpha)}\Big)
          = \sum_{\alpha \in E} \phi(\alpha) \, Y^{d(\alpha)},
\end{equation}
at least formally, as in (\ref{f(g(y)) = ... = sum_{alpha in E}
  phi(alpha) y^{d(alpha)}}).  Similarly, if $A_n$ and $c_n$ are as in
(\ref{A_n = {alpha in E : d(alpha) = n} = ...}) and (\ref{c_n =
  sum_{alpha in A_n} phi(alpha)}), respectively, then we get that
\begin{equation}
\label{f(g(Y)) = ... = sum_{n = 0}^infty c_n Y^n}
 f(g(Y)) = \sum_{n = 0}^\infty \Big(\sum_{\alpha \in A_n} \phi(\alpha)
                                                       \, Y^{d(\alpha)}\Big)
         = \sum_{n = 0}^\infty c_n \, Y^n,
\end{equation}
at least formally, as in (\ref{f(g(y)) = ... = sum_{n = 0}^infty c_n
  y^n}).  More precisely, note that
\begin{equation}
\label{A_n = bigcup_{j = 1}^infty A_{j, n}}
        A_n = \bigcup_{j = 1}^\infty A_{j, n}
\end{equation}
when $n \ge 1$, and that
\begin{equation}
\label{A_0 = E_0 cup (bigcup_{j = 1}^infty A_{j, 0})}
        A_0 = E_0 \cup \Big(\bigcup_{j = 1}^\infty A_{j, 0}\Big),
\end{equation}
where $E_0$ is as in the previous section.  This implies that
\begin{equation}
\label{c_n = sum_{j = 1}^infty a_j c_{j, n}}
        c_n = \sum_{j = 1}^\infty a_j \, c_{j, n}
\end{equation}
when $n \ge 1$, and that
\begin{equation}
\label{c_0 = a_0 + sum_{j = 1}^infty a_j c_{j, 0}}
        c_0 = a_0 + \sum_{j = 1}^\infty a_j \, c_{j, 0},
\end{equation}
at least formally.  Thus (\ref{f(g(Y)) = ... = sum_{n = 0}^infty c_n
  Y^n}) is basically the same as saying that
\begin{equation}
\label{f(g(Y)) = ... = sum_{n = 0}^infty c_n Y^n, 2}
 \quad  f(g(Y)) = a_0 + \sum_{j = 1}^\infty a_j \, g(Y)^j
            = a_0 + \sum_{j = 1}^\infty \sum_{n = 0}^\infty a_j \, c_{j, n} \, Y^n
                = \sum_{n = 0}^\infty c_n \, Y^n,
\end{equation}
at least formally again, using (\ref{g(Y)^j = sum_{n = 0}^infty c_{j,
    n} Y^n}) in the second step, and interchanging the order of
summation in the third step.  If $a_j = 0$ for all but finitely many
$j$, so that $f(X)$ is actually a formal polynomial in $X$, then
(\ref{f(g(Y)) = sum_{j = 0}^infty a_j g(Y)^j}) reduces to a finite sum
of products of formal power series, as in Section \ref{formal power
  series}.  In this case, (\ref{c_n = sum_{j = 1}^infty a_j c_{j, n}})
and (\ref{c_0 = a_0 + sum_{j = 1}^infty a_j c_{j, 0}}) reduce to
finite sums in $k$, and there is no problem with (\ref{f(g(Y)) = ... =
  sum_{n = 0}^infty c_n Y^n, 2}).

        Observe that
\begin{equation}
\label{c_{j, 0} = b_0^j}
        c_{j, 0} = b_0^j
\end{equation}
for each $j \ge 1$, so that (\ref{c_0 = a_0 + sum_{j = 1}^infty a_j c_{j,
    0}}) becomes
\begin{equation}
\label{c_0 = sum_{j = 0}^infty a_j b_0^j}
        c_0 = \sum_{j = 0}^\infty a_j \, b_0^j.
\end{equation}
As in (\ref{f(g(Y)) = ... = sum_{n = 0}^infty c_n Y^n, 2}), this is
the expected constant term in $f(g(Y))$, at least formally.  If $a_j
\ne 0$ for infinitely many $j$, and $b_0 \ne 0$, then one would
normally need some additional convergence hypotheses to make sense of
(\ref{c_0 = sum_{j = 0}^infty a_j b_0^j}).  Similarly, $c_n$ may
involve sums of infinitely many nonzero terms when $a_j \ne 0$ for
infinitely many $j$, $b_0 \ne 0$, and $n \ge 1$.  However, if $b_0 =
0$, then it is easy to see that (\ref{f(g(Y)) = sum_{j = 0}^infty a_j
  g(Y)^j}) makes sense as a formal power series in $Y$.  In this case,
\begin{equation}
\label{beta_j(alpha) = 0}
        \beta_j(\alpha) = 0
\end{equation}
for every $\alpha \in E_j$ such that $d_j(\alpha) < j$, because at
least one of the coordinates of $\alpha$ has to be equal to $0$.  This
implies that
\begin{equation}
\label{c_{j, n} = 0}
        c_{j, n} = 0
\end{equation}
when $j > n$, so that the sums in (\ref{c_n = sum_{j = 1}^infty a_j
  c_{j, n}}) and (\ref{c_0 = a_0 + sum_{j = 1}^infty a_j c_{j, 0}})
have only finitely many nonzero terms.  It is a bit simpler to take
\begin{equation}
\label{E_j = {bf Z}_+^j}
        E_j = {\bf Z}_+^j
\end{equation}
in this situation, instead of (\ref{E_j = ({bf Z}_+ cup {0})^j}),
which amounts to throwing away the terms that are automatically equal
to $0$ when $b_0 = 0$.  With this definition of $E_j$, we have that
$d_j(\alpha) \ge j$ on $E_j$, $A_{j, n} = \emptyset$ when $j > n$, and
that $A_n$ has only finitely many elements for each $n \ge 0$.

        Now let $k_0$ be a field, let $T$ be an indeterminate, and
let $k_0((T))$ be the corresponding field of formal Laurent series
with coefficients in $k_0$ and poles of finite order in $T$, as in
Section \ref{formal laurent series}.  Also let $r$ be a positive real
number strictly less than $1$, and let $|\cdot|$ be the corresponding
absolute value function on $k_0((T))$, as before.  If $f(X)$ is a
formal power series in an indeterminate $X$ with coefficients $a_j \in
k_0$, as in (\ref{f(X) = sum_{j = 0}^infty a_j X^j, 2}), and if
\begin{equation}
\label{g(T) = sum_{l = 0}^infty b_l T^l}
        g(T) = \sum_{l = 0}^\infty b_l \, T^l
\end{equation}
is a formal power series in $T$ with coefficients $b_l \in k$, then
\begin{equation}
\label{f(g(T)) = sum_{j = 0}^infty a_j g(T)^j}
        f(g(T)) = \sum_{j = 0}^\infty a_j \, g(T)^j
\end{equation}
may be considered as an infinite series with terms in $k_0((T))$.  Of
course, if $a_j = 0$ for all but finitely many $j$, so that $f(X)$ is
a formal polynomial in $X$, then (\ref{f(g(T)) = sum_{j = 0}^infty a_j
  g(T)^j}) reduces to a finite sum, which makes sense for every $g(T)
\in k_0((T))$.  Otherwise, if $a_j \ne 0$ for infinitely many $j$,
then (\ref{f(g(T)) = sum_{j = 0}^infty a_j g(T)^j}) converges in
$k_0((T))$ with respect to the ultrametric associated to the
absolute value function $|\cdot|$ when $|g(T)| < 1$, which means
that $g(T) \in k_0[[T]]$ and $b_0 = 0$.

\section{Changing centers}
\label{changing centers}

        Let $k$ be a field, let
\begin{equation}
\label{f(x) = sum_{j = 0}^infty a_j x^j, 2}
        f(x) = \sum_{j = 0}^\infty a_j \, x^j
\end{equation}
be a power series with coefficients in $k$, and let $b_0$ be an
element of $k$.  We would like to consider
\begin{equation}
\label{f(b_0 + y) = sum_{j = 0}^infty a_j (b_0 + y)^j}
        f(b_0 + y) = \sum_{j = 0}^\infty a_j \, (b_0 + y)^j
\end{equation}
as a power series in $y$, at least formally, which corresponds to
(\ref{f(g(y)) = sum_{j = 0}^infty a_j g(y)^j}) in Section
\ref{compositions} with $g(y) = b_0 + y$.  Using the binomial
theorem, we get that
\begin{equation}
\label{f(b_0 + y) = sum_j sum_{l = 0}^j a_j {j choose l} b_0^{j - l} y^l = ...}
 f(b_0 + y) = \sum_{j = 0}^\infty \sum_{l = 0}^j a_j \, {j \choose l} \cdot
                                                        b_0^{j - l} \, y^l
      = \sum_{l = 0}^\infty \Big(\sum_{j = l}^\infty a_j \, {j \choose l} \cdot
                                                        b_0^{j - l}\Big) \, y^l
\end{equation}
at least formally again.  As usual, there is no problem with this when
$a_j = 0$ for all but finitely many $j$.

        Suppose that $k = {\bf R}$ or ${\bf C}$, with the standard
absolute value function, and that
\begin{equation}
\label{sum_{j = 0}^infty |a_j| r^j, 3}
        \sum_{j = 0}^\infty |a_j| \, r^j
\end{equation}
converges for some $r > 0$.  Of course, this implies that the series
in (\ref{f(x) = sum_{j = 0}^infty a_j x^j, 2}) converges absolutely
when $x \in k$ satisfies $|x| \le r$.  Suppose also that
\begin{equation}
\label{|b_0| + t le r}
        |b_0| + t \le r
\end{equation}
for some $t > 0$, so that the series in (\ref{f(b_0 + y) = sum_{j =
    0}^infty a_j (b_0 + y)^j}) converges absolutely when $y \in k$
satisfies $|y| \le t$.  Using the binomial theorem again, we get that
\begin{equation}
\label{sum_j sum_{l le j} |a_j| {j choose l} |b_0|^{j - l} t^l = ...}
 \sum_{j = 0}^\infty \sum_{l = 0}^j |a_j| \, {j \choose l} |b_0|^{j - l} \, t^l
        = \sum_{j = 0}^\infty |a_j| (|b_0| + t)^j
        \le \sum_{j = 0}^\infty |a_j| \, r^j < \infty.
\end{equation}
It follows that
\begin{equation}
\label{sum_l (sum_{j ge l} |a_j| {j choose l} |b_0|^{j - l}) t^l = ...}
 \sum_{l = 0}^\infty \Big(\sum_{j = l}^\infty |a_j| \, {j \choose l} \,
                                             |b_0|^{j - l}\Big) \, t^l
 = \sum_{j = 0}^\infty \sum_{l = 0}^j |a_j| \, {j \choose l} \,
                                              |b_0|^{j - l} \, t^l < \infty.
\end{equation}
In particular,
\begin{equation}
\label{sum_{j = l}^infty |a_j| {j choose l} |b_0|^{j - l} < infty}
 \sum_{j = l}^\infty |a_j| \, {j \choose l} \, |b_0|^{j - l} < \infty
\end{equation}
for each $l \ge 0$, which means that the sum in $j$ on the right side
of (\ref{f(b_0 + y) = sum_j sum_{l = 0}^j a_j {j choose l} b_0^{j - l}
  y^l = ...}) converges absolutely for each $l$.  The finiteness of
(\ref{sum_l (sum_{j ge l} |a_j| {j choose l} |b_0|^{j - l}) t^l =
  ...})  implies that the sums in (\ref{f(b_0 + y) = sum_j sum_{l =
    0}^j a_j {j choose l} b_0^{j - l} y^l = ...}) converge absolutely
for every $y \in k$ with $|y| \le t$, and permits the interchange of
summation in the second step in (\ref{f(b_0 + y) = sum_j sum_{l = 0}^j
  a_j {j choose l} b_0^{j - l} y^l = ...}), as in Section \ref{sums of
  sums}.

        Now let $k$ be an arbitrary field with an ultrametric
absolute value function $|\cdot|$ such that $k$ is complete
with respect to the associated ultrametric.  Suppose that
\begin{equation}
\label{lim_{j to infty} |a_j| r^j = 0, 3}
        \lim_{j \to \infty} |a_j| r^j = 0
\end{equation}
for some $r > 0$, so that the series in (\ref{f(x) = sum_{j = 0}^infty
  a_j x^j, 2}) converges in $k$ when $x \in k$ satisfies $|x| \le r$.
Suppose also that
\begin{equation}
\label{|b_0| le r}
        |b_0| \le r,
\end{equation}
which implies that the series in (\ref{f(b_0 + y) = sum_{j = 0}^infty
  a_j (b_0 + y)^j}) converges in $k$ for every $y \in k$ with $|y| \le
r$.  Observe that
\begin{equation}
\label{|a_j {j choose l} cdot b_0^{j - l}| r^l le ... le |a_j| r^j}
        \biggl|a_j \, {j \choose l} \cdot b_0^{j - l}\biggr| \, r^l
                    \le |a_j| \, |b_0|^{j - l} \, r^l \le |a_j| \, r^j
\end{equation}
for every $j \ge l \ge 0$, using the ultrametric version of the
triangle inequality and the fact that the binomial coefficients are
integers in the first step.  Thus
\begin{equation}
\label{|a_j {j choose l} cdot b_0^{j - l} y^l| = ... le |a_j| r^j}
        \biggl|a_j \, {j \choose l} \cdot b_0^{j - l} \, y^l\biggr|
          = \biggl|a_j \, {j \choose l} \cdot b_0^{j - l}\biggr| \, |y|^l
          \le |a_j| \, r^j
\end{equation}
for every $j \ge l \ge 0$ when $y \in k$ satisfies $|y| \le r$.

        Put
\begin{equation}
\label{widetilde{a}_l = sum_{j = l}^infty a_j {j choose l} cdot b_0^{j - l}}
 \widetilde{a}_l = \sum_{j = l}^\infty a_j \, {j \choose l} \cdot b_0^{j - l}
\end{equation}
for each $l \ge 0$, where the convergence of the series in $k$ follows
from (\ref{lim_{j to infty} |a_j| r^j = 0, 3}) and (\ref{|a_j {j
    choose l} cdot b_0^{j - l}| r^l le ... le |a_j| r^j}).
The ultrametric version of the triangle inequality implies that
\begin{equation}
\label{|widetilde{a}_l| le ... le max_{j ge l} (|a_j| r^{j - l})}
 |\widetilde{a}_l|
 \le \max_{j \ge l} \biggl|a_j \, {j \choose l} \cdot b_0^{j - l}\biggr|
 \le \max_{j \ge l} (|a_j| \, |b_0|^{j - l}) \le \max_{j \ge l} (|a_j| \, r^{j - l})
\end{equation}
for each $l \ge 0$, and hence that
\begin{equation}
\label{|widetilde{a}_l| r^l le max_{j ge l} (|a_j| r^j)}
        |\widetilde{a}_l| \, r^l \le \max_{j \ge l} (|a_j| \, r^j)
\end{equation}
for each $l \ge 0$.  It follows that
\begin{equation}
\label{sum_{l = 0}^infty widetilde{a}_l y^l}
        \sum_{l = 0}^\infty \widetilde{a}_l \, y^l
\end{equation}
coverges in $k$ for every $y \in k$ with $|y| \le r$, because
(\ref{|widetilde{a}_l| r^l le max_{j ge l} (|a_j| r^j)}) tends to $0$
as $l \to \infty$, by (\ref{lim_{j to infty} |a_j| r^j = 0, 3}).  Of
course, (\ref{sum_{l = 0}^infty widetilde{a}_l y^l}) is the same as
the right side of (\ref{f(b_0 + y) = sum_j sum_{l = 0}^j a_j {j choose
    l} b_0^{j - l} y^l = ...}), and one can check that (\ref{f(b_0 +
  y) = sum_j sum_{l = 0}^j a_j {j choose l} b_0^{j - l} y^l = ...})
holds for every $y \in k$ with $|y| \le r$ under these conditions.
More precisely, this uses the fact that (\ref{|a_j {j choose l} cdot
  b_0^{j - l} y^l| = ... le |a_j| r^j}) tends to $0$ as $j \to
\infty$, by (\ref{lim_{j to infty} |a_j| r^j = 0, 3}), in order to
interchange the order of summation in the second step in (\ref{f(b_0 +
  y) = sum_j sum_{l = 0}^j a_j {j choose l} b_0^{j - l} y^l = ...}),
as in Section \ref{sums of sums}.

\section{The residue field}
\label{residue field}

        Let $k$ be a field with an ultrametric absolute value
function $|\cdot|$.  Observe that the closed unit ball
\begin{equation}
\label{overline{B}(0, 1) = {x in k : |x| le 1}}
        \overline{B}(0, 1) = \{x \in k : |x| \le 1\}
\end{equation}
in $k$ is a subring of $k$, and that the open unit ball
\begin{equation}
\label{B(0, 1) = {x in k : |x| < 1}}
        B(0, 1) = \{x \in k : |x| < 1\}
\end{equation}
in $k$ is an ideal in $\overline{B}(0, 1)$.  Thus the quotient
\begin{equation}
\label{overline{B}(0, 1) / B(0, 1)}
        \overline{B}(0, 1) / B(0, 1)
\end{equation}
is defined as a commutative ring, and in fact it is a field, known as
the \emph{residue field}\index{residue field} associated to $|\cdot|$
on $k$.  More precisely, the multiplicative identity element $1$ in
$k$ satisfies $|1| = 1$, so that its image in the quotient is nonzero,
which is the multiplicative identity element in the quotient.  An
element $x$ of $\overline{B}(0, 1)$ has a multiplicative inverse in
$\overline{B}(0, 1)$ exactly when $|x| = 1$, which implies that
nonzero elements of the quotient have multiplicative inverses in the
quotient.

         If $|\cdot|$ is the trivial absolute value function on $k$,
then $\overline{B}(0, 1) = k$, $B(0, 1) = \{0\}$, and hence the
residue field is the same as $k$.  If $k = {\bf Q}_p$ equipped with
the $p$-adic absolute value function for some prime number $p$, then
$\overline{B}(0, 1)$ is the ring ${\bf Z}_p$ of $p$-adic integers,
$B(0, 1) = p \, {\bf Z}_p$, and the residue field is isomorphic to
${\bf Z} / p \, {\bf Z}$, as in Section \ref{p-adic integers}.

        If $|x|$ is an ultrametric absolute value function on any
field $k$, then $|x|^a$ is also an ultrametric absolute value function
on $k$ for every positive real number $a$, as in Section
\ref{quasimetric absolute value functions}.  The open and closed unit
balls in $k$ with respect to $|x|^a$ are the same as for $|x|$ for
each $a > 0$, which implies that the residue field associated to
$|x|^a$ is the same as the residue field associated to $|x|$.

        Let $k$ be any field with an ultrametric absolute value
function $|\cdot|$ again, and let $k_1$ be a subfield of $k$.
The restriction of $|\cdot|$ to $k_1$ is an absolute value function
on $k_1$, and it is easy to see that there is a natural induced
injective homomorphism from the residue field associated to $k_1$
into the residue field associated to $k$.  If $k_1$ is dense in $k$
with respect to the ultrametric corresponding to $|\cdot|$, then
one can check that the induced homomorphism between the residue
fields is surjective.  In particular, the residue field associated
to the completion of a field with an unltrametric absolute value function
is isomorphic to the residue field associated to the original field
in a natural way.

        Suppose that $k$ is a field with characteristic $p$ for
some prime number $p$, and equipped with an ultrametric absolute value
function $|\cdot|$.  Thus $p \cdot 1 = 0$ in $k$, which implies that
the analogous statement holds in the associated residue field, so that
the residue field has characteristic $p$ too.  Alternatively, if $k$
has characteristic $p$, then there is a natural embedding of ${\bf Z}
/ p \, {\bf Z}$ into $k$.  Let $k_1$ be the image of ${\bf Z} / p \,
{\bf Z}$ in $k$ under this embedding, and note that the restriction
of $|\cdot|$ to $k_1$ is trivial, by (\ref{|x|^n = |x^n| = |1| = 1})
in Section \ref{absolute value functions}.  This implies that the
residue field associated to $k_1$ is isomorphic to ${\bf Z} / p \,
{\bf Z}$ as well, which leads to an embedding of ${\bf Z} / p \, {\bf
  Z}$ into the residue field associated to $k$, by the remarks in the
preceding paragraph.

        Let $k_0$ be a field, let $T$ be an indeterminate, and
let $|f(T)|$ be the absolute value function on $k_0((T))$ associated
to some $r \in (0, 1)$, as in Section \ref{formal laurent series}.
The corresponding closed unit in $k_0((T))$ is equal to $k_0[[T]]$,
the open unit ball is equal to $T \, k_0[[T]]$, and the residue
field is isomorphic to $k_0$.  Of course, $k_0$ can also be identified
with a subfield of $k_0((T))$.

        Suppose for the moment that $|\cdot|$ is a nontrivial discrete
ultrametric absolute value function on a field $k$, and that the
associated residue field has exactly $N$ elements for some integer $N
\ge 2$.  As in Section \ref{discrete absolute value functions}, there
is a $\rho_1 \in (0, 1)$ such that the nonzero values of $|\cdot|$ on
$k$ are the same as the integer powers of $\rho_1$.  Thus open balls
in $k$ of radius $1$ are the same as closed balls of radius $\rho_1$,
so that $\overline{B}(0, 1)$ can be expressed as the union of $N$
pairwise-disjoint closed balls of radius $\rho_1$.  Using this, one
can check that each closed ball in $k$ of radius $\rho_1^j$ for some
$j \in {\bf Z}$ can be expressed as the union of $N$ pairwise-disjoint
closed balls of radius $\rho_1^{j + 1}$.  Repeating the process, we
get that each closed ball in $k$ of radius $\rho_1^j$ for some $j \in
{\bf Z}$ can be expressed as the union of $N^l$ pairwise-disjoint
closed balls of radius $\rho_1^{j + l}$ for every $l \in {\bf Z}_+$.
In particular, this implies that bounded subsets of $k$ are totally
bounded.  If $k$ is complete with respect to the ultrametric
associated to $|\cdot|$, then it follows that closed and bounded
subsets of $k$ are compact.

        Let $k$ be any field with an ultrametric absolute value function
$|\cdot|$ again.  If the closed unit ball in $k$ is totally bounded,
then it is easy to see that the associated residue field is finite,
and one can also check that $|\cdot|$ has to be discrete on $k$ in
this case.  More precisely, the residue field is finite exactly when
the closed unit ball can be covered by finitely many open balls of
radius $1$, which can be taken to be centered at points in
$\overline{B}(0, 1)$.  Similarly, if the open unit ball can be covered
by finitely many closed balls of radius less than $1$, which can be
taken to be centered at points in $B(0, 1)$, then one can verify that
$|\cdot|$ is discrete on $k$.  If $k$ is locally compact with respect
to the topology determined by the metric associated to $|\cdot|$, and
if $|\cdot|$ is not the trivial absolute value function on $k$, then
the closed unit ball in $k$ is compact.  In particular, this implies
that the closed unit ball in $k$ is totally bounded.  Remember too
that $k$ is complete with respect to the metric associated to
$|\cdot|$ when $k$ is locally compact.

\chapter{Geometry of mappings}
\label{geometry of mappings}

\section{Differentiation}
\label{differentiation}

        Let $k$ be a field, and let $|\cdot|$ be an absolute value
function on $k$.  Also let $E$ be a subset of $k$, and let $x$ be an
element of $E$ that is a limit point of $E$ with respect to the metric
associated to $|\cdot|$.  Note that any interior point of $E$ is a
limit point of $E$ when $|\cdot|$ is not the trivial absolute value
function on $k$, and that $k$ has no limit points when $|\cdot|$ is
the trivial absolute value function on $k$.  As usual, a $k$-valued
function $f$ on $E$ is said to be
\emph{differentiable}\index{differentiability} at $x$ if the limit of
\begin{equation}
\label{frac{f(y) - f(x)}{y - x}}
        \frac{f(y) - f(x)}{y - x}
\end{equation}
as $y \in E$ approaches $x$ exists in $k$.  In this case, the
\emph{derivative}\index{derivative of a function} $f'(x)$ of $f$ at
$x$ is defined to be the value of this limit of (\ref{frac{f(y) -
    f(x)}{y - x}}).  Equivalently, this means that
\begin{equation}
\label{lim_{y to x atop y in E} frac{f(y) - f(x) - f'(x) (y - x)}{y - x} = 0}
 \lim_{y \to x \atop y \in E} \frac{f(y) - f(x) - f'(x) \, (y - x)}{y - x} = 0.
\end{equation}
In particular, this implies that
\begin{equation}
\label{lim_{y to x atop y in E} (f(y) - f(x) - f'(x) (y - x)) = 0}
        \lim_{y \to x \atop y \in E} (f(y) - f(x) - f'(x) \, (y - x)) = 0,
\end{equation}
and hence that
\begin{equation}
\label{lim_{y to x atop y in E} (f(y) - f(x)) = 0}
        \lim_{y \to x \atop y \in E} (f(y) - f(x)) = 0,
\end{equation}
so that $f$ is continuous at $x$ as a $k$-valued function on $E$.

        Let $g$ be another $k$-valued function on $E$, and suppose
that $f$ and $g$ are both differentiable at $x$.  Under these
conditions, it is easy to see that $f + g$ is differentiable at $x$
too, with
\begin{equation}
\label{(f + g)'(x) = f'(x) + g'(x)}
        (f + g)'(x) = f'(x) + g'(x).
\end{equation}
Similarly, one can check that $f \, g$ is differentiable at $x$, with
\begin{equation}
\label{(f g)'(x) = f'(x) g(x) + f(x) g'(x)}
        (f \, g)'(x) = f'(x) \, g(x) + f(x) \, g'(x),
\end{equation}
as in the classical product rule\index{product rule} for derivatives.
More precisely, we have that
\begin{equation}
\label{frac{f(y) g(y) - f(x) g(x)}{y - x} = ...}
        \frac{f(y) \, g(y) - f(x) \, g(x)}{y - x}
           = \Big(\frac{f(y) - f(x)}{y - x}\Big) \, g(y)
                    + f(x) \, \Big(\frac{g(y) - g(x)}{y - x}\Big)
\end{equation}
for every $y \in E$ with $y \ne x$.  This implies (\ref{(f g)'(x) =
  f'(x) g(x) + f(x) g'(x)}), by taking the limit as $y \to x$ of both
sides of (\ref{frac{f(y) g(y) - f(x) g(x)}{y - x} = ...}), and using
the fact that $g$ is continuous at $x$ as a $k$-valued function on $E$
to deal with the factor of $g(y)$ in the first term on the right side
of (\ref{frac{f(y) g(y) - f(x) g(x)}{y - x} = ...}).  In particular,
if $\alpha \in k$, then $\alpha \, f$ is a $k$-valued function on $E$
that is differentiable at $x$, with
\begin{equation}
\label{(alpha f)'(x) = alpha f'(x)}
        (\alpha f)'(x) = \alpha \, f'(x).
\end{equation}
This corresponds to the case where $g$ is the constant function on $E$
equal to $\alpha$ at every point, so that $g'(x) = 0$, although one
can verify (\ref{(alpha f)'(x) = alpha f'(x)}) more directly too.

        In order to formulate the chain rule\index{chain rule} in this
setting, let $g$ be a $k$-valued function defined on a set $A
\subseteq k$, and let $f$ be a $k$-valued function defined on a set $E
\subseteq k$.  Thus the composition $f \circ g$ is defined on the set
\begin{equation}
\label{A cap g^{-1}(E)}
        A \cap g^{-1}(E).
\end{equation}
Let $x$ be an element of (\ref{A cap g^{-1}(E)}), so that $x \in A$
and $g(x) \in E$.  Suppose also that $x$ is a limit point of (\ref{A
  cap g^{-1}(E)}), which implies in particular that $x$ is a limit
point of $A$.  In the other direction, if $x$ is a limit point of $A$,
$g(x)$ is an element of the interior of $E$, and $g$ is continuous at
$x$ as a $k$-valued function on $A$, then $x$ is a limit point of
(\ref{A cap g^{-1}(E)}) too.  In addition, we would like $g(x)$ to be
a limit point of $E$.  This follows from the condition that $x$ be a
limit point of (\ref{A cap g^{-1}(E)}) when $g$ is continuous at $x$
as a $k$-valued function on (\ref{A cap g^{-1}(E)}) and $g$ is not
constant on the intersection of (\ref{A cap g^{-1}(E)}) with any
neighborhood of $x$ in $k$.  

        Suppose that $g$ is differentiable at $x$ as a $k$-valued
function on $A$, or at least on (\ref{A cap g^{-1}(E)}).  This implies
that $g$ is continuous at $x$, as before, which is relevant for some
of the remarks in the preceding paragraph.  If $f$ is differentiable
at $g(x)$ as a $k$-valued function on $E$, then one can verify that $f
\circ g$ is differentiable at $x$ as a $k$-valued function on (\ref{A
  cap g^{-1}(E)}), with
\begin{equation}
\label{(f circ g)'(x) = f'(g(x)) g'(x)}
        (f \circ g)'(x) = f'(g(x)) \, g'(x).
\end{equation}
More precisely, the differentiability of $f$ at $g(x)$ implies that
\begin{equation}
\label{f(g(y)) - f(g(x)) - f'(g(x)) (g(y) - g(x))}
        f(g(y)) - f(g(x)) - f'(g(x)) \, (g(y) - g(x))
\end{equation}
is small compared to $g(y) - g(x)$ when $g(y)$ is close to $g(x)$, for
$g(y) \in E$ and hence for $y$ in (\ref{A cap g^{-1}(E)}).  Similarly,
the differentiability of $g$ at $x$ means that
\begin{equation}
\label{g(y) - g(x) - g'(x) (y - x)}
        g(y) - g(x) - g'(x) \, (y - x)
\end{equation}
is small compared to $y - x$, for $y$ in (\ref{A cap g^{-1}(E)}) close
to $x$.  Combining these two statements, one can get that
\begin{equation}
\label{f(g(y)) - f(g(x)) - f'(g(x)) g'(x) (x - y)}
        f(g(y)) - f(g(x)) - f'(g(x)) \, g'(x) \, (x - y)
\end{equation}
is small compared to $x - y$, for $y$ in (\ref{A cap g^{-1}(E)}) close
to $x$, as desired.  This uses the differentiability of $g$ at $x$ to
get that $|g(y) - g(x)|$ is bounded by a constant times $|y - x|$ when
$y$ in (\ref{A cap g^{-1}(E)}) is close to $x$, so that (\ref{f(g(y))
  - f(g(x)) - f'(g(x)) (g(y) - g(x))}) is small compared to $y - x$.

\section{Mappings between metric spaces}
\label{mappings between metric spaces}

        Let $(M_1, d_1(x, y))$ and $(M_2, d_2(u, v))$ be metric spaces,
and let $f$ be a mapping from $M_1$ into $M_2$.  Put
\begin{equation}
\label{D_r(f)(x) = r^{-1} sup {d_2(f(x), f(y)) : y in M_1, d_1(x, y) le r}}
  D_r(f)(x) = r^{-1} \, \sup \{d_2(f(x), f(y)) : y \in M_1, \ d_1(x, y) \le r\}
\end{equation}
for each $x \in M_1$ and $r > 0$, where the supremum is defined as a
nonnegative extended real number.\index{D_r(f)(x)@$D_r(f)(x)$}
Similarly, put
\begin{equation}
\label{widetilde{D}_t(f)(x) = sup_{0 < r le t} D_r(f)(x)}
        \widetilde{D}_t(f)(x) = \sup_{0 < r \le t} D_r(f)(x)
\end{equation}
for each $x \in M_1$ and $t > 0$, which is also defined as an extended
real number.\index{D_t~(f)(x)@$\widetilde{D}_t(f)(x)$}  Equivalently,
\begin{equation}
\label{widetilde{D}_t(f)(x) = sup {frac{d_2(f(x), f(y))}{d_1(x, y)} : ...}}
 \widetilde{D}_t(f)(x) = \sup \bigg\{\frac{d_2(f(x), f(y))}{d_1(x, y)} :
                             y \in M_1, \ 0 < d_1(x, y) \le t\bigg\}
\end{equation}
when there is a $y \in M_1$ such that $0 < d_1(x, y) \le t$, and
otherwise $\widetilde{D}_t(f)(x) = 0$.  This can be seen by taking $r
= d_1(x, y)$ in (\ref{widetilde{D}_t(f)(x) = sup_{0 < r le t}
  D_r(f)(x)}). 

        By construction, $\widetilde{D}_t(f)(x)$ increases monotonically
in $t$, and we put\index{D(f)(x)@$D(f)(x)$}
\begin{equation}
\label{D(f)(x) = limsup_{r to 0} D_r(f)(x) = inf_{t > 0} widetilde{D}_t(f)(x)}
 D(f)(x) = \limsup_{r \to 0} D_r(f)(x) = \inf_{t > 0} \widetilde{D}_t(f)(x)
\end{equation}
for each $x \in M_1$.  This is defined as a nonnegative extended real
number as well, which may be considered as the limit of
$\widetilde{D}_t(f)(x)$ as $t \to 0$.  If $x \in M_1$ is a limit point
of $M_1$, then there are $y \in M_1$ as in (\ref{widetilde{D}_t(f)(x)
  = sup {frac{d_2(f(x), f(y))}{d_1(x, y)} : ...}}) for each $t > 0$.
In this case, $D(f)(x)$ may be expressed equivalently by
\begin{equation}
\label{D(f)(x) = limsup_{y to x} frac{d_2(f(x), f(y))}{d_1(x, y)}}
        D(f)(x) = \limsup_{y \to x} \frac{d_2(f(x), f(y))}{d_1(x, y)}.
\end{equation}
Otherwise, if $x$ is an isolated point in $M_1$, then $D(f)(x) = 0$.

        If $D(f)(x) < A$ for some $x \in M_1$ and $A \in {\bf R}$, then
$\widetilde{D}_t(f)(x) < A$ for some $t > 0$, and hence
\begin{equation}
\label{d_2(f(x), f(y)) le A d_1(x, y)}
        d_2(f(x), f(y)) \le A \, d_1(x, y)
\end{equation}
for every $y \in M_1$ with $d_1(x, y) \le t$.  Conversely, if
(\ref{d_2(f(x), f(y)) le A d_1(x, y)}) holds for some $x \in M_1$, $A
\ge 0$, $t > 0$, and every $y \in M_1$ with $d_1(x, y) \le t$, then we
get that
\begin{equation}
\label{D(f)(x) le widetilde{D}_t(f)(x) le A}
        D(f)(x) \le \widetilde{D}_t(f)(x) \le A.
\end{equation}
Thus $D(f)(x)$ may be described as the infimum of the nonnegative real
numbers $A$ for which there is a $t > 0$ such that (\ref{d_2(f(x),
  f(y)) le A d_1(x, y)}) holds for every $y \in M_1$ with $d_1(x, y)
\le t$, at least when there is such an $A$.  Otherwise, if there is no
such $A$, then $D(f)(x) = +\infty$.  Note that $f$ is continuous at
$x$ when $D(f)(x) < \infty$, by (\ref{d_2(f(x), f(y)) le A d_1(x,
  y)}).

        Let $(M_3, d_3(w, z))$ be another metric space, let $f_1$
be a mapping from $M_1$ into $M_2$, and let $f_2$ be a mapping from
$M_2$ into $M_3$.  Thus the composition $f_2 \circ f_1$ is defined
as a mapping from $M_1$ to $M_3$, and we would like to show that
\begin{equation}
\label{D(f_2 circ f_1)(x) le D(f_2)(f_1(x)) D(f_1)(x)}
        D(f_2 \circ f_1)(x) \le D(f_2)(f_1(x)) \, D(f_1)(x)
\end{equation}
for every $x \in M_1$ such that
\begin{equation}
\label{D(f_1)(x), D(f_2)(f_1(x)) < +infty}
        D(f_1)(x), \, D(f_2)(f_1(x)) < +\infty.
\end{equation}
More precisely, $D(f_1)(x)$ is defined in exactly the same way as
before, while $D(f_2 \circ f_1)(x)$ and $D(f_2)(f_1(x))$ are defined
analogously for mappings from $M_1$ and $M_2$ into $M_3$,
respectively.  To do this, let such a point $x \in M_1$ be given, and
suppose that $A_1, A_2 \in {\bf R}$ satisfy
\begin{equation}
\label{D(f_1)(x) < A_1}
        D(f_1)(x) < A_1
\end{equation}
and
\begin{equation}
\label{D(f_2)(f_1(x)) < A_2}
        D(f_2)(f_1(x)) < A_2.
\end{equation}
As in the preceding paragraph, (\ref{D(f_1)(x) < A_1}) implies that
there is a $t_1 > 0$ such that
\begin{equation}
\label{d_2(f_1(x), f_1(y)) le A_1 d_1(x, y)}
        d_2(f_1(x), f_1(y)) \le A_1 \, d_1(x, y)
\end{equation}
for every $y \in M_1$ with $d_1(x, y) \le t_1$.  Similarly,
(\ref{D(f_2)(f_1(x)) < A_2}) implies that there is a $t_2 > 0$ such
that
\begin{equation}
\label{d_3(f_2(f_1(x)), f_2(v)) le A_2 d_2(f_1(x), v)}
        d_3(f_2(f_1(x)), f_2(v)) \le A_2 \, d_2(f_1(x), v)
\end{equation}
for every $v \in M_2$ with $d_2(f_1(x), v) \le t_2$.  Put
\begin{equation}
\label{t_3 = min(t_1, A_1^{-1} t_2)}
        t_3 = \min(t_1, A_1^{-1} \, t_2).
\end{equation}
If $y \in M_1$ satisfies $d_1(x, y) \le t_3$, then (\ref{d_2(f_1(x),
  f_1(y)) le A_1 d_1(x, y)}) implies that
\begin{equation}
\label{d_2(f_1(x), f_1(y)) le A_1 d_1(x, y) le A_1 t_3 le t_2}
        d_2(f_1(x), f_1(y)) \le A_1 \, d_1(x, y) \le A_1 \, t_3 \le t_2,
\end{equation}
so that (\ref{d_3(f_2(f_1(x)), f_2(v)) le A_2 d_2(f_1(x), v)}) holds
with $v = f_1(y)$.  It follows that
\begin{equation}
\label{d_3(f_2(f_1(x)), f_2(f_1(y))) le ... le A_2 A_1 d_1(x, y)}
        d_3(f_2(f_1(x)), f_2(f_1(y))) \le A_2 \, d_2(f_1(x), f_1(y))
                                      \le A_2 \, A_1 \, d_1(x, y)
\end{equation}
for every $y \in M_1$ with $d_1(x, y) \le t_3$, by (\ref{d_2(f_1(x),
  f_1(y)) le A_1 d_1(x, y)}) and (\ref{d_3(f_2(f_1(x)), f_2(v)) le A_2
  d_2(f_1(x), v)}).  This shows that
\begin{equation}
\label{D(f_2 circ f_1)(x) le widetilde{D}_{t_3}(f_2 circ f_1)(x) le A_2 A_1}
 D(f_2 \circ f_1)(x) \le \widetilde{D}_{t_3}(f_2 \circ f_1)(x) \le A_2 \, A_1,
\end{equation}
as in (\ref{D(f)(x) le widetilde{D}_t(f)(x) le A}), where
$\widetilde{D}_{t_3}(f_2 \circ f_1)(x)$ is defined in the same way as
before, but for mappings from $M_1$ into $M_3$.  It is easy to get
(\ref{D(f_2 circ f_1)(x) le D(f_2)(f_1(x)) D(f_1)(x)}) from
(\ref{D(f_2 circ f_1)(x) le widetilde{D}_{t_3}(f_2 circ f_1)(x) le A_2
  A_1}), by taking the infimum over $A_1, A_2 \in {\bf R}$ that
satisfy (\ref{D(f_1)(x) < A_1}) and (\ref{D(f_2)(f_1(x)) < A_2}),
respectively.

\section{$k$-Valued functions}
\label{k-valued functions}

        Let $k$ be a field, and let $|\cdot|$ be an absolute value
function on $k$.  Also let $E$ be a subset of $k$, and let $x$ be an
element of $E$ that is a limit point of $E$ with respect to the metric
associated to $|\cdot|$.  If $f$ is a $k$-valued function on $E$ that
is differentiable at $x$, then
\begin{equation}
\label{lim_{y to x atop y in M_1} frac{|f(y) - f(x)|}{|x - y|} = |f'(x)|}
        \lim_{y \to x \atop y \in M_1} \frac{|f(y) - f(x)|}{|x - y|} = |f'(x)|.
\end{equation}
In particular, this implies that
\begin{equation}
\label{D(f)(x) = |f'(x)|}
        D(f)(x) = |f'(x)|
\end{equation}
in the notation of the previous section, where $M_1 = E$ and $M_2 = k$
are equipped with the metric associated to $|\cdot|$.  In the other
direction, if any $k$-valued function $f$ on $E$ satisfies
\begin{equation}
\label{D(f)(x) = 0}
        D(f)(x) = 0,
\end{equation}
then $f$ is differentiable at $x$, with $f'(x) = 0$.

        Now let $(M_1, d_1(x, y))$ be any metric space again,
and let us take $M_2 = k$, equipped with the metric associated to
$|\cdot|$.  If $f$ is any $k$-valued function on $M_1$ and $a \in k$,
then it is easy to see that
\begin{equation}
\label{D(a f)(x) = |a| D(f)(x)}
        D(a \, f)(x) = |a| \, D(f)(x)
\end{equation}
for every $x \in M_1$.  More precisely, the right side of (\ref{D(a
  f)(x) = |a| D(f)(x)}) should be interpreted as being equal to
$+\infty$ when $D(f)(x) = \infty$ and $a \ne 0$, and the right side of
(\ref{D(a f)(x) = |a| D(f)(x)}) may be interpreted as being $0$ when
$a = 0$ and $D(f)(x) = +\infty$.  If $g$ is another $k$-valued
function on $M_1$, then
\begin{equation}
\label{D(f + g)(x) le D(f)(x) + D(g)(x)}
        D(f + g)(x) \le D(f)(x) + D(g)(x)
\end{equation}
for every $x \in M_1$, with the usual interpretations when $D(f)(x)$
or $D(g)(x)$ is infinite.  Similarly,
\begin{equation}
\label{D(f + g)(x) le max(D(f)(x), D(g)(x))}
        D(f + g)(x) \le \max(D(f)(x), D(g)(x))
\end{equation}
for every $x \in M_1$ when $|\cdot|$ is an ultrametric absolute
value function on $k$.

        Of course,
\begin{equation}
\label{f(y) g(y) - f(x) g(x) = (f(y) - f(x)) g(y) + f(x) (g(y) - g(x))}
 f(y) \, g(y) - f(x) \, g(x) = (f(y) - f(x)) \, g(y) + f(x) \, (g(y) - g(x))
\end{equation}
for every $x, y \in M_1$, which implies that
\begin{equation}
\label{|f(y) g(y) - f(x) g(x)| le |f(y) - f(x)| |g(y)| + |f(x)| |g(y) - g(x)|}
        |f(y) \, g(y) - f(x) \, g(x)| \le |f(y) - f(x)| \, |g(y)|
                                             + |f(x)| \, |g(y) - g(x)|.
\end{equation}
If $D(f)(x), D(g)(x) < \infty$, then one can use (\ref{|f(y) g(y) -
  f(x) g(x)| le |f(y) - f(x)| |g(y)| + |f(x)| |g(y) - g(x)|}) to show that
\begin{equation}
\label{D(f g)(x) le D(f)(x) |g(x)| + |f(x)| D(g)(x)}
        D(f \, g)(x) \le D(f)(x) \, |g(x)| + |f(x)| \, D(g)(x).
\end{equation}
More precisely, if $D(g)(x) < \infty$, then $g$ is continuous at $x$,
which permits one to approximate $|g(y)|$ in the first term on the
right side of (\ref{|f(y) g(y) - f(x) g(x)| le |f(y) - f(x)| |g(y)| +
  |f(x)| |g(y) - g(x)|}) by $|g(x)|$.  Similarly, if $|\cdot|$ is an
ultrametric absolute value function on $k$, then (\ref{f(y) g(y) -
  f(x) g(x) = (f(y) - f(x)) g(y) + f(x) (g(y) - g(x))}) implies that
\begin{equation}
\label{|f(y) g(y) - f(x) g(x)| le ...}
 \qquad |f(y) \, g(y) - f(x) \, g(x)|
            \le \max(|f(y) - f(x)| \, |g(y)|, |f(x)| \, |g(y) - g(x)|)
\end{equation}
for every $x, y \in M_1$.  Using this, one can check that
\begin{equation}
\label{D(f g)(x) le max(D(f)(x) |g(x)|, |f(x)| D(g)(x))}
        D(f \, g)(x) \le \max(D(f)(x) \, |g(x)|, |f(x)| \, D(g)(x))
\end{equation}
when $D(f)(x), D(g)(x) < \infty$.

        Suppose for the moment that $x$ is a limit point of $M_1$,
since otherwise $D$ of any function on $M_1$ is equal to $0$ at $x$,
and (\ref{|f(y) g(y) - f(x) g(x)| le ...}) and (\ref{D(f g)(x) le
  max(D(f)(x) |g(x)|, |f(x)| D(g)(x))}) are trivial.  Thus
\begin{equation}
\label{limsup_{y to x} |g(y)|}
        \limsup_{y \to x} |g(y)|
\end{equation}
is defined, which is equal to $|g(x)|$ when $|g|$ is continuous at
$x$.  If (\ref{limsup_{y to x} |g(y)|}) is less than or equal to
$|g(x)|$, then $|g|$ is said to be upper semicontinuous at $x$.  If
$f(x) = 0$, then the computations in the previous paragraph can be
simplified, and it is easy to see that
\begin{equation}
\label{D(f g)(x) le D(f)(x) limsup_{y to x} |g(y)|}
        D(f \, g)(x) \le D(f)(x) \, \limsup_{y \to x} |g(y)|
\end{equation}
when $D(f)(x)$ and (\ref{limsup_{y to x} |g(y)|}) are finite, even if
$D(g)(x) = +\infty$.  Of course, there is an analogous statement when
$g(x) = 0$.

\section{Lipschitz mappings}
\label{lipschitz mappings}

        Let $(M_1, d_1(x, y))$ and $(M_2, d_2(u, v))$ be metric spaces,
and let $a$ be a positive real number.  A mapping $f : M_1 \to M_2$
is said to be \emph{Lipschitz of order $a$}\index{Lipschitz mappings}
if there is a nonnegative real number $C$ such that
\begin{equation}
\label{d_2(f(x), f(y)) le C d_1(x, y)^a}
        d_2(f(x), f(y)) \le C \, d_1(x, y)^a
\end{equation}
for every $x, y \in M_1$.  Note that $f$ satisfies (\ref{d_2(f(x),
  f(y)) le C d_1(x, y)^a}) with $C = 0$ if and only if $f$ is
constant, and that Lipschitz mappings of any positive order are
uniformly continuous.  One sometimes simply says that $f$ is a
Lipschitz mapping when $f$ is Lipschitz of order $a = 1$.

        If $f : M_1 \to M_2$ is Lipschitz of order $a > 0$ with constant
$C \ge 0$, then
\begin{equation}
\label{D_r(f)(x) le C r^{a - 1}}
        D_r(f)(x) \le C \, r^{a - 1}
\end{equation}
for every $x \in M_1$ and $r > 0$, where $D_r(f)(x)$ is as in
(\ref{D_r(f)(x) = r^{-1} sup {d_2(f(x), f(y)) : y in M_1, d_1(x, y) le
    r}}) in Section \ref{mappings between metric spaces}.  This
implies that
\begin{equation}
\label{widetilde{D}_t(f)(x) le C t^{a - 1}}
        \widetilde{D}_t(f)(x) \le C \, t^{a - 1}
\end{equation}
for every $t > 0$ when $a \ge 1$, where $\widetilde{D}_t(f)(x)$ is as
in (\ref{widetilde{D}_t(f)(x) = sup_{0 < r le t} D_r(f)(x)}).  Of
course, (\ref{widetilde{D}_t(f)(x) le C t^{a - 1}}) can also be
derived from (\ref{widetilde{D}_t(f)(x) = sup {frac{d_2(f(x),
      f(y))}{d_1(x, y)} : ...}}).  It follows that
\begin{equation}
\label{D(f)(x) le C}
        D(f)(x) \le C
\end{equation}
when $a = 1$, and that
\begin{equation}
\label{D(f)(x) = 0, 2}
        D(f)(x) = 0
\end{equation}
when $a > 1$, where $D(f)(x)$ is as in (\ref{D(f)(x) = limsup_{r to 0}
  D_r(f)(x) = inf_{t > 0} widetilde{D}_t(f)(x)}).

        Let $(M_3, d_3(w, z))$ be another metric space, and suppose that
$f_1 : M_1 \to M_2$ is Lipschitz of order $a_1 > 0$ with constant
$C_1 \ge 0$, and that $f_2 : M_2 \to M_3$ is Lipschitz of order $a_2 > 0$
with constant $C_2 \ge 0$.  This implies that
\begin{equation}
\label{d_3(f_2(f_1(x)), f_2(f_1(y))) le... le C_2 C_1^{a_2} d_1(x, y)^{a_1 a_2}}
 \qquad d_3(f_2(f_1(x)), f_2(f_1(y))) \le C_2 \, d_2(f_1(x), f_1(y))^{a_2}
                          \le C_2 \, C_1^{a_2} \, d_1(x, y)^{a_1 \, a_2}
\end{equation}
for every $x, y \in M_1$.  Thus the composition $f_2 \circ f_1$ is
Lipschitz of order $a_1 \, a_2$ as a mapping from $M_1$ into $M_3$,
with constant $C_2 \, C_1^{a_2}$.

        In some situations, we may have a mapping $f : M_1 \to M_2$
that satisfies
\begin{equation}
\label{C^{-1} d_1(x, y)^a le d_2(f(x), f(y)) le C d_1(x, y)^a}
        C^{-1} \, d_1(x, y)^a \le d_2(f(x), f(y)) \le C \, d_1(x, y)^a
\end{equation}
for some $a > 0$ and $C \ge 1$, and every $x, y \in M_1$.  If $a = 1$,
then $f$ is said to be \emph{bilipschitz}\index{bilipschitz mappings}
with constant $C$.  Note that $f$ is bilipschitz with constant $C = 1$
if and only if $f$ is an isometric embedding.  If $a$ is any positive
real number, then (\ref{C^{-1} d_1(x, y)^a le d_2(f(x), f(y)) le C
  d_1(x, y)^a}) is equivalent to saying that $f$ is Lipschitz of order
$a$ with constant $C$, and that the inverse mapping $f^{-1}$ is defined
and Lipschitz of order $1/a$ on $f(M_1)$, with constant $C^{1/a}$.

        Let $k$ be a field with an absolute value function $|\cdot|$,
and let us take $M_2 = k$, with the metric associated to $|\cdot|$.
Also let $f_1$ and $f_2$ be $k$-valued functions on $M_1$ that are
Lipschitz of order $a > 0$ with constants $C_1$ and $C_2$,
respectively.  It is easy to see that $f_1 + f_2$ is Lipschitz of
order $a$ on $M_1$ too, with constant
\begin{equation}
\label{C_1 + C_2}
        C_1 + C_2.
\end{equation}
Similarly, if $|\cdot|$ is an ultrametric absolute value function on
$k$, then $f_1 + f_2$ is Lipschitz of order $a$ on $M_1$ with constant
\begin{equation}
\label{max(C_1, C_2)}
        \max(C_1, C_2).
\end{equation}
If $\alpha \in k$, then $\alpha \, f_1$ is Lipschitz of order $a$ on
$M_1$ with constant
\begin{equation}
\label{|alpha| C_1}
        |\alpha| \, C_1.
\end{equation}
If $f_1$ and $f_2$ are also bounded functions on $M_1$, then one can
check that $f_1 \, f_2$ is Lipschitz of order $a$ on $M_1$, with
constant
\begin{equation}
\label{C_1 (sup_{x in M_1} |f_2(x)|) + (sup_{x in M_1} |f_1(x)|) C_2}
        C_1 \, \Big(\sup_{x \in M_1} |f_2(x)|\Big)
                 + \Big(\sup_{x \in M_1} |f_1(x)|\Big) \, C_2,
\end{equation}
using (\ref{|f(y) g(y) - f(x) g(x)| le |f(y) - f(x)| |g(y)| + |f(x)|
  |g(y) - g(x)|}) in Section \ref{k-valued functions}.  In this case,
if $|\cdot|$ is an ultrametric absolute value function on $k$, then
$f_1 \, f_2$ is Lipschitz of order $a$ on $M_1$ with constant
\begin{equation}
\label{max(C_1 (sup_{x in M_1} |f_2(x)|), (sup_{x in M_1} |f_1(x)|) C_2)}
        \max\Big(C_1 \, \Big(\sup_{x \in M_1} |f_2(x)|\Big),
                  \Big(\sup_{x \in M_1} |f_1(x)|\Big) \, C_2\Big),
\end{equation}
because of (\ref{|f(y) g(y) - f(x) g(x)| le ...}).

        Let us now take $M_2 = {\bf R}$, equipped with the standard metric,
and let $f$ be a real-valued function on $M_1$.  If
\begin{equation}
\label{f(x) le f(y) + C d_1(x, y)^a}
        f(x) \le f(y) + C \, d_1(x, y)^a
\end{equation}
for some $a > 0$ and $C \ge 0$, and for every $x, y \in M_1$, then we
also have that
\begin{equation}
\label{f(y) le f(x) + C d_1(x, y)^a}
        f(y) \le f(x) + C \, d_1(x, y)^a
\end{equation}
for every $x, y \in M_1$, by interchanging the roles of $x$ and $y$.
It follows that
\begin{equation}
\label{|f(x) - f(y)| = max(f(x) - f(y), f(y) - f(x)) le C d_1(x, y)^}
        |f(x) - f(y)| = \max(f(x) - f(y), f(y) - f(x)) \le C \, d_1(x, y)^a
\end{equation}
for every $x, y \in M_1$, so that $f$ is Lipschitz of order $a$ with
constant $C$.  If
\begin{equation}
\label{d_1(x, y)^a}
        d_1(x, y)^a
\end{equation}
is a metric on $M_1$, then
\begin{equation}
\label{f_{p, a}(x) = d_1(x, p)^a}
        f_{p, a}(x) = d_1(x, p)^a
\end{equation}
satisfies (\ref{f(x) le f(y) + C d_1(x, y)^a}) with $C = 1$ for every
$p \in M_1$, because of the triangle inequality.  In particular, this
holds for $0 < a \le 1$ when $d_1(x, y)$ is a metric on $M_1$, and for
every $a > 0$ when $d_1(x, y)$ is an ultrametric on $M_1$, as in
Section \ref{quasimetrics}.

        If $d(x, y)$ is any metric on $M_1$ and $0 < b \le 1$, then
\begin{equation}
\label{d_1(x, y) = d(x, y)^b}
        d_1(x, y) = d(x, y)^b
\end{equation}
is also a metric on $M_1$, as in Section \ref{quasimetrics}.  In this
case, (\ref{d_1(x, y)^a}) is a metric on $M_1$ when $0 < a \le 1/b$,
so that (\ref{f_{p, a}(x) = d_1(x, p)^a}) is Lipschitz of order $a$
with constant $C = 1$ when $0 < a \le 1/b$.  If $b < 1$ and $M_1$ has
at least two elements, then this leads to nonconstant real-valued
functions on $M_1$ that are Lipschitz of order $a > 1$.  Note that a
locally constant mapping $f$ on any metric space $M_1$ satisfies
(\ref{D(f)(x) = 0, 2}) for every $x \in M_1$.  If $M_1$ is not
connected, then there are locally constant mappings on $M_1$ that are
not constant on $M_1$.

\section{Lipschitz mappings, continued}
\label{lipschitz mappings, continued}

        Let $k$ be a field, and let $|\cdot|$ be an absolute value
function on $k$.  If $x, y \in k$ and $j$ is a positive integer, then
\begin{equation}
\label{(x - y) sum_{l = 0}^{j - 1} x^l y^{j - l - 1} = x^j - y^j}
        (x - y) \, \sum_{l = 0}^{j - 1} x^l \, y^{j - l - 1} = x^j - y^j,
\end{equation}
by a standard computation.  This implies that
\begin{equation}
\label{|x^j - y^j| le ... le j |x - y| (max(|x|, |y|))^{j - 1}}
 |x^j - y^j| \le |x - y| \, \sum_{l = 0}^{j - 1} |x|^l \, |y|^{j - l - 1}
              \le j \, |x - y| \, \Big(\max(|x|, |y|)\Big)^{j - 1}.
\end{equation}
If $|\cdot|$ is an ultrametric absolute value function on $k$,
then we get that
\begin{equation}
\label{|x^j - y^j| le ... le |x - y| (max(|x|, |y|))^{j - 1}}
 \qquad |x^j - y^j|
        \le |x - y| \, \max_{0 \le l \le j - 1} (|x|^l \, |y|^{j - l - 1})
         \le |x - y| \, \Big(\max(|x|, |y|)\Big)^{j - 1}.
\end{equation}  

        Let $a_0, a_1, a_2, a_3, \ldots$ be a sequence of elements of $k$,
and suppose for the moment that the series
\begin{equation}
\label{sum_{j = 1}^infty j |a_j| r^{j - 1}}
        \sum_{j = 1}^\infty j \, |a_j| \,  r^{j - 1}
\end{equation}
converges for some positive real number $r$.  This implies that the series
\begin{equation}
\label{sum_{j = 0}^infty |a_j| r^j, 4}
        \sum_{j = 0}^\infty |a_j| \, r^j
\end{equation}
converges, so that the corresponding power series
\begin{equation}
\label{sum_{j = 0}^infty a_j x^j, 2}
        \sum_{j = 0}^\infty a_j \, x^j
\end{equation}
converges absolutely when $|x| \le r$.  If $k$ is complete with
respect to the metric associated to $|\cdot|$, then it follows that
(\ref{sum_{j = 0}^infty a_j x^j, 2}) converges in $k$, and we let
$f(x)$ denote the value of the sum.  Of course, if $a_j = 0$ for all
but finitely many $j$, then the completeness of $k$ is not needed, and
we still let $f(x)$ denote the value of the sum (\ref{sum_{j =
    0}^infty a_j x^j, 2}).  In both cases, if $x, y \in k$ and $|x|,
|y| \le r$, then we get that
\begin{equation}
\label{|f(x) - f(y)| le ... le (sum_{j = 1}^infty j |a_j| r^{j - 1}) |x - y}
 |f(x) - f(y)| \le \sum_{j = 1}^\infty |a_j| \, |x^j - y^j|
         \le \Big(\sum_{j = 1}^\infty j \, |a_j| \, r^{j - 1}\Big) \, |x - y|,
\end{equation}
by (\ref{|x^j - y^j| le ... le j |x - y| (max(|x|, |y|))^{j - 1}}).

        Suppose now that $|\cdot|$ is an ultrametric absolute value
function on $k$, and that the $a_j$'s satisfy
\begin{equation}
\label{lim_{j to infty} |a_j| r^j = 0, 4}
        \lim_{j \to \infty} |a_j| \, r^j = 0
\end{equation}
for some $r > 0$, instead of (\ref{sum_{j = 1}^infty j |a_j| r^{j -
    1}}).  If $k$ is complete, then it follows that the power series
(\ref{sum_{j = 0}^infty a_j x^j, 2}) converges in $k$ for each $x \in
k$ with $|x| \le r$, and we let $f(x)$ denote the value of the sum
again.  As before, the completeness of $k$ is not needed when $a_j =
0$ for all but finitely many $j$.  In both cases, if $x, y \in k$ and
$|x|, |y| \le r$, then we have that
\begin{equation}
\label{|f(x) - f(y)| le ... le max_{j ge 1} (|a_j| r^{j - 1}) |x - y|}
 |f(x) - f(y)| \le \max_{j \ge 1} (|a_j| \, |x^j - y^j|)
                \le \max_{j \ge 1} (|a_j| \, r^{j - 1}) \, |x - y|,
\end{equation}
by (\ref{|x^j - y^j| le ... le |x - y| (max(|x|, |y|))^{j - 1}}).

\section{Differentiation, continued}
\label{differentiation, continued}

        Let $k$ be a field, and let
\begin{equation}
\label{f(X) = sum_{j = 0}^infty a_j X^j, 3}
        f(X) = \sum_{j = 0}^\infty a_j \, X^j
\end{equation}
be a formal power series with coefficients in $k$.  The formal
derivative of $f(X)$ is defined to be the formal power series
\begin{equation}
\label{f'(X) = sum_{j = 1}^infty j cdot a_j X^{j - 1}}
        f'(X) = \sum_{j = 1}^\infty j \cdot a_j \, X^{j - 1}.
\end{equation}
If
\begin{equation}
\label{g(X) = sum_{j = 0}^infty b_j X^j}
        g(X) = \sum_{j = 0}^\infty b_j \, X^j
\end{equation}
is another formal power series with coefficients in $k$, then the sum
$(f + g)(X) = f(X) + g(X)$ is a formal power series too, whose
derivative is given by
\begin{equation}
\label{(f + g)'(X) = f'(X) + g'(X)}
        (f + g)'(X) = f'(X) + g'(X).
\end{equation}
Similarly, if $\alpha \in k$, then $(\alpha \, f)(X) = \alpha \, f(X)$
is a formal power series, whose derivative is given by
\begin{equation}
\label{(alpha f)'(X) = alpha f'(X)}
        (\alpha \, f)'(X) = \alpha \, f'(X).
\end{equation}

        Suppose for the moment that $k = {\bf R}$ or ${\bf C}$,
with the standard absolute value function.  Suppose also that the
series (\ref{sum_{j = 1}^infty j |a_j| r^{j - 1}}) converges for some
positive real number $r$, so that the power series
\begin{equation}
\label{sum_{j = 1}^infty j cdot a_j x^{j - 1}}
        \sum_{j = 1}^\infty j \cdot a_j \, x^{j - 1}
\end{equation}
converges absolutely for every $x \in k$ with $|x| \le r$.  The
convergence of the series (\ref{sum_{j = 1}^infty j |a_j| r^{j - 1}})
implies that the series (\ref{sum_{j = 0}^infty |a_j| r^j, 4})
converges as well, and hence that $f(x)$ can be defined as a
$k$-valued function on the closed ball
\begin{equation}
\label{overline{B}(0, r) = {x in k : |x| le r}, 2}
        \overline{B}(0, r) = \{x \in k : |x| \le r\}
\end{equation}
by the power series (\ref{sum_{j = 0}^infty a_j x^j, 2}).  If $x \in
k$ and $|x| \le r$, then one can check that
\begin{equation}
\label{frac{f(y) - f(x)}{y - x}, 2}
        \frac{f(y) - f(x)}{y - x}
\end{equation}
tends to (\ref{sum_{j = 1}^infty j cdot a_j x^{j - 1}}) as $y \in k$
with $|y| \le r$ tends to $x$.  Thus the derivative $f'(x)$ of $f$ at
$x$ exists and is equal to (\ref{sum_{j = 1}^infty j cdot a_j x^{j -
    1}}), for $f$ as a $k$-valued function defined on $\overline{B}(0,
r)$.  This is elementary when $f$ is a polynomial, and otherwise one
should be careful about interchanging the order of the limits for the
infinite sum.  The main point is that the errors in the relevant
approximations can be estimated as in (\ref{|f(x) - f(y)| le ... le
  (sum_{j = 1}^infty j |a_j| r^{j - 1}) |x - y}).

        If (\ref{f(X) = sum_{j = 0}^infty a_j X^j, 3}) has radius of
convergence $\rho > 0$, then the series (\ref{sum_{j = 0}^infty |a_j|
  r^j, 4}) converges when $r < \rho$.  It is well known that the
series (\ref{sum_{j = 1}^infty j |a_j| r^{j - 1}}) converges when $r <
\rho$ too, which can be derived from the fact that the series
\begin{equation}
\label{sum_{j = 0}^infty |a_j| t^j}
        \sum_{j = 0}^\infty |a_j| \, t^j
\end{equation}
converges when $r < t < \rho$.  Thus $f(x)$ can be defined as a
$k$-valued function on the open ball
\begin{equation}
\label{B(0, rho) = {x in k : |x| < rho}, 2}
        B(0, \rho) = \{x \in k : |x| < \rho\}
\end{equation}
by the power series (\ref{sum_{j = 0}^infty a_j x^j, 2}), and the
power series (\ref{sum_{j = 1}^infty j cdot a_j x^{j - 1}}) converges
absolutely at every point in $B(0, \rho)$ as well.  If $x \in k$ and
$|x| < \rho$, then (\ref{frac{f(y) - f(x)}{y - x}, 2}) tends to
(\ref{sum_{j = 1}^infty j cdot a_j x^{j - 1}}) as $y \in k$ tends to
$x$.  This follows from the analogous statement in the preceding
paragraph applied to $r$ such that $|x| < r < \rho$.  In this case, we
do not need to explicitly restrict our attention to $|y| \le r$, since
this holds automatically when $y$ is sufficiently close to $x$.  As
before, this means that the derivative $f'(x)$ of $f$ at $x$ exists
and is equal to (\ref{sum_{j = 1}^infty j cdot a_j x^{j - 1}}), for
$f$ as a $k$-valued function defined on $B(0, \rho)$.

        Suppose now that $k$ is any field equipped with an ultrametric
absolute value function $|\cdot|$, and that the $a_j$'s satisfy
(\ref{lim_{j to infty} |a_j| r^j = 0, 4}) for some $r > 0$.  This
implies that
\begin{equation}
\label{lim_{j to infty} |j cdot a_j| r^{j - 1} = 0}
        \lim_{j \to \infty} |j \cdot a_j| \, r^{j - 1} = 0,
\end{equation}
since $|j \cdot a_j| \le |a_j|$ for each $j \in {\bf Z}_+$, by the
ultrametric version of the triangle inequality.  If $k$ is complete
with respect to the ultrametric associated to $|\cdot|$, then it
follows that the power series (\ref{sum_{j = 0}^infty a_j x^j, 2}) and
(\ref{sum_{j = 1}^infty j cdot a_j x^{j - 1}}) converge in $k$ for
every $x \in k$ with $|x| \le r$.  As usual, these sums also make
sense when $a_j = 0$ for all but finitely many $j$, even if $k$ is not
complete.  In both cases, $f(x)$ can be defined as a $k$-valued
function on the closed ball $\overline{B}(0, r)$ 
by the power series (\ref{sum_{j = 0}^infty a_j x^j, 2}).
If $|\cdot|$ is not the trivial absolute value function on $k$, then
(\ref{frac{f(y) - f(x)}{y - x}, 2}) tends to (\ref{sum_{j = 1}^infty j
  cdot a_j x^{j - 1}}) as $y \in k$ tends to $x \in k$ with $|x| \le
r$, for essentially the same reasons as before.  Thus the derivative
$f'(x)$ of $f$ at $x$ exists and is given by the power series
(\ref{sum_{j = 1}^infty j cdot a_j x^{j - 1}}) under these conditions.

\section{Derivative $0$}
\label{derivative 0}

        Let $k$ be a field, and let $f(X)$ be a formal power series
with coefficients in $k$, as in (\ref{f(X) = sum_{j = 0}^infty a_j X^j, 3}).
Thus the formal derivative $f'(X)$ is equal to $0$ as a formal power
series if and only if
\begin{equation}
\label{j cdot a_j = 0}
        j \cdot a_j = 0
\end{equation}
for every positive integer $j$.  If $k$ has characteristic $0$, then
this is the same as saying that
\begin{equation}
\label{a_j = 0}
        a_j = 0
\end{equation}
for every $j \ge 1$.  Otherwise, if $k$ has characteristic $p$ for
some prime number $p$, then (\ref{j cdot a_j = 0}) holds for every
$j \ge 1$ if and only if (\ref{a_j = 0}) holds when $j$ is not a
multiple of $p$.

        Let $|\cdot|$ be an absolute value function on $k$, and
suppose that $k$ is complete with respect to the associated metric.
Suppose also that $f$ has positive radius of convergence, and that
$a_j \ne 0$ for some $j \in {\bf Z}_+$.  Let $j_0$ be the smallest
positive integer with this property, and note that the power series
\begin{equation}
\label{sum_{j = j_0}^infty a_j x^{j - j_0}}
        \sum_{j = j_0}^\infty a_j \, x^{j - j_0}
\end{equation}
has the same radius of convergence as $f$.  Thus (\ref{sum_{j =
    j_0}^infty a_j x^{j - j_0}}) is defined and continuous as a
$k$-valued function on a ball in $k$ centered at $0$ with positive
radius.  This implies that (\ref{sum_{j = j_0}^infty a_j x^{j - j_0}})
is not equal to $0$ when $x \in k$ and $|x|$ is sufficiently small,
because (\ref{sum_{j = j_0}^infty a_j x^{j - j_0}}) is equal to
$a_{j_0} \ne 0$ when $x = 0$.  It follows that
\begin{equation}
\label{f(x) = a_0 + x^{j_0} sum_{j = j_0}^infty a_j x^{j - j_0}}
        f(x) = a_0 + x^{j_0} \, \sum_{j = j_0}^\infty a_j \, x^{j - j_0}
\end{equation}
is different from $f(0) = a_0$ when $|x|$ is sufficiently small and $x
\ne 0$.  This shows that $f(x)$ is not constant on any neighborhood of
$0$ when $|\cdot|$ is not the trivial absolute value function on $k$.

        Let $k$ be a field with characteristic $p$ for some prime
number $p$.  It is well known that
\begin{equation}
\label{(x + y)^p = x^p + y^p}
        (x + y)^p = x^p + y^p
\end{equation}
for every $x, y \in k$, by the binomial theorem, because ${p \choose j}$
is divisible by $p$ when $1 \le j < p$.  Similarly,
\begin{equation}
\label{(x - y)^p = x^p - y^p}
        (x - y)^p = x^p - y^p
\end{equation}
for every $x, y \in k$, because $(-1)^p = -1$ automatically when $p$
is odd, and $(-1)^2 = 1 = -1$ when $p = 2$.  If $|\cdot|$ is any
absolute value function on $k$, then it follows that
\begin{equation}
\label{|x^p - y^p| = |x - y|^p}
        |x^p - y^p| = |x - y|^p
\end{equation}
for every $x, y \in k$.  Thus $x \mapsto x^p$ defines a Lipschitz
mapping of order $p$ from $k$ into itself, with respect to the metric
on $k$ associated to the absolute value function.  More precisely,
this corresponds to (\ref{C^{-1} d_1(x, y)^a le d_2(f(x), f(y)) le C
  d_1(x, y)^a}) in Section \ref{lipschitz mappings}, with $a = p$ and
$C = 1$.  Remember that any absolute value function on $k$ is
non-archimedian, as in Section \ref{archimedian property}.

        If $f(X)$ is a formal power series with coefficients in $k$
such that $f'(X) = 0$ as a formal power series, then $f(X)$ can be
expressed as
\begin{equation}
\label{f(X) = g(X^p)}
        f(X) = g(X^p)
\end{equation}
for some other formal power series $g$, by the remarks at the
beginning of the section.  Suppose again that $k$ is complete with
respect to the metric associated to $|\cdot|$, and that $f$ has
positive radius of convergence.  This implies that $g$ has positive
radius of convergence too, which is equal to the $p$th power of the
radius of convergence of $f$.  As in Section \ref{lipschitz mappings,
  continued}, the function corresponding to $g$ satisfies Lipschitz
conditions of order $1$ on closed balls in $k$ centered at $0$ with
suitable radii.  This leads to Lipschitz conditions of order $p$ for
the function corresponding to $f$ on closed balls in $k$ centered at
$0$ with suitable radii, because of (\ref{|x^p - y^p| = |x - y|^p}).
If $g' = 0$ as a formal power series, then one can repeat the process.
Of course, the process can only be repeated finitely many times,
unless $a_j = 0$ for every $j \in {\bf Z}_+$.

\section{Some related estimates}
\label{some related estimates}

        Let $k$ be a field, and let $j$ be an integer with $j \ge 2$.
Observe that
\begin{eqnarray}
\label{x^j - y^j - j cdot (x - y) y^{j - 1} = ...}
 \qquad x^j - y^j - j \cdot (x - y) \, y^{j - 1} 
 & = & (x - y) \sum_{l = 0}^{j - 1} x^l \, y^{j - l - 1} 
                                     - j \cdot (x - y) \, y^{j - 1} \\
 & = & (x - y) \, \sum_{l = 0}^{j - 1} (x^l - y^l) \, y^{j - l - 1}, \nonumber
\end{eqnarray}
for each $x, y \in k$, using (\ref{(x - y) sum_{l = 0}^{j - 1} x^l
  y^{j - l - 1} = x^j - y^j}) in Section \ref{lipschitz mappings,
  continued} in the first step.  If $|\cdot|$ is an absolute value
function on $k$, then we get that
\begin{equation}
\label{|x^j - y^j - j cdot (x - y) y^{j - 1}| le ...}
 |x^j - y^j - j \cdot (x - y) \, y^{j - 1}|
   \le |x - y| \, \sum_{l = 0}^{j - 1} |x^l - y^l| \, |y|^{j - l - 1}
\end{equation}
for every $x, y \in k$.  Remember that
\begin{equation}
\label{|x^l - y^l| le l |x - y| (max(|x|, |y|))^{l - 1}}
        |x^l - y^l| \le l \, |x - y| \, \Big(\max(|x|, |y|)\Big)^{l - 1}
\end{equation}
for every $x, y \in k$ and $l \in {\bf Z}_+$, as in (\ref{|x^j - y^j|
  le ... le j |x - y| (max(|x|, |y|))^{j - 1}}) in Section
\ref{lipschitz mappings, continued}.  Combining (\ref{|x^j - y^j - j
  cdot (x - y) y^{j - 1}| le ...}) and (\ref{|x^l - y^l| le l |x - y|
  (max(|x|, |y|))^{l - 1}}), we obtain that
\begin{eqnarray}
\label{|x^j - y^j - j cdot (x - y) y^{j - 1}| le ..., 2}
 \qquad |x^j - y^j - j \cdot (x - y) \, y^{j - 1}|
 & \le & |x - y|^2 \, \sum_{l = 0}^{j - 1} l \, \Big(\max(|x|, |y|)\Big)^{j - 2} \\
 & = & \frac{j \, (j - 1)}{2} \, |x - y|^2 \Big(\max(|x|, |y|)\Big)^{j - 2}
                                                           \nonumber
\end{eqnarray}
for every $x, y \in k$.  Similarly, if $|\cdot|$ is an ultrametric
absolute value function on $k$, then (\ref{x^j - y^j - j cdot (x - y)
  y^{j - 1} = ...})  implies that
\begin{equation}
\label{|x^j - y^j - j cdot (x - y) y^{j - 1}| le ..., 3}
        |x^j - y^j - j \cdot (x - y) \, y^{j - 1}|
          \le |x - y| \, \max_{0 \le l \le j - 1} (|x^l - y^l| \, |y|^{j - l - 1})
\end{equation}
for every $x, y \in k$.  In this case, we also have that
\begin{equation}
\label{|x^l - y^l| le |x - y| (max(|x|, |y|))^{l - 1}}
        |x^l - y^l| \le |x - y| \, \Big(\max(|x|, |y|)\Big)^{l - 1}
\end{equation}
for every $x, y \in k$ and $l \in {\bf Z}_+$, as in (\ref{|x^j - y^j|
  le ... le |x - y| (max(|x|, |y|))^{j - 1}}) in Section
\ref{lipschitz mappings, continued}.  It follows that
\begin{equation}
\label{|x^j - y^j - j cdot (x - y) y^{j - 1}| le ..., 4}
        |x^j - y^j - j \cdot (x - y) \, y^{j - 1}|
                        \le |x - y|^2 \, \Big(\max(|x|, |y|)\Big)^{j - 2}
\end{equation}
for every $x, y \in k$ under these conditions, by combining (\ref{|x^j
  - y^j - j cdot (x - y) y^{j - 1}| le ..., 3}) and (\ref{|x^l - y^l|
  le |x - y| (max(|x|, |y|))^{l - 1}}).

        Suppose for the moment that $k = {\bf R}$ or ${\bf C}$,
with the standard absolute value function.  Let $a_0, a_1, a_2, a_3,
\ldots$ be a sequence of elements of $k$ such that the series
\begin{equation}
\label{sum_{j = 2}^infty j (j - 1) |a_j| r^{j - 2}}
        \sum_{j = 2}^\infty j \, (j - 1) \, |a_j| \, r^{j - 2}
\end{equation}
converges for some positive real number $r$.  Thus the series
(\ref{sum_{j = 1}^infty j |a_j| r^{j - 1}}) and (\ref{sum_{j =
    0}^infty |a_j| r^j, 4}) in Section \ref{lipschitz mappings,
  continued} converge too, which implies that the power series
(\ref{sum_{j = 0}^infty a_j x^j, 2}) in Section \ref{lipschitz
  mappings, continued} and (\ref{sum_{j = 1}^infty j cdot a_j x^{j -
    1}}) in Section \ref{differentiation, continued} converge
absolutely when $x \in k$ satisfies $|x| \le r$.  Let $f(x)$ be the
$k$-valued function defined on the closed ball $\overline{B}(0, r)$ by
the power series (\ref{sum_{j = 0}^infty a_j x^j, 2}), whose
derivative $f'(x)$ is given by the power series (\ref{sum_{j =
    1}^infty j cdot a_j x^{j - 1}}), as in Section
\ref{differentiation, continued}.  Note that
\begin{equation}
\label{f(x) - f(y) - f'(y) (x - y) = ...}
        f(x) - f(y) - f'(y) \, (x - y)
 = \sum_{j = 2}^\infty a_j(x^j - y^j - j \, (x - y) \, y^{j - 1})
\end{equation}
for every $x, y \in k$ with $|x|, |y| \le r$, where the contributions
from the $j = 0$ and $j = 1$ terms automatically cancel.  Hence
\begin{equation}
\label{|f(x) - f(y) - f'(y) (x - y)| le ...}
        |f(x) - f(y) - f'(y) \, (x - y)|
 \le \sum_{j = 2}^\infty |a_j| \, |x^j - y^j - j \, (x - y) \, y^{j - 1}|
\end{equation}
for every $x, y \in k$ with $|x|, |y| \le r$.  Combining this with
(\ref{|x^j - y^j - j cdot (x - y) y^{j - 1}| le ..., 2}), we get that
\begin{equation}
\label{|f(x) - f(y) - f'(y) (x - y)| le ..., 2}
  \quad |f(x) - f(y) - f'(y) \, (x - y)|
 \le \Big(\sum_{j = 2}^\infty \frac{j \, (j - 1)}{2} \, |a_j| \, r^{j - 2}\Big)
                                                             \, |x - y|^2
\end{equation}
for every $x, y \in k$ with $|x|, |y| \le r$.

        Now let $|\cdot|$ be an ultrametric absolute value function on
any field $k$, and let $a_0, a_1, a_2, a_3, \ldots$ be a sequence of
elements of $k$ that satisfies (\ref{lim_{j to infty} |a_j| r^j = 0,
  4}) in Section \ref{lipschitz mappings, continued} for some $r > 0$.
If $k$ is complete with respect to the ultrametric associated to
$|\cdot|$, then the power series (\ref{sum_{j = 0}^infty a_j x^j, 2})
and (\ref{sum_{j = 1}^infty j cdot a_j x^{j - 1}}) converge in $k$ for
every $x \in k$ with $|x| \le r$, as before.  Let the sums of these
series be denoted $f(x)$ and $f'(x)$, respectively, which can also be
defined when $a_j = 0$ for all but finitely many $j$, even if $k$ is
not complete, as usual.  As in Section \ref{differentiation,
  continued}, $f'(x)$ is the derivative of $f(x)$ as a $k$-valued
function $f(x)$ on the closed ball $\overline{B}(0, r)$ when $|\cdot|$
is nontrivial on $k$.  Otherwise, if $|\cdot|$ is the trivial absolute
value function on $k$, then one can still define $f'(x)$ by the power
series (\ref{sum_{j = 1}^infty j cdot a_j x^{j - 1}}), but the
estimates that follow would be trivial.

        Of course, (\ref{f(x) - f(y) - f'(y) (x - y) = ...}) holds
in this situation too, so that
\begin{equation}
\label{|f(x) - f(y) - f'(y) (x - y)| le ..., 3}
 \qquad |f(x) - f(y) - f'(y) \, (x - y)|
  \le \max_{j \ge 2} (|a_j| \, |x^j - y^j - j \cdot (x - y) \, y^{j - 1}|)
\end{equation}
for every $x, y \in k$ with $|x|, |y| \le r$, by the ultrametric
version of the triangle inequality.  Combining this with (\ref{|x^j -
  y^j - j cdot (x - y) y^{j - 1}| le ..., 4}), we obtain that
\begin{equation}
\label{|f(x) - f(y) - f'(y) (x - y)| le ..., 4}
        |f(x) - f(y) - f'(y) \, (x - y)|
                    \le \max_{j \ge 2} (|a_j| \, r^{j - 2}) \, |x - y|^2
\end{equation}
for every $x, y \in k$ with $|x|, |y| \le r$.  In particular, it
follows that
\begin{eqnarray}
\label{|f(x) - f(y)| le ...}
 \qquad |f(x) - f(y)| & \le & \max(|f(x) - f(y) - f'(y) \, (x - y)|,
                                                  |f'(y)| \, |x - y|) \\
 & \le & \max\Big(\max_{j \ge 2} (|a_j| \, r^{j - 2}) \, |x - y|, |f'(y)|\Big)
                                                    \, |x - y| \nonumber
\end{eqnarray}
for every $x, y \in k$ with $|x|, |y| \le r$.  Note that
\begin{equation}
\label{|f'(x)| le ... le max_{j ge 1} (|a_j| r^{j - 1})}
        |f'(x)| \le \max_{j \ge 1} (|j \cdot a_j| \, r^{j - 1})
                 \le \max_{j \ge 1} (|a_j| \, r^{j - 1})
\end{equation}
for every $x \in k$ with $|x| \le r$, by the definition (\ref{sum_{j =
    1}^infty j cdot a_j x^{j - 1}}) of $f'(x)$ and the ultrametric
version of the triangle inequality.  One can also check that
\begin{equation}
\label{|f'(x) - f'(y)| le ... le max_{j ge 2} (|a_j| r^{j - 2}) |x - y|}
 \qquad |f'(x) - f'(y)|
             \le \max_{j \ge 2} (|j \cdot a_j| \, r^{j - 2}) \, |x - y|
                  \le \max_{j \ge 2} (|a_j| \, r^{j - 2}) \, |x - y|
\end{equation}
for every $x, y \in k$ with $|x|, |y| \le r$, by applying (\ref{|f(x)
  - f(y)| le ... le max_{j ge 1} (|a_j| r^{j - 1}) |x - y|}) in
Section \ref{lipschitz mappings, continued} to $f'$ instead of $f$.

\section{Some related estimates, continued}
\label{some related estimates, continued}

        Let us continue with the same notation and hypotheses as at the
end of the preceding section.  Also let $x_0 \in k$ and $t > 0$ be
given, with $|x_0| \le r$ and $t \le r$, so that
\begin{equation}
\label{overline{B}(x_0, t) subseteq overline{B}(x_0, r) = overline{B}(0, r)}
 \overline{B}(x_0, t) \subseteq \overline{B}(x_0, r) = \overline{B}(0, r),
\end{equation}
by the ultrametric version of the triangle inequality.  If $y \in
\overline{B}(x_0, t)$, then we get that
\begin{equation}
\label{|f'(y) - f'(x_0)| le max_{j ge 2} (|a_j| r^{j - 2}) t}
        |f'(y) - f'(x_0)| \le \max_{j \ge 2} (|a_j| \, r^{j - 2}) \, t,
\end{equation}
by (\ref{|f'(x) - f'(y)| le ... le max_{j ge 2} (|a_j| r^{j - 2}) |x -
  y|}) applied to $x = x_0$.  In particular,
\begin{equation}
\label{|f'(y)| le max(|f'(x_0)|, max_{j ge 2} (|a_j| r^{j - 2}) t)}
 |f'(y)| \le \max\Big(|f'(x_0)|, \max_{j \ge 2} (|a_j| \, r^{j - 2}) \, t\Big)
\end{equation}
when $y \in \overline{B}(x_0, t)$.  Combining this with (\ref{|f(x) -
  f(y)| le ...}), we obtain that
\begin{equation}
\label{|f(x) - f(y)| le ..., 2}
 |f(x) - f(y)| \le \max\Big(|f'(x_0)|,
                     \max_{j \ge 2} (|a_j| \, r^{j - 2}) \, t\Big) \, |x - y|
\end{equation}
for every $x, y \in \overline{B}(x_0, t)$, since $|x - y| \le t$ in
this case too.

        Put
\begin{equation}
\label{g_0(x) = f(x) - f(x_0) - f'(x_0) (x - x_0)}
        g_0(x) = f(x) - f(x_0) - f'(x_0) \, (x - x_0)
\end{equation}
for each $x \in \overline{B}(0, r)$, so that
\begin{equation}
\label{|g_0(x)| le max_{j ge 2} (|a_j| r^{j - 2}) |x - x_0|^2}
        |g_0(x)| \le \max_{j \ge 2} (|a_j| \, r^{j - 2}) \, |x - x_0|^2
\end{equation}
for every $x \in \overline{B}(0, r)$, by (\ref{|f(x) - f(y) - f'(y) (x
  - y)| le ..., 4}).  We also have that
\begin{eqnarray}
\label{g_0(x) - g_0(y) = f(x) - f(y) - f'(x_0) (x - y) = ...}
        g_0(x) - g_0(y) & = & f(x) - f(y) - f'(x_0) \, (x - y) \nonumber \\
  & = & f(x) - f(y) - f'(y) \, (x - y) + (f'(y) - f'(x_0)) \, (x - y)
\end{eqnarray}
for every $x, y \in \overline{B}(0, r)$, and hence
\begin{eqnarray}
\label{|g_0(x) - g_0(y)| le ...}
\lefteqn{|g_0(x) - g_0(y)|} \\
 & \le & \max(|f(x) - f(y) - f'(y) \, (x - y)|, |f'(y) - f'(x_0)| \, |x - y|).
                                                               \nonumber
\end{eqnarray}
It follows that
\begin{equation}
\label{|g_0(x) - g_0(y)| le ..., 2}
 \quad |g_0(x) - g_0(y)| \le \Big(\max_{j \ge 2} (|a_j| \, r^{j - 2})\Big)
                              \, \max(|x - y|, |y - x_0|) \, |x - y|
\end{equation}
for every $x, y \in \overline{B}(0, r)$, by (\ref{|f(x) - f(y) - f'(y)
  (x - y)| le ..., 4}) and (\ref{|f'(x) - f'(y)| le ... le max_{j ge
    2} (|a_j| r^{j - 2}) |x - y|}).  If $x, y \in \overline{B}(x_0,
t)$, so that $|x - y| \le t$ too, then (\ref{|g_0(x) - g_0(y)| le ...,
  2}) implies that that
\begin{equation}
\label{|g_0(x) - g_0(y)| le (max_{j ge 2} (|a_j| r^{j - 2})) t |x - y}
        |g_0(x) - g_0(y)| \le \Big(\max_{j \ge 2} (|a_j| \, r^{j - 2})\Big)
                                                         \, t \, |x - y|.
\end{equation}

        As in Section \ref{changing centers}, we can reexpress
$f(x) = f(x_0 + (x - x_0))$ as a power series in $x - x_0$, which
corresponds to taking $b_0 = x_0$ and $y = x - x_0$ in the earlier
notation.  More precisely, we have that
\begin{equation}
\label{f(x) = sum_{j = 0}^infty a_j (x_0 + (x - x_0))^j = ...}
  \quad f(x) = \sum_{j = 0}^\infty a_j \, (x_0 + (x - x_0))^j 
 = \sum_{j = 0}^\infty \sum_{l = 0}^j a_j {j \choose l} \cdot x_0^{j - l} \, 
                                                    (x - x_0)^l
\end{equation}
for each $x \in k$ with $|x| \le r$, using the binomial theorem in the
third step.  Put
\begin{equation}
\label{widetilde{a}_l = sum_{j = l}^infty a_j {j choose l} cdot x_0^{j - l}}
 \widetilde{a}_l = \sum_{j = l}^\infty a_j \, {j \choose l} \cdot x_0^{j - l}
\end{equation}
for each nonnegative integer $l$, where the series converges in $k$
for the same reasons as before.  We have also seen that
\begin{equation}
\label{|widetilde{a}_l| r^l le max_{j ge l} (|a_j| r^j), 2}
        |\widetilde{a}_l| \, r^l \le \max_{j \ge l} (|a_j| \, r^j)
\end{equation}
for each $l \ge 0$, and that
\begin{equation}
\label{f(x) = sum_{l = 0}^infty widetilde{a}_l (x - x_0)^l}
        f(x) = \sum_{l = 0}^\infty \widetilde{a}_l \, (x - x_0)^l
\end{equation}
for every $x \in k$ with $|x| \le r$.  Of course, $\widetilde{a}_0 =
f(x_0)$ and $\widetilde{a}_1 = f'(x_0)$, so that
\begin{equation}
\label{g_0(x) = sum_{l = 2}^infty widetilde{a}_l (x - x_0)^l}
        g_0(x) = \sum_{l = 2}^\infty \widetilde{a}_l \, (x - x_0)^l
\end{equation}
for every $x \in k$ with $|x| \le r$.

        If $x, y \in \overline{B}(x_0, t)$, then one can check that
\begin{equation}
\label{|f(x) - f(y)| le max_{l ge 1} (|widetilde{a}_l| t^{l - 1}) |x - y|}
 |f(x) - f(y)| \le \max_{l \ge 1} (|\widetilde{a}_l| \, t^{l - 1}) \, |x - y|,
\end{equation}
using the same type of estimate as in (\ref{|f(x) - f(y)| le ... le
  max_{j ge 1} (|a_j| r^{j - 1}) |x - y|}) in Section \ref{lipschitz
  mappings, continued}, applied to the expansion (\ref{f(x) = sum_{l =
    0}^infty widetilde{a}_l (x - x_0)^l}).  Equivalently,
\begin{equation}
\label{|f(x) - f(y)| le ..., 3}
 |f(x) - f(y)| \le \max\Big(|f'(x_0)|,
                 \max_{l \ge 2} (|\widetilde{a}_l| \, t^{l - 1})\Big) \, |x - y|
\end{equation}
for every $x, y \in \overline{B}(x_0, t)$, since $\widetilde{a}_1 =
f'(x_0)$.  We also have that
\begin{equation}
\label{max_{l ge 2} (|widetilde{a}_l| t^{l - 1}) le ...}
        \max_{l \ge 2} (|\widetilde{a}_l| \, t^{l - 1})
          \le \max_{l \ge 2} (|\widetilde{a}_l| \, r^{l - 2}) \, t
           \le \max_{j \ge 2} (|a_j| \, r^{j - 2}) \, t,
\end{equation}
using the fact that $t \le r$ in the first step, and
(\ref{|widetilde{a}_l| r^l le max_{j ge l} (|a_j| r^j), 2}) in the
second step.  This gives another way to look at (\ref{|f(x) - f(y)| le
  ..., 2}), by combining (\ref{|f(x) - f(y)| le ..., 3}) and
(\ref{max_{l ge 2} (|widetilde{a}_l| t^{l - 1}) le ...}).  Similarly,
one can verify that
\begin{equation}
\label{|g_0(x) - g_0(y)| le max_{l ge 2} (|widetilde{a}_l| t^{l - 1}) |x - y|}
 |g_0(x) - g_0(y)| \le \max_{l \ge 2} (|\widetilde{a}_l| \, t^{l - 1}) \, |x - y|
\end{equation}
for every $x, y \in \overline{B}(x_0, t)$.  As before, this uses the
same type of estimate as in (\ref{|f(x) - f(y)| le ... le max_{j ge 1}
  (|a_j| r^{j - 1}) |x - y|}) in Section \ref{lipschitz mappings,
  continued}, applied to the expansion (\ref{g_0(x) = sum_{l =
    2}^infty widetilde{a}_l (x - x_0)^l}).  This gives another way to
look at (\ref{|g_0(x) - g_0(y)| le (max_{j ge 2} (|a_j| r^{j - 2})) t
  |x - y}), by combining (\ref{|g_0(x) - g_0(y)| le max_{l ge 2}
  (|widetilde{a}_l| t^{l - 1}) |x - y|}) and (\ref{max_{l ge 2}
  (|widetilde{a}_l| t^{l - 1}) le ...}).

\section{Hensel's lemma}
\label{hensel's lemma}

        Let us continue with the same notation and hypotheses as in
the preceding section again.  Suppose for the moment that $x, y \in k$
satisfy $|x|, |y| \le r$ and
\begin{equation}
\label{max_{j ge 2} (|a_j| r^{j - 2}) |x - y| < |f'(y)|}
        \max_{j \ge 2} (|a_j| \, r^{j - 2}) \, |x - y| < |f'(y)|.
\end{equation}
This implies that
\begin{equation}
\label{|f'(x) - f'(y)| < |f'(y)|}
        |f'(x) - f'(y)| < |f'(y)|,
\end{equation}
by (\ref{|f'(x) - f'(y)| le ... le max_{j ge 2} (|a_j| r^{j - 2}) |x -
  y|}) in Section \ref{some related estimates}.  It follows that
\begin{equation}
\label{|f'(x)| = |f'(y)|}
        |f'(x)| = |f'(y)|,
\end{equation}
as in (\ref{|x + y| = |y|}) in Section \ref{absolute value functions}.
Similarly, let us check that
\begin{equation}
\label{|f(x) - f(y)| = |f'(y)| |x - y|}
        |f(x) - f(y)| = |f'(y)| \, |x - y|
\end{equation}
when $x$, $y$ satisfy (\ref{max_{j ge 2} (|a_j| r^{j - 2}) |x - y| <
  |f'(y)|}).  Of course, (\ref{|f(x) - f(y)| = |f'(y)| |x - y|}) is
trivial when $x = y$, and so we may suppose that $x \ne y$.  In this
case, we can multiply (\ref{max_{j ge 2} (|a_j| r^{j - 2}) |x - y| <
  |f'(y)|}) by $|x - y|$, to get that
\begin{equation}
\label{max_{j ge 2} (|a_j| r^{j - 2}) |x - y|^2 < |f'(y)| |x - y|}
        \max_{j \ge 2} (|a_j| \, r^{j - 2}) \, |x - y|^2 < |f'(y)| \, |x - y|.
\end{equation}
Combining this with (\ref{|f(x) - f(y) - f'(y) (x - y)| le ..., 4}) in
Section \ref{some related estimates}, we obtain that
\begin{equation}
\label{|f(x) - f(y) - f'(y) (x - y)| < |f'(y)| |x - y|}
        |f(x) - f(y) - f'(y) \, (x - y)| < |f'(y)| \, |x - y|.
\end{equation}
This implies (\ref{|f(x) - f(y)| = |f'(y)| |x - y|}), as in (\ref{|x +
  y| = |y|}) in Section \ref{absolute value functions} again.

        Let $x_0 \in k$ and $t > 0$ be as in the previous section,
and suppose from now on in this section that
\begin{equation}
\label{t max_{j ge 2} (|a_j| r^{j - 2}) < |f'(x_0)|}
        t \, \max_{j \ge 2} (|a_j| \, r^{j - 2}) < |f'(x_0)|.
\end{equation}
This implies that every $x \in \overline{B}(x_0, t)$ satisfies
(\ref{max_{j ge 2} (|a_j| r^{j - 2}) |x - y| < |f'(y)|}) with $y = x_0$.
It follows that
\begin{equation}
\label{|f'(x)| = |f'(x_0)|}
        |f'(x)| = |f'(x_0)|
\end{equation}
for every $x \in \overline{B}(x_0, t)$, by (\ref{|f'(x)| = |f'(y)|})
applied to $y = x_0$.  Of course, this is the same as saying that
\begin{equation}
\label{|f'(y)| = |f'(x_0)|}
        |f'(y)| = |f'(x_0)|
\end{equation}
for every $y \in \overline{B}(x_0, t)$.  If $x, y \in
\overline{B}(x_0, t)$, then $|x - y| \le t$, and (\ref{t max_{j ge 2}
  (|a_j| r^{j - 2}) < |f'(x_0)|}) implies that (\ref{max_{j ge 2}
  (|a_j| r^{j - 2}) |x - y| < |f'(y)|}) holds, because of
(\ref{|f'(y)| = |f'(x_0)|}).  Thus (\ref{|f(x) - f(y)| = |f'(y)| |x -
  y|}) implies that
\begin{equation}
\label{|f(x) - f(y)| = |f'(x_0)| |x - y|}
        |f(x) - f(y)| = |f'(x_0)| \, |x - y|,
\end{equation}
using (\ref{|f'(y)| = |f'(x_0)|}) again.

        It follows that
\begin{equation}
\label{f(overline{B}(x_0, t)) subseteq overline{B}(f(x_0), |f'(x_0)| t)}
        f(\overline{B}(x_0, t)) \subseteq \overline{B}(f(x_0), |f'(x_0)| \, t)
\end{equation}
under these conditions, and in fact we have that
\begin{equation}
\label{f(overline{B}(x_0, t)) = overline{B}(f(x_0), |f'(x_0)| t)}
        f(\overline{B}(x_0, t)) = \overline{B}(f(x_0), |f'(x_0)| \, t)
\end{equation}
when $k$ is complete with respect to the ultrametric associated to
$|\cdot|$.  This is basically \emph{Hensel's lemma},\index{Hensel's
  lemma} as in \cite{c, fg}.  To prove (\ref{f(overline{B}(x_0, t)) =
  overline{B}(f(x_0), |f'(x_0)| t)}), we shall use the contraction
mapping theorem.  Note that the completeness of $k$ is important here,
even when $a_j = 0$ for all but finitely many $j$.

        Let $z$ be an element of the right side of
(\ref{f(overline{B}(x_0, t)) = overline{B}(f(x_0), |f'(x_0)| t)}),
so that $z \in k$ satisfies
\begin{equation}
\label{|f(x_0) - z| le |f'(x_0)| t}
        |f(x_0) - z| \le |f'(x_0)| \, t.
\end{equation}
Put
\begin{equation}
\label{h_z(x) = x_0 + f'(x_0)^{-1}(z - f(x_0)) - f'(x_0)^{-1} g_0(x)}
        h_z(x) = x_0 + f'(x_0)^{-1}(z - f(x_0)) - f'(x_0)^{-1} \, g_0(x)
\end{equation}
for each $x \in \overline{B}(x_0, t)$, where $g_0(x)$ is as in
(\ref{g_0(x) = f(x) - f(x_0) - f'(x_0) (x - x_0)}).  Thus
\begin{equation}
\label{f'(x_0) (x - h_z(x)) = ... = f(x) - z}
 \qquad f'(x_0) \, (x - h_z(x)) = f'(x_0) \, (x - x_0) - z + f(x_0) + g_0(x)
                                = f(x) - z
\end{equation}
for every $x \in \overline{B}(x_0, t)$, by the definition (\ref{g_0(x)
  = f(x) - f(x_0) - f'(x_0) (x - x_0)}) of $g_0(x)$.  It is easy to
see that $h_z$ is Lipschitz of order $1$ with constant
\begin{equation}
\label{|f'(x_0)|^{-1} (max_{j ge 2} (|a_j| r^{j - 2})) t}
        |f'(x_0)|^{-1} \, \Big(\max_{j \ge 2} (|a_j| \, r^{j - 2})\Big) \, t
\end{equation}
on $\overline{B}(x_0, t)$, by (\ref{|g_0(x) - g_0(y)| le (max_{j ge 2}
  (|a_j| r^{j - 2})) t |x - y}).  The hypothesis (\ref{t max_{j ge 2}
  (|a_j| r^{j - 2}) < |f'(x_0)|}) says exactly that
(\ref{|f'(x_0)|^{-1} (max_{j ge 2} (|a_j| r^{j - 2})) t}) is strictly
less than $1$.  We also have that
\begin{equation}
\label{|h_z(x) - x_0| le |f'(x_0)|^{-1} max(|z - f(x_0)|, |g_0(x)|)}
        |h_z(x) - x_0| \le |f'(x_0)|^{-1} \, \max(|z - f(x_0)|, |g_0(x)|)
\end{equation}
for every $x \in \overline{B}(x_0, t)$, by the definition (\ref{h_z(x)
  = x_0 + f'(x_0)^{-1}(z - f(x_0)) - f'(x_0)^{-1} g_0(x)}) of
$h_z(x)$.  Observe that
\begin{equation}
\label{|g_0(x)| le max_{j ge 2} (|a_j| r^{j - 2}) t^2 < t |f'(x_0)|}
 |g_0(x)| \le \max_{j \ge 2} (|a_j| \, r^{j - 2}) \, t^2 < t \, |f'(x_0)|
\end{equation}
for every $x \in \overline{B}(x_0, t)$, using (\ref{|g_0(x)| le max_{j
    ge 2} (|a_j| r^{j - 2}) |x - x_0|^2}) in the first step, and
(\ref{t max_{j ge 2} (|a_j| r^{j - 2}) < |f'(x_0)|}) in the second
step.  This implies that
\begin{equation}
\label{|h_z(x) - x_0| le t}
        |h_z(x) - x_0| \le t
\end{equation}
for every $x \in \overline{B}(x_0, t)$, by (\ref{|f(x_0) - z| le
  |f'(x_0)| t}) and (\ref{|h_z(x) - x_0| le |f'(x_0)|^{-1} max(|z -
  f(x_0)|, |g_0(x)|)}).  Equivalently,
\begin{equation}
\label{h_z(overline{B}(x_0, t)) subseteq overline{B}(x_0, t)}
        h_z(\overline{B}(x_0, t)) \subseteq \overline{B}(x_0, t).
\end{equation}
If $k$ is complete with respect to the metric associated to $|\cdot|$,
then $\overline{B}(x_0, t)$ is also complete as a metric space with
respect to the restriction of this metric to $\overline{B}(x_0, t)$,
because $\overline{B}(x_0, t)$ is a closed set in $k$.  This permits
us to apply the contraction mapping theorem to $h_z$ on
$\overline{B}(x_0, t)$, to get that $h_z$ has a fixed point in
$\overline{B}(x_0, t)$.  This is the same as saying that $f(x) = z$
for some $x \in \overline{B}(x_0, t)$, by (\ref{f'(x_0) (x - h_z(x)) =
  ... = f(x) - z}), which implies (\ref{f(overline{B}(x_0, t)) =
  overline{B}(f(x_0), |f'(x_0)| t)}).

\section{Some variants}
\label{some variants}

        Let us go back to the same notation and hypotheses as at
the beginning of Section \ref{some related estimates, continued}.
Suppose for the moment that $x, y \in k$ satisfy $|x|, |y| \le r$ and
\begin{equation}
\label{max_{j ge 2} (|a_j| r^{j - 2}) |x - y| le |f'(y)|}
        \max_{j \ge 2} (|a_j| \, r^{j - 2}) \, |x - y| \le |f'(y)|,
\end{equation}
instead of (\ref{max_{j ge 2} (|a_j| r^{j - 2}) |x - y| < |f'(y)|}).
This implies that
\begin{equation}
\label{|f'(x) - f'(y)| le |f'(y)|}
        |f'(x) - f'(y)| \le |f'(y)|,
\end{equation}
by (\ref{|f'(x) - f'(y)| le ... le max_{j ge 2} (|a_j| r^{j - 2}) |x -
  y|}) in Section \ref{some related estimates}, and hence
\begin{equation}
\label{|f'(x)| le |f'(y)|}
        |f'(x)| \le |f'(y)|,
\end{equation}
by the ultrametric version of the triangle inequality.  Similarly,
(\ref{max_{j ge 2} (|a_j| r^{j - 2}) |x - y| le |f'(y)|}) implies that
\begin{equation}
\label{|f(x) - f(y)| le |f'(y)| |x - y|}
        |f(x) - f(y)| \le |f'(y)| \, |x - y|,
\end{equation}
by (\ref{|f(x) - f(y)| le ...}) in Section \ref{some related estimates}.

        Let $x_0 \in k$ be given, with $|x_0| \le r$, as in Section
\ref{some related estimates, continued}.  Also let $t_0 > 0$ be given,
with $t_0 \le r$, so that
\begin{equation}
\label{overline{B}(x_0, t_0) subseteq overline{B}(x_0, r) = overline{B}(0, r)}
 \overline{B}(x_0, t_0) \subseteq \overline{B}(x_0, r) = \overline{B}(0, r),
\end{equation}
as before.  Let us suppose from now on in this section that
\begin{equation}
\label{t_0 max_{j ge 2} (|a_j| r^{j - 2}) le |f'(x_0)|}
        t_0 \, \max_{j \ge 2} (|a_j| \, r^{j - 2}) \le |f'(x_0)|,
\end{equation}
instead of (\ref{t max_{j ge 2} (|a_j| r^{j - 2}) < |f'(x_0)|}).  Of
course, if $f'(x_0) = 0$, then (\ref{t_0 max_{j ge 2} (|a_j| r^{j -
    2}) le |f'(x_0)|}) implies that $a_j = 0$ for every $j \ge 2$,
which implies in turn that $a_1 = f'(x_0) = 0$ too.  Thus we also ask that
\begin{equation}
\label{f'(x_0) ne 0}
        f'(x_0) \ne 0.
\end{equation}
If $x \in \overline{B}(x_0, t_0)$, then (\ref{t_0 max_{j ge 2} (|a_j|
  r^{j - 2}) le |f'(x_0)|}) implies that $x$ satisfies (\ref{max_{j ge
    2} (|a_j| r^{j - 2}) |x - y| le |f'(y)|}) with $y = x_0$, so that
\begin{equation}
\label{|f'(x)| le |f'(x_0)|}
        |f'(x)| \le |f'(x_0)|,
\end{equation}
by (\ref{|f'(x)| le |f'(y)|}).  If $x, y \in \overline{B}(x_0, t_0)$,
then $|x - y| \le t_0$ too, and we get that
\begin{equation}
\label{|f(x) - f(y)| le |f'(x_0)| |x - y|}
        |f(x) - f(y)| \le |f'(x_0)| \, |x - y|,
\end{equation}
by (\ref{|f(x) - f(y)| le ...}) in Section \ref{some related
  estimates}, (\ref{t_0 max_{j ge 2} (|a_j| r^{j - 2}) le |f'(x_0)|}),
and (\ref{|f'(x)| le |f'(x_0)|}) applied to $y$ instead of $x$.

        However, if $x \in B(x_0, t_0)$, then (\ref{t_0 max_{j ge 2} 
(|a_j| r^{j - 2}) le |f'(x_0)|}) implies that $x$ satisfies
(\ref{max_{j ge 2} (|a_j| r^{j - 2}) |x - y| < |f'(y)|}) in the
previous section with $y = x_0$.  It follows that (\ref{|f'(x)| =
  |f'(x_0)|}) holds for every $x$ in $B(x_0, t_0)$, as before.
Similarly, if $x, y \in B(x_0, t_0)$, then (\ref{|f(x) - f(y)| =
  |f'(x_0)| |x - y|}) holds, for the same reasons as before.  If $0 <
t < t_0$, then $t \le r$ in particular, as in Section \ref{some
  related estimates, continued}, and $t$ also satisfies (\ref{t max_{j
    ge 2} (|a_j| r^{j - 2}) < |f'(x_0)|}) in the previous section,
because of (\ref{t_0 max_{j ge 2} (|a_j| r^{j - 2}) le |f'(x_0)|}).
This implies that (\ref{f(overline{B}(x_0, t)) = overline{B}(f(x_0),
  |f'(x_0)| t)}) holds for every $t \in (0, t_0)$ when $k$ is
complete, as before, and hence that
\begin{equation}
\label{f(B(x_0, t_0)) = B(f(x_0), |f'(x_0)| t_0)}
        f(B(x_0, t_0)) = B(f(x_0), |f'(x_0)| \, t_0).
\end{equation}

        Now let $r_1 > 0$ be given, and suppose that our sequence
$a_0, a_1, a_2, a_3, \ldots$ of coefficients has the property that
\begin{equation}
\label{|a_j| r_1^j}
        |a_j| \, r_1^j
\end{equation}
is bounded.  This implies that the convergence condition (\ref{lim_{j
    to infty} |a_j| r^j = 0, 4}) in Section \ref{lipschitz mappings,
  continued} holds for every $r > 0$ with $r < r_1$.  Thus $f(x)$ is
defined by the power series (\ref{sum_{j = 0}^infty a_j x^j, 2}) in
Section \ref{lipschitz mappings} for every $x \in k$ with $|x| < r_1$,
and similarly $f'(x)$ is defined by the power series (\ref{sum_{j =
    1}^infty j cdot a_j x^{j - 1}}) in Section \ref{differentiation,
  continued} for every $x \in k$ with $|x| < r_1$.  Also let $x_0 \in
k$ be given, with $|x_0| < r_1$, and suppose in addition that
\begin{equation}
\label{r_1 max_{j ge 2} (|a_j| r_1^{j - 2}) le ... le |f'(x_0)|}
        r_1 \, \max_{j \ge 2} (|a_j| \, r_1^{j - 2})
              = \max_{j \ge 2} (|a_j| \, r_1^{j - 1}) \le |f'(x_0)|.
\end{equation}
As before, we ask that $f'(x_0) \ne 0$ too, since otherwise $a_j = 0$
for every $j \ge 1$.  Under these conditions, we can apply the
discussion in the previous section to $r = t > 0$ such that $r < r_1$
and $|x_0| \le r$, in which case (\ref{t max_{j ge 2} (|a_j| r^{j -
    2}) < |f'(x_0)|}) follows from (\ref{r_1 max_{j ge 2} (|a_j|
  r_1^{j - 2}) le ... le |f'(x_0)|}).  This implies that (\ref{|f'(x)|
  = |f'(x_0)|}) and (\ref{|f(x) - f(y)| = |f'(x_0)| |x - y|}) hold for
every $x$, $y$ in
\begin{equation}
\label{overline{B}(x_0, t) = overline{B}(x_0, r) = overline{B}(0, r)}
        \overline{B}(x_0, t) = \overline{B}(x_0, r) = \overline{B}(0, r).
\end{equation}
In particular, (\ref{|f'(x)| = |f'(x_0)|}) holds when $x = 0$, which
means that we could have taken $x_0 = 0$.  It follows that
(\ref{|f'(x)| = |f'(x_0)|}) and (\ref{|f(x) - f(y)| = |f'(x_0)| |x -
  y|}) hold for every $x, y \in B(0, r_1)$, by taking $r$ close to
$r_1$.  If $k$ is complete, then
\begin{equation}
\label{f(B(x_0, r_1)) = B(f(x_0), |f'(x_0)| r_1}
        f(B(x_0, r_1)) = B(f(x_0), |f'(x_0)| \, r_1),
\end{equation}
by (\ref{f(overline{B}(x_0, t)) = overline{B}(f(x_0), |f'(x_0)| t)})
with $t = r$ close to $r_1$ again.

\section{Some variants, continued}
\label{some variants, continued}

        Let us return to the same notation and hypotheses as at
the beginning of Section \ref{some related estimates, continued}
again.  In particular, the sequence of coefficients $a_0, a_1, a_2,
a_3, \ldots$ is supposed to satisfy the convergence condition
(\ref{lim_{j to infty} |a_j| r^j = 0, 4}) in Section \ref{lipschitz
  mappings, continued} for the given $r > 0$, which means that the
analogous condition holds for any smaller value of $r$ too.  Of
course, the conditions (\ref{t max_{j ge 2} (|a_j| r^{j - 2}) <
  |f'(x_0)|}) and (\ref{t_0 max_{j ge 2} (|a_j| r^{j - 2}) le
  |f'(x_0)|}) in Sections \ref{hensel's lemma} and \ref{some
  variants}, respectively, are more easily satisfied with smaller
values of $r$.  However, one is also supposed to have $|x_0| \le r$,
and either $t \le r$, as in Sections \ref{some related estimates,
  continued} and \ref{hensel's lemma}, or $|t_0| \le r$, as in Section
\ref{some variants}.

        The condition that $|x_0| \le r$ is trivial when $x_0 = 0$,
and one can basically reduce to that case by expanding $f(x)$ into a
power series centered at $x_0$, as in (\ref{f(x) = sum_{l = 0}^infty
  widetilde{a}_l (x - x_0)^l}) in Section \ref{some related estimates,
  continued}.  Thus one can get the same conclusions as in Section
\ref{hensel's lemma} when $0 < t \le r$ and
\begin{equation}
\label{t max_{l ge 2} (|widetilde{a}_l| r^{l - 2}) < |f'(x_0)|}
        t \, \max_{l \ge 2} (|\widetilde{a}_l| \, r^{l - 2}) < |f'(x_0)|,
\end{equation}
where $\widetilde{a}_l$ is as in (\ref{widetilde{a}_l = sum_{j =
    l}^infty a_j {j choose l} cdot x_0^{j - l}}) in Section \ref{some
  related estimates, continued}.  More precisely, (\ref{t max_{l ge 2}
  (|widetilde{a}_l| r^{l - 2}) < |f'(x_0)|}) replaces (\ref{t max_{j
    ge 2} (|a_j| r^{j - 2}) < |f'(x_0)|}) in Section \ref{hensel's
  lemma}, and it is easy to see that (\ref{t max_{j ge 2} (|a_j| r^{j
    - 2}) < |f'(x_0)|}) implies (\ref{t max_{l ge 2} (|widetilde{a}_l|
  r^{l - 2}) < |f'(x_0)|}), because of (\ref{|widetilde{a}_l| r^l le
  max_{j ge l} (|a_j| r^j), 2}) in Section \ref{some related
  estimates, continued}.  If $0 < t \le r$, then one might as well
replace $r$ with $t$, as in the previous paragraph, so that (\ref{t
  max_{l ge 2} (|widetilde{a}_l| r^{l - 2}) < |f'(x_0)|}) becomes
\begin{equation}
\label{t max_{l ge 2} (|widetilde{a}_l| t^{l - 2}) = ... < |f'(x_0)}
        t \, \max_{l \ge 2} (|\widetilde{a}_l| \, t^{l - 2}) 
            = \max_{l \ge 2} (|\widetilde{a}_l| \, t^{l - 1}) < |f'(x_0)|.
\end{equation}

        Similarly, one can get the same conclusions as in the
preceding section when $0 < t_0 \le r$, $f'(x_0) \ne 0$, and
\begin{equation}
\label{t_0 max_{l ge 2} (|widetilde{a}_l| r^{l - 2}) le |f'(x_0)}
        t_0 \, \max_{l \ge 2} (|\widetilde{a}_l| \, r^{l - 2}) \le |f'(x_0)|,
\end{equation}
instead of (\ref{t_0 max_{j ge 2} (|a_j| r^{j - 2}) le |f'(x_0)|}).
As in the previous paragraph, (\ref{t_0 max_{j ge 2} (|a_j| r^{j - 2})
  le |f'(x_0)|}) implies (\ref{t_0 max_{l ge 2} (|widetilde{a}_l| r^{l
    - 2}) le |f'(x_0)}), because of (\ref{|widetilde{a}_l| r^l le
  max_{j ge l} (|a_j| r^j), 2}) in Section \ref{some related
  estimates, continued}.  If $0 < t_0 \le r$ and $f'(x_0) \ne 0$, then
one might as well replace $r$ with $t_0$, as before, so that (\ref{t_0
  max_{l ge 2} (|widetilde{a}_l| r^{l - 2}) le |f'(x_0)}) becomes
\begin{equation}
\label{t_0 max_{l ge 2} (|widetilde{a}_l| t_0^{l - 2}) = ... le |f'(x_0)|}
        t_0 \, \max_{l \ge 2} (|\widetilde{a}_l| \, t_0^{l - 2})
              = \max_{l \ge 2} (|\widetilde{a}_l| \, t_0^{l - 1}) \le |f'(x_0)|.
\end{equation}

\section{A basic situation}
\label{a basic situation}

        Let $k$ be a field, and let $|\cdot|$ be an ultrametric
absolute value function on $k$.  Also let $a_0, a_1, a_2, a_3, \ldots$
be a sequence of elements of $k$, and suppose that
\begin{equation}
\label{|a_j| le 1}
        |a_j| \le 1
\end{equation}
for each $j \ge 0$, and that
\begin{equation}
\label{lim_{j to infty} |a_j| = 0}
        \lim_{j \to \infty} |a_j| = 0.
\end{equation}
As before, we would like to put
\begin{equation}
\label{f(x) = sum_{j = 0}^infty a_j x^j, 3}
        f(x) = \sum_{j = 0}^\infty a_j \, x^j
\end{equation}
and
\begin{equation}
\label{f'(x) = sum_{j = 1}^infty j cdot a_j x^{j - 1}}
        f'(x) = \sum_{j = 1}^\infty j \cdot a_j \, x^{j - 1}
\end{equation}
for each $x \in k$ with $|x| \le 1$.  Of course, these series reduce
to finite sums when $a_j = 0$ for all but finitely many $j$, and
otherwise these series converge in $k$ when $k$ is complete with
respect to the ultrametric associated to $|\cdot|$, because of
(\ref{lim_{j to infty} |a_j| = 0}).  If $|\cdot|$ is not the trivial
absolute value function on $k$, then (\ref{f'(x) = sum_{j = 1}^infty j
  cdot a_j x^{j - 1}}) is indeed the derivative of (\ref{f(x) = sum_{j
    = 0}^infty a_j x^j, 3}), as in Section \ref{differentiation,
  continued}.

        Under these conditions, we have that
\begin{equation}
\label{|f(x)| le max_{j ge 0} (|a_j| |x|^j) le 1}
        |f(x)| \le \max_{j \ge 0} (|a_j| \, |x|^j) \le 1
\end{equation}
for every $x \in k$ with $|x| \le 1$, by the ultrametric version
of the triangle inequality.  Similarly,
\begin{equation}
\label{|f'(x)| le max_{j ge 1} (|j cdot a_j| |x|^{j - 1}) le 1}
        |f'(x)| \le \max_{j \ge 1} (|j \cdot a_j| \, |x|^{j - 1}) \le 1
\end{equation}
for every $x \in k$ with $|x| \le 1$, which also corresponds to
(\ref{|f'(x)| le ... le max_{j ge 1} (|a_j| r^{j - 1})}) in Section
\ref{some related estimates}, with $r = 1$.  Note that $f$ is
Lipschitz of order $1$ with constant equal to $1$ as a mapping from
the closed unit ball $\overline{B}(0, 1)$ in $k$ into $k$, by
(\ref{|f(x) - f(y)| le ... le max_{j ge 1} (|a_j| r^{j - 1}) |x - y|})
in Section \ref{lipschitz mappings, continued}, with $r = 1$ again.

        Let $x_0 \in k$ be given, with $|x_0| \le 1$.  Suppose for the
moment that
\begin{equation}
\label{max_{j ge 2} |a_j| < |f'(x_0)|}
        \max_{j \ge 2} |a_j| < |f'(x_0)|,
\end{equation}
which is the same as saying that (\ref{t max_{j ge 2} (|a_j| r^{j -
    2}) < |f'(x_0)|}) in Section \ref{hensel's lemma} holds with $r =
t = 1$.  This implies that (\ref{|f'(x)| = |f'(x_0)|}) and (\ref{|f(x)
  - f(y)| = |f'(x_0)| |x - y|}) hold for every $x$, $y$ in
\begin{equation}
\label{overline{B}(x_0, 1) = overline{B}(0, 1)}
        \overline{B}(x_0, 1) = \overline{B}(0, 1),
\end{equation}
as before.  In particular, (\ref{|f'(x)| = |f'(x_0)|}) holds with $x =
0$, so that one could have taken $x_0 = 0$ here.  If $k$ is complete,
then
\begin{equation}
\label{f(overline{B}(x_0, 1)) = overline{B}(f(x_0), |f'(x_0)|)}
        f(\overline{B}(x_0, 1)) = \overline{B}(f(x_0), |f'(x_0)|),
\end{equation}
by (\ref{f(overline{B}(x_0, t)) = overline{B}(f(x_0), |f'(x_0)| t)}).

        Suppose now that $f'(x_0) \ne 0$, and that
\begin{equation}
\label{|f'(x_0)| le max_{j ge 2} |a_j|}
        |f'(x_0)| \le \max_{j \ge 2} |a_j|.
\end{equation}
Let $t_0$ be a positive real number such that
\begin{equation}
\label{t_0 max_{j ge 2} |a_j| le |f'(x_0)|}
        t_0 \, \max_{j \ge 2} |a_j| \le |f'(x_0)|,
\end{equation}
which implies that $t_0 \le 1$, by (\ref{|f'(x_0)| le max_{j ge 2}
  |a_j|}).  Note that
\begin{equation}
\label{t_0 = |f'(x_0)|}
        t_0 = |f'(x_0)|
\end{equation}
satisfies (\ref{t_0 max_{j ge 2} |a_j| le |f'(x_0)|}), because of
(\ref{|a_j| le 1}).  Of course, (\ref{t_0 max_{j ge 2} |a_j| le
  |f'(x_0)|}) is the same as (\ref{t_0 max_{j ge 2} (|a_j| r^{j - 2})
  le |f'(x_0)|}), with $r = 1$.  This implies that (\ref{|f'(x)| =
  |f'(x_0)|}) and (\ref{|f(x) - f(y)| = |f'(x_0)| |x - y|}) hold for
every $x, y \in B(x_0, t_0)$, as in Section \ref{some variants}.  If
$k$ is complete, then (\ref{f(B(x_0, t_0)) = B(f(x_0), |f'(x_0)|
  t_0)}) holds too, as before.  As in the previous section, one can
also consider analogous arguments for the expansion of $f$ as a power
series around $x_0$.

\section{Some examples}
\label{some examples}

        Let $k$ be a field with an ultrametric absolute value function
$|\cdot|$, as before.  Also let $n$ be an integer with $n \ge 2$, and
put
\begin{equation}
\label{f(x) = x^n}
        f(x) = x^n
\end{equation}
for each $x \in k$ with $|x| \le 1$.  Thus $f$ is Lipschitz of order
$1$ with constant equal to $1$ as a mapping from the closed unit ball
$\overline{B}(0, 1)$ in $k$ into $k$, as in Section \ref{lipschitz
  mappings, continued}, and as mentioned in the preceding section.  In
this case, it is easy to see that $f$ cannot be Lipschitz of order $1$
with constant strictly less than $1$ on $\overline{B}(0, 1)$, because
$f(0) = 0$ and $f(1) = 1$.

        Put
\begin{equation}
\label{f'(x) = n x^{n - 1}}
        f'(x) = n \, x^{n - 1}
\end{equation}
for each $x \in k$ with $|x| \le 1$, which corresponds to the formal
derivative of $f$, and which is the derivative of $f$ as a function on
$\overline{B}(0, 1)$ when $|\cdot|$ is nontrivial, as usual.  It
follows that
\begin{equation}
\label{|f'(x)| = |n cdot 1| |x|^{n - 1}}
        |f'(x)| = |n \cdot 1| \, |x|^{n - 1},
\end{equation}
where $n \cdot 1$ is the sum of $n$ $1$'s in $k$.  Suppose for the
moment that
\begin{equation}
\label{|n cdot 1| = 1}
        |n \cdot 1| = 1,
\end{equation}
and let $x_0 \in k$ be given, with $|x_0| = 1$, so that
\begin{equation}
\label{|f'(x_0)| = 1}
        |f'(x_0)| = 1.
\end{equation}
In this case, we have equality in (\ref{|f'(x_0)| le max_{j ge 2}
  |a_j|}), and we can take $t_0 = 1$, as in (\ref{t_0 = |f'(x_0)|}).
As before, the restriction of $f$ to $B(x_0, 1)$ is an isometry, and
$f$ maps $B(x_0, 1)$ onto $B(f(x_0, 1))$ when $k$ is complete.

        More precisely, if $k$ has characteristic $0$, then there is
a natural embedding of ${\bf Q}$ into $k$, so that $|\cdot|$ induces
an ultrametric absolute value function on ${\bf Q}$.  If the induced
absolute value function on ${\bf Q}$ is trivial, then (\ref{|n cdot 1|
  = 1}) holds automatically.  Otherwise, Ostrowski's theorem implies
that the induced absolute value function on ${\bf Q}$ is equivalent to
the $p$-adic absolute value function on ${\bf Q}$ for some prime
number $p$, as in Section \ref{ostrowski's theorems}.  In this case,
(\ref{|n cdot 1| = 1}) holds when $n$ is not divisible by $p$.
Similarly, if $k$ has characteristic $p$ for some prime number $p$,
then (\ref{|n cdot 1| = 1}) holds exactly when $n$ is not divisible by
$p$.

        If $x, y \in k$ and $|x| \ne |y|$, then
\begin{equation}
\label{|x - y| = max(|x|, |y|)}
        |x - y| = \max(|x|, |y|),
\end{equation}
as in (\ref{|x + y| = |y|}) in Section \ref{absolute value functions}.
In this case, $|x^n| = |x|^n \ne |y|^n = |y^n|$ for each positive
integer $n$, so that
\begin{equation}
\label{|x^n - y^n| = max(|x|^n, |y|^n)}
        |x^n - y^n| = \max(|x|^n, |y|^n),
\end{equation}
for the same reasons as in (\ref{|x - y| = max(|x|, |y|)}).  Note that
this is consistent with the fact that (\ref{f(x) = x^n}) is Lipschitz
of order $1$ with constant equal to $1$ on $\overline{B}(0, 1)$, as
before.

\section{Some examples, continued}
\label{some examples, continued}

        Let $k$ be a field with an ultrametric absolute value function
$|\cdot|$ again, and consider $f(x) = x^n$, as in (\ref{f(x) = x^n}).
Suppose now that $n = p$ for some prime number $p$, and that (\ref{|n
  cdot 1| = 1}) does not hold.  The case where $k$ has characteristic
$p$ was discussed in Section \ref{derivative 0}, and so we suppose
that $k$ has characteristic $0$.  Thus Ostrowski's theorem implies
that the induced absolute value function on ${\bf Q}$ is equivalent to
the $p$-adic absolute value function, as before.  Let us ask that the
induced absolute value function on ${\bf Q}$ actually be equal to the
$p$-adic absolute value function, which can be arranged by replacing
the given absolute value function on $k$ by an appropriate positive
power of itself.

        Let $x_0 \in k$ be given again, with $|x_0| = 1$, so that
\begin{equation}
\label{|f'(x_0)| = 1/p}
        |f'(x_0)| = 1/p.
\end{equation}
The right side of (\ref{|f'(x_0)| le max_{j ge 2} |a_j|}) in Section
\ref{a basic situation} is equal to $1$ in this situation, which
implies that (\ref{|f'(x_0)| le max_{j ge 2} |a_j|}) holds with strict
inequality.  Similarly, the maximal value of $t_0 > 0$ that satisfies
(\ref{t_0 max_{j ge 2} |a_j| le |f'(x_0)|}) is given by (\ref{t_0 =
  |f'(x_0)|}).  Let us now consider the analogous arguments for the
expansion of $f$ around $x_0$, as in Section \ref{some variants,
  continued}.

        Using the binomial theorem, we get that
\begin{equation}
\label{f(x) = ... = sum_{l = 0}^p {p choose l} cdot x_0^{p - l} (x - x_0)^l}
        f(x) = (x_0 + (x - x_0))^p
             = \sum_{l = 0}^p {p \choose l} \cdot x_0^{p - l} \, (x - x_0)^l,
\end{equation}
which corresponds to (\ref{f(x) = sum_{j = 0}^infty a_j (x_0 + (x -
  x_0))^j = ...}) in Section \ref{some related estimates, continued}.
Thus the coefficients $\widetilde{a}_l$ in (\ref{widetilde{a}_l =
  sum_{j = l}^infty a_j {j choose l} cdot x_0^{j - l}}) reduce to
\begin{equation}
\label{widetilde{a}_l = {p choose l} cdot x_0^{p - l}}
        \widetilde{a}_l = {p \choose l} \cdot x_0^{p - l}
\end{equation}
when $l \le p$, and $\widetilde{a}_l = 0$ when $l > p$.  This
implies that
\begin{equation}
\label{|widetilde{a}_l| = 1/p}
        |\widetilde{a}_l| = 1/p
\end{equation}
for $l = 1, \ldots, p - 1$, while $|\widetilde{a}_0| =
|\widetilde{a}_p| = 1$.  As in Section \ref{some variants, continued},
we would like to choose $t_0 > 0$ as large as possible so that $t_0
\le 1$ and (\ref{t_0 max_{l ge 2} (|widetilde{a}_l| t_0^{l - 2}) =
  ... le |f'(x_0)|}) holds.  If $p = 2$, then (\ref{t_0 max_{l ge 2}
  (|widetilde{a}_l| t_0^{l - 2}) = ... le |f'(x_0)|}) reduces to
\begin{equation}
\label{t_0 le |f'(x_0)| = 1/2}
        t_0 \le |f'(x_0)| = 1/2,
\end{equation}
so that the maximal value of $t_0$ is $1/2$.

        Otherwise, if $p > 2$, then
\begin{equation}
\label{max_{l ge 2} (|widetilde{a}_l| t_0^{l - 1}) = max (t_0/p, t_0^{p - 1})}
        \max_{l \ge 2} (|\widetilde{a}_l| \, t_0^{l - 1})
          = \max (t_0/p, t_0^{p - 1})
\end{equation}
for each $0 < t_0 \le 1$, which is to say that the maximum on the left
side occurs with either $l = 2$ or $l = p$.  In this case, (\ref{t_0
  max_{l ge 2} (|widetilde{a}_l| t_0^{l - 2}) = ... le |f'(x_0)|})
reduces to
\begin{equation}
\label{max (t_0/p, t_0^{p - 1}) le |f'(x_0)| = 1/p}
        \max (t_0/p, t_0^{p - 1}) \le |f'(x_0)| = 1/p.
\end{equation}
Of course, $t_0/p \le 1/p$ automatically when $t_0 \le 1$, and so we
can take
\begin{equation}
\label{t_0 = p^{-1/(p - 1)}}
        t_0 = p^{-1/(p - 1)}
\end{equation}
in the context of Section \ref{some variants, continued}.

        We can also use the binomial theorem to look at the behavior
of $f$ more directly.  Let $x, y \in k$ be given, with $|x| = |y|$,
since otherwise we already have (\ref{|x - y| = max(|x|, |y|)}) and
(\ref{|x^n - y^n| = max(|x|^n, |y|^n)}).  Put
\begin{equation}
\label{R = |x| = |y|}
        R = |x| = |y|,
\end{equation}
and note that
\begin{equation}
\label{|x - y| le R}
        |x - y| \le R,
\end{equation}
by the ultrametric version of the triangle inequality.  Using the
binomial theorem again, we have that
\begin{eqnarray}
\label{f(x) - f(y) = x^p - y^p = (y + (x - y))^p - y^p = ...}
        f(x) - f(y) = x^p - y^p & = & (y + (x - y))^p - y^p \\
 & = & \sum_{l = 0}^p {p \choose l} \cdot y^{p - l} \, (x - y)^l - y^p
                                                  \nonumber \\
 & = & \sum_{l = 1}^p {p \choose l} \cdot y^{p - l} \, (x - y)^l. \nonumber
\end{eqnarray}
Observe that
\begin{equation}
\label{|{p choose l} cdot y^{p - l} (x - y)^l|}
        \biggl|{p \choose l} \cdot y^{p - l} \, (x - y)^l\biggr|
\end{equation}
is equal to
\begin{equation}
\label{(1/p) R^{p - l} |x - y|^l}
         (1/p) \, R^{p - l} \, |x - y|^l
\end{equation}
when $ 1 \le l \le p - 1$, and to
\begin{equation}
\label{|x - y|^p}
        |x - y|^p
\end{equation}
when $l = p$.  In particular, (\ref{|{p choose l} cdot y^{p - l} (x -
  y)^l|}) is equal to
\begin{equation}
\label{(1/p) R^{p - 1} |x - y|}
        (1/p) \, R^{p - 1} \, |x - y|
\end{equation}
when $l = 1$, and (\ref{|{p choose l} cdot y^{p - l} (x - y)^l|}) is
less than or equal to (\ref{(1/p) R^{p - 1} |x - y|}) when $2 \le l
\le p - 1$, because of (\ref{|x - y| le R}).  Similarly, if $|x - y| <
R$ and $x \ne y$, then (\ref{|{p choose l} cdot y^{p - l} (x - y)^l|})
is strictly less than (\ref{(1/p) R^{p - 1} |x - y|}) when $2 \le l
\le p - 1$.

        Suppose for the moment that
\begin{equation}
\label{|x - y| < p^{-1/(p - 1)} R}
        |x - y| < p^{-1/(p - 1)} \, R,
\end{equation}
so that
\begin{equation}
\label{|x - y|^{p - 1} < (1/p) R^{p - 1}}
        |x - y|^{p - 1} < (1/p) \, R^{p - 1}.
\end{equation}
If $x \ne y$, then this implies that
\begin{equation}
\label{|x - y|^p < (1/p) R^{p - 1} |x - y|}
        |x - y|^p < (1/p) \, R^{p - 1} \, |x - y|.
\end{equation}
Note that (\ref{|x - y| < p^{-1/(p - 1)} R}) implies that $|x - y| <
R$, so that (\ref{|{p choose l} cdot y^{p - l} (x - y)^l|}) is
strictly less than (\ref{(1/p) R^{p - 1} |x - y|}) when $2 \le l \le p
-1$ and $x \ne y$, as in the preceding paragraph.  It follows that
\begin{equation}
\label{|f(x) - f(y)| = (1/p) R^{p - 1} |x - y|}
        |f(x) - f(y)| = (1/p) \, R^{p - 1} \, |x - y|
\end{equation}
when $x, y \in k$ satisfy (\ref{R = |x| = |y|}), (\ref{|x - y| <
  p^{-1/(p - 1)} R}), and $x \ne y$, and of course (\ref{|f(x) - f(y)|
  = (1/p) R^{p - 1} |x - y|}) holds trivially when $x = y$.  More
precisely, if $x \ne y$ satisfy (\ref{R = |x| = |y|}) and (\ref{|x -
  y| < p^{-1/(p - 1)} R}), then the absolute values of the terms on
the right side of (\ref{f(x) - f(y) = x^p - y^p = (y + (x - y))^p -
  y^p = ...})  corresponding to $l \ge 2$ are strictly less than
(\ref{(1/p) R^{p - 1} |x - y|}), which is the absolute value of the $l
= 1$ term on the right side of (\ref{f(x) - f(y) = x^p - y^p = (y + (x
  - y))^p - y^p = ...}).  This implies that the absolute value of the
left side of (\ref{f(x) - f(y) = x^p - y^p = (y + (x - y))^p - y^p =
  ...})  is equal to the absolute value of the $l = 1$ term on the
right side of (\ref{f(x) - f(y) = x^p - y^p = (y + (x - y))^p - y^p =
  ...}), which is the same as (\ref{|f(x) - f(y)| = (1/p) R^{p - 1} |x
  - y|}).  This uses the same type of argument as in (\ref{|x + y| =
  |y|}) in Section \ref{absolute value functions}.  The same
conclusion follows from the earlier discussion when $R = 1$, and it is
easy to reduce to that case in this situation anyway, or to deal with
it in the same way as before.

        Similarly, if
\begin{equation}
\label{|x - y| = p^{-1/(p - 1)} R}
        |x - y| = p^{-1/(p - 1)} \, R,
\end{equation}
then
\begin{equation}
\label{|x - y|^{p - 1} = (1/p) R^{p - 1}}
        |x - y|^{p - 1} = (1/p) \, R^{p - 1},
\end{equation}
and hence
\begin{equation}
\label{|x - y|^p = (1/p) R^{p - 1} |x - y|}
        |x - y|^p = (1/p) \, R^{p - 1} \, |x - y|.
\end{equation}
This implies that
\begin{equation}
\label{|f(x) - f(y)| le (1/p) R^{p - 1} |x - y|}
        |f(x) - f(y)| \le (1/p) \, R^{p - 1} \, |x - y|,
\end{equation}
because (\ref{|{p choose l} cdot y^{p - l} (x - y)^l|}) is less than
or equal to (\ref{(1/p) R^{p - 1} |x - y|}) when $2 \le l \le p - 1$,
as before.  More precisely, all of the terms on the right side of
(\ref{f(x) - f(y) = x^p - y^p = (y + (x - y))^p - y^p = ...}) have
absolute value less than or equal to (\ref{(1/p) R^{p - 1} |x - y|})
under these conditions, which implies (\ref{|f(x) - f(y)| le (1/p)
  R^{p - 1} |x - y|}).  As in the previous paragraph, one could get
the same conclusion from the earlier discussion when $R = 1$, and the
restriction to $R = 1$ is not very serious anyway.

        Now suppose that
\begin{equation}
\label{|x - y| > p^{-1/(p - 1)} R}
        |x - y| > p^{-1/(p - 1)} \, R,
\end{equation}
so that
\begin{equation}
\label{|x - y|^{p - 1} > (1/p) R^{p - 1}}
        |x - y|^{p - 1} > (1/p) \, R^{p - 1},
\end{equation}
and thus
\begin{equation}
\label{|x - y|^p > (1/p) R^{p - 1} |x - y|}
        |x - y|^p > (1/p) \, R^{p - 1} \, |x - y|.
\end{equation}
This means that (\ref{|{p choose l} cdot y^{p - l} (x - y)^l|}) is
strictly less than (\ref{|x - y|^p}) when $1 \le l \le p - 1$, while
(\ref{|{p choose l} cdot y^{p - l} (x - y)^l|}) is equal to (\ref{|x -
  y|^p}) when $l = p$.  It follows that
\begin{equation}
\label{|f(x) - f(y)| = |x - y|^p}
        |f(x) - f(y)| = |x - y|^p
\end{equation}
in this case, using (\ref{f(x) - f(y) = x^p - y^p = (y + (x - y))^p -
  y^p = ...}) and the same type of argument as in (\ref{|x + y| =
  |y|}) in Section \ref{absolute value functions} again.

\chapter{Some additional topics}
\label{some additional topics}

\section{Sums of functions}
\label{sums of functions}

        Let $k$ be a field with an absolute value function $|\cdot|$,
and suppose that $k$ is complete with respect to the metric associated
to $|\cdot|$.  Also let $E$ and $M$ be nonempty sets, and let
$a_\alpha(x)$ be a $k$-valued function on $M$ for each $\alpha \in E$.
Let us say that
\begin{equation}
\label{sum_{alpha in E} a_alpha}
        \sum_{\alpha \in E} a_\alpha
\end{equation}
converges \emph{pointwise}\index{pointwise convergence} on $M$ if for
each $x \in M$,
\begin{equation}
\label{sum_{alpha in E} a_alpha(x)}
        \sum_{\alpha \in E} a_\alpha(x)
\end{equation}
converges as a sum in $k$ in the sense discussed in Section
\ref{generalized convergence}.  By definition, this means that the
corresponding net of finite sums
\begin{equation}
\label{sum_{alpha in A} a_alpha(x)}
        \sum_{\alpha \in A} a_\alpha(x)
\end{equation}
converges in $k$ for each $x \in M$, where $A$ is a finite subset of
$E$.  Remember that this happens exactly when (\ref{sum_{alpha in E}
  a_alpha(x)}) satisfies the generalized Cauchy criterion as a sum in
$k$ for each $x \in M$, as in Sections \ref{generalized convergence}
and \ref{generalized convergence, continued}, because $k$ is complete.

        Let $\ell^\infty(M, k)$ be the vector space of bounded
$k$-valued functions on $M$, equipped with the supremum norm,
as in Section \ref{supremum norm}.  If $a_\alpha$ is bounded as a
$k$-valued function on $M$ for each $\alpha \in E$, then
(\ref{sum_{alpha in E} a_alpha}) may be considered as a sum in
$\ell^\infty(M, k)$, as in Section \ref{generalized convergence}.
Thus (\ref{sum_{alpha in E} a_alpha}) converges in $\ell^\infty(M, k)$
if the corresponding net of finite sums
\begin{equation}
\label{sum_{alpha in A} a_alpha}
        \sum_{\alpha \in A} a_\alpha
\end{equation}
converges in $\ell^\infty(M, k)$ with respect to the supremum norm,
where $A$ is a finite subset of $E$ again.  In this case, this means
that the finite sums (\ref{sum_{alpha in A} a_alpha(x)}) converge
uniformly on $M$.  In particular, this implies that the finite sums
(\ref{sum_{alpha in A} a_alpha(x)}) converge pointwise on $M$.
Remember that $\ell^\infty(M, k)$ is complete with respect to the
metric associated to the supremum norm, because $k$ is complete, by
hypothesis.  It follows that (\ref{sum_{alpha in E} a_alpha})
converges in $\ell^\infty(M, k)$ with respect to the supremum norm if
and only if (\ref{sum_{alpha in E} a_alpha}) satisfies the generalized
Cauchy criterion with respect to the supremum norm, as in Sections
\ref{generalized convergence} and \ref{generalized convergence,
  continued}.

        Let us continue to ask that $a_\alpha$ be a bounded $k$-valued
function on $M$ for each $\alpha \in E$, and let $a$ denote the mapping
from $E$ into $\ell^\infty(M, k)$ defined by
\begin{equation}
\label{alpha mapsto a_alpha}
        \alpha \mapsto a_\alpha.
\end{equation}
Similarly, for each $x \in M$, let $a(x)$ denote the mapping from $E$
into $k$ defined by
\begin{equation}
\label{alpha mapsto a_alpha(x)}
        \alpha \mapsto a_\alpha(x).
\end{equation}
Suppose that $a$ has bounded finite sums on $E$ with respect to the
supremum norm on $\ell^\infty(M, k)$, as in Section \ref{bounded
  finite sums}.  This implies that $a(x)$ has bounded finite sums on
$E$ for each $x \in M$.  More precisely, for each $x \in M$, we have
that
\begin{equation}
\label{||a(x)||_{BFS(E, k)} le ||a||_{BFS(E, ell^infty(M, k))}}
        \|a(x)\|_{BFS(E, k)} \le \|a\|_{BFS(E, \ell^\infty(M, k))},
\end{equation}
where the BFS norms are as defined in Section \ref{bounded finite
  sums}.  Conversely, if $a(x)$ has bounded finite sums on $E$ for
each $x \in M$, and if the BFS norm of $a(x)$ on $E$ is uniformly
bounded over $x \in M$, then $a$ has bounded finite sums on $E$ with
respect to the supremum norm on $\ell^\infty(M, k)$.  In this
situation, it is easy to see that
\begin{equation}
\label{||a||_{BFS(E, ell^infty(M, k))} = sup_{x in M} ||a(x)||_{BFS(E, k)}}
 \|a\|_{BFS(E, \ell^\infty(M, k))} = \sup_{x \in M} \|a(x)\|_{BFS(E, k)}.
\end{equation}

        If (\ref{sum_{alpha in E} a_alpha}) converges in $\ell^\infty(M, k)$,
then $a$ has bounded finite sums on $E$ , as in Section \ref{bounded
  finite sums}.  Similarly, if (\ref{sum_{alpha in E} a_alpha(x)})
converges in $k$ for some $x \in M$, then $a(x)$ has bounded finite
sums on $E$, and
\begin{equation}
\label{|sum_{alpha in E} a_alpha(x)| le ||a(x)||_{BFS(E, k)}}
 \biggl|\sum_{\alpha \in E} a_\alpha(x)\biggr| \le \|a(x)\|_{BFS(E, k)}.
\end{equation}
If $a$ has bounded finite sums on $E$, and if (\ref{sum_{alpha in E}
  a_alpha(x)}) converges in $k$ for every $x \in M$, then
(\ref{sum_{alpha in E} a_alpha(x)}) defines a bounded $k$-valued
function on $M$, with
\begin{equation}
\label{|sum_{alpha in E} a_alpha(x)| le ||a||_{BFS(E, ell^infty(M, k))}}
 \biggl|\sum_{\alpha \in E} a_\alpha(x)\biggr| \le \|a\|_{BFS(E, \ell^\infty(M, k))}
\end{equation}
for each $x \in M$, by (\ref{||a(x)||_{BFS(E, k)} le ||a||_{BFS(E,
    ell^infty(M, k))}}) and (\ref{|sum_{alpha in E} a_alpha(x)| le
  ||a(x)||_{BFS(E, k)}}).

        Suppose for the moment that $k = {\bf R}$, with the standard
absolute value function.  If $a(x)$ has bounded finite sums on $E$
for some $x \in M$, then $a(x)$ is summable on $E$, as mentioned
in Section \ref{bounded finite sums}.  More precisely, we have that
\begin{equation}
\label{sum_{alpha in E} |a_alpha(x)| le 2 ||a(x)||_{BFS(E, {bf R})}}
        \sum_{\alpha \in E} |a_\alpha(x)| \le 2 \, \|a(x)\|_{BFS(E, {\bf R})},
\end{equation}
as in (\ref{sum_{x in X} |f(x)| le 2 C}) in Section \ref{bounded
  finite sums}.  In particular, this implies that (\ref{sum_{alpha in
    E} a_alpha(x)}) converges in ${\bf R}$, as in Section
\ref{generalized convergence, continued}.  If $a$ has bounded finite
sums on $E$, then $a(x)$ has bounded finite sums on $E$ for each $x
\in M$, as before.  This implies that
\begin{equation}
\label{sum_{alpha in E} |a_alpha(x)| le 2 ||a||_{BFS(E, ell^infty(M, {bf R}))}}
 \sum_{\alpha \in E} |a_\alpha(x)| \le 2 \, \|a\|_{BFS(E, \ell^\infty(M, {\bf R}))}
\end{equation}
for every $x \in M$, by combining (\ref{||a(x)||_{BFS(E, k)} le
  ||a||_{BFS(E, ell^infty(M, k))}}) and (\ref{sum_{alpha in E}
  |a_alpha(x)| le 2 ||a(x)||_{BFS(E, {bf R})}}).  Thus
(\ref{sum_{alpha in E} a_alpha(x)}) converges in ${\bf R}$ for each $x
\in M$ under these conditions, and defines a bounded real-valued
function on $M$, as in the preceding paragraph.

        Similarly, if $k = {\bf C}$ with the standard absolute value
function, then we can apply the remarks in the previous paragraph to
the real and imaginary parts of $a$.  If $a(x)$ has bounded finite
sums on $E$ for some $x \in M$, then it follows that $a(x)$ is
summable on $E$, with
\begin{equation}
\label{sum_{alpha in E} |a_alpha(x)| le 4 ||a(x)||_{BSF(E, {bf C})}}
        \sum_{\alpha \in E} |a_\alpha(x)| \le 4 \, \|a(x)\|_{BSF(E, {\bf C})},
\end{equation}
by applying (\ref{sum_{alpha in E} |a_alpha(x)| le 2 ||a(x)||_{BFS(E,
    {bf R})}}) to the real and imaginary parts of $a(x)$.  If $a$ has
bounded finite sums on $E$, then $a(x)$ has bounded finite sums on $E$
for every $x \in M$, and we get that
\begin{equation}
\label{sum_{alpha in E} |a_alpha(x)| le 4 ||a||_{BFS(E, ell^infty(M, {bf C}))}}
 \sum_{\alpha \in E} |a_\alpha(x)| \le 4 \, \|a\|_{BFS(E, \ell^\infty(M, {\bf C}))}
\end{equation}
for every $x \in M$, by combining (\ref{||a(x)||_{BFS(E, k)} le
  ||a||_{BFS(E, ell^infty(M, k))}}) and (\ref{sum_{alpha in E}
  |a_alpha(x)| le 4 ||a(x)||_{BSF(E, {bf C})}}).  This implies that
(\ref{sum_{alpha in E} a_alpha(x)}) converges in ${\bf C}$ for each $x
\in M$, and defines a bounded complex-valued function on $M$, as
before.

        Suppose now that $k$ is a field with an ultrametric
absolute value function $|\cdot|$, and that $k$ is still complete
with respect to the associated ultrametric.  In this case, the
supremum norm on $\ell^\infty(M, k)$ is an ultranorm, as in Section
\ref{supremum norm}.  This implies that the BFS norm of a function
on $E$ with values in $k$ or $\ell^\infty(M, k)$ is the same as the
corresponding supremum norm on $E$, as in Section \ref{bounded finite
sums}.  Similarly, the sum over $E$ of a function on $E$ with values
in $k$ or in $\ell^\infty(M, k)$ satisfies the generalized Cauchy
criterion if and only if the function vanishes at infinity on $E$,
as in Sections \ref{generalized convergence} and \ref{generalized
convergence, continued}.

        Let $k$ be any field with an absolute value function $|\cdot|$
again, and where $k$ is still complete with respect to the associated
metric.  Suppose that $M$ is now also equipped with a topology, and
let $C_b(M, k)$\index{C_b(M, k)@$C_b(M, k)$} be the space of bounded
continuous $k$-valued functions on $M$.  Thus $C_b(M, k)$ is a
subalgebra of $\ell^\infty(M, k)$ with respect to pointwise addition
and multiplication, and a closed set in $\ell^\infty(M, k)$ with
respect to the supremum norm.  If $a_\alpha \in C_b(M, k)$ for each
$\alpha \in E$, and if the sum (\ref{sum_{alpha in E} a_alpha})
converges in $\ell^\infty(M, k)$ with respect to the supremum norm,
then the sum is a continuous function on $M$ too.

\section{Lipschitz seminorms}
\label{lipschitz seminorms}

        Let $k$ be a field, let $|\cdot|$ be an absolute value function
on $k$, and let $V$ be a vector space over $k$.  A nonnegative
real-valued function $N$ on $V$ is said to be a
\emph{seminorm}\index{seminorms} on $V$ if
\begin{equation}
\label{N(t v) = |t| N(v), 2}
        N(t \, v) = |t| \, N(v)
\end{equation}
for every $v \in V$ and $t \in k$, and
\begin{equation}
\label{N(v + w) le N(v) + N(w), 2}
        N(v + w) \le N(v) + N(w)
\end{equation}
for every $v, w \in V$.  Note that (\ref{N(t v) = |t| N(v), 2})
implies that $N(0) = 0$, by taking $t = 0$, and that a seminorm $N$ on
$V$ is a norm on $V$ when $N(v) > 0$ for every $v \in V$ with $v \ne
0$.  A seminorm $N$ on $V$ may be called an
\emph{ultra-seminorm}\index{ultra-seminorm} if
\begin{equation}
\label{N(v + w) le max(N(v), N(w)), 2}
        N(v + w) \le \max(N(v), N(w))
\end{equation}
for every $v, w \in V$, which automatically implies (\ref{N(v + w) le
  N(v) + N(w), 2}).

        Let $(M, d(x, y))$ be a (nonempty) metric space, and let
$\gamma$ be a positive real number.  Consider the space
$\Lip_\gamma(M, k)$\index{Lip_{gamma}(M, k)@$\Lip_\gamma(M, k)$}
of $k$-valued functions on $M$ that are Lipschitz of order $\gamma$,
as in Section \ref{lipschitz mappings}.  More precisely, this uses the
metric on $k$ associated to the absolute value function.  It is easy
to see that $\Lip_\gamma(M, k)$ is a vector space over $k$ with
respect to pointwise addition and scalar multiplication, as in Section
\ref{lipschitz mappings} again.  If $f \in \Lip_\gamma(M, k)$, then
put
\begin{equation}
\label{||f||_{Lip_gamma(M, k)} = ...}
 \|f\|_{\Lip_\gamma(M, k)} = \sup \bigg\{\frac{|f(x) - f(y)|}{d(x, y)^\gamma} :
                                   x, y \in M, \ x \ne y\bigg\}
\end{equation}
when $M$ has at least two elements, and otherwise put
$\|f\|_{\Lip_\gamma(M, k)} = 0$.  If $M$ has at least two elements,
then the right side of (\ref{||f||_{Lip_gamma(M, k)} = ...}) is the
supremum of a bounded nonempty set of nonnegative real numbers, and
hence is a nonnegative real number.  Of course, if $f$ is Lipschitz of
order $\gamma$ on $M$ with constant $C$, then
(\ref{||f||_{Lip_gamma(M, k)} = ...}) is less than or equal to $C$.
It is easy to see that $f \in \Lip_\gamma(M, k)$ is Lipschitz of order
$\gamma$ with constant equal to (\ref{||f||_{Lip_gamma(M, k)} = ...}),
so that (\ref{||f||_{Lip_gamma(M, k)} = ...}) may be characterized
equivalently as the smallest nonnegative real number $C$ such that $f$
is Lipschitz of order $\gamma$ with constant $C$.  Observe that
(\ref{||f||_{Lip_gamma(M, k)} = ...}) is equal to $0$ if and only if
$f$ is constant on $M$.  One can check that (\ref{||f||_{Lip_gamma(M,
    k)} = ...}) defines a seminorm on $\Lip_\gamma(M, k)$, as in
Section \ref{lipschitz mappings}.  If $|\cdot|$ is an ultrametric absolute
value function on $k$, then (\ref{||f||_{Lip_gamma(M, k)} = ...})
defines an ultra-seminorm on $\Lip_\gamma(M, k)$.

        Let
\begin{equation}
\label{{Lip}_{b, gamma}(M, k) = {Lip}_gamma(M, k) cap ell^infty(M, k)}
        {\Lip}_{b, \gamma}(M, k) = {\Lip}_\gamma(M, k) \cap \ell^\infty(M, k)
\end{equation}
be the space of bounded $k$-valued functions on $M$ that are Lipschitz
of order $\gamma$.\index{Lip_{b, gamma}(M, k)@$\Lip_{b, \gamma}(M,
  k)$} Of course, this is a linear subspace of both $\Lip_\gamma(M,
k)$ and $\ell^\infty(M, k)$, which are both linear subspaces of the
space of all $k$-valued functions on $M$.  If $M$ is bounded, then
every $k$-valued Lipschitz function on $M$ of any positive order is
bounded on $M$, so that (\ref{{Lip}_{b, gamma}(M, k) = {Lip}_gamma(M,
  k) cap ell^infty(M, k)}) is the same as $\Lip_\gamma(M, k)$.

        There are two particularly simple ways in which to define a norm on
(\ref{{Lip}_{b, gamma}(M, k) = {Lip}_gamma(M, k) cap ell^infty(M, k)}).
The first possibility is to put
\begin{equation}
\label{||f||_{Lip_{b, gamma}(M, k)} = sum of norms}
 \|f\|_{\Lip_{b, \gamma}(M, k)} = \|f\|_{\Lip_\gamma(M, k)} + \|f\|_{\ell^\infty(M, k)}
\end{equation}
for every $f \in \Lip_{b, \gamma}(M, k)$, where $\|f\|_{\ell^\infty(M, k)}$
denotes the supremum norm of $f$, as in Section \ref{supremum norm}.
The second possibility is to put
\begin{equation}
\label{||f||_{Lip_{b, gamma}(M, k)} = max of norms}
 \|f\|_{\Lip_{b, \gamma}(M, k)} = \max(\|f\|_{\Lip_\gamma(M, k)},
                                         \|f\|_{\ell^\infty(M, k)})
\end{equation}
for every $f \in \Lip_{b, \gamma}(M, k)$.  It is easy to see that both
(\ref{||f||_{Lip_{b, gamma}(M, k)} = sum of norms}) and
(\ref{||f||_{Lip_{b, gamma}(M, k)} = max of norms}) define norms on
$\Lip_{b, \gamma}(M, k)$, because $\|f\|_{\Lip_\gamma(M, k)}$ is a
seminorm on $\Lip_\gamma(M, k)$, and $\|f\|_{\ell^\infty(M, k)}$ is a
norm on $\ell^\infty(M, k)$.  Note that (\ref{||f||_{Lip_{b, gamma}(M,
    k)} = max of norms}) is less than or equal to (\ref{||f||_{Lip_{b,
      gamma}(M, k)} = sum of norms}), and that (\ref{||f||_{Lip_{b,
      gamma}(M, k)} = sum of norms}) is less than or equal to $2$
times (\ref{||f||_{Lip_{b, gamma}(M, k)} = max of norms}).  This
implies that (\ref{||f||_{Lip_{b, gamma}(M, k)} = sum of norms}) and
(\ref{||f||_{Lip_{b, gamma}(M, k)} = max of norms}) determine the same
topology on $\Lip_{b, \gamma}(M, k)$.  If $|\cdot|$ is an ultrametric
absolute value function on $k$, then (\ref{||f||_{Lip_{b, gamma}(M,
    k)} = max of norms}) has the advantage of being an ultranorm on
$\Lip_{b, \gamma}(M, k)$, because $\|f\|_{\Lip_\gamma(M, k)}$ is an
ultra-seminorm on $\Lip_\gamma(M, k)$, and $\|f\|_{\ell^\infty(M, k)}$
is an ultranorm on $\ell^\infty(M, k)$.

        If $f$, $g$ are bounded $k$-valued functions on $M$, then
their product $f \, g$ is bounded on $M$ too, and satisfies
\begin{equation}
\label{||f g||_{ell^infty(M, k)} le ...}
        \|f \, g\|_{\ell^\infty(M, k)} \le \|f\|_{\ell^\infty(M, k)} \,
                                         \|g\|_{\ell^\infty(M, k)}.
\end{equation}
If $f$ and $g$ are also both Lipschitz of order $\gamma$ on $M$, then
one can check that $f \, g$ is Lipschitz of order $\gamma$ on $M$ as
well, as in Section \ref{lipschitz mappings}.  More precisely, we have
that
\begin{equation}
\label{||f g||_{Lip_gamma(M, k)} le ...}
 \qquad \|f \, g\|_{\Lip_\gamma(M, k)}
           \le \|f\|_{\Lip_\gamma(M, k)} \, \|g\|_{\ell^\infty(M, k)}
              + \|f\|_{\ell^\infty(M, k)} \, \|g\|_{\Lip_\gamma(M, k)},
\end{equation}
as in (\ref{C_1 (sup_{x in M_1} |f_2(x)|) + (sup_{x in M_1} |f_1(x)|)
  C_2}) in Section \ref{lipschitz mappings}.  If the norm on $\Lip_{b,
  \gamma}(M, k)$ is defined as in (\ref{||f||_{Lip_{b, gamma}(M, k)} =
  sum of norms}), then it follows that
\begin{equation}
\label{||f g||_{Lip_{b, gamma}(M, k)} le ...}
        \|f \, g\|_{\Lip_{b, \gamma}(M, k)}
         \le \|f\|_{\Lip_{b, \gamma}(M, k)} \, \|g\|_{\Lip_{b, \gamma}(M, k)}.
\end{equation}
This is an advantage of (\ref{||f||_{Lip_{b, gamma}(M, k)} = sum of
  norms}) in the archimedian case.  If $|\cdot|$ is an ultrametric
absolute value function on $k$, then we have that
\begin{equation}
\label{||f g||_{Lip_gamma(M, k)} le ..., 2}
 \qquad \|f \, g\|_{\Lip_\gamma(M, k)}
          \le \max(\|f\|_{Lip_\gamma(M, k)} \, \|g\|_{\ell^\infty(M, k)},
                    \|f\|_{\ell^\infty(M, k)} \, \|g\|_{\Lip_\gamma(M, k)}),
\end{equation}
as in (\ref{max(C_1 (sup_{x in M_1} |f_2(x)|), (sup_{x in M_1}
  |f_1(x)|) C_2)}) in Section \ref{lipschitz mappings}.  In this case,
(\ref{||f g||_{Lip_{b, gamma}(M, k)} le ...}) still holds when the
$\Lip_{b, \gamma}(M, k)$ norm is defined as in (\ref{||f||_{Lip_{b,
      gamma}(M, k)} = max of norms}).

        If $k$ is complete with respect to the metric associated to
$|\cdot|$, then $\Lip_{b, \gamma}(M, k)$ is complete with respect to
the norms (\ref{||f||_{Lip_{b, gamma}(M, k)} = sum of norms}) and
(\ref{||f||_{Lip_{b, gamma}(M, k)} = max of norms}), which are
essentially equivalent for this purpose.  To see this, let $\{f_j\}_{j
  = 1}^\infty$ be a Cauchy sequence of elements of $\Lip_{b,
  \gamma}(M, k)$ with respect to either of these norms.  This
basically means that $\{f_j\}_{j = 1}^\infty$ is a Cauchy sequence
with respect to both the supremum norm and the Lipschitz seminorm
(\ref{||f||_{Lip_gamma(M, k)} = ...}).  The completeness of
$\ell^\infty(M, k)$ implies that $\{f_j\}_{j = 1}^\infty$ converges to
a bounded $k$-valued function on $M$ with respect to the supremum
norm.  It is easy to see that $f$ is also Lipschitz of order $\gamma$
on $M$ under these conditions, because the corresponding Lipschitz seminorms
of the $f_j$'s are bounded.  Similarly, one can check that
$\{f_j\}_{j = 1}^\infty$ converges to $f$ with respect to the Lipschitz
seminorm (\ref{||f||_{Lip_gamma(M, k)} = ...}) under these conditions,
using the Cauchy condition for $\{f_j\}_{j = 1}^\infty$ with respect
to the Lipschitz seminorm.  This implies that $\{f_j\}_{j = 1}^\infty$
converges to $f$ in $\Lip_{b, \gamma}(M, k)$, as desired.

        Suppose now that $M$ is a nonempty subset of $k$, equipped with
the metric that is the restriction to $M$ of the metric associated to
$|\cdot|$ on $M$, and let us take $\gamma = 1$.  Also let $\{f_j\}_{j
  = 1}^\infty$ be a sequence of elements of $\Lip_1(M, k)$ that
converges to some $f \in \Lip_1(M, k)$ with respect to the $\Lip_1(M,
k)$ seminorm, in the sense that
\begin{equation}
\label{lim_{j to infty} ||f_j - f||_{Lip_1(M, k)} = 0}
        \lim_{j \to \infty} \|f_j - f\|_{\Lip_1(M, k)} = 0.
\end{equation}
Let $x$ be an element of $M$ that is a limit point of $M$ too, and
suppose that $f_j$ is differentiable at $x$ for each $j \ge 1$, as in
Section \ref{differentiation}.  It is easy to see that
\begin{equation}
\label{|f_j'(x) - f_l'(x)| le ||f_j - f_l||_{Lip_1(M, k)}}
        |f_j'(x) - f_l'(x)| \le \|f_j - f_l\|_{\Lip_1(M, k)}
\end{equation}
for every $j, l \ge 1$, since the diffence quotients of a $k$-valued
function on $M$ are bounded by the $\Lip_1(M, k)$ seminorm of the
function.  Note that
\begin{equation}
\label{lim_{j, l to infty} ||f_j - f_l||_{Lip_1(M, k)} = 0}
        \lim_{j, l \to \infty} \|f_j - f_l\|_{\Lip_1(M, k)} = 0,
\end{equation}
because of (\ref{lim_{j to infty} ||f_j - f||_{Lip_1(M, k)} = 0}),
so that
\begin{equation}
\label{lim_{j, l to infty} |f_j'(x) - f_l'(x)| = 0}
        \lim_{j, l \to \infty} |f_j'(x) - f_l'(x)| = 0,
\end{equation}
by (\ref{|f_j'(x) - f_l'(x)| le ||f_j - f_l||_{Lip_1(M, k)}}).  Of
course, this says exactly that $\{f_j'(x)\}_{j = 1}^\infty$ is a
Cauchy sequence in $k$.  If $k$ is complete with respect to the metric
associated to $|\cdot|$, then $\{f_j'(x)\}_{j = 1}^\infty$ converges
to an element of $k$.  Under these conditions, one can check that $f$
is differentiable at $x$ as well, with
\begin{equation}
\label{f'(x) = lim_{j to infty} f_j'(x)}
        f'(x) = \lim_{j \to \infty} f_j'(x).
\end{equation}
More precisely, the difference quotients for $f_j$ at $x$ converge
uniformly to the corresponding difference quotients for $f$ at $x$,
because of (\ref{lim_{j to infty} ||f_j - f||_{Lip_1(M, k)} = 0}).
This permits one to interchange the order of the limits, by standard
arguments.

\section{The product rule}
\label{product rule}

        Let $k$ be a field, and let
\begin{equation}
\label{f(X) = sum_{j = 0}^infty a_j X^j, g(X) = sum_{j = 0}^infty b_j X^j}
        f(X) = \sum_{j = 0}^\infty a_j \, X^j, \quad
        g(X) = \sum_{j = 0}^\infty b_j \, X^j
\end{equation}
be formal power series with coefficients in $k$.  The derivatives
$f'(X)$, $g'(X)$ of $f(X)$, $g(X)$ are defined as formal power series
by
\begin{equation}
\label{f'(X) = sum j cdot a_j X^{j - 1}, g'(X) = sum j cdot b_j X^{j - 1}}
        f'(X) = \sum_{j = 1}^\infty j \cdot a_j \, X^{j - 1}, \quad
        g'(X) = \sum_{j = 1}^\infty j \cdot b_j \, X^{j - 1},
\end{equation}
as in (\ref{f'(X) = sum_{j = 1}^infty j cdot a_j X^{j - 1}}) in
Section \ref{differentiation, continued}. 
The product of $f(X)$ and $g(X)$ is given by
\begin{equation}
\label{(f g)(X) = f(X) g(X) = sum_{n = 0}^infty c_n X^n}
        (f \, g)(X) = f(X) \, g(X) = \sum_{n = 0}^\infty c_n \, X^n,
\end{equation}
where
\begin{equation}
\label{c_n = sum_{j = 0}^n a_j b_{n - j}, 3}
        c_n = \sum_{j = 0}^n a_j \, b_{n - j}
\end{equation}
is the Cauchy product of the coefficients of $f(X)$ and $g(X)$, as in
(\ref{f(T) g(T) = sum_{n = 0}^infty (f g)_n T^n}) and (\ref{(f g)_n =
  sum_{j = 0}^n f_j g_{n - j}}) in Section \ref{formal power series}.
Thus
\begin{equation}
\label{(f g)'(X) = sum n cdot c_n X^{n - 1}}
        (f \, g)'(X) = \sum_{n = 1}^\infty n \cdot c_n \, X^{n - 1},
\end{equation}
and one can check that
\begin{equation}
\label{(f g)'(X) = f'(X) g(X) + f(X) g'(X)}
        (f \, g)'(X) = f'(X) \, g(X) + f(X) \, g'(X),
\end{equation}
as in the usual product rule\index{product rule} for derivatives.
More precisely, it is easy to see that
\begin{eqnarray}
\label{n cdot c_n = sum_{j = 0}^n n cdot a_j b_{n - j} = ...}
        n \cdot c_n = \sum_{j = 0}^n n \cdot a_j \, b_{n - j}
   & = & \sum_{j = 0}^n (j \cdot a_j) \, b_{n - j} 
           + \sum_{j = 0}^n a_j \, ((n - j) \, b_{n - j}) \nonumber \\
   & = & \sum_{j = 1}^n (j \cdot a_j) \, b_{n - j}
           + \sum_{j = 0}^{n - 1} a_j \, ((n - j) \cdot b_{n - j})
\end{eqnarray}
 for each $n \ge 1$.  The two sums on the right side of (\ref{n cdot
   c_n = sum_{j = 0}^n n cdot a_j b_{n - j} = ...}) correspond to the
 coefficients of $X^{n - 1}$ in the two terms on the right side of
 (\ref{n cdot c_n = sum_{j = 0}^n n cdot a_j b_{n - j} = ...}), which
 are given by Cauchy products.  Of course, these expressions for the
 Cauchy products are slightly different from the usual ones, because
 of the shifts in the indices for the derivatives.

        Suppose for the moment that $k = {\bf R}$ or ${\bf C}$, with
the standard absolute value function, and that
\begin{equation}
\label{sum_{j = 0}^infty |a_j| r^j, sum_{j = 0}^infty |b_j| r^j}
        \sum_{j = 0}^\infty |a_j| \, r^j, \quad \sum_{j = 0}^\infty |b_j| \, r^j
\end{equation}
converge for some positive real number $r$.  Put
\begin{equation}
\label{f(x) = sum_{j = 0}^infty a_j x^j, g(x) = sum_{j = 0}^infty b_j x^j}
 f(x) = \sum_{j = 0}^\infty a_j \, x^j, \quad g(x) = \sum_{j = 0}^\infty b_j \, x^j
\end{equation}
for every $x \in k$ with $|x| \le r$, where the series in (\ref{f(x) =
  sum_{j = 0}^infty a_j x^j, g(x) = sum_{j = 0}^infty b_j x^j})
converge absolutely by the comparison test.  If $c_n$ is as in
(\ref{c_n = sum_{j = 0}^n a_j b_{n - j}, 3}), then it is easy to see
that
\begin{equation}
\label{|c_n| r^n le sum_{j = 0}^n (|a_j| r^j) (|b_{n - j}| r^{n - j})}
 |c_n| \, r^n \le \sum_{j = 0}^n (|a_j| \, r^j) \, (|b_{n - j}| \, r^{n - j})
\end{equation}
for each $n \ge 0$.  The right side of (\ref{|c_n| r^n le sum_{j =
    0}^n (|a_j| r^j) (|b_{n - j}| r^{n - j})}) corresponds exactly to
the Cauchy product of the series in (\ref{sum_{j = 0}^infty |a_j| r^j,
  sum_{j = 0}^infty |b_j| r^j}), so that
\begin{equation}
\label{sum_{n = 0}^infty |c_n| r^n le ...}
 \sum_{n = 0}^\infty |c_n| \, r^n \le \Big(\sum_{j = 0}^\infty |a_j| \, r^j\Big)
                                  \, \Big(\sum_{l = 0}^\infty |b_l| \, r^l\Big),
\end{equation}
as in Section \ref{cauchy products}.  In particular, the left side of
(\ref{sum_{n = 0}^infty |c_n| r^n le ...}) converges, and we have that
\begin{equation}
\label{f(x) g(x) = sum_{n = 0}^infty c_n x^n}
        f(x) \, g(x) = \sum_{n = 0}^\infty c_n \, x^n
\end{equation}
for every $x \in k$ with $|x| \le r$, as in Sections \ref{cauchy
  products} and \ref{radius of convergence}.

        Let us ask that
\begin{equation}
\label{sum j |a_j| r^{j - 1}, sum j |b_j| r^{j - 1}}
        \sum_{j = 1}^\infty j \, |a_j| \, r^{j - 1}, \quad
        \sum_{j = 1}^\infty j \, |b_j| \, r^{j - 1}
\end{equation}
converge, which implies the convergence of the series in (\ref{sum_{j
    = 0}^infty |a_j| r^j, sum_{j = 0}^infty |b_j| r^j}).  As in
Section \ref{differentiation, continued}, this also implies that $f$
and $g$ are differentiable as functions defined on the closed ball
$\overline{B}(0, r)$ in $k$, with
\begin{equation}
\label{f'(x) = sum j cdot a_j x^{j - 1}, g'(x) = sum j cdot b_j x^{j - 1}}
        f'(x) = \sum_{j = 1}^\infty j \cdot a_j \, x^{j - 1}, \quad
        g'(x) = \sum_{j = 1}^\infty j \cdot b_j \, x^{j - 1}
\end{equation}
for every $x \in k$ with $|x| \le r$.  Observe that
\begin{eqnarray}
\label{n |c_n| r^{n - 1} le ...}
 n \, |c_n| \, r^{n - 1} 
 & \le & \sum_{j = 1}^n (j \, |a_j| \, r^{j - 1}|) \, (|b_{n - j}| \, r^{n -j}) \\
   & & + \sum_{j = 0}^{n - 1} (|a_j| \, r^j) \,
                           ((n - j) \, |b_{n - j}| \, r^{n - j - 1})  \nonumber
\end{eqnarray}
for each $n \ge 1$, using (\ref{n cdot c_n = sum_{j = 0}^n n cdot a_j
  b_{n - j} = ...}) and the triangle inequality.  The two sums on the
right side of (\ref{n |c_n| r^{n - 1} le ...}) may be considered as
Cauchy products, which can be summed over $n$ to obtain that
\begin{eqnarray}
\label{sum_{n = 1}^infty n |c_n| r^{n - 1} le ...}
 \sum_{n = 1}^\infty n \, |c_n| \, r^{n - 1} 
 & \le & \Big(\sum_{j = 1}^\infty j \, |a_j| \, r^{j - 1}\Big) \,
                                  \Big(\sum_{l = 0}^\infty |b_l| \, r^l\Big) \\
    & & + \Big(\sum_{j = 0}^\infty |a_j| \, r^j\Big) \,
           \Big(\sum_{l = 1}^\infty l \, |b_l| \, r^{l - 1}\Big). \nonumber
\end{eqnarray}
The sums on the right side of (\ref{sum_{n = 1}^infty n |c_n| r^{n -
    1} le ...}) are finite by hypothesis, so that the sum on the left
side of (\ref{sum_{n = 1}^infty n |c_n| r^{n - 1} le ...}) is finite
too.  Thus the discussion in Section \ref{differentiation, continued}
implies that $f \, g$ is differentiable as a function defined on
$\overline{B}(0, r)$, with
\begin{equation}
\label{(f g)'(x) = sum n cdot c_n x^{n - 1}}
        (f \, g)'(x) = \sum_{n = 1}^\infty n \cdot c_n \, x^{n - 1}
\end{equation}
for every $x \in k$ with $|x| \le r$, and where the series converges
absolutely by the comparison test.  This implies that
\begin{equation}
\label{(f g)'(x) = f'(x) g(x) + f(x) g'(x), 2}
        (f \, g)'(x) = f'(x) \, g(x) + f(x) \, g'(x)
\end{equation}
for every $x \in k$ with $|x| \le r$, by treating the products on the
right side of (\ref{(f g)'(x) = f'(x) g(x) + f(x) g'(x), 2}) as Cauchy
products, as in (\ref{(f g)'(X) = f'(X) g(X) + f(X) g'(X)}).  Note
that the discussion of the product rule in Section
\ref{differentiation} also implies that $f \, g$ is differentiable as
a function defined on $\overline{B}(0, r)$, with derivative as in
(\ref{(f g)'(x) = f'(x) g(x) + f(x) g'(x), 2}).

        Now let $k$ be a field with an ultrametric absolute value
function $|\cdot|$, and suppose that $k$ is complete with respect to
the ultrametric associated to $|\cdot|$.  Suppose also that
\begin{equation}
\label{lim_{j to infty} |a_j| r^j = lim_{j to infty} |b_j| r^j = 0}
        \lim_{j \to \infty} |a_j| \, r^j = \lim_{j \to \infty} |b_j| \, r^j = 0
\end{equation}
for some positive real number $r$, which implies that the series in
(\ref{f(x) = sum_{j = 0}^infty a_j x^j, g(x) = sum_{j = 0}^infty b_j
  x^j}) converge in $k$ for every $x \in k$ with $|x| \le r$.  If
$c_n$ is as in (\ref{c_n = sum_{j = 0}^n a_j b_{n - j}, 3}), then we
get that
\begin{equation}
\label{|c_n| r^n le max_{0 le j le n} (|a_j| |b_{n - j}|) r^n = ...}
 |c_n| \, r^n \le \max_{0 \le j \le n} (|a_j| \, |b_{n - j}|) \, r^n
               = \max_{0 \le j \le n} ((|a_j| \, r^j) \, (|b_{n - j}| \, r^{n - j}))
\end{equation}
for each $n \ge 0$.  It follows from this and (\ref{lim_{j to infty}
  |a_j| r^j = lim_{j to infty} |b_j| r^j = 0}) that
\begin{equation}
\label{lim_{n to infty} |c_n| r^n = 0}
        \lim_{n \to \infty} |c_n| \, r^n = 0,
\end{equation}
and hence that the series on the right side of (\ref{f(x) g(x) =
  sum_{n = 0}^infty c_n x^n}) converges in $k$ for every $x \in k$
with $|x| \le r$.  The right side of (\ref{f(x) g(x) = sum_{n =
    0}^infty c_n x^n}) is the same as the Cauchy product of the series
in (\ref{f(x) = sum_{j = 0}^infty a_j x^j, g(x) = sum_{j = 0}^infty
  b_j x^j}), so that (\ref{f(x) g(x) = sum_{n = 0}^infty c_n x^n})
holds for every $x \in k$ with $|x| \le r$, as in Section \ref{cauchy
  products}.  If $|\cdot|$ is not the trivial absolute value function
on $k$, then the discussion in Section \ref{differentiation,
  continued} implies that $f$, $g$, and $f \, g$ are differentiable as
functions on $\overline{B}(0, r)$ in $k$, with derivatives given by
the series in (\ref{f'(x) = sum j cdot a_j x^{j - 1}, g'(x) = sum j
  cdot b_j x^{j - 1}}) and (\ref{(f g)'(x) = sum n cdot c_n x^{n -
    1}}).  In this case, the convergence of these series for $x \in k$
with $|x| \le r$ follows from (\ref{lim_{j to infty} |a_j| r^j =
  lim_{j to infty} |b_j| r^j = 0}) and (\ref{lim_{n to infty} |c_n|
  r^n = 0}), as in Section \ref{differentiation, continued}.  As
before, one can derive (\ref{(f g)'(x) = f'(x) g(x) + f(x) g'(x), 2})
from (\ref{(f g)'(x) = sum n cdot c_n x^{n - 1}}), by treating the
product on the right side of (\ref{(f g)'(x) = f'(x) g(x) + f(x)
  g'(x), 2}) as Cauchy products.  One can also use the discussion of
the product rule in Section \ref{differentiation} to get that $f \, g$
is differentiable as a function defined on $\overline{B}(0, r)$, with
derivative as in (\ref{(f g)'(x) = f'(x) g(x) + f(x) g'(x), 2}), as in
the previous situation.

\section{The chain rule}
\label{chain rule}

        Let $k$ be a field again, and let
\begin{equation}
\label{f(X) = sum_{j = 0}^infty a_j X^j, g(Y) = sum_{l = 0}^infty b_l Y^l}
        f(X) = \sum_{j = 0}^\infty a_j \, X^j, \quad
        g(Y) = \sum_{l = 0}^\infty b_l \, Y^l
\end{equation}
be formal power series with coefficients in $k$.  As before, the
derivatives $f'(X)$, $g'(Y)$ of $f(X)$, $g(Y)$ are defined as formal
power series by
\begin{equation}
\label{f'(X) = sum_j j cdot a_j X^{j - 1}, g'(Y) = sum_l l cdot b_l Y^{l - 1}}
        f'(X) = \sum_{j = 1}^\infty j \cdot a_j \, X^{j - 1}, \quad
        g'(Y) = \sum_{l = 1}^\infty l \cdot b_l \, Y^{l - 1},
\end{equation}
as in (\ref{f'(X) = sum_{j = 1}^infty j cdot a_j X^{j - 1}}) in
Section \ref{differentiation, continued}.  Suppose for the moment that
$a_j = 0$ for all but finitely many $j$, or that $b_0 = 0$.  In both
cases, the composition
\begin{equation}
\label{(f circ g)(Y) = f(g(Y)) = sum_{j = 0}^infty a_j g(Y)^j}
        (f \circ g)(Y) = f(g(Y)) = \sum_{j = 0}^\infty a_j \, g(Y)^j
\end{equation}
of $f(X)$ and $g(Y)$ can be defined as a formal power series with
coefficients in $k$ too, as in Section \ref{compositions, continued}.
Similarly, the composition
\begin{equation}
\label{(f' circ g)(Y) = f'(g(Y)) = sum_{j = 1}^infty j cdot a_j g(Y)^{j - 1}}
 (f' \circ g)(Y) = f'(g(Y)) = \sum_{j = 1}^\infty j \cdot a_j \, g(Y)^{j - 1}
\end{equation}
of $f'(X)$ with $g(Y)$ can be defined as a formal power series with
coefficients in $k$ as well.  Put
\begin{equation}
\label{(g^j)(Y) = g(Y)^j}
        (g^j)(Y) = g(Y)^j
\end{equation}
for each positive integer $j$, and observe that
\begin{equation}
\label{(g^j)'(Y) = j cdot g(Y)^{j - 1} g'(Y)}
        (g^j)'(Y) = j \cdot g(Y)^{j - 1} \, g'(Y)
\end{equation}
by the product rule, where $(g^j)'(Y)$ is the derivative of
$(g^j)(Y)$.  Using this, one can verify that
\begin{equation}
\label{(f circ g)'(Y) = f'(g(Y)) g'(Y)}
        (f \circ g)'(Y) = f'(g(Y)) \, g'(Y),
\end{equation}
where $(f \circ g)'(Y)$ is the derivative of $(f \circ g)(Y)$.  Of
course, this is a version of the chain rule\index{chain rule} for
formal power series.

        Suppose now for the moment that $a_j = 0$ for all but
finitely many $j$, and that $b_l = 0$ for all but finitely many $l$.
Thus
\begin{equation}
\label{f(x) = sum_{j = 0}^infty a_j x^j, g(y) = sum_{l = 0}^infty b_l y^l}
        f(x) = \sum_{j = 0}^\infty a_j \, x^j, \quad
        g(y) = \sum_{l = 0}^\infty b_l \, y^l
\end{equation}
are defined for all $x, y \in k$, and hence their composition
\begin{equation}
\label{(f circ g)(y) = f(g(y)) = sum_{j = 0}^infty a_j g(y)^j}
        (f \circ g)(y) = f(g(y)) = \sum_{j = 0}^\infty a_j \, g(y)^j
\end{equation}
is defined for every $y \in k$ too.  If $k$ is equipped with a
nontrivial absolute value function $|\cdot|$, then the derivatives of
$f$ and $g$ can be defined as in Section \ref{differentiation}, and
are given by
\begin{equation}
\label{f'(x) = sum_j j cdot a_j x^{j - 1}, g'(y) = sum_l l cdot b_l y^{l - 1}}
        f'(x) = \sum_{j = 1}^\infty j \cdot a_j \, x^{j - 1}, \quad
        g'(y) = \sum_{l = 1}^\infty l \cdot b_l \, y^{l - 1}
\end{equation}
for every $x, y \in k$.  Similarly, the composition
\begin{equation}
\label{(f' circ g)(y) = f'(g(y)) = sum_{j = 1}^infty j cdot a_j g(y)^{j - 1}}
 (f' \circ g)(y) = f'(g(y)) = \sum_{j = 1}^\infty j \cdot a_j \, g(y)^{j - 1}
\end{equation}
of $f'$ with $g$ is defined for every $y \in k$ under these
conditions.  The usual version of the chain rule implies that
\begin{equation}
\label{(f circ g)'(y) = f'(g(y)) g'(y)}
        (f \circ g)'(y) = f'(g(y)) \, g'(y)
\end{equation}
for every $y \in k$, where $(f \circ g)'(y)$ is the derivative of $f
\circ g$ at $y$.  More precisely, the discussion of the chain rule in
Section \ref{differentiation} implies that the value of the derivative
of $f \circ g$ at any point $y \in k$ is given as in (\ref{(f circ
  g)'(y) = f'(g(y)) g'(y)}).  Alternatively, $(f \circ g)(y)$ may
be expressed as a polynomial in $y$, whose derivative is the same as
the right side of (\ref{(f circ g)'(y) = f'(g(y)) g'(y)}) as a product
of polynomials in $y$, as in the previous paragraph.

        Now suppose that $|\cdot|$ is a nontrivial ultrametric
absolute value function on $k$, and that $k$ is complete with respect
to the ultrametric associated to $|\cdot|$.  Suppose also that $a_0,
a_1, a_2, a_3, \ldots$ and $b_0, b_1, b_2, b_3, \ldots$ are sequences
of elements of $k$ that satisfy (\ref{lim_{j to infty} |a_j| r^j = 0,
  2}), (\ref{lim_{l to infty} |b_l| t^l = 0}), and (\ref{max_{l ge 0}
  |b_l| t^l le r}) in Section \ref{compositions} for some $r, t > 0$.
This implies that $f(x)$ and $g(y)$ can be defined as in (\ref{f(x) =
  sum_{j = 0}^infty a_j x^j, g(y) = sum_{l = 0}^infty b_l y^l}) for
every $x, y \in k$ with $|x| \le r$ and $|y| \le t$, and that $|g(y)|
\le r$ for all such $y$.  This permits one to define $f(g(y))$ for
every $y \in k$ with $|y| \le t$, and the discussion in Section
\ref{compositions} shows that $f(g(y))$ is given by a power series in
$y$, with suitable convergence properties.  As in Section
\ref{differentiation, continued}, $f$ and $g$ are differentiable as
$k$-valued functions on the closed balls $\overline{B}(0, r)$ and
$\overline{B}(0, t)$ in $k$, respectively, with derivatives given as
in (\ref{f'(x) = sum_j j cdot a_j x^{j - 1}, g'(y) = sum_l l cdot b_l
  y^{l - 1}}).  These series for the derivatives have convergence
properties like those for $f$ and $g$, and in particular the discussion
in Section \ref{compositions} can be applied to $f'$ instead of $f$,
to get that $f'(g(y))$ can be expressed by a power series in $y$
with suitable convergence properties as well.  As before, the
discussion of the chain rule in Section \ref{differentiation}
implies that (\ref{(f circ g)'(y) = f'(g(y)) g'(y)}) holds
for every $y \in k$ with $|y| \le t$.  Note that the right side
of (\ref{(f circ g)'(y) = f'(g(y)) g'(y)}) is given by the
Cauchy product of the power series for $f'(g(y))$ and $g'(y)$,
which is a power series with the same type of convergence properties.
Similarly, the left side of (\ref{(f circ g)'(y) = f'(g(y)) g'(y)})
can be obtained by differentiating the power series expansion for
$(f \circ g)(y)$, as in Section \ref{differentiation, continued}.
One can check that (\ref{(f circ g)'(y) = f'(g(y)) g'(y)}) holds as an
equality between power series in $y$ as well, which is to say that the
power series on both sides of the equation have the same coefficients.
More precisely, this can be verified directly, but we shall not go
through the details here.  Otherwise, this can be obtained from the
fact that (\ref{(f circ g)'(y) = f'(g(y)) g'(y)}) holds for every $y
\in k$ with $|y| \le t$, since $t > 0$ and $|\cdot|$ is nontrivial on
$k$, and since both sides of (\ref{(f circ g)'(y) = f'(g(y)) g'(y)})

        Suppose now that $k = {\bf R}$ or ${\bf C}$, with the standard
absolute value function, and that $a_0, a_1, a_2, a_3, \ldots$ and
$b_0, b_1, b_2, b_3, \ldots$ are sequences of elements of $k$ that
satisfy (\ref{sum_{j = 0}^infty |a_j| r^j, 2}) and (\ref{sum_{l =
    0}^infty |b_l| t^l le r}) in Section \ref{compositions} for some
$r, t >0$.  As in the previous situation, this implies that $f(x)$ and
$g(y)$ can be defined as in (\ref{f(x) = sum_{j = 0}^infty a_j x^j,
  g(y) = sum_{l = 0}^infty b_l y^l}) for every $x, y \in k$ with $|x|
\le r$ and $|y| \le t$, and that $|g(y)| \le r$ for all such $y$.
Thus $f(g(y))$ is also defined for every $y \in k$ with $|y| \le t$,
and the discussion in Section \ref{compositions} shows that $f(g(y))$
is given by an absolutely convergent power series.  In order to deal
with derivatives, let us ask in addition that
\begin{equation}
\label{sum_j j |a_j| r^{j - 1}, sum_l l |b_l| t^{l - 1}}
        \sum_{j = 1}^\infty j \, |a_j| \, r^{j - 1}, \quad
         \sum_{l = 1}^\infty l \, |b_l| \, t^{j - 1}
\end{equation}
converge.  This implies that $f$ and $g$ are differentiable on the
closed balls $\overline{B}(0, r)$ and $\overline{B}(0, t)$ in $k$,
respectively, as in Section \ref{differentiation, continued}, with
derivatives given as in (\ref{f'(x) = sum_j j cdot a_j x^{j - 1},
  g'(y) = sum_l l cdot b_l y^{l - 1}}).  As before, the discussion in
Section \ref{compositions} can be applied to $f'$ instead of $f$, to
get that $f'(g(y))$ is given by an absolutely convergent power series
in $y$ too.  The discussion of the chain rule in Section
\ref{differentiation, continued} implies that $f \circ g$ is
differentiable on $\overline{B}(0, t)$, and that (\ref{(f circ g)'(y)
  = f'(g(y)) g'(y)}) holds for every $y \in k$ with $|y| \le t$.
However, one should be a bit more careful about the convergence of the
power series for the derivative of $f \circ g$ in this case.
Let $c_0, c_1, c_2, c_3, \ldots$ be the coefficients of the power
series expansion for $(f \circ g)(y)$, as in Section \ref{compositions}.
One can check that
\begin{equation}
\label{sum_{n = 1}^infty n |c_n| t^{n - 1}}
        \sum_{n = 1}^\infty n \, |c_n| \, t^{n - 1}
\end{equation}
converges under these conditions, using the convergence of the series
in (\ref{sum_j j |a_j| r^{j - 1}, sum_l l |b_l| t^{l - 1}}).  This can
be done directly, along with showing that (\ref{(f circ g)'(y) =
  f'(g(y)) g'(y)}) holds as an equality between power series.
Alternatively, the convergence of the power series for $(f \circ
g)(y)$ when $|y| \le t$ implies that the power series for the
derivative of $(f \circ g)(y)$ converges when $|y| < t$.  It follows
that the derivative of $(f \circ g)(y)$ is given by this power series
when $|y| < t$, as in Section \ref{differentiation, continued}.  Both
factors on the right side of (\ref{(f circ g)'(y) = f'(g(y)) g'(y)})
are already given by power series that converge absolutely when $|y|
\le t$, and so their product has the same property.  As before, the
power series on both sides of (\ref{(f circ g)'(y) = f'(g(y)) g'(y)})
have to have the same coefficients, because (\ref{(f circ g)'(y) =
  f'(g(y)) g'(y)}) holds for all $y \in k$ with $|y| < t$ and $t > 0$.
This permits the convergence of (\ref{sum_{n = 1}^infty n |c_n| t^{n -
    1}}) to be derived from the absolute convergence of the series on
the right side of (\ref{(f circ g)'(y) = f'(g(y)) g'(y)}) when $|y| =
t$.  Of course, the convergence of (\ref{sum_{n = 1}^infty n |c_n|
  t^{n - 1}}) implies that the derivative of $(f \circ g)(y)$ on
$\overline{B}(0, t)$ is given by the corresponding power series for
every $y \in k$ with $|y| \le t$, as in Section \ref{differentiation,
  continued}.

\section{Functions of sums}
\label{functions of sums}

        Let $k$ be a field, and let
\begin{equation}
\label{f(x) = sum_{j = 0}^infty a_j x^j, 4}
        f(x) = \sum_{j = 0}^\infty a_j \, x^j
\end{equation}
be a power series with coefficients in $k$.  Also let $B$ be a
nonempty set, and let $b$ be a $k$-valued function on $B$.  We would
like to consider
\begin{equation}
\label{f(sum_{l in B} b(l)) = sum_{j = 0}^infty a_j (sum_{l in B} b(l))^j}
        f\Big(\sum_{l \in B} b(l)\Big)
          = \sum_{j = 0}^\infty a_j \, \Big(\sum_{l \in B} b(l)\Big)^j,
\end{equation}
at least formally at first.  Of course, there is no problem with this
when $a_j = 0$ for all but finitely many $j \ge 0$, and $b(l) = 0$ for
all but finitely many $l \in B$.  This is analogous to the discussion
in Section \ref{compositions}, which corresponds to the case where $B
= {\bf Z}_+ \cup \{0\}$ and $b(l) = b_l \, y^l$ for some $b_l, y \in
k$.

        As in Section \ref{compositions}, let
\begin{equation}
\label{E_j = B^j}
        E_j = B^j
\end{equation}
be the $j$th Cartesian power of $B$ for each positive integer $j$,
consisting of $j$-tuples $\alpha = (\alpha_1, \ldots, \alpha_j)$ of
elements of $B$.  Put
\begin{equation}
\label{beta_j(alpha) = b(alpha_1) b(alpha_2) cdots b(alpha_j)}
        \beta_j(\alpha) = b(\alpha_1) \, b(\alpha_2) \cdots b(\alpha_j)
\end{equation}
for each $j \in {\bf Z}_+$ and $\alpha \in E_j$, so that
\begin{equation}
\label{(sum_{l in B} b(l))^j = sum_{alpha in E_j} beta_j(alpha)}
 \Big(\sum_{l \in B} b(l)\Big)^j = \sum_{\alpha \in E_j} \beta_j(\alpha),
\end{equation}
at least formally again.  As before, there is no problem with this
when $b(l) = 0$ for all but finitely many $l \in B$.  Note that the
sets $E_j$ are considered to be pairwise disjoint.

        As in Section \ref{compositions} again, we let $E_0$ be a set
with exactly one element, not contained in $E_j$ for any $j \ge 1$.
Thus the sets $E_j$ are pairwise disjoint for all $j \ge 0$, and we put
\begin{equation}
\label{E = bigcup_{j = 0}^infty E_j, 2}
        E = \bigcup_{j = 0}^\infty E_j.
\end{equation}
Let $\phi$ be the $k$-valued function on $E$ defined by
\begin{equation}
\label{phi(alpha) = a_j beta_j(alpha), 2}
        \phi(\alpha) = a_j \, \beta_j(\alpha)
\end{equation}
for each $\alpha \in E_j$ when $j \ge 1$, and $\phi = a_0$ on $E_0$.
Combining (\ref{f(sum_{l in B} b(l)) = sum_{j = 0}^infty a_j (sum_{l
    in B} b(l))^j}) and (\ref{(sum_{l in B} b(l))^j = sum_{alpha in
    E_j} beta_j(alpha)}), we get that
\begin{equation}
\label{f(sum_{l in B} b(l)) = ... = sum_{alpha in E} phi(alpha)}
        f\Big(\sum_{l \in B} b(l)\Big)
  = \sum_{j = 0}^\infty a_j \, \Big(\sum_{\alpha \in E_j} \beta_j(\alpha)\Big)
             = \sum_{\alpha \in E} \phi(\alpha),
\end{equation}
at least formally again.  As usual, if $a_j = 0$ for all but finitely
many $j \ge 0$, and $b(l) = 0$ for all but finitely many $l \in B$,
then $\phi \in c_{00}(E, k)$, and there is no problem with
(\ref{f(sum_{l in B} b(l)) = ... = sum_{alpha in E} phi(alpha)}).

        Suppose for the moment that $k = {\bf R}$ or ${\bf C}$,
with the standard absolute value function.  Suppose also that
\begin{equation}
\label{sum_{j = 0}^infty |a_j| r^j, 5}
        \sum_{j = 0}^\infty |a_j| \, r^j
\end{equation}
converges for some nonnegative real number $r$, and that
\begin{equation}
\label{sum_{l in B} |b(l)| le r}
        \sum_{l \in B} |b(l)| \le r,
\end{equation}
where the sum on the left side of (\ref{sum_{l in B} |b(l)| le r}) can
be defined as in Section \ref{nonnegative sums}.  The convergence of
(\ref{sum_{j = 0}^infty |a_j| r^j, 5}) implies that the right side of
(\ref{f(x) = sum_{j = 0}^infty a_j x^j, 4}) converges absolutely for
every $x \in k$ with $|x| \le r$.  The finiteness of the sum on the
left side of (\ref{sum_{l in B} |b(l)| le r}) means that $b(l)$ is
summable as a $k$-valued function on $B$, as in Section \ref{ell^r
  norms}.  This implies that $\sum_{l \in B} b(l)$ can be defined as
in Sections \ref{generalized convergence} and \ref{generalized
  convergence, continued}, and that
\begin{equation}
\label{|sum_{l in B} b(l)| le sum_{l in B} |b(l)| le r}
        \biggl|\sum_{l \in B} b(l)\biggr| \le \sum_{l \in B} |b(l)| \le r,
\end{equation}
as in (\ref{N(sum_{x in X} f(x)) le sum_{x in X} N(f(x))}) in Section
\ref{generalized convergence, continued}.

        It follows that the series in $j$ on the right side of
(\ref{f(sum_{l in B} b(l)) = sum_{j = 0}^infty a_j (sum_{l in B} b(l))^j})
converges absolutely under these conditions, which can be used to
define the left side of (\ref{f(sum_{l in B} b(l)) = sum_{j = 0}^infty
  a_j (sum_{l in B} b(l))^j}).  We also have that
\begin{equation}
\label{sum_{alpha in E_j} |beta_j(alpha)| = ... le r^j}
        \sum_{\alpha \in E_j} |\beta_j(\alpha)|
           = \Big(\sum_{l \in B} |b(l)|\Big)^j \le r^j
\end{equation}
for each $j \ge 0$, so that
\begin{equation}
\label{sum_{alpha in E} |phi(alpha)| = ... le sum_{j = 0}^infty |a_j| r^j}
        \sum_{\alpha \in E} |\phi(\alpha)|
          = |a_0| + \sum_{j = 1}^\infty \Big(\sum_{\alpha \in E_j} |a_j| \,
                                                      |\beta_j(\alpha)|\Big)
           \le \sum_{j = 0}^\infty |a_j| \, r^j.
\end{equation}
In particular, (\ref{sum_{alpha in E_j} |beta_j(alpha)| = ... le r^j})
implies that $\beta_j$ is summable on $E_j$ for each $j \ge 0$, and
(\ref{sum_{alpha in E} |phi(alpha)| = ... le sum_{j = 0}^infty |a_j|
  r^j}) implies that $\phi$ is summable on $E$.  The summability of
$\beta_j$ on $E_j$ means that the right side of (\ref{(sum_{l in B}
  b(l))^j = sum_{alpha in E_j} beta_j(alpha)}) can be defined as in
Sections \ref{generalized convergence} and \ref{generalized
  convergence, continued}, and this sum can be evaluated as an
iterated sum to get (\ref{(sum_{l in B} b(l))^j = sum_{alpha in E_j}
  beta_j(alpha)}).  Similarly, the summability of $\phi$ on $E$ means
that the right side of (\ref{f(sum_{l in B} b(l)) = ... = sum_{alpha
    in E} phi(alpha)}) can be defined as in Sections \ref{generalized
  convergence} and \ref{generalized convergence, continued}, and the
second step in (\ref{f(sum_{l in B} b(l)) = ... = sum_{alpha in E}
  phi(alpha)}) can be obtained as in Section \ref{sums of sums}.

        Now let $k$ be any field with an ultrametric absolute value
function $|\cdot|$, where $k$ is complete with respect to the
corresponding ultrametric.  Suppose that
\begin{equation}
\label{lim_{j to infty} |a_j| r^j = 0, 5}
        \lim_{j \to \infty} |a_j| \, r^j = 0
\end{equation}
for some nonnegative real number $r$,
\begin{equation}
\label{b vanishes at infinity on B}
        b \hbox{ vanishes at infinity on } B,
\end{equation}
and
\begin{equation}
\label{max_{l in B} |b(l)| le r}
        \max_{l \in B} |b(l)| \le r.
\end{equation}
Thus the series on the right side of (\ref{f(x) = sum_{j = 0}^infty
  a_j x^j, 4}) converges in $k$ for every $x \in k$ with $|x| \le r$,
because of (\ref{lim_{j to infty} |a_j| r^j = 0, 5}) and the
completeness of $k$.  Similarly, $\sum_{l \in B} b(l)$ can be defined
as an element of $k$ as in Sections \ref{generalized convergence} and
\ref{generalized convergence, continued}, because of (\ref{b vanishes
  at infinity on B}) and the completeness of $k$.  We also have that
\begin{equation}
\label{|sum_{l in B} b(l)| le max_{l in B} |b(l)| le r}
        \biggl|\sum_{l \in B} b(l)\biggr| \le \max_{l \in B} |b(l)| \le r,
\end{equation}
as in (\ref{N(sum_{x in X} f(x)) le max_{x in X} N(f(x))}) in Section
\ref{generalized convergence, continued}.

        As before, the series on the right side of
(\ref{f(sum_{l in B} b(l)) = sum_{j = 0}^infty a_j (sum_{l in B} b(l))^j})
converges in $k$ in this situation, which can be used to define the
left side of (\ref{f(sum_{l in B} b(l)) = sum_{j = 0}^infty a_j
  (sum_{l in B} b(l))^j}).  It is easy to see that
\begin{equation}
\label{beta_j vanishes at infinity on E_j}
        \beta_j \hbox{ vanishes at infinity on } E_j
\end{equation}
for each $j \ge 1$, because of (\ref{b vanishes at infinity on B}),
and that
\begin{equation}
\label{max_{alpha in E_j} |beta_j(alpha)| = (max_{l in B} |b(l)|)^j le r^j}
        \max_{\alpha \in E_j} |\beta_j(\alpha)|
              = \Big(\max_{l \in B} |b(l)|\Big)^j \le r^j,
\end{equation}
by (\ref{max_{l in B} |b(l)| le r}).  This implies that
\begin{equation}
\label{max_{alpha in E_j} |phi(alpha)| = ... le |a_j| r^j}
        \max_{\alpha \in E_j} |\phi(\alpha)|
 = |a_j| \, \max_{\alpha \in E_j} |\beta_j(\alpha)| \le |a_j| \, r^j
\end{equation}
for each $j \ge 1$, which tends to $0$ as $j \to \infty$,
by (\ref{lim_{j to infty} |a_j| r^j = 0, 5}).  It follows that
\begin{equation}
\label{phi vanishes at infinity on E}
        \phi \hbox{ vanishes at infinity on } E,
\end{equation}
since (\ref{beta_j vanishes at infinity on E_j}) implies that the
restriction of $\phi$ to $E_j$ vanishes at infinity for each $j \ge
1$.  As in Section \ref{generalized convergence}, (\ref{beta_j
  vanishes at infinity on E_j}) implies that
\begin{equation}
\label{sum_{alpha in E_j} beta_j(alpha) satisfies the gcc}
        \sum_{\alpha \in E_j} \beta_j(\alpha)
               \hbox{ satisfies the generalized Cauchy criterion}
\end{equation}
for each $j \ge 1$, and similarly (\ref{phi vanishes at infinity on
  E}) implies that
\begin{equation}
\label{sum_{alpha in E} phi(alpha) satisfies the gcc}
        \sum_{\alpha \in E} \phi(\alpha)
               \hbox{ satisfies the generalized Cauchy criterion}.
\end{equation}
Thus these sums can be defined in $k$, because $k$ is complete, as in
Section \ref{generalized convergence, continued}.  One can also check
that (\ref{(sum_{l in B} b(l))^j = sum_{alpha in E_j} beta_j(alpha)})
and (\ref{f(sum_{l in B} b(l)) = ... = sum_{alpha in E} phi(alpha)})
hold under these conditions, using the remarks in Section \ref{sums of
  sums}.

\section{The logarithm}
\label{logarithm}

        Let $k$ be a field of characteristic $0$, and consider
\begin{equation}
\label{log (1 + X) = sum_{j = 1}^infty frac{(-1)^{j + 1}}{j} X^j}
        \log (1 + X) = \sum_{j = 1}^\infty \frac{(-1)^{j + 1}}{j} \, X^j
\end{equation}
as a formal power series with coefficients in $k$, where $X$ is an
indeterminate.  If $k$ is equipped with an absolute value function
$|\cdot|$, then we may put
\begin{equation}
\label{log (1 + x) = sum_{j = 1}^infty frac{(-1)^{j + 1}}{j} x^j}
        \log (1 + x) = \sum_{j = 1}^\infty \frac{(-1)^{j + 1}}{j} \, x^j
\end{equation}
for every $x \in k$ such that the series on the right side of
(\ref{log (1 + x) = sum_{j = 1}^infty frac{(-1)^{j + 1}}{j} x^j})
converges.  In particular, this series converges when $x = 0$, so that
\begin{equation}
\label{log 1 = 0}
        \log 1 = 0.
\end{equation}
Of course, this is the usual power series of the
logarithm,\index{logarithm} so that (\ref{log (1 + x) = sum_{j =
    1}^infty frac{(-1)^{j + 1}}{j} x^j}) may be considered as the
definition of a logarithm function on a subset of $k$.  Note that the
formal derivative of (\ref{log (1 + X) = sum_{j = 1}^infty
  frac{(-1)^{j + 1}}{j} X^j}) is given by
\begin{equation}
\label{sum_{j = 1}^infty (-1)^{j + 1} X^{j - 1} = sum_{j = 0}^infty (-1)^j X^j}
 \sum_{j = 1}^\infty (-1)^{j + 1} \, X^{j - 1} = \sum_{j = 0}^\infty (-1)^j \, X^j,
\end{equation}
which is the power series associated to
\begin{equation}
\label{frac{1}{1 + X}}
        \frac{1}{1 + X}.
\end{equation}

        Suppose for the moment that $k = {\bf R}$ or ${\bf C}$, with
the standard absolute value function.  In this case, it is well known
that radius of convergence of (\ref{log (1 + X) = sum_{j = 1}^infty
  frac{(-1)^{j + 1}}{j} X^j}) is equal to $1$.  More precisely, if $x
\in k$ and $|x| < 1$, then the series on the right side of (\ref{log
  (1 + x) = sum_{j = 1}^infty frac{(-1)^{j + 1}}{j} x^j}) converges
absolutely, by comparison with the convergent series
\begin{equation}
\label{sum_{j = 1}^infty |x|^j}
        \sum_{j = 1}^\infty |x|^j
\end{equation}
of nonnegative real numbers.  However, the series on the right side of
(\ref{log (1 + x) = sum_{j = 1}^infty frac{(-1)^{j + 1}}{j} x^j}) does
not converge absolutely when $|x| = 1$, which is the same in this case
as saying that it does not converge when $x = -1$.  If $x = 1$, then
the right side of (\ref{log (1 + x) = sum_{j = 1}^infty frac{(-1)^{j +
      1}}{j} x^j}) does converge, by Leibniz' alternating series test.
Similarly, if $x \in {\bf C}$, $|x| = 1$, and $x \ne -1$, then one can
show that the right side of (\ref{log (1 + x) = sum_{j = 1}^infty
  frac{(-1)^{j + 1}}{j} x^j}) converges in ${\bf C}$, as in Theorem
3.44 on p71 of \cite{r1}.  Of course, the logarithm can be defined in
other ways for all positive real numbers, and it can be extended
holomorphically to suitable domains in the complex plane.  These
extensions satisfy
\begin{equation}
\label{frac{d}{dz} log z = frac{1}{z}}
        \frac{d}{dz} \, \log z = \frac{1}{z},
\end{equation}
and indeed this can be used to define the logarithm, together with
$\log 1 = 0$.

        Now let $k$ be a field of characteristic $0$ with an ultrametric
absolute value function $|\cdot|$, and suppose that $k$ is complete
with respect to the ultrametric associated to $|\cdot|$.  Thus the
series on the right side of (\ref{log (1 + x) = sum_{j = 1}^infty
  frac{(-1)^{j + 1}}{j} x^j}) converges in $k$ exactly when the terms
of the series converge to $0$ in $k$.  Remember that $|\cdot|$ induces
an ultrametric absolute value function on ${\bf Q}$, using the natural
embedding of ${\bf Q}$ in $k$.  If the induced absolute value function
on ${\bf Q}$ is trivial, then
\begin{equation}
\label{|(-1)^{j + 1} x^j / j| = |x|^j}
        |(-1)^{j + 1} \, x^j / j| = |x|^j
\end{equation}
for every $x \in k$ and $j \in {\bf Z}_+$.  In this case, the series
on the right side of (\ref{log (1 + x) = sum_{j = 1}^infty
  frac{(-1)^{j + 1}}{j} x^j}) converges exactly when $|x| < 1$.
Suppose for the moment that $|\cdot|$ is not the trivial absolute
value function on $k$, so that every point in the open unit ball $B(0,
1)$ in $k$ is a limit point of $B(0, 1)$.  As in Section
\ref{differentiation, continued}, the derivative of $\log (1 + x)$ as
a $k$-valued function defined on $B(0, 1)$ exists at every point in
$B(0, 1)$, and is given by the corresponding power series for the
derivative.  Thus the derivative of $\log (1 + x)$ is equal to
\begin{equation}
\label{frac{1}{1 + x}}
        \frac{1}{1 + x}
\end{equation}
for every $x \in k$ with $|x| < 1$, since the power series for the
derivative of $\log (1 + x)$ corresponds to the usual power series for
(\ref{frac{1}{1 + x}}), as before.  Of course, if $|\cdot|$ is the
trivial absolute value function on $k$, then $\log (1 + x)$ is defined
only for $x = 0$, and so the derivative is not defined.

        If the induced absolute value function on ${\bf Q}$ is not trivial,
then there is a prime number $p$ such that the induced absolute value
function on ${\bf Q}$ is equivalent to the $p$-adic absolute value
function, by Ostrowski's theorem, as in Section \ref{ostrowski's
  theorems}.  In this case, we may as well ask that the induced
absolute value function on ${\bf Q}$ be equal to the $p$-adic absolute
value, since this can be arranged by replacing the given absolute
value function $|\cdot|$ on $k$ by a suitable positive power of
itself.  This implies that
\begin{equation}
\label{|(-1)^{j + 1} x^j / j| = |x|^j / |j|_p}
        |(-1)^{j + 1} \, x^j / j| = |x|^j / |j|_p
\end{equation}
for every $x \in k$ and $j \in {\bf Z}_+$, where $|j|_p$ is the
$p$-adic absolute value of $j$.  It is easy to see that
\begin{equation}
\label{1/j le |j|_p le 1}
        1/j \le |j|_p \le 1
\end{equation}
for every $j \in {\bf Z}_+$, so that
\begin{equation}
\label{|x|^j le |(-1)^{j + 1} x^j / j| le j |x|^j}
        |x|^j \le |(-1)^{j + 1} \, x^j / j| \le j \, |x|^j
\end{equation}
for every $x \in k$ and $j \in {\bf Z}_+$.  It follows that
(\ref{|(-1)^{j + 1} x^j / j| = |x|^j / |j|_p}) tends to $0$ as $j \to
\infty$ exactly when $|x| < 1$ in this situation.  Equivalently, this
means that the right side of (\ref{log (1 + x) = sum_{j = 1}^infty
  frac{(-1)^{j + 1}}{j} x^j}) converges in $k$ exactly when $|x| < 1$
under these conditions.  As in the preceding paragraph, the derivative
of $\log (1 + x)$ as a $k$-valued function on the open unit ball $B(0,
1)$ in $k$ exists at every point in $B(0, 1)$, and is equal to
(\ref{frac{1}{1 + x}}), by the discussion in Section
\ref{differentiation, continued}.

\section{The usual identity}
\label{usual identity}

        If $u$, $v$ are positive real numbers, then
\begin{equation}
\label{log (u v) = log u + log v}
        \log (u \, v) = \log u + \log v,
\end{equation}
where the logarithm refers to the standard real-valued logarithm
function on the positive half-line.  There are analogous statements
for complex numbers, but one should be careful about the conditions
under which they hold.

        Now let $k$ be a field with characteristic $0$, and let
$Y$ and $Z$ be commuting indeterminates.  Observe that
\begin{equation}
\label{(1 + Y) (1 + Z) = 1 + Y + Z + Y Z}
        (1 + Y) \, (1 + Z) = 1 + Y + Z + Y \, Z,
\end{equation}
so that
\begin{eqnarray}
\label{log ((1 + Y) (1 + Z)) = log (1 + Y + Z + Y Z) = ...}
        \log ((1 + Y) \, (1 + Z)) & = & \log (1 + Y + Z + Y \, Z) \\
 & = & \sum_{j = 1}^\infty \frac{(-1)^{j + 1}}{j} \, (Y + Z + Y \, Z)^j, \nonumber
\end{eqnarray}
at least formally.  More precisely, the right side of (\ref{log ((1 +
  Y) (1 + Z)) = log (1 + Y + Z + Y Z) = ...}) can be defined as
a formal power series in $Y$ and $Z$, by expanding
\begin{equation}
\label{(Y + Z + Y Z)^j}
        (Y + Z + Y \, Z)^j
\end{equation}
into a finite sum of monomials for each $j$.  Each of the monomials
that occurs in the sum has total degree in $Y$ and $Z$ greater than or
equal to $j$, and less than or equal to $2 \, j$.  Any particular
monomial in $Y$ and $Z$ can occur in (\ref{(Y + Z + Y Z)^j}) for only
finitely many $j$, and hence the coefficient of such a monomial in the
right side of (\ref{log ((1 + Y) (1 + Z)) = log (1 + Y + Z + Y Z) =
  ...})  is given by a finite sum in $k$.

        In analogy with (\ref{log (u v) = log u + log v}), we have that
\begin{equation}
\label{log ((1 + Y) (1 + Z)) = log (1 + Y) + log (1 + Z)}
        \log ((1 + Y) \, (1 + Z)) = \log (1 + Y) + \log (1 + Z),
\end{equation}
as an equality between formal power series in $Y$ and $Z$.  One way to
look at this is to start with the case where $k = {\bf R}$ or ${\bf
  C}$, where the identity for the corresponding functions near $1$
implies the appropriate identities between the power series
coefficients.  These identities between the power series coefficients
are simply statements about finite sums of rational numbers, which
carry over to any field of characteristic $0$.  Alternatively, one can
check that the formal derivatives of both sides of (\ref{log ((1 + Y)
  (1 + Z)) = log (1 + Y) + log (1 + Z)}) in $Y$ and $Z$ are the same,
as formal power series in $Y$ and $Z$.  This implies (\ref{log ((1 +
  Y) (1 + Z)) = log (1 + Y) + log (1 + Z)}), because the constant
terms on both sides of the equation are equal to $0$, and because
$k$ has characteristic zero, so that the non-constant terms
can be recovered from their first derivatives.

        If $|\cdot|$ is an absolute value function on $k$, then
\begin{equation}
\label{{w in k : |w| = 1}}
        \{w \in k : |w| = 1\}
\end{equation}
is a group with respect to multiplication.  Suppose now that $|\cdot|$
is an ultrametric absolute value function on $k$, and let us check that
\begin{equation}
\label{{w in k : |w - 1| < 1}}
        \{w \in k : |w - 1| < 1\}
\end{equation}
is a subgroup of (\ref{{w in k : |w| = 1}}).  If $w \in k$ and $|w -
1| < 1$, then it is easy to see that $|w| = 1$, as in (\ref{|x + y| =
  |y|}) in Section \ref{absolute value functions}.  We also have that
\begin{equation}
\label{1/w - 1 = (1 - w)/w}
        1/w - 1 = (1 - w)/w,
\end{equation}
so that
\begin{equation}
\label{|1/w - 1| = |1 - w|/|w| = |1 - w| < 1}
        |1/w - 1| = |1 - w|/|w| = |1 - w| < 1,
\end{equation}
which means that $1/w$ is an element of (\ref{{w in k : |w - 1| < 1}})
too.  Note that
\begin{equation}
\label{(1 + y) (1 + z) = 1 + y + z + y z}
        (1 + y) \, (1 + z) = 1 + y + z + y \, z
\end{equation}
for every $y, z \in k$, as in (\ref{(1 + Y) (1 + Z) = 1 + Y + Z + Y
  Z}), so that
\begin{equation}
\label{|(1 + y) (1 + z) - 1| = |y + z + y z| le max(|y|, |z|, |y| |z|)}
 |(1 + y) \, (1 + z) - 1| = |y + z + y \, z| \le \max(|y|, |z|, |y| \, |z|).
\end{equation}
If $|y|, |z| < 1$, then it follows that
\begin{equation}
\label{|(1 + y) (1 + z) - 1| < 1}
        |(1 + y) \, (1 + z) - 1| < 1.
\end{equation}
This implies that (\ref{{w in k : |w - 1| < 1}}) is closed under
multiplication, as desired.

        Let us continue to suppose that $|\cdot|$ is an ultrametric
absolute value function on $k$, and let us also ask that $k$ be
complete with respect to the ultrametric associated to $|\cdot|$.  If
$x \in k$ and $|x| < 1$, then $\log (1 + x)$ can be defined in $k$ by
(\ref{log (1 + x) = sum_{j = 1}^infty frac{(-1)^{j + 1}}{j} x^j}), as
in the previous section.  Similarly, if $y, z \in k$ satisfy $|y|, |z|
< 1$, then
\begin{eqnarray}
\label{log ((1 + y) (1 + z)) = log (1 + y + z + y z) = ...}
        \log ((1 + y) \, (1 + z)) & = & \log (1 + y + z + y \, z) \\
 & = & \sum_{j = 1}^\infty \frac{(-1)^{j + 1}}{j} \, (y + z + y \, z)^j \nonumber
\end{eqnarray}
can be defined in $k$ as before, because of (\ref{|(1 + y) (1 + z) -
  1| < 1}).  Under these conditions, one can show that
\begin{equation}
\label{log ((1 + y) (1 + z)) = log (1 + y) + log(1 + z)}
        \log ((1 + y) \, (1 + z)) = \log (1 + y) + \log(1 + z).
\end{equation}
Let us mention two ways to look at this, following \cite{c, fg}.

        In the first approach, one can begin by expanding
\begin{equation}
\label{(y + z + y z)^j}
        (y + z + y \, z)^j
\end{equation}
into a sum of monomials in $y$ and $z$, as we did before for (\ref{(Y
  + Z + Y Z)^j}).  If one plugs the resulting sums into the right side
of (\ref{log ((1 + y) (1 + z)) = log (1 + y + z + y z) = ...}), then
one would like to rearrange the terms to get a sum that corresponds to
the formal power series expansion in $y$ and $z$.  More precisely,
this can be done using the discussions in Sections \ref{sums of sums}
and \ref{functions of sums}.  This permits (\ref{log ((1 + y) (1 + z))
  = log (1 + y) + log(1 + z)}) to be derived from the analogous
statement (\ref{log ((1 + Y) (1 + Z)) = log (1 + Y) + log (1 + Z)})
for formal power series, as on p281 of \cite{c}.

        Alternatively, let $y \in k$ with $|y| < 1$ be given, and
consider
\begin{eqnarray}
\label{log((1 + y) (1 + z)) = log (1 + y + (1 + y) z) = ...}
        \log((1 + y) \, (1 + z)) & = & \log (1 + y + (1 + y) \, z) \\
 & = & \sum_{j = 1}^\infty \frac{(-1)^{j + 1}}{j} \, (y + (1 + y) \, z)^j \nonumber
\end{eqnarray}
as a function of $z$.  This can be converted into a power series in
$z$ that converges for $|z| < 1$, as in Section \ref{changing
  centers}.  Thus (\ref{log ((1 + y) (1 + z)) = log (1 + y) + log(1 +
  z)}) may be treated as an equation relating power series in $z$,
with $y$ as a constant.  Note that (\ref{log ((1 + y) (1 + z)) = log
  (1 + y) + log(1 + z)}) holds trivially for every $y, z \in k$ with
$|y|, |z| < 1$ when $|\cdot|$ is the trivial absolute value function
on $k$, since $y = z = 0$.  Thus we may as well suppose that $|\cdot|$
is not the trivial absolute value function on $k$, so that every
element of the open unit ball $B(0, 1)$ in $k$ is a limit point of
$B(0, 1)$.  In this case, the derivative of $\log (1 + z)$ as a
$k$-valued function on $B(0, 1)$ is equal to
\begin{equation}
\label{frac{1}{1 + z}}
        \frac{1}{1 + z}
\end{equation}
for every $z \in B(0, 1)$, as in the previous section.  This implies
that the derivative of the right side of (\ref{log ((1 + y) (1 + z)) =
  log (1 + y) + log(1 + z)}) as a $k$-valued function of $z$ on $B(0,
1)$ is equal to (\ref{frac{1}{1 + z}}) too.  Similarly, one can check
that the derivative of the left side of (\ref{log ((1 + y) (1 + z)) =
  log (1 + y) + log(1 + z)}) as a $k$-valued function of $z$ on $B(0,
1)$ is equal to (\ref{frac{1}{1 + z}}), using the chain rule.  Both
sides of (\ref{log ((1 + y) (1 + z)) = log (1 + y) + log(1 + z)}) can
be expressed as power series in $z$, as before, and so their
derivatives are given by the corresponding differentiated power series
in $z$, as in Section \ref{differentiation, continued}.  It follows
that the coefficients of the differentiated power series in $z$
corresponding to both sides of (\ref{log ((1 + y) (1 + z)) = log (1 +
  y) + log(1 + z)}) are the same, since the values of the derivatives
are the same on $B(0, 1)$, and $|\cdot|$ is not the trivial absolute
value function on $k$.  Thus the coefficients of the power series in
$z$ corresponding to both sides of (\ref{log ((1 + y) (1 + z)) = log
  (1 + y) + log(1 + z)}) are the same, except perhaps for the constant
terms, because $k$ has characteristic $0$.  Of course, (\ref{log ((1 +
  y) (1 + z)) = log (1 + y) + log(1 + z)}) obviously holds when $z =
0$, which means that the constant terms of the power series in $z$
corresponding to both sides of (\ref{log ((1 + y) (1 + z)) = log (1 +
  y) + log(1 + z)}) are the same as well.  This implies that (\ref{log
  ((1 + y) (1 + z)) = log (1 + y) + log(1 + z)}) holds for every $z
\in k$ with $|z| < 1$, as in the proof of Proposition 4.5.3 on p110 of
\cite{fg}.

\section{Some additional properties}
\label{some additional properties}

        Let $k$ be a field, and let $|\cdot|$ be an absolute value
function on $k$.  If $u, v \in k$ and $|u| = |v| = 1$, then
\begin{equation}
\label{|(u/v) - 1| = |u - v|/|v| = |u - v|}
        |(u/v) - 1| = |u - v|/|v| = |u - v|.
\end{equation}
Suppose now that $|\cdot|$ is an ultrametric absolute value function
on $k$, and that $y, z \in k$ satisfy $|y|, |z| < 1$.  Remember that
$u = 1 + y$ and $v = 1 + z$ satisfy $|u| = |v| = 1$, as in the
previous section.  In this case, (\ref{|(u/v) - 1| = |u - v|/|v| = |u
  - v|}) implies that
\begin{equation}
\label{|(1 + y)/(1 + z) - 1| = |y - z|}
        |(1 + y)/(1 + z) - 1| = |y - z|.
\end{equation}

        Let us suppose for the rest of the section that $k$ is a field
of characteristic $0$ with an ultrametric absolute value function
$|\cdot|$, and that $k$ is complete with respect to the ultrametric
associated to $|\cdot|$.  Let us also suppose that the induced
absolute value function on ${\bf Q}$ is trivial.  If $x \in k$ and
$|x| < 1$, then
\begin{equation}
\label{|log (1 + x)| = |x|}
        |\log (1 + x)| = |x|.
\end{equation}
Of course, this is trivial when $x = 0$, and so it suffices to verify
that (\ref{|log (1 + x)| = |x|}) holds when $x \ne 0$.  Observe that
\begin{equation}
\label{|log (1 + x) - x| = ... le max_{j ge 2} |x^j/j|}
        |\log (1 + x) - x|
 = \biggl|\sum_{j = 2}^\infty \frac{(-1)^{j + 1}}{j} \, x^j\biggr|
                \le \max_{j \ge 2} |x^j/j|,
\end{equation}
using the definition (\ref{log (1 + x) = sum_{j = 1}^infty
  frac{(-1)^{j + 1}}{j} x^j}) of $\log (1 + x)$ in the first step, and
the ultrametric version of the triangle inequality in the second step.
It follows that
\begin{equation}
\label{|log (1 + x) - x| le max_{j ge 2} |x|^j = |x|^2}
        |\log (1 + x) - x| \le \max_{j \ge 2} |x|^j = |x|^2
\end{equation}
when the induced absolute value function on ${\bf Q}$ is trivial, and
hence that
\begin{equation}
\label{|log (1 + x) - x| < |x|}
        |\log (1 + x) - x| < |x|
\end{equation}
when $x \ne 0$ and $|x| < 1$.  This implies (\ref{|log (1 + x)| =
  |x|}), using the ultrametric version of the triangle inequality
again, as in (\ref{|x + y| = |y|}) in Section \ref{absolute value
  functions}.

        Let us check that
\begin{equation}
\label{|log(1 + y) - log (1 + z)| = |y - z|}
        |\log(1 + y) - \log (1 + z)| = |y - z|
\end{equation}
for every $y, z \in k$ with $|y|, |z| < 1$ under these conditions.
Note that
\begin{equation}
\label{|(1 + y)/(1 + z) - 1| < 1}
        |(1 + y)/(1 + z) - 1| < 1,
\end{equation}
by (\ref{|(1 + y)/(1 + z) - 1| = |y - z|}), or the fact that (\ref{{w
    in k : |w - 1| < 1}}) is a group with respect to multiplication.
Thus the logarithm of $(1 + y)/(1 + z)$ is defined, and 
\begin{equation}
\label{log ((1 + y)/(1 + z)) = log (1 + y) - log (1 + z)}
        \log ((1 + y)/(1 + z)) = \log (1 + y) - \log (1 + z),
\end{equation}
by (\ref{log ((1 + y) (1 + z)) = log (1 + y) + log(1 + z)}).  It
follows that
\begin{eqnarray}
\label{|log (1 + y) - log (1 + z)| = ... = |y - z|}
 |\log (1 + y) - \log (1 + z)| & = & |\log ((1 + y)/(1 + z))| \\
                     & = & |(1 + y)/(1 + z) - 1| = |y - z|, \nonumber
\end{eqnarray}
as desired, using (\ref{|log (1 + x)| = |x|}) in the second step, and
(\ref{|(1 + y)/(1 + z) - 1| = |y - z|}) in the third step.  One could
also get (\ref{|log(1 + y) - log (1 + z)| = |y - z|}) as in
(\ref{|f(x) - f(y)| = |f'(y)| |x - y|}) or (\ref{|f(x) - f(y)| =
  |f'(x_0)| |x - y|}) in Section \ref{hensel's lemma}.

        It will be convenient to put
\begin{equation}
\label{f(x) = log (1 + x)}
        f(x) = \log (1 + x)
\end{equation}
for every $x \in k$ with $|x| < 1$, so that $f$ defines a $k$-valued
function on the open unit ball $B(0, 1)$ in $k$.  Of course, $f$ maps
$B(0, 1)$ into itself, by (\ref{|log (1 + x)| = |x|}), and in fact one
can check that
\begin{equation}
\label{f(B(0, 1)) = B(0, 1)}
        f(B(0, 1)) = B(0, 1),
\end{equation}
using Hensel's lemma.  More precisely, one can take $x_0 = 0$ in
Section \ref{hensel's lemma}, and $0 < r = t < 1$.  Thus $f$ maps
$\overline{B}(0, t)$ onto itself for every $t \in (0, 1)$, as in
(\ref{f(overline{B}(x_0, t)) = overline{B}(f(x_0), |f'(x_0)| t)}),
since $f(0) = 0$ and $f'(0) = 1$.  This implies (\ref{f(B(0, 1)) =
  B(0, 1)}), which also corresponds to (\ref{f(B(x_0, r_1)) =
  B(f(x_0), |f'(x_0)| r_1}) in Section \ref{some variants}, with $r_1
= 1$.

\section{Some additional properties, continued}
\label{some additional properties, continued}

        Let $k$ be a field with an ultrametric absolute value function
$|\cdot|$.  It is easy to see that
\begin{equation}
\label{{w in k : |w - 1| < r}}
        \{w \in k : |w - 1| < r\}
\end{equation}
is a group with respect to multiplication when $0 < r \le 1$, and
similarly that
\begin{equation}
\label{{w in k : |w - 1| le r}}
        \{w \in k : |w - 1| \le r\}
\end{equation}
is a group with respect to multiplication when $0 \le r < 1$.  More
precisely, (\ref{{w in k : |w - 1| < r}}) and (\ref{{w in k : |w - 1|
    le r}}) are subgroups of the group
\begin{equation}
\label{{w in k : |w| = 1}, 2}
        \{w \in k : |w| = 1\}
\end{equation}
with respect to multiplication.  If $r = 1$, then (\ref{{w in k : |w -
    1| < r}}) is the same as (\ref{{w in k : |w - 1| < 1}}) in Section
\ref{usual identity}, which we have already seen is a subgroup of
(\ref{{w in k : |w| = 1}, 2}).  Essentially the same argument works in
the other cases, using (\ref{|1/w - 1| = |1 - w|/|w| = |1 - w| < 1})
and (\ref{|(1 + y) (1 + z) - 1| = |y + z + y z| le max(|y|, |z|, |y|
  |z|)}).

        Let us suppose for the rest of this section that $k$ is a
field of characteristic $0$ with an ultrametric absolute value
function $|\cdot|$, and that $k$ is complete with respect to the
associated ultrametric.  If the induced absolute value function on
${\bf Q}$ is not trivial, then it is equivalent to the $p$-adic absolute
value function on ${\bf Q}$ for some prime number $p$, by Ostrowski's
theorem, as in Section \ref{ostrowski's theorems}.  Let us suppose also
that the induced absolute value function on ${\bf Q}$ is equal to
the $p$-adic absolute value function, which can always be arranged
by replacing $|\cdot|$ on $k$ by a suitable positive power of itself.

        Observe that
\begin{equation}
\label{p^l - 1 = (p - 1) sum_{m = 0}^{l - 1} p^m ge (p - 1) l}
        p^l - 1 = (p - 1) \, \sum_{m = 0}^{l - 1} p^m \ge (p - 1) \, l
\end{equation}
for every positive integer $l$, which also holds when $l = 0$, with
the sum interpreted as being equal to $0$.  If $r$ is a nonnegative
real number such that
\begin{equation}
\label{r le p^{-1/(p - 1)}}
        r \le p^{-1/(p - 1)},
\end{equation}
then it follows that
\begin{equation}
\label{p^l r^{p^l - 1} le p^l r^{(p - 1) l} le 1}
        p^l \, r^{p^l - 1} \le p^l \, r^{(p - 1) \, l} \le 1
\end{equation}
for each nonnegative integer $l$.  Remember that the $p$-adic absolute
value $|j|_p$ of $j \in {\bf Z}_+$ is given by
\begin{equation}
\label{|j|_p = p^{-l(j)}}
        |j|_p = p^{-l(j)},
\end{equation}
where $l(j)$ is the largest nonnegative integer such that $j$ is an
integer multiple of $p^{l(j)}$.  If $r \ge 0$ satisfies (\ref{r le
  p^{-1/(p - 1)}}), then it follows that
\begin{equation}
\label{r^{j - 1}/|j|_p = p^{l(j)} r^{j - 1} le p^{l(j)} r^{p^{l(j)} - 1} le 1}
 r^{j - 1}/|j|_p = p^{l(j)} \, r^{j - 1} \le p^{l(j)} \, r^{p^{l(j)} - 1} \le 1
\end{equation}
for every $j \ge 1$.  More precisely, this uses the facts that $j \ge
p^{l(j)}$ and $r \le 1$ in the first inequality, and (\ref{p^l r^{p^l
    - 1} le p^l r^{(p - 1) l} le 1}) in the second inequality.

        If $x \in k$ and $|x| < 1$, then $\log (1 + x)$ is defined
in $k$ as in Section \ref{logarithm}, and
\begin{equation}
\label{|log (1 + x)| le ... = max_{j ge 1} (p^{l(j)} |x|^j)}
 |\log (1 + x)| \le \max_{j \ge 1} |x^j/j| = \max_{j \ge 1} (|x|^j/|j|_p)
                                          = \max_{j \ge 1} (p^{l(j)} \, |x|^j),
\end{equation}
using the ultrametric version of the triangle inequality in the
first step.  If
\begin{equation}
\label{|x| le p^{-1/(p - 1)}}
        |x| \le p^{-1/(p - 1)},
\end{equation}
then we have that
\begin{equation}
\label{p^{l(j)} |x|^{j - 1} le 1}
        p^{l(j)} \, |x|^{j - 1} \le 1
\end{equation}
for every $j \ge 1$, by (\ref{r^{j - 1}/|j|_p = p^{l(j)} r^{j - 1} le
  p^{l(j)} r^{p^{l(j)} - 1} le 1}).  Combining this with (\ref{|log (1
  + x)| le ... = max_{j ge 1} (p^{l(j)} |x|^j)}), we get that
\begin{equation}
\label{|log (1 + x)| le |x|}
        |\log (1 + x)| \le |x|
\end{equation}
for every $x \in k$ that satisfies (\ref{|x| le p^{-1/(p - 1)}}).

        Suppose for the moment that $y, z \in k$ satisfy $|y|, |z| < 1$ and
\begin{equation}
\label{|(1 + y)/(1 + z) - 1| le p^{-1/(p - 1)}}
        |(1 + y)/(1 + z) - 1| \le p^{-1/(p - 1)}.
\end{equation}
In particular, this holds when
\begin{equation}
\label{|y|, |z| le p^{-1/(p - 1)}}
        |y|, |z| \le p^{-1/(p - 1)},
\end{equation}
because (\ref{{w in k : |w - 1| le r}}) is a group with respect to
mulitplication when
\begin{equation}
\label{r = p^{-1/(p - 1)} < 1}
        r = p^{-1/(p - 1)} < 1.
\end{equation}
Under these conditions, we get that
\begin{eqnarray}
\label{|log (1 + y) - log (1 + z)| = ... = |y - z|, 2}
 |\log (1 + y) - \log (1 + z)| & = & |\log ((1 + y)/(1 + z))| \\
                & \le & |(1 + y)/(1 + z) - 1| = |y - z|, \nonumber
\end{eqnarray}
using (\ref{log ((1 + y) (1 + z)) = log (1 + y) + log(1 + z)}) in
Section \ref{usual identity} in the first step, (\ref{|log (1 + x)| le
  |x|}) in the second step, and (\ref{|(1 + y)/(1 + z) - 1| = |y -
  z|}) in the previous section in the third step.

        Note that
\begin{equation}
\label{r max_{j ge 2} (r^{j - 2}/|j|_p) = max_{j ge 2} (r^{j - 1}/|j|_p) le 1}
 r \, \max_{j \ge 2} (r^{j - 2}/|j|_p) = \max_{j \ge 2} (r^{j - 1}/|j|_p) \le 1
\end{equation}
when $r \ge 0$ satisfies (\ref{r le p^{-1/(p - 1)}}), using (\ref{r^{j
    - 1}/|j|_p = p^{l(j)} r^{j - 1} le p^{l(j)} r^{p^{l(j)} - 1} le
  1}) in the second step.  This implies that
\begin{equation}
\label{r max_{j ge 2} (r^{j - 2}/|j|_p) = max_{j ge 2} (r^{j - 1}/|j|_p) < 1}
 r \, \max_{j \ge 2} (r^{j - 2}/|j|_p) = \max_{j \ge 2} (r^{j - 1}/|j|_p) < 1
\end{equation}
when
\begin{equation}
\label{0 le r < p^{-1/(p - 1)}}
        0 \le r < p^{-1/(p - 1)}.
\end{equation}

        Let us check that
\begin{equation}
\label{|log (1 + x)| = |x|, 2}
        |\log (1 + x)| = |x|
\end{equation}
for every $x \in k$ with
\begin{equation}
\label{|x| < p^{-1/(p - 1)}}
        |x| < p^{-1/(p - 1)}.
\end{equation}
As in (\ref{|log (1 + x) - x| = ... le max_{j ge 2} |x^j/j|}) in the
previous section, we have that
\begin{equation}
\label{|log (1 + x) - x| = ... = max_{j ge 2} (|x|^j/|j|_p)}
 \qquad |\log (1 + x) - x|
 = \biggl|\sum_{j = 2}^\infty \frac{(-1)^{j + 1}}{j} \, x^j\biggr|
                \le \max_{j \ge 2} |x^j/j| = \max_{j \ge 2} (|x|^j/|j|_p),
\end{equation}
for every $x \in k$ with $|x| < 1$, using the current hypothesis about
the absolute value function on $k$ in the last step.  We also have
that
\begin{equation}
\label{max_{j ge 2} (|x|^{j - 1}/|j|_p) < 1}
        \max_{j \ge 2} (|x|^{j - 1}/|j|_p) < 1
\end{equation}
when $x \in k$ satisfies (\ref{|x| < p^{-1/(p - 1)}}), by (\ref{r
  max_{j ge 2} (r^{j - 2}/|j|_p) = max_{j ge 2} (r^{j - 1}/|j|_p) <
  1}) with $r = |x|$.  Combining this with (\ref{|log (1 + x) - x| =
  ... = max_{j ge 2} (|x|^j/|j|_p)}), we get that that
\begin{equation}
\label{|log (1 + x) - x| < |x|, 2}
        |\log (1 + x) - x| < |x|
\end{equation}
when $x \in k$ satisfies (\ref{|x| < p^{-1/(p - 1)}}) and $x \ne 0$.
As before, (\ref{|log (1 + x)| = |x|, 2}) follows from (\ref{|log (1 +
  x) - x| < |x|, 2}) when $x \in k$ satisfies (\ref{|x| < p^{-1/(p -
    1)}}) and $x \ne 0$, as in (\ref{|x + y| = |y|}) in Section
\ref{absolute value functions}, and of course (\ref{|log (1 + x)| =
  |x|, 2}) is trivial when $x = 0$.

        Suppose for the moment that $y, z \in k$ satisfy $|y|, |z| < 1$ and
\begin{equation}
\label{|(1 + y)/(1 + z) - 1| < p^{-1/(p - 1)}}
        |(1 + y)/(1 + z) - 1| < p^{-1/(p - 1)}
\end{equation}
which holds in particular when
\begin{equation}
\label{|y|, |z| < p^{-1/(p - 1)}}
        |y|, |z| < p^{-1/(p - 1)},
\end{equation}
because (\ref{{w in k : |w - 1| < r}}) is a group when $r$ is as in
(\ref{r = p^{-1/(p - 1)} < 1}).  Under these conditions, we get that
\begin{eqnarray}
\label{|log (1 + y) - log (1 + z)| = ... = |(1 + y)/(1 + z) - 1| = |y - z|}
 |\log (1 + y) - \log (1 + z)| & = & |\log ((1 + y)/(1 + z))| \\
                          & = & |(1 + y)/(1 + z) - 1| = |y - z|, \nonumber
\end{eqnarray}
using (\ref{log ((1 + y) (1 + z)) = log (1 + y) + log(1 + z)}) in
Section \ref{usual identity} in the first step, (\ref{|log (1 + x)| =
  |x|, 2}) in the second step, and (\ref{|(1 + y)/(1 + z) - 1| = |y -
  z|}) in the third step.

        Put
\begin{equation}
\label{f(x) = log (1 + x), 2}
        f(x) = \log (1 + x)
\end{equation}
for every $x \in k$ with $|x| < 1$ again, as in the previous section.
Using Hensel's lemma, one can check that
\begin{equation}
\label{f(B(0, p^{-1/(p - 1)})) = B(0, p^{-1/(p - 1)})}
        f(B(0, p^{-1/(p - 1)})) = B(0, p^{-1/(p - 1)}).
\end{equation}
More precisely, one can show that $f$ maps $\overline{B}(0, t)$
onto itself when
\begin{equation}
\label{0 < t < p^{-1/(p - 1)}}
        0 < t < p^{-1/(p - 1)},
\end{equation}
as in (\ref{f(overline{B}(x_0, t)) = overline{B}(f(x_0), |f'(x_0)|
  t)}) in Section \ref{hensel's lemma}, with $x_0 = 0$.  This also
corresponds to (\ref{f(B(x_0, r_1)) = B(f(x_0), |f'(x_0)| r_1}) in
Section \ref{some variants}, with $r_1 = p^{-1/(p - 1)}$.

\section{Outer measures}
\label{outer measures}

        Let $X$ be a set, and let $\mu$ be a function defined on
the collection of all subsets of $X$ with values in the set of
nonnegative extended real numbers.  If $\mu$ satisfies the following
three conditions, then $\mu$ is said to be an \emph{outer
  measure}\index{outer measures} on $X$.  First,
\begin{equation}
\label{mu(emptyset) = 0}
        \mu(\emptyset) = 0.
\end{equation}
Second, if $A \subseteq B \subseteq X$, then
\begin{equation}
\label{mu(A) le mu(B)}
        \mu(A) \le \mu(B).
\end{equation}
Third, if $A_1, A_2, A_3, \ldots$ is any infinite sequence of subsets
of $X$, then
\begin{equation}
\label{mu(bigcup_{j = 1}^infty A_j) le sum_{j = 1}^infty mu(A_j)}
 \mu\Big(\bigcup_{j = 1}^\infty A_j\Big) \le \sum_{j = 1}^\infty \mu(A_j),
\end{equation}
where the sum on the right side is interpreted as in Section
\ref{nonnegative sums}.  In particular, the sum on the right side of
(\ref{mu(bigcup_{j = 1}^infty A_j) le sum_{j = 1}^infty mu(A_j)}) is
automatically interpreted as being $+\infty$ when $\mu(A_j) = +\infty$
for any $j$.  Otherwise, if $\mu(A_j) < \infty$ for each $j$, then the
sum on the right side of (\ref{mu(bigcup_{j = 1}^infty A_j) le sum_{j
    = 1}^infty mu(A_j)}) may be considered as an ordinary infinite
series of nonnegative real numbers, which is equal to $+\infty$ when
the series does not converge in the usual sense in ${\bf R}$.  Of
course, (\ref{mu(bigcup_{j = 1}^infty A_j) le sum_{j = 1}^infty
  mu(A_j)}) holds trivially when the sum on the right side is equal to
$+\infty$.

        The second condition (\ref{mu(A) le mu(B)}) may be described as
monotonicity of $\mu$, and the third condition (\ref{mu(bigcup_{j =
    1}^infty A_j) le sum_{j = 1}^infty mu(A_j)}) is known as countable
subadditivity.\index{countable subadditivity}  If $A_1, A_2, \ldots,
A_n$ is any finite sequence of subsets of $X$, then (\ref{mu(emptyset)
  = 0}) and (\ref{mu(bigcup_{j = 1}^infty A_j) le sum_{j = 1}^infty
  mu(A_j)}) imply that
\begin{equation}
\label{mu(bigcup_{j = 1}^n A_j) le sum_{j = 1}^n mu(A_j)}
        \mu\Big(\bigcup_{j = 1}^n A_j\Big) \le \sum_{j = 1}^n \mu(A_j),
\end{equation}
by taking $A_j = \emptyset$ when $j > n$.  This may be described as
finite subadditivity of $\mu$.  As before, the sum on the right side
of (\ref{mu(bigcup_{j = 1}^n A_j) le sum_{j = 1}^n mu(A_j)}) is
interpreted as being equal to $+\infty$ when $\mu(A_j) = +\infty$ for
any $j = 1, 2, \ldots, n$, in which case (\ref{mu(bigcup_{j = 1}^n
  A_j) le sum_{j = 1}^n mu(A_j)}) holds trivially.

        If $A$ is a subset of $X$, and $A_1, A_2, A_3, \ldots$ is an
infinite sequence of subsets of $X$ such that
\begin{equation}
\label{A subseteq bigcup_{j = 1}^infty A_j}
        A \subseteq \bigcup_{j = 1}^\infty A_j,
\end{equation}
then (\ref{mu(A) le mu(B)}) and (\ref{mu(bigcup_{j = 1}^infty A_j) le
  sum_{j = 1}^infty mu(A_j)}) imply that
\begin{equation}
\label{mu(A) le sum_{j = 1}^infty mu(A_j)}
        \mu(A) \le \sum_{j = 1}^\infty \mu(A_j),
\end{equation}
with $B = \bigcup_{j = 1}^\infty A_j$.  Of course, this property
implies (\ref{mu(bigcup_{j = 1}^infty A_j) le sum_{j = 1}^infty
  mu(A_j)}), by taking $A = \bigcup_{j = 1}^\infty A_j$.  One can also
get (\ref{mu(A) le mu(B)}) from (\ref{mu(A) le sum_{j = 1}^infty
  mu(A_j)}) and (\ref{mu(emptyset) = 0}), by taking $A_1 = B$ and $A_j
= \emptyset$ when $j > 1$.  Thus the definition of an outer measure
can be equivalently formulated in terms of (\ref{mu(emptyset) = 0})
and (\ref{mu(A) le sum_{j = 1}^infty mu(A_j)}), instead of
(\ref{mu(emptyset) = 0}), (\ref{mu(A) le mu(B)}), and
(\ref{mu(bigcup_{j = 1}^infty A_j) le sum_{j = 1}^infty mu(A_j)}).

        Let $I$ be a countably infinite set, and suppose that for
each $j \in I$, $A_j$ is a subset of $X$.  The countable subadditivity
property (\ref{mu(bigcup_{j = 1}^infty A_j) le sum_{j = 1}^infty mu(A_j)})
can be reformulated as saying that
\begin{equation}
\label{mu(bigcup_{j in I} A_j) le sum_{j in I} mu(A_j)}
        \mu\Big(\bigcup_{j \in I} A_j\Big) \le \sum_{j \in I} \mu(A_j)
\end{equation}
under these conditions, where the sum on the right side of
(\ref{mu(bigcup_{j in I} A_j) le sum_{j in I} mu(A_j)}) is defined as
in Section \ref{nonnegative sums}.  Similarly, (\ref{mu(A) le sum_{j =
    1}^infty mu(A_j)}) can be reformulated as saying that
\begin{equation}
\label{mu(A) le sum_{j in I} mu(A_j)}
        \mu(A) \le \sum_{j \in I} \mu(A_j)
\end{equation}
when $A \subseteq \bigcup_{j \in I} A_j$.  In both cases, one may as
well allow $I$ to be a nonempty set with only finitely or countably
many elements, using (\ref{mu(emptyset) = 0}), as before.  One could
even allow $I$ to be the empty set, and interpret any sum over $I$ as
being equal to $0$, and any union over $I$ as being the empty set.
With these interpretations, (\ref{mu(bigcup_{j in I} A_j) le sum_{j in
    I} mu(A_j)}) or (\ref{mu(A) le sum_{j in I} mu(A_j)}) may be
considered to imply that $\mu(\emptyset) \le 0$, and hence
(\ref{mu(emptyset) = 0}), since $\mu$ is supposed to be nonnegative by
hypothesis.  Thus the notion of an outer measure may be defined in
terms of (\ref{mu(A) le sum_{j in I} mu(A_j)}), where $I$ is allowed
to be any set with only finitely or countably many elements, including
the empty set.  Strictly speaking, it is better to consider finite or
countable collections of subsets of $X$, so that an auxiliary set $I$
of indices is not needed.

        A subset $B$ of $X$ is said to be \emph{measurable}\index{measurable
sets} with respect to an outer measure $\mu$ on $X$ if
\begin{equation}
\label{mu(A) = mu(A cap B) + mu(A setminus B)}
        \mu(A) = \mu(A \cap B) + \mu(A \setminus B)
\end{equation}
for every subset $A$ of $X$.  It is well known that the collection of
subsets of $X$ that are measuable with respect to $\mu$ forms a
$\sigma$-algebra, and that the restriction of $\mu$ to this
$\sigma$-algebra of measurable sets is countably additive.  Note that
$B \subseteq X$ is measurable when $\mu(B) = 0$.  An outer measure
$\mu$ on a topological space $X$ is said to be a \emph{Borel outer
  measure}\index{Borel outer measures}\index{outer measures!Borel} if
every Borel subset of $X$ is measurable with respect to $\mu$.

        Suppose for the moment that $(X, d(x, y))$ is a metric space.
If $A$ and $B$ are nonempty subsets of $X$, then the
\emph{distance}\index{distance between sets} between $A$ and $B$
is defined by
\begin{equation}
\label{dist(A, B) = inf {d(x, y) : x in A, y in B}}
        \dist(A, B) = \inf \{d(x, y) : x \in A, y \in B\}.
\end{equation}
Thus $\dist(A, B) > 0$ if and only if there is an $\eta > 0$ such that
\begin{equation}
\label{d(x, y) ge eta}
        d(x, y) \ge \eta
\end{equation}
for every $x \in A$ and $y \in B$, which implies in particular that
$A$ and $B$ are disjoint.  An outer measure $\mu$ on $X$ is said to be
a \emph{metric outer measure}\index{metric outer measures} if
\begin{equation}
\label{mu(A cup B) = mu(A) + mu(B)}
        \mu(A \cup B) = \mu(A) + \mu(B)
\end{equation}
for every pair of nonempty subsets $A$, $B$ of $X$ with $\dist(A, B) >
0$.  Of course, it suffices to show that
\begin{equation}
\label{mu(A cup B) ge mu(A) + mu(B)}
        \mu(A \cup B) \ge \mu(A) + \mu(B),
\end{equation}
since the opposite inequality follows from finite subadditivity, as in
(\ref{mu(bigcup_{j = 1}^n A_j) le sum_{j = 1}^n mu(A_j)}).  If $\mu$
is a metric outer measure on $X$, then it is well known that the Borel
sets in $X$ are measurable with respect to $\mu$.  This is called
\emph{Carath\'eodory's criterion}.\index{Carath\'eodory's criterion}

        Let $E$ be a nonempty proper subset of $X$ such that
\begin{equation}
\label{dist(E, X setminus E) > 0}
        \dist(E, X \setminus E) > 0,
\end{equation}
which implies that $E$ is both open and closed in $X$, so that $X$ is
not connected.  If $\mu$ is a metric doubling measure on $X$, then it
is easy to see directly that $E$ is measurable with respect to $\mu$,
which is a special case of Carath\'eodory's criterion.  Note that
(\ref{dist(E, X setminus E) > 0}) holds when $E$ is a nonempty proper
compact open subset of $X$.  If $d(x, y)$ is an ultrametric on $X$ and
$E$ is an open or closed ball in $X$ of positive radius which is a
proper subset of $X$, then $E$ satisfies (\ref{dist(E, X setminus E) >
  0}).  Similarly, if $d(x, y)$ is an ultrametric on $X$ and $E$ is a
proper nonempty subset of $X$ that can be expressed as the union of a
family of balls of a fixed positive radius, then $E$ satisfies
(\ref{dist(E, X setminus E) > 0}).

        An outer measure $\mu$ on a set $X$ is said to be 
\emph{regular}\index{regular outer measures}\index{outer measures!regular}
if for each subset $A$ of $X$ there is a subset $B$ of $X$ that
is measurable with respect to $\mu$ and satisfies
\begin{equation}
\label{A subseteq B}
        A \subseteq B
\end{equation}
and
\begin{equation}
\label{mu(A) = mu(B)}
        \mu(A) = \mu(B).
\end{equation}
Similarly, an outer measure $\mu$ on a topological space $X$ is said
to be \emph{Borel regular}\index{Borel regular outer
  measures}\index{outer measures!Borel regular} if $\mu$ is a Borel
outer measure on $X$, and if for each $A \subseteq X$ there is a Borel
set $B \subseteq X$ that satisfies (\ref{A subseteq B}) and
(\ref{mu(A) = mu(B)}).  Of course, (\ref{A subseteq B}) implies that
$\mu(A) \le \mu(B)$, as in (\ref{mu(A) le mu(B)}), and so it suffices
to check that the opposite inequality holds to get (\ref{mu(A) =
  mu(B)}).  In both cases, if $\mu(A) = +\infty$, then one can simply
take $B = X$.  Some texts use the term ``measure'' for what is called
an outer measure here, and then use the adjectives ``Borel'',
``regular'', and ``Borel regular'' as defined here.  Otherwise, the
term ``measure'' is often used for a countably-additive nonnegative
extended-real-valued function defined an a $\sigma$-algebra of
``measurable'' subsets of a set $X$, for which the measure of the
empty set is equal to $0$.  In this terminology, a Borel measure is a
measure defined on the $\sigma$-algebra of Borel subsets of a
topological space $X$, and somewhat different regularity properties
are typically considered, especially for Borel measures on locally
compact Hausdorff topological spaces.

        Let $\mu$ be an outer measure on a set $X$, and let $A_1, A_2, A_3,
\ldots$ be a sequence of subsets of $X$ such that
\begin{equation}
\label{A_j subseteq A_{j + 1}}
        A_j \subseteq A_{j + 1}
\end{equation}
for each $j$.  Thus
\begin{equation}
\label{mu(A_j) le mu(A_{j + 1}) le mu(bigcup_{l = 1}^infty A_l)}
        \mu(A_j) \le \mu(A_{j + 1}) \le \mu\Big(\bigcup_{l = 1}^\infty A_l\Big)
\end{equation}
for each $j \ge 1$, and hence
\begin{equation}
\label{sup_{j ge 1} mu(A_j) le mu(bigcup_{l = 1}^infty A_l)}
        \sup_{j \ge 1} \mu(A_j) \le \mu\Big(\bigcup_{l = 1}^\infty A_l\Big).
\end{equation}
Note that $\mu(A_j)$ tends to the supremum as $j \to \infty$, because
$\mu(A_j)$ increases monotonically in $j$, with the usual
interpretations for extended real numbers.  Suppose that $A_j$ is
measurable with respect to $\mu$ for each $j \ge 1$, which implies
that $A_j \setminus A_{j - 1}$ is measurable with respect to $\mu$ for
each $j \ge 2$.  The sets $A_j \setminus A_{j - 1}$ with $j \ge 2$ are
pairwise disjoint, because of (\ref{A_j subseteq A_{j + 1}}), and
disjoint from $A_1$.  We also have that
\begin{equation}
\label{A_j = A_1 cup (bigcup_{l = 2}^j (A_l setminus A_{l - 1})}
        A_j = A_1 \cup \Big(\bigcup_{l = 2}^j (A_l \setminus A_{l - 1})\Big)
\end{equation}
for each $j \ge 1$, with suitable interpretations when $j = 1$, and that
\begin{equation}
\label{bigcup_{j = 1}^infty A_j = ...}
 \bigcup_{j = 1}^\infty A_j
         = A_1 \cup \Big(\bigcup_{l = 2}^\infty (A_l \setminus A_{l - 1})\Big).
\end{equation}
This implies that
\begin{equation}
\label{mu(A_j) = mu(A_1) + sum_{l = 2}^j mu(A_l setminus A_{l - 1})}
        \mu(A_j) = \mu(A_1) + \sum_{l = 2}^j \mu(A_l \setminus A_{l - 1})
\end{equation}
for each $j \ge 1$, with suitable interpretations when $j = 1$, and that
\begin{equation}
\label{mu(bigcup_{j = 1}^infty A_j) = ...}
 \mu\Big(\bigcup_{j = 1}^\infty A_j\Big)
                 = \mu(A_1) + \sum_{l = 2}^\infty \mu(A_l \setminus A_{l - 1}),
\end{equation}
because $\mu$ is countably additive on measurable sets.  It follows that
\begin{equation}
\label{lim_{j to infty} mu(A_j) = mu(bigcup_{j = 1}^infty A_j)}
        \lim_{j \to \infty} \mu(A_j) = \mu\Big(\bigcup_{j = 1}^\infty A_j\Big).
\end{equation}
Of course, this is a standard fact about countably-additive measures
on $\sigma$-algebras of measurable sets.

        If $\mu$ is a regular outer measure on $X$, then it is well known
that (\ref{lim_{j to infty} mu(A_j) = mu(bigcup_{j = 1}^infty A_j)})
holds even when the $A_j$'s are not asked to be measurable with
respect to $\mu$.  We already have (\ref{sup_{j ge 1} mu(A_j) le
  mu(bigcup_{l = 1}^infty A_l)}), and so the point is to show that the
opposite inequality holds as well.  The regularity of $\mu$ implies
that each $A_j$ is contained in a measurable set with the same
measure, and so one would like to use this to reduce to the previous
case of measurable sets.  However, one should be a bit careful to
choose these measurable sets so that they are also monotonically
increasing with respect to inclusion, which is not too difficult to
do.

\section{Hausdorff measures}
\label{hausdorff measures}

        Let $(M, d(x, y))$ be a metric space, and let $\alpha$ be a
nonnegative real number.  Remember that a subset of $M$ is said to be
\emph{bounded}\index{bounded sets} if it is contained in a ball of
finite radius.  The \emph{diameter}\index{diameters of sets} of a
nonempty bounded set $A \subseteq M$ is defined by
\begin{equation}
\label{diam A = sup {d(x, y) : x, y in A}}
        \diam A = \sup \{d(x, y) : x, y \in A\},
\end{equation}
in which case
\begin{equation}
\label{(diam A)^alpha}
        (\diam A)^\alpha
\end{equation}
can be defined in the usual way for each $\alpha > 0$.  Let us
interpret (\ref{(diam A)^alpha}) as being equal to $1$ when $A$ is
bounded and nonempty, even if $A$ has only one element, so that $\diam
A = 0$.  If $A = \emptyset$, then we interpret (\ref{(diam A)^alpha})
as being equal to $0$ for every $\alpha \ge 0$, and we interpret
(\ref{(diam A)^alpha}) as being equal to $+\infty$ for every $\alpha
\ge 0$ when $A$ is unbounded.

        The \emph{$\alpha$-dimensional Hausdorff
content}\index{Hausdorff content} $H^\alpha_{con}(E)$ of a set $E
        \subseteq M$ is defined to be the infimum of the sums
\begin{equation}
\label{sum_j (diam A_j)^alpha}
        \sum_j (\diam A_j)^\alpha
\end{equation}
over all collections $\{A_j\}_j$ of finitely or countably many subsets
of $M$ such that
\begin{equation}
\label{E subseteq bigcup_j A_j}
        E \subseteq \bigcup_j A_j.
\end{equation}
More precisely, the sum (\ref{sum_j (diam A_j)^alpha}) can be defined
as a nonnegative extended real number as in Section \ref{nonnegative
  sums}, and $H^\alpha_{con}(E)$ is also a nonnegative extended real
number.  Note that
\begin{equation}
\label{H^alpha_{con}(E) le (diam E)^alpha}
        H^\alpha_{con}(E) \le (\diam E)^\alpha
\end{equation}
for every $E \subseteq M$, since one can use the covering of $E$ by
itself in the previous definition.  It is not difficult to verify that
$H^\alpha_{con}$ is an outer measure on $M$, by standard arguments.
Basically, $H^\alpha_{con}$ is the largest possible outer measure on
$M$ that satisfies (\ref{H^alpha_{con}(E) le (diam E)^alpha}).

        Similarly, $H^\alpha_\delta(E)$ is defined for $0 < \delta
\le +\infty$ as the infimum of the sums (\ref{sum_j (diam A_j)^alpha})
over all collections $\{A_j\}_j$ of finitely or countably many subsets
of $M$ that satisfy (\ref{E subseteq bigcup_j A_j}) and
\begin{equation}
\label{diam A_j < delta}
        \diam A_j < \delta
\end{equation}
for each $j$, when there are such coverings of $E$ in $M$.  Otherwise,
if there is no such covering of $E$ in $M$, then we simply put
$H^\alpha_\delta(E) = +\infty$.  Note that
\begin{equation}
\label{H^alpha_infty(E) = H^alpha_{con}(E)}
        H^\alpha_\infty(E) = H^\alpha_{con}(E)
\end{equation}
for every $E \subseteq M$, and that $H^\alpha_\delta(E)$ decreases
monotonically in $\delta$.  One can also check that $H^\alpha_\delta$
is an outer measure on $M$ for every $\delta > 0$, as before.

        Remember that $M$ is said to be
\emph{separable}\index{separable metric spaces} if there is a dense
subset of $M$ with only finitely or countably many elements.  This
implies that $M$ can be covered by finitely or countably many subsets
with diameter less than $\delta$ for each $\delta > 0$, and in fact
this condition is equivalent to separability.  It follows that every
subset of $M$ can be covered in this way when $M$ is separable.

        The \emph{$\alpha$-dimensional Hausdorff
measure}\index{Hausdorff measure} of $E \subseteq M$ is defined by
\begin{equation}
\label{H^alpha(E) = sup_{delta > 0} H^alpha_delta(E)}
        H^\alpha(E) = \sup_{\delta > 0} H^\alpha_\delta(E),
\end{equation}
which can also be considered as the limit of $H^\alpha_\delta(E)$ as
$\delta \to 0$, because of the monotonicity of $H^\alpha_\delta(E)$ in
$\delta$.  It is easy to see that $H^\alpha$ is an outer measure on
$M$, since $H^\alpha_\delta(E)$ is an outer measure on $M$ for each
$\delta > 0$.  If $H^\alpha(E) = 0$ for some $E \subseteq M$, then
$H^\alpha_\delta(E) = 0$ for each $\delta > 0$, and in particular
$H^\alpha_{con}(E) = 0$.  Conversely, if $H^\alpha_{con}(E) = 0$, then
the coverings of $E$ for which the sums (\ref{sum_j (diam A_j)^alpha})
are small involve subsets $A_j$ of $M$ with small diameter.  This
implies that $H^\alpha_\delta(E) = 0$ for every $\delta > 0$, and
hence that $H^\alpha(E) = 0$.

        Suppose that $E_1$ and $E_2$ are nonempty subsets of $M$
such that
\begin{equation}
\label{d(x, y) ge eta, 2}
        d(x, y) \ge \eta
\end{equation}
for some $\eta > 0$ and every $x \in E_1$ and $y \in E_2$.  If $A$ is
any subset of $M$ with diameter less than $\eta$, then it follows that
$A$ cannot intersect both $E_1$ and $E_2$.  Using this, one can check
that
\begin{equation}
\label{H^alpha_delta(E_1 cup E_2) ge H^alpha_delta(E_1) + H^alpha_delta(E_2)}
 H^\alpha_\delta(E_1 \cup E_2) \ge H^\alpha_\delta(E_1) + H^\alpha_\delta(E_2)
\end{equation}
when $0 < \delta \le \eta$.  More precisely, if $\delta \le \eta$,
then any covering of $E_1 \cup E_2$ by finitely or countably many
subsets of $M$ with diameter less than $\delta$ can be split into
coverings of $E_1$ and $E_2$ separately.  This leads to a splitting of
the corresponding sums (\ref{sum_j (diam A_j)^alpha}), which can be
used to obtain (\ref{H^alpha_delta(E_1 cup E_2) ge H^alpha_delta(E_1)
  + H^alpha_delta(E_2)}).  Under these conditions, we get that
\begin{equation}
\label{H^alpha(E_1 cup E_2) ge H^alpha(E_1) + H^alpha(E_2)}
        H^\alpha(E_1 \cup E_2) \ge H^\alpha(E_1) + H^\alpha(E_2),
\end{equation}
by taking the limit as $\delta \to 0$.  Thus $H^\alpha$ satisfies
Carath\'eodory's criterion, which implies that Borel subsets of $M$
are measurable with respect to $H^\alpha$.

        One can check that $H^\alpha$ reduces to counting measure on $M$
when $\alpha = 0$, using the conventions for defining (\ref{(diam
  A)^alpha}) when $\alpha = 0$ mentioned at the beginning of the
section.  In particular, to get $H^0(\emptyset) = 0$, one can cover
the empty set by itself.  Otherwise, one can let the empty set be
covered by the empty collection of subsets of $M$, and interpret an
empty sum as being equal to $0$, as before.

        It is easy to see that the diameter of a set $A \subseteq M$
is equal to the diameter of the closure $\overline{A}$ of $A$ in $M$.
One can also show that every set $A \subseteq M$ is contained in open
subsets of $M$ with diameter arbitrarily close to the diameter of $A$.
This implies that one can restrict one's attention to coverings of a
set $E \subseteq M$ by open or closed subsets of $M$ in the
definitions of $H^\alpha_{con}(E)$ and $H^\alpha_\delta(E)$, and get
the same result as before.  If $E$ is compact, then it follows that
one can restrict one's attention to coverings of $E$ by finitely many
subsets of $M$ in the definitions of $H^\alpha_{con}(E)$ and
$H^\alpha_\delta(E)$.

        Note that $H^\alpha_\delta$ is often defined using the condition
\begin{equation}
\label{diam A_j le delta}
        \diam A_j \le \delta
\end{equation}
instead of (\ref{diam A_j < delta}), which leads to an equivalent
definition of $H^\alpha$.  An advantage of using (\ref{diam A_j <
  delta}) is that one can more easily restrict one's attention to
coverings by open subsets of $M$, as in the previous paragraph.
Otherwise, if one uses (\ref{diam A_j le delta}) and restricts one's
attention to covering by open sets, then one gets the same result for
$H^\alpha$ in the limit as $\delta \to 0$, but not necessarily for
each $\delta > 0$ individually.

        If $M = {\bf R}$ with the standard metric, then one can
restrict one's attention to coverings of $E \subseteq {\bf R}$ by
intervals in the definitions of $H^\alpha_{con}(E)$ and
$H^\alpha_\delta(E)$.  This is especially nice when $\alpha = 1$, for
which one can check that $H^1_\delta = H^1_{con}$ for each $\delta >
0$, by subdividing intervals into smaller pieces.  It follows that
$H^1 = H^1_{con}$ on ${\bf R}$, which is the same as Lebesgue outer
measure on ${\bf R}$.  If $\alpha = 0$, then it is helpful to consider
the empty set as an interval in ${\bf R}$, to deal with the case where
$E = \emptyset$.  Otherwise, one can avoid the problem by allowing
the empty set to be covered by the empty collection of subsets of ${\bf R}$.

        If $A$ is a bounded subset of any (nonempty) metric space $M$,
then $A$ is contained in a closed ball in $M$ with radius equal to the
diameter of $A$.  Any closed ball $B$ in $M$ with radius $r$ has
diameter less than or equal to $2 \, r$, by the triangle inequality.
If $d(x, y)$ is an ultrametric on $M$, then the diameter of a closed
ball in $M$ of radius $r$ is less than or equal to $r$.  In this case,
it follows that one can restrict one's attention to coverings of a set
$E \subseteq M$ by closed balls in $M$ in the definitions of
$H^\alpha_{con}(E)$ and $H^\alpha_\delta(E)$ , and get the same result
as before, at least if $\alpha > 0$.  If $\alpha = 0$, then one should
allow the empty set to be covered by itself, or by the empty collection
of subsets of $M$, as usual.

        Let $(M, d(x, y))$ be an arbitrary metric space again,
and let $Y$ be a subset of $M$.  Of course, $Y$ can also be considered
as a metric space, using the restriction of $d(x, y)$ to $x, y \in Y$.
If $E \subseteq Y$, then one can restrict one's attention to coverings
of $E$ by subsets of $Y$ in the definition of $H^\alpha_{con}(E)$ and
$H^\alpha_\delta(E)$, and get the same result as for $E$ as a subset
of $M$.  More precisely, every covering of $E$ in $Y$ can be
considered as a covering of $E$ in $M$, and every covering of $E$ in
$M$ leads to a covering of $E$ in $Y$, by taking the intersections of
the subsets of $M$ in the covering of $E$ with $Y$.  It follows that
the definitions of $H^\alpha_{con}(E)$, $H^\alpha_\delta(E)$, and
$H^\alpha(E)$ for $E$ as a subset of $Y$ are equivalent to the
analogous definitions for $E$ as a subset of $M$.

        Let $E$ be a subset of $M$ such that $H^\alpha_\delta(E) < \infty$
for some $\delta > 0$.  By definition of $H^\alpha_\delta(E)$, for
each $\epsilon > 0$ there is a collection $\{A_j\}_j$ of finitely or
countably many subsets of $M$ such that $E$ is contained in the union
of the $A_j$'s, the diameter of $A_j$ is less than $\delta$ for each $j$,
and
\begin{equation}
\label{sum_j (diam A_j)^alpha < H^alpha_delta(E) + epsilon}
        \sum_j (\diam A_j)^\alpha < H^\alpha_\delta(E) + \epsilon.
\end{equation}
As before, we can also choose the $A_j$'s to be open subsets of $M$,
so that
\begin{equation}
\label{U = U(epsilon, delta) = bigcup_j A_j}
        U = U(\epsilon, \delta) = \bigcup_j A_j
\end{equation}
is an open set in $M$ as well.  By construction,
\begin{equation}
\label{H^alpha_delta(U) le sum_j (diam A_j)^alpha}
        H^\alpha_\delta(U) \le \sum_j (\diam A_j)^\alpha,
\end{equation}
since the $A_j$'s can be used to cover $U$ in the definition of
$H^\alpha_\delta(U)$ too.

        Suppose now that $H^\alpha(E) < \infty$, so that
$H^\alpha_\delta(E) < \infty$ for each $\delta > 0$.  Thus
we can apply the remarks in the previous paragraph to $\epsilon = \delta
= 1/n$ for each positive integer $n$, to get an open set $U_n \subseteq M$
such that $E \subseteq U_n$ and
\begin{equation}
\label{H^alpha_{1/n}(U_n) < H^alpha_{1/n}(E) + 1/n le H^alpha(E) + 1/n}
        H^\alpha_{1/n}(U_n) < H^\alpha_{1/n}(E) + 1/n \le H^\alpha(E) + 1/n.
\end{equation}
Put
\begin{equation}
\label{widetilde{E} = bigcap_{n = 1}^infty U_n}
        \widetilde{E} = \bigcap_{n = 1}^\infty U_n,
\end{equation}
so that $\widetilde{E}$ is a $G_\delta$ set in $M$, and hence a Borel
set, and $E \subseteq \widetilde{E}$.  Of course, $\widetilde{E}
\subseteq U_n$ for each $n$, which implies that
\begin{equation}
\label{H^alpha_{1/n}(widetilde{E}) le H^alpha_{1/n}(U_n)}
        H^\alpha_{1/n}(\widetilde{E}) \le H^\alpha_{1/n}(U_n)
\end{equation}
for each $n$.  Combining this with (\ref{H^alpha_{1/n}(U_n) <
  H^alpha_{1/n}(E) + 1/n le H^alpha(E) + 1/n}) and taking the limit as
$n \to \infty$, we get that
\begin{equation}
\label{H^alpha(widetilde{E}) le H^alpha(E)}
        H^\alpha(\widetilde{E}) \le H^\alpha(E).
\end{equation}
It follows that
\begin{equation}
\label{H^alpha(widetilde{E}) = H^alpha(E)}
        H^\alpha(\widetilde{E}) = H^\alpha(E),
\end{equation}
since the opposite inequality holds automatically, because $E
\subseteq \widetilde{E}$.  This shows that $H^\alpha$ is Borel
regular as an outer measure on $X$.

\section{Hausdorff measures, continued}
\label{hasudorff measures, continued}

        Let $k$ be a field, and let $|\cdot|$ be an absolute value
function on $k$ which is nontrivial and discrete.  As in Section
\ref{discrete absolute value functions}, this implies that there is a
real number $\rho_1$ such that $0 < \rho_1 < 1$ and the positive
values of $|\cdot|$ on $k$ are the same as the integer powers of
$\rho_1$.  Remember that every closed ball in $k$ of radius $r \ge 0$
with respect to the ultrametric associated to $|\cdot|$ has diameter
less than or equal to $r$, as mentioned in the previous section.  In
this situation, if $r$ is an integer power of $\rho_1$, then the
diameter of a closed ball in $k$ of radius $r$ is equal to $r$.
Of course, every closed ball in $k$ of positive radius is the same
as a closed ball with radius equal to an integer power of $\rho_1$.

        Suppose in addition that the residue field associated to
$|\cdot|$ on $k$ as in Section \ref{residue field} has exactly $N$
elements for some integer $N \ge 2$.  Let $\alpha$ be the real number
determined by
\begin{equation}
\label{rho_1^alpha = 1/N}
        \rho_1^\alpha = 1/N,
\end{equation}
and observe that $\alpha > 0$ under these conditions.  Also let
$H^\alpha_{con}$, $H^\alpha_\delta$, and $H^\alpha$ be the outer
measures on $k$ corresponding to the ultrametric on $k$ associated to
$|\cdot|$ as in the previous section.  If $j$, $l$ are integers and $l
\ge 0$, then every closed ball in $k$ of radius $\rho_1^j$ can be
expressed as the union of $N^l$ pairwise-disjoint closed balls of
radius $\rho_1^{j + l}$, as in Section \ref{residue field} again.
This implies that
\begin{equation}
\label{H^alpha_delta(overline{B}(x, rho_1^j)) le ... = rho_1^{alpha j}}
 H^\alpha_\delta(\overline{B}(x, \rho_1^j)) \le N^l \, (\rho_1^{j + l})^\alpha
                                            = \rho_1^{\alpha \, j}
\end{equation}
for every $x \in k$ when $\rho_1^{j + l} < \delta$, using the
definition of $H^\alpha_\delta$ in the first step, and
(\ref{rho_1^alpha = 1/N}) in the second step.  It follows that
\begin{equation}
\label{H^alpha(overline{B}(x, rho_1^j)) le rho_1^{alpha j}}
        H^\alpha(\overline{B}(x, \rho_1^j)) \le \rho_1^{\alpha \, j}
\end{equation}
for every $x \in k$ and $j \in {\bf Z}$, by taking the supremum of the
left side of (\ref{H^alpha_delta(overline{B}(x, rho_1^j)) le ... =
  rho_1^{alpha j}}) over $\delta > 0$.  Let us check that
\begin{equation}
\label{H^alpha_delta(E) = H^alpha_{con}(E)}
        H^\alpha_\delta(E) = H^\alpha_{con}(E)
\end{equation}
for every $E \subseteq k$ and $\delta > 0$, so that
\begin{equation}
\label{H^alpha(E) = H^alpha_{con}(E)}
        H^\alpha(E) = H^\alpha_{con}(E).
\end{equation}
Of course, $H^\alpha_{con}(E)$ is automatically less than or equal to
$H^\alpha_\delta(E)$ for each $\delta > 0$, and so it suffices to
verify the opposite inequality.  One way to do this is to observe that
any covering of $E$ by finitely or countably many closed balls can be
replaced by a covering of $E$ by closed balls with arbitrarily small
radius, by covering each ball by balls of smaller radius, as before.
This keeps the same sum as in (\ref{sum_j (diam A_j)^alpha}), because
of the way that $\alpha$ was chosen.  Alternatively, one can use the
definition of $H^\alpha_{con}(E)$, (\ref{H^alpha_delta(overline{B}(x,
  rho_1^j)) le ... = rho_1^{alpha j}}), and countable subadditivity of
$H^\alpha_\delta$.

        Let us now ask that $k$ be complete with respect to the
ultrametric associated to $|\cdot|$, in addition to the hypotheses
mentioned earlier.  It follows that closed balls in $k$ are compact,
because closed totally bounded subsets of a complete metric space are
compact.  In this case, it is not too difficult to show that
\begin{equation}
\label{H^alpha_{con}(overline{B}(x, rho_1^j)) ge rho_1^{alpha j}}
        H^\alpha_{con}(\overline{B}(x, \rho_1^j)) \ge \rho_1^{\alpha \, j}
\end{equation}
for every $x \in k$ and $j \in {\bf Z}$.  More precisely, because
$\overline{B}(x, \rho_1^j)$ is compact, it suffices to consider
coverings of $\overline{B}(x, \rho_1^j)$ by finitely many open subsets
of $k$ in the definition of $H^\alpha_{con}(\overline{B}(x,
\rho_1^j))$, as in the previous section.  In fact, it suffices to
consider coverings of $\overline{B}(x, \rho_1^j)$ by finitely many
closed balls of positive radius, because the metric on $k$ is an
ultrametric.  In the present situation, one can reduce further to
coverings of $\overline{B}(x, \rho_1^j)$ by finitely many closed balls
of the same radius $\rho_1^{j + l}$ for some $l \in {\bf Z}_+$, since
one can cover balls of arbitrary radius by balls with smaller radius
as before.  This also uses the way that $\alpha$ was chosen, to ensure
that sums like (\ref{sum_j (diam A_j)^alpha}) are maintained when one
reduces to covering by balls of smaller radius.  However, one can
check that $\overline{B}(x, \rho_1^j)$ cannot be covered by fewer that
$N^l$ closed balls of radius $\rho_1^{j + l}$ for any $l \in {\bf
  Z}_+$ under these conditions.  This implies
(\ref{H^alpha_{con}(overline{B}(x, rho_1^j)) ge rho_1^{alpha j}}), as
desired.

        Combining (\ref{H^alpha(overline{B}(x, rho_1^j)) le rho_1^{alpha j}})
and (\ref{H^alpha_{con}(overline{B}(x, rho_1^j)) ge rho_1^{alpha j}}),
we get that
\begin{equation}
\label{H^alpha(overline{B}(x, rho_1^j)) = rho_1^{alpha j}}
        H^\alpha(\overline{B}(x, \rho_1^j)) = \rho_1^{\alpha \, j}
\end{equation}
for every $x \in k$ and $j \in {\bf Z}$.  Of course, Hausdorff measure
of any dimension is invariant under isometries on any metric space, by
construction.  In particular, Hausdorff measure of any dimension is
invariant under translations on $k$.  This implies that $H^\alpha$
satisfies the requirements of Haar measure on $k$, because $H^\alpha$
is finite on bounded subsets of $k$, and positive on nonempty open
subsets of $k$, by (\ref{H^alpha(overline{B}(x, rho_1^j)) =
  rho_1^{alpha j}}).

\section{Lipschitz mappings, revisited}
\label{lipschitz mappings, revisited}

        Let $(M_1, d_1(x, y))$ and $(M_2, d_2(u, v))$ be metric spaces,
and suppose that $f$ is a Lipschitz mapping of order $a > 0$ from $M_1$
into $M_2$ with constant $C$.  If $A$ is a bounded subset of $M_1$,
then $f(A)$ is a bounded subset of $M_2$, and
\begin{equation}
\label{diam f(A) le C (diam A)^a}
        \diam f(A) \le C \, (\diam A)^a,
\end{equation}
where the diameters in (\ref{diam f(A) le C (diam A)^a}) are taken in
the appropriate metric space.  Using this, one can check that
\begin{equation}
\label{H^alpha_{con}(f(E)) le C^alpha H^{alpha a}_{con}(E)}
        H^\alpha_{con}(f(E)) \le C^\alpha \, H^{\alpha \, a}_{con}(E)
\end{equation}
for every $E \subseteq M_1$ and $\alpha \ge 0$, where the Hausdorff
content on each side is taken in the appropriate metric space.  More
precisely, $C^\alpha$ should be interpreted as being equal to $1$ for
every $C \ge 0$ when $\alpha = 0$.  If $\alpha > 0$ and $C = 0$, then
the right side of (\ref{H^alpha(f(E)) le C^alpha H^{alpha a}(E)}) may
be interpreted as being equal to $0$, even when $H^{\alpha \, a}(E) =
\infty$.

        Let $\delta > 0$ be given, and suppose that
\begin{equation}
\label{delta' > 0 and delta' ge C delta^a}
        \delta' > 0 \quad\hbox{and}\quad \delta' \ge C \, \delta^a.
\end{equation}
In analogy with (\ref{H^alpha_{con}(f(E)) le C^alpha H^{alpha a}_{con}(E)}),
we have that
\begin{equation}
\label{H^alpha_{delta'}(f(E)) le C^alpha H^{alpha a}_delta(E)}
        H^\alpha_{\delta'}(f(E)) \le C^\alpha \, H^{\alpha \, a}_\delta(E)
\end{equation}
for every $E \subseteq M_1$ and $\alpha \ge 0$, and with the same
conventions as before.  If $C > 0$, then one might as well take
\begin{equation}
\label{delta' = C delta^a}
        \delta' = C \, \delta^a,
\end{equation}
and otherwise one can take any $\delta' > 0$ when $C = 0$.  It follows
that
\begin{equation}
\label{H^alpha(f(E)) le C^alpha H^{alpha a}(E)}
        H^\alpha(f(E)) \le C^\alpha \, H^{\alpha \, a}(E)
\end{equation}
for every $E \subseteq M_1$ and $\alpha \ge 0$, since $\delta > 0$ is
arbitrary.

        Suppose now that
\begin{equation}
\label{C^{-1} d_1(x, y)^a le d_2(f(x), f(y)) le C d_1(x, y)^a, 2}
        C^{-1} \, d_1(x, y)^a \le d_2(f(x), f(y)) \le C \, d_1(x, y)^a
\end{equation}
for some $a > 0$ and $C \ge 1$, and for every $x, y \in M_1$.  This
implies that
\begin{equation}
\label{C^{-1} (diam A)^a le diam f(A) le C (diam A)^a}
        C^{-1} \, (\diam A)^a \le \diam f(A) \le C \, (\diam A)^a
\end{equation}
for every $A \subseteq M_1$, as in (\ref{diam f(A) le C (diam A)^a}).
As in Section \ref{lipschitz mappings}, the inverse of $f$ is
Lipschitz of order $1/a$ with constant $C^{1/a}$ as a mapping from
$f(M_1)$ into $M_1$, because of the first inequality in (\ref{C^{-1}
  d_1(x, y)^a le d_2(f(x), f(y)) le C d_1(x, y)^a, 2}).  Remember that
the Hausdorff content of $f(E)$ as a subset of $M_2$ is the same as
the Hausdorff content of $f(E)$ as a subset of $f(M_1)$, as in Section
\ref{hausdorff measures}, because $f(E) \subseteq f(M_1) \subseteq
M_2$.  Applying (\ref{H^alpha_{con}(f(E)) le C^alpha H^{alpha
    a}_{con}(E)}) to the inverse of $f$, we get that
\begin{equation}
\label{H^{alpha a}_{con}(E) le ... = C^alpha H^alpha_{con}(f(E))}
        H^{\alpha \, a}_{con}(E) 
               \le (C^{1/a})^{\alpha \, a} \, H^{(\alpha \, a)/a}_{con}(f(E))
                  = C^\alpha \, H^\alpha_{con}(f(E))
\end{equation}
for every $E \subseteq M_1$ and $\alpha \ge 0$.  More precisely, we
are applying (\ref{H^alpha_{con}(f(E)) le C^alpha H^{alpha
    a}_{con}(E)}) to $f(E)$ instead of $E$, to $\alpha \, a$ instead
of $\alpha$, to $1/a$ instead of $a$, and to $C^{1/a}$ instead of $C$.
It follows that
\begin{equation}
\label{C^{-alpha} H^{alpha a}_{con}(E) le H^alpha_{con}(f(E)) le ...}
        C^{-\alpha} \, H^{\alpha \, a}_{con}(E) \le H^\alpha_{con}(f(E))
                                      \le C^\alpha \, H^{\alpha \, a}_{con}(E)
\end{equation}
for every $E \subseteq M_1$ and $\alpha \ge 0$, where the first
inequality is equivalent to (\ref{H^{alpha a}_{con}(E) le ... =
  C^alpha H^alpha_{con}(f(E))}), and the second inequality is the
same as (\ref{H^alpha_{con}(f(E)) le C^alpha H^{alpha a}_{con}(E)}).

        Similarly, let $\delta_2 > 0$ be given, and put
\begin{equation}
\label{delta_1 = C^{1/a} delta_2^{1/a}}
        \delta_1 = C^{1/a} \, \delta_2^{1/a},
\end{equation}
which is the same as saying that
\begin{equation}
\label{delta_2 = C^{-1} delta_1^a}
        \delta_2 = C^{-1} \, \delta_1^a.
\end{equation}
One can check that
\begin{equation}
\label{H^{alpha a}_{delta_1}(E) le ... = C^alpha H^alpha_{delta_2}(f(E))}
        H^{\alpha \, a}_{\delta_1}(E)
               \le (C^{1/a})^{\alpha \, a} \, H^{(\alpha \, a)/a}_{\delta_2}(f(E))
                        = C^\alpha \, H^\alpha_{\delta_2}(f(E))
\end{equation}
for every $E \subseteq M_1$ and $\alpha \ge 0$, by applying
(\ref{H^alpha_{delta'}(f(E)) le C^alpha H^{alpha a}_delta(E)}) to the
inverse of $f$, with $\delta_2$, $\delta_1$ in place of $\delta$,
$\delta'$, respectively.  As before, we are also applying
(\ref{H^alpha_{delta'}(f(E)) le C^alpha H^{alpha a}_delta(E)}) to
$f(E)$ instead of $E$, to $\alpha \, a$ instead of $\alpha$, to $1/a$
instead of $a$, and to $C^{1/a}$ instead of $C$.  This implies that
\begin{equation}
\label{H^{alpha a}(E) le C^alpha H^alpha(f(E))}
        H^{\alpha \, a}(E) \le C^\alpha \, H^\alpha(f(E))
\end{equation}
for every $E \subseteq M_1$ and $\alpha \ge 0$, which could be derived
from (\ref{H^alpha(f(E)) le C^alpha H^{alpha a}(E)}) applied to the
inverse of $f$ as well, with the same substitutions as before.  It
follows that
\begin{equation}
\label{C^{-alpha} H^{alpha a}(E) le H^alpha(f(E)) le C^alpha H^{alpha a}(E)}
        C^{-\alpha} \, H^{\alpha \, a}(E) \le H^\alpha(f(E))
                                    \le C^\alpha \, H^{\alpha \, a}(E)
\end{equation}
for every $E \subseteq M$ and $\alpha \ge 0$, using (\ref{H^{alpha
    a}(E) le C^alpha H^alpha(f(E))}) in the first step, and
(\ref{H^alpha(f(E)) le C^alpha H^{alpha a}(E)}) in the second step.

\section{Local Lipschitz conditions}
\label{local lipschitz conditions}

        Let $(M_1, d_1(x, y))$, $(M_2, d_2(u, v))$ be metric spaces
again, let $a$, $\eta$ be positive real numbers, and let $C$ be a
nonnegative real number.  A mapping $f : M_1 \to M_2$ is said to be
\emph{locally Lipschitz}\index{locally Lipschitz mappings} of order
$a$ at the scale of $\eta$ with constant $C$ if
\begin{equation}
\label{d_2(f(x), f(y)) le C d_1(x, y)^a, 2}
        d_2(f(x), f(y)) \le C \, d_1(x, y)^a
\end{equation}
for every $x, y \in M_1$ with
\begin{equation}
\label{d_1(x, y) le eta}
        d_1(x, y) \le \eta.
\end{equation}
Alternatively, one might prefer to ask that (\ref{d_2(f(x), f(y)) le C
  d_1(x, y)^a, 2}) hold for every $x, y \in M_1$ with
\begin{equation}
\label{d_1(x, y) < eta}
        d_1(x, y) < \eta,
\end{equation}
instead of (\ref{d_1(x, y) le eta}).  This would imply that
(\ref{d_2(f(x), f(y)) le C d_1(x, y)^a, 2}) holds for all $x, y \in
M_1$ that can be approximated by elements of $M_1$ at distance
strictly less than $\eta$, by continuity.  If $x, y \in M_1$ can be
approximated by elements of $M_1$ at distance strictly less than
$\eta$, then $x$, $y$ satisfy (\ref{d_1(x, y) le eta}), but the
converse does not always hold.

        Suppose that $f : M_1 \to M_2$ is locally Lipschitz of
order $a > 0$ at the scale of $\eta > 0$ and with constant $C \ge 0$.
It is easy to see that (\ref{diam f(A) le C (diam A)^a}) still holds
for any bounded set $A \subseteq M_1$ with
\begin{equation}
\label{diam A le eta}
        \diam A \le \eta.
\end{equation}
If
\begin{equation}
\label{0 < delta le eta}
        0 < \delta \le \eta
\end{equation}
and $\delta'$ satisfies (\ref{delta' > 0 and delta' ge C delta^a}),
then (\ref{H^alpha_{delta'}(f(E)) le C^alpha H^{alpha a}_delta(E)})
still holds for every $E \subseteq M_1$ and $\alpha > 0$, for
essentially the same reasons as before.  This implies that
(\ref{H^alpha(f(E)) le C^alpha H^{alpha a}(E)}) holds for every $E
\subseteq M_1$ and $\alpha > 0$, basically by taking $\delta$ and
$\delta'$ arbitrarily small, as before.  If the local Lipschitz
condition is defined in terms of (\ref{d_1(x, y) < eta}) instead of
(\ref{d_1(x, y) le eta}), then (\ref{diam f(A) le C (diam A)^a}) holds
when
\begin{equation}
\label{diam A < eta}
        \diam A < \eta,
\end{equation}
instead of (\ref{diam A le eta}).  In this case,
(\ref{H^alpha_{delta'}(f(E)) le C^alpha H^{alpha a}_delta(E)}) still
holds for every $E \subseteq M_1$ and $\alpha > 0$ when $\delta$ and
$\delta'$ satisfy (\ref{0 < delta le eta}) and (\ref{delta' > 0 and
  delta' ge C delta^a}), because of the strict inequality in
(\ref{diam A_j < delta}) in Section \ref{hausdorff measures}.  If one
were to also use the non-strict inequality (\ref{diam A_j le delta})
instead of the strict inequality (\ref{diam A_j < delta}) in the
definition of these outer measures, then one should ask that
\begin{equation}
\label{0 < delta < eta}
        0 < \delta < \eta
\end{equation}
instead of (\ref{0 < delta le eta}), in order to get
(\ref{H^alpha_{delta'}(f(E)) le C^alpha H^{alpha a}_delta(E)}).  In
each variant, one still gets (\ref{H^alpha(f(E)) le C^alpha H^{alpha
    a}(E)}) for every $E \subseteq M_1$ and $\alpha > 0$, by taking
$\delta$ and $\delta'$ to be arbitrarily small.

        Let us now restrict our attention to $a = 1$, in which case
we may simply say that a mapping is locally Lipschitz at the scale of
$\eta$ with constant $C$.  Let $f$ be a mapping from $M_1$ into $M_2$
again, and let $x$ be an element of $M_1$.  Using the notation in
Section \ref{mappings between metric spaces}, we have that
\begin{equation}
\label{D_r(f)(x) le C}
        D_r(f)(x) \le C
\end{equation}
for $0 < r \le \eta$ if and only if
\begin{equation}
\label{widetilde{D}_eta(f)(x) le C}
        \widetilde{D}_\eta(f)(x) \le C,
\end{equation}
which happens if and only if (\ref{d_2(f(x), f(y)) le C d_1(x, y)^a,
  2}) holds for every $y \in M_1$ that satisfies (\ref{d_1(x, y) le
  eta}).  In particular, $f$ is locally Lipschitz at the scale of
$\eta$ with constant $C$ on $M_1$ if and only if (\ref{D_r(f)(x) le
  C}) holds for every $x \in M_1$ and $0 < r \le \eta$, which is the
same as saying that (\ref{widetilde{D}_eta(f)(x) le C}) holds for
every $x \in M_1$.  It follows that the restriction of $f$ to
\begin{equation}
\label{{x in M_1 : widetilde{D}_eta(f)(x) le C}}
        \{x \in M_1 : \widetilde{D}_\eta(f)(x) \le C\}
\end{equation}
is automatically locally Lipschitz at the scale of $\eta$ with
constant $C$.  If one defines local Lipschitz conditions in terms of
(\ref{d_1(x, y) < eta}) instead of (\ref{d_1(x, y) le eta}), then
it is better to consider $0 < r < \eta$ in (\ref{D_r(f)(x) le C}),
and
\begin{equation}
\label{widetilde{D}_t(f)(x) le C}
        \widetilde{D}_t(f)(x) \le C
\end{equation}
for $0 < t < \eta$ instead of (\ref{widetilde{D}_eta(f)(x) le C}).
Alternatively, one could modify the definition (\ref{D_r(f)(x) =
  r^{-1} sup {d_2(f(x), f(y)) : y in M_1, d_1(x, y) le r}}) of
$D_r(f)(x)$, by taking the supremum over $y \in M_1$ with $d_1(x, y) <
r$.  One could still define $\widetilde{D}_t(f)(x)$ as in
(\ref{widetilde{D}_t(f)(x) = sup_{0 < r le t} D_r(f)(x)}), using the
modified definition of $D_r(f)(x)$, or by taking the supremum over $0
< r < t$ of either definition of $D_r(f)(x)$.  In the second
characterization (\ref{widetilde{D}_t(f)(x) = sup {frac{d_2(f(x),
      f(y))}{d_1(x, y)} : ...}}) of $\widetilde{D}_t(f)(x)$, one
should then take the supremum over $y \in M_1$ such that $0 < d_1(x,
y) < t$, when there is such a point $y$.  This second characterization
of $\widetilde{D}_t(f)(x)$ corresponds to taking $r > d_1(x, y)$ with
$r$ close to $d_1(x, y)$ in the definition of $\widetilde{D}_t(f)(x)$,
instead of $r = d_1(x, y)$, as before.

        Note that
\begin{equation}
\label{D(f)(x) < C'}
        D(f)(x) < C'
\end{equation}
for some $x \in M_1$ and $C' > 0$ if and only if there is a $t > 0$
such that
\begin{equation}
\label{widetilde{D}_t(f)(x) < C'}
        \widetilde{D}_t(f)(x) < C',
\end{equation}
by the definition (\ref{D(f)(x) = limsup_{r to 0} D_r(f)(x) = inf_{t >
    0} widetilde{D}_t(f)(x)}) of $D(f)(x)$.  Thus
\begin{equation}
\label{{x in M_1 : D(f)(x) < C'}}
        \{x \in M_1 : D(f)(x) < C'\}
\end{equation}
is the same as the union of
\begin{equation}
\label{{x in M_1 : widetilde{D}_t(f)(x) < C'}}
        \{x \in M_1 : \widetilde{D}_t(f)(x) < C'\}
\end{equation}
over $t > 0$.  More precisely, it suffices to use a sequence of
positive real numbers $t$ that converges to $0$.  One can also check
that $D(f)(x)$ would not be affected by modifying the definitions of
$D_r(f)(x)$ or $\widetilde{D}_t(f)(x)$ as in the preceding paragraph,
using strict inequalities instead of non-strict inequalities in the
relevant suprema.  This would correspond to slightly different
interpretations of (\ref{widetilde{D}_t(f)(x) < C'}) and (\ref{{x in
    M_1 : widetilde{D}_t(f)(x) < C'}}).

\backmatter

\newpage

\addcontentsline{toc}{chapter}{Index}

\printindex

\end{document}